%% file: AAA_main_manuscript_v26.tex
\documentclass[a4paper,12pt]{book}

\usepackage{amssymb,amsmath,amsthm,amsfonts}
\usepackage{bookmark}
\usepackage{mathrsfs}
\usepackage{multirow}
\usepackage{graphicx}
\usepackage{authblk}
\usepackage{indentfirst}
\usepackage{multicol}
\usepackage{tabu}
\usepackage{url}
\usepackage{fancyhdr}
\usepackage[numbers]{natbib}
\usepackage[all]{xy}
\usepackage{mathtools}
\usepackage[english]{babel}
\usepackage{colortbl}
\usepackage{caption}
\usepackage{hyperref}\usepackage{subcaption}
\usepackage{tikz}\usepackage{tikz-cd}\usepackage{float}\usepackage{enumitem}\usetikzlibrary{arrows.meta}
\usepackage{algorithm}
\usepackage{algpseudocode}
\usepackage{booktabs}
\usepackage{xcolor}
\definecolor{axisgray}{HTML}{808080}
\definecolor{mainblue}{HTML}{004488}
\definecolor{symred}{HTML}{CC3311}
\definecolor{arcangle}{HTML}{C86400}

\newtheorem{theorem}{Theorem}[section]
\newtheorem{proposition}[theorem]{Proposition}
\newtheorem{lemma}[theorem]{Lemma}
\newtheorem{corollary}[theorem]{Corollary}

\newtheorem{definition}{Definition}[section]
\newtheorem{example}{Example}[section]
\newtheorem{remark}{Remark}[section]

\DeclareMathOperator{\tr}{tr}
\DeclareMathOperator{\Iso}{Iso}
\DeclareMathOperator{\Sym}{Sym}
\DeclareMathOperator{\Conf}{Conf}
\DeclareMathOperator{\Span}{Span}
\DeclareMathOperator{\Aut}{Aut}

\usepackage{xr}
\makeatletter
    \newcommand*{\addFileDependency}[1]{
    \typeout{(#1)}
    \@addtofilelist{#1}
    \IfFileExists{#1}{}{\typeout{No file #1.}}
    }
\makeatother

\title{Parametrization of Symmetry in Data}

\author[1,2]{Jian Liu}
\author[2]{Dong Chen}
\author[2,3,4]{Guo-Wei Wei \thanks{Corresponding author: weig@msu.edu}}

\affil[1]{Mathematical Science Research Center, Chongqing University of Technology, Chongqing 400054, China}
\affil[2]{Department of Mathematics, Michigan State University, MI 48824, USA}
\affil[3]{Department of Electrical and Computer Engineering, Michigan State University, MI 48824, USA}
\affil[4]{Department of Biochemistry and Molecular Biology, Michigan State University, MI 48824, USA}

\makeatletter
    \renewcommand*{\@fnsymbol}[1]{\ensuremath{\ifcase#1\or \dagger\or *\or *\or
   \mathsection\or \else\@ctrerr\fi}}
\makeatother
\date{}

\begin{document}
    \maketitle


		\paragraph{Abstract}

        Symmetry plays a fundamental role in understanding natural phenomena and mathematical structures. This work develops a comprehensive theory for studying the persistent symmetries and degree of asymmetry of finite point configurations under parameterization in metric spaces. Leveraging category theory and span categories, we define persistent symmetry groups and introduce novel invariants called symmetry barcodes and polybarcodes that capture the birth, death, persistence, and reappearance of symmetries during parameter evolution. Metrics and stability theorems are established for these invariants. The concept of symmetry types is formalized via the action of isometry groups on configuration spaces. To quantitatively characterize symmetry and asymmetry, measures such as the degree of symmetry and symmetry defect are introduced; the latter leads to the study of approximate subgroups of orthogonal groups in Euclidean settings. Moreover, a theory of persistence representations of persistence groups is developed, generalizing the classical decomposition theorem of persistence modules. Persistent Fourier analysis on persistence groups is further proposed to characterize dynamic phenomena, including symmetry breaking and phase transitions. We also investigate equivariant persistent cohomology for persistence $G$-spaces, introducing a stratified barcode to visualize symmetry degradation encoded in torsion modules and establishing stability together with a finite algorithm based on Bredon cohomology. Algorithms for computing symmetry groups, barcodes, and symmetry defect in low-dimensional spaces are presented, complemented by discussions on extending symmetry analysis beyond geometric contexts. This work thus bridges geometric group theory, topological data analysis, representation theory, and machine learning, providing novel tools for the analysis of the parametrized symmetry of data.



\section*{Acknowledgments}

This work was supported in part by NIH grants R01GM126189, R01AI164266, and R35GM148196, National Science Foundation grants DMS2052983, DMS-1761320, and IIS-1900473, NASA  grant 80NSSC21M0023,   Michigan State University Research Foundation, and  Bristol-Myers Squibb  65109. Jian was also supported by the Natural Science Foundation of China (NSFC Grant No. 12401080) and the start-up research fund from Chongqing University of Technology.

\tableofcontents

\include{symmetry_introduction}

\chapter{Persistence symmetries}
\include{preliminaries}

\include{topological_dynamics}

\include{categorification_symmetries}
\include{metrics_symmetries}
\include{stability_symmetries}

\chapter{Symmetry analysis}
\include{symmetry_types}
\include{symmetry_degree_analysis}
\include{symmetry_defect_analysis}

\chapter{Persistence theory for group actions}
\include{persistent_representation}
\include{persistent_Fourier_analysis}

\include{equivariant_persistence}

\chapter{Computational aspects of symmetries}
\include{computation_aspects}

\bibliographystyle{plain}  
\bibliography{Reference}

\end{document}

%% file: symmetry_introduction.tex
\section{Introduction}

Symmetry is one of the fundamental concepts through which humans have come to understand nature and the world around them. It manifests in various forms throughout the natural world, influencing physical structures, biological systems, and even fundamental laws of physics, such as parity. Recognizing and analyzing symmetry not only enhances our comprehension of natural phenomena but also plays a crucial role in fields ranging from mathematics and physics to art and engineering.
There are numerous systematic and influential studies on symmetry across various fields, including mathematics and theoretical physics \cite{gromov1983filling,noether1983invariante,weyl2015symmetry}, chemistry \cite{hargittai1995symmetry,jaffe2002symmetry}, computer science \cite{forsyth2002computer,liu2010computational}, biology \cite{alberch1979size,thomson1917growth}, and art \cite{kubler2008shape}.

The concept of symmetry is intimately connected to one of the fundamental pillars of modern mathematics---group theory. A wide range of problems concerning symmetry can be rigorously formulated, analyzed, and resolved within the framework of group theory \cite{noether1921idealtheorie,schur2024neue}. In modern mathematics, many profound theories and work in geometry and topology are deeply rooted in symmetry, including geometric invariant theory \cite{mumford1994geometric}, linear algebraic groups \cite{borel2012linear}, fixed point theory \cite{atiyah1968lefschetz}, and the Weil conjectures \cite{deligne1974conjecture}.

With the rapid growth of big data and the surge of interest in artificial intelligence, the symmetry of data has attracted increasing attention. Most research on data symmetry is primarily focused on image and geometric data, with a particular emphasis on the long-standing and extensively explored problem of detecting symmetries in images and shapes \cite{kazhdan2004symmetry,marola2002detection,mitra2006partial,reisfeld1995context,sun2011reflection}. In contrast, symmetry detection in general data forms, such as point clouds, polygonal meshes, or other unstructured geometric datasets, remains relatively underdeveloped. 
Recently, symmetry has been used in neural networks to improve representation, generalization, and reduce model complexity by leveraging data invariances \cite{keck2025impact,keck2023learning}. Current symmetry analysis is primarily qualitative, often limited to categorical statements, from symmetric or asymmetric objects to quasiconformal and quasiregular mappings \cite{vuorinen2006conformal}. 
For instance, it typically lacks quantitative measures to assess the degree of asymmetry when symmetry is broken. Furthermore, there is little exploration of how symmetry or asymmetry persists under time evolution, systematic perturbation, or other parametric variation. Challenges stem not only from irregular parameterizations and noisy data but also from the need to define symmetry in a representation-independent, mathematically rigorous, and computationally efficient manner.

Persistent homology \cite{carlsson2009topology,edelsbrunner2002topological,zomorodian2004computing}, a recently emerging algebraic topology tool in data science, has demonstrated remarkable success in fields such as computer science, materials science, image science \cite{carlsson2008local}, chemistry \cite{townsend2020representation},  and molecular biology \cite{papamarkou2024position, pun2022persistent}. Homology groups characterize properties of data structures such as connected components, loops, and cavities, which capture certain global features of the data and exhibit strong robustness. However, these structural descriptors are often relatively coarse and limited \cite{su2025topological}. The key to the success of persistent homology lies in introducing ``persistence'', which enriches  structural characterization by providing a wealth of topological features over scales for complex and high-dimensional data. Recently, this idea has been extended to commutative algebra \cite{suwayyid2025persistent}. The fundamental principle of ``persistence'' is to study the dynamic, multi-scale, or multi-level information inherent in the data.

Symmetry, as another structural feature of data, has significant potential for applications, particularly in fields such as chemistry, materials science, and molecular biology, where symmetric objects are abundant in complex data. The core idea of this work is to introduce quantitative measures to assess the degree of asymmetry (asymmetricity) or the persistence of symmetry by incorporating dynamic information and/or parametrization into the study of data symmetry and asymmetry, thereby extracting sufficiently rich features for comprehensive multiscale analysis of complex data. These data can be dynamic or static. These features will be applied to machine learning and data-driven discovery. 

The symmetry of a geometric object can be understood as a transformation that preserves the object's rigidity---its shape. More precisely, let $M$ be a metric space and $X \subseteq M$ be a geometric object. A symmetry is an isometry $\pi: M \to M$ such that $\pi(X) = X$. The collection of all such symmetries forms the symmetry group of $X$, denoted by $\Sym(X)$. In this work, the geometric objects considered are typically finite point sets, which we call finite configurations. Since $M$ is a metric space, the set of all $n$-point configurations forms an $n$-point configuration space that describes all possible states of a system of $n$ points. This perspective motivates us to introduce the dynamics of finite configurations. However, in this work, time can be regarded as a parameter, and time evolution can be regarded as a filtration process. 

In Section \ref{section:dynamics_symmetries}, we introduce the dynamics of symmetries by examining how the symmetries of an $n$-point configuration (or simply an $n$-configuration) evolve over time. 
Certain symmetries may appear and disappear as time progresses, analogous to the birth and death of topological invariants observed in persistent homology filtration. Moreover, the space generated by these symmetries forms a persistent module, giving rise to the notion of persistent symmetries.

Section \ref{section:categorification} provides a categorical formalization of persistent symmetry. We construct the category $\mathcal{S}_n(M)$ of $n$-configurations on $M$, and a persistent $n$-configuration
\[
  \mathcal{F} : (\mathbb{R}, \leq) \to \mathcal{S}_n(M)
\]
is a functor encoding the dynamics of $n$-configurations. For real numbers $a \leq b$, a natural idea to define persistence is to construct a group homomorphism $\Sym(\mathcal{F}_a) \to \Sym(\mathcal{F}_b)$. Unfortunately, such a homomorphism generally does not exist. To address this, we introduce the \textit{path 2-category} $\mathbb{S}_n(M)$, where morphisms are finite sequences of configuration maps. This construction is essential for capturing the trajectory of symmetries: it distinguishes symmetries that survive continuously through a path from those that merely exist at the endpoints. Within this construction, we have a span of groups
\[
  \Sym(\mathcal{F}_a) \xleftarrow{f_{a,b}^\flat} \Sym_{f_{a,b}}(\mathcal{F}_a) \xrightarrow{f_{a,b}^\sharp} \Sym(\mathcal{F}_b).
\]
The assignment $\Sym : \mathbb{S}_n(M) \to \mathrm{Span}(\mathbf{Grp})$ then becomes a lax functor (see Theorem \ref{theorem:lax_functor}). Consequently, we define the \emph{persistent symmetry group} as the image of $f_{a,b}^\sharp$. It is worth noting that persistent symmetry groups are generally non-abelian, distinguishing them from persistent homology groups. Nonetheless, this does not hinder the definition of a \emph{symmetry barcode}, which records the birth and death information of these symmetries. On the other hand, by taking vector spaces generated by these persistent group elements, we obtain a persistence module $\mathcal{M} : \mathbb{S}_n(M) \to \mathbf{Vec}_{\mathbb{K}}$, which is a strict functor (see Proposition \ref{proposition:functor}). Finally, since each symmetry corresponds to a concrete isometry, we focus on the time intervals during which a particular symmetry exists. This leads to the concept of a \emph{polybarcode}, which may be a union of disjoint intervals as a symmetry that dies may later reappear, opening a new perspective for studying symmetry representations.

In Section \ref{section:polybarcodes_metrics}, we introduce several metrics on polybarcodes, including the symmetric difference distance, the expansion distance, and the interleaving distance. Section \ref{section:stability} develops stability results for both symmetry barcodes and polybarcodes. For the study of symmetry barcodes, we formalize the interleaving distance between pseudofunctors by lifting the symmetry construction to the category of truncated dynamical systems, providing the basis for their stability analysis. In the case of polybarcodes, the interleaving distance directly on persistence configurations may be insufficient to capture geometric disparity. To address this, we employ persistent dynamical systems to define the \emph{indexed interleaving distance} between persistence configurations. This refined approach yields effective bounds for the distances between polybarcodes. Specifically, we show that the left expansion distance coincides with the interleaving distance on polybarcodes and is bounded above by the indexed interleaving distance, leading to the following stability result (see Theorem~\ref{theorem:polybarcodes_stability2}):
\begin{theorem}
Let $\mathcal{F}, \mathcal{G}: (\mathbb{R}, \leq) \to \mathcal{S}_{n}(M)$ be two persistence $n$-configurations. Then
\[
  d_{L}(\mathcal{B}(\mathcal{F}), \mathcal{B}(\mathcal{G}))
  = d_{I}(\mathcal{B}(\mathcal{F}), \mathcal{B}(\mathcal{G}))
  \leq d_{II}(\mathcal{F}, \mathcal{G}),
\]
where $d_{L}$, $d_{I}$, and $d_{II}$ denote the left expansion distance, the interleaving distance, and the indexed interleaving distance, respectively.
\end{theorem}
It is worth noting that the expansion distance is often more intuitive and computationally accessible, whereas the interleaving distance is more abstract and well-suited for the formulation of stability theorems.

We then turn to the analysis of symmetries. In Section~\ref{section:symmetry_types}, the notion of \emph{symmetry type} is introduced. Specifically, let $\Conf_n(M)$ denote the configuration space and $\Iso(M)$ the isometry group of $M$. The natural action of $\Iso(M)$ on $\Conf_n(M)$ induces the orbit space
\[
\mathcal{M}_n = \Conf_n(M)/\Iso(M).
\]
Two $n$-configurations $X, Y \in \Conf_n(M)$ are said to have the same symmetry type if their symmetry groups are conjugate in $\Iso(M)$. Based on this notion, we construct a bicategory of symmetry types, and show that the assignment
\[
\Theta : \mathbb{S}_n(M) \longrightarrow \mathsf{SymType},
\]
is a pseudofunctor (see Propositions~\ref{proposition:symmetry_type_functor}). Symmetry types are further employed to investigate topological properties of point sets related to symmetries. For example, in Euclidean space $\mathbb{R}^{k}$, let $H(k) \subset \Conf_n(\mathbb{R}^k)$ denote the space of all configurations with centroid at the origin. For a fixed finite subgroup $V \leq \Iso(\mathbb{R}^{k})$, the set of $n$-configurations whose symmetry group is exactly $V$ is either empty or locally closed in $H(k)$ (see Proposition~\ref{proposition:locally_closed}). We also establish a decomposition of the configuration space based on symmetry types and analyze their closure relations (see Theorem~\ref{theorem:symmetry_type}). Finally, the notion of \emph{persistent symmetry type} is introduced. We define the symmetry type trajectory and use a digraph of groups to model the evolution of symmetry types over time.

Section \ref{section:symmetry_degree} introduces quantitative invariants to characterize the symmetry groups of configurations. We define the \emph{degree of symmetry} as a numerical measure of symmetric extent, alongside \emph{symmetry entropy} and the \emph{symmetry degree polynomial} to capture the distributional complexity of symmetry information. For persistence configurations, we extend these concepts to define the \emph{persistent degree of symmetry}, and visualize the evolution of symmetry degrees via weighted paths and analyze the structural relationships between symmetry groups using Cayley graphs constructed from the spans of persistent symmetry groups.

In contrast to degrees of symmetry, Section \ref{section:symmetry_defect} studies the degree of asymmetry of finite configurations. Asymmetry has been studied since early times \cite{klingenberg2002shape,kovesi1997symmetry}. To quantify how asymmetric a finite configuration is, we introduce the concept of \emph{symmetry defect}. Specifically, for an $n$-configuration $X$ in $\mathcal{S}_n(M)$ and a given isometry $\pi \in \Iso(M)$, the symmetry defect of $X$ with respect to $\pi$ is defined as
\[
  \mu(X,\pi) = \left( \inf_{\gamma \neq \pi|_X} \sum_{x \in X} d(x, \gamma(x))^p \right)^{1/p},
\]
where $d$ is the metric on the space $M$, $\pi|_X$ denotes the restriction of the isometry $\pi: M \to M$ to the set $X$, and the infimum is taken over all bijections $\gamma: X \to \pi(X)$ different from $\pi|_X$. The overall symmetry defect of $X$ is then defined by
\[
  \mu(X) = \inf_{\pi \in \Iso(M)} \mu(X,\pi),
\]
where the infimum is taken over all isometries $\pi$ of $M$. When $M$ is a Euclidean space, the symmetry defect is proved to be continuous (see Proposition~\ref{proposition:continuous_defect}). An interesting phenomenon occurs in Euclidean space: for a fixed configuration, the set of $\varepsilon$-approximate symmetries forms a compact \emph{approximate subgroup} of the orthogonal group. Approximate subgroups were introduced by Tao and others as a powerful tool in additive combinatorics and number theory \cite{breuillard2012structure,tao2008product,tao2006additive}. Finally, we propose a \emph{symmetry measure} ranging from 0 to 1 to quantify the fuzziness of asymmetry from completely asymmetric to fully symmetric configurations. The stability of the symmetry defect is also established.

A natural direction following the study of persistent symmetries is the investigation of \emph{persistence groups}.
The representation theory of persistence groups provides a tool to transform intricate non-abelian structures into linear, and hence commutative, algebraic objects.
In Section~\ref{section:persistent_representation}, we introduce the notion of a \emph{persistence representation} and analyze the behavior of its irreducible subrepresentations under the structural morphisms. As shown in Proposition~\ref{proposition:irreducible}, an irreducible subrepresentation of a persistence representation is mapped by the structure maps either to another irreducible representation or to the zero space.
This observation suggests that irreducible subrepresentations exhibit a \emph{birth-death} behavior analogous to that in persistent homology.
Moreover, an irreducible persistence representation supported in an interval corresponds precisely to the \emph{bar} associated with that irreducible representation.
We then establish the following decomposition theorem for persistence representations (see Theorem~\ref{theorem:decomposition_representation}):
\begin{theorem}
Let $\mathcal{R} : (\mathbb{Z}, \leq) \to \mathbf{Rep}_{\mathbb{K}}^{G}$ be a finite persistence representation.
Then there exists a finite collection of irreducible persistence representations such that
\[
\mathcal{R} \cong \bigoplus_{i \in I} \mathcal{R}_{i},
\]
where each $\mathcal{R}_{i}$ is an irreducible persistence representation supported on the interval $[s_i, t_i]$.
\end{theorem}
This result may be regarded as a generalization of the structure theorem for persistence modules.
In fact, when $G$ is the trivial group, each irreducible representation is simply a one-dimensional vector space generated by a single element, and the persistence of this one-dimensional space corresponds exactly to the persistence of a generator in a persistence module.
Furthermore, persistence representations can be interpreted as graded modules on a polynomial ring generated by the shift operator induced by the persistence morphisms (see Section~\ref{section:module_representations}).
Finally, we study the regular representations of persistence groups to understand how the irreducible subrepresentations of the persistent regular representation evolve with respect to the filtration parameter.

In Section~\ref{section:Fourier_analysis}, we introduce the theory of persistent Fourier analysis for persistence groups as a tool to characterize and understand the structure of persistence groups.
In contrast to classical Fourier analysis, the persistent setting is capable of capturing essential dynamical phenomena such as symmetry breaking, phase transitions, and stabilization.
Key quantities arising in persistent Fourier analysis, including spectral energy curves, persistent entropy, and dominant frequency ratios, provide stable and interpretable features that can be effectively applied to tasks such as classification, clustering, and model analysis.
Within this setting, we develop and analyze several fundamental notions, including the Fourier inversion formula, persistent convolution, persistent Fourier spectrum, and persistent Laplacian operator, together with their associated results.

In Section \ref{section:equivariant_persistenc}, we investigate the equivariant persistent cohomology of persistence $G$-spaces. Let $\mathcal{H}_G^\ast$ denote the cohomology ring of the classifying space $BG$. For a persistence $G$-space, the corresponding equivariant persistent cohomology is a persistence $\mathcal{H}_G^\ast$-module. The $\mathcal{H}_G^\ast$-module torsion encodes intrinsic symmetry information, such as the existence of non-free orbits, distinguishing it from the birth–death intervals in classical persistence. To visualize this algebraic stratification, we introduce a \textit{stratified barcode} representation that explicitly tracks symmetry degradation, such as transitions from free to torsion modules. We establish the stability of this invariant with respect to the interleaving distance. Finally, to address computational challenges, we develop a finite algorithm based on Bredon cohomology and a minimal regular subdivision procedure, enabling the practical computation of equivariant topological signatures.

Section \ref{section:computation} focuses on the computational aspects of symmetry analysis. We present algorithms and strategies for determining and representing symmetries in low-dimensional spaces, and for computing the symmetry groups of finite configurations. We detail procedures for extracting symmetry barcodes and polybarcodes, and develop approaches for quantifying symmetry defects. Furthermore, we introduce parametrization techniques for point cloud data to facilitate multiscale analysis. Finally, we propose strategies for the analysis of local symmetry.

For a brief outlook on symmetry analysis in data science, beyond the study of geometric symmetry, one can also develop tools based on other types of symmetry, such as combinatorial symmetry and topological symmetry, to enhance data processing and analysis. Moreover, although groups constitute the classical setting for characterizing symmetry, it is worthwhile to explore and develop alternative algebraic or geometric tools for investigating data symmetry. Additionally, current computational methods and algorithms for symmetry representation and analysis can be further optimized and extended. 

%% file: preliminaries.tex
\section{Preliminaries}

Symmetry is a fundamental concept across mathematics, physics, and many scientific disciplines. In modern mathematics, group theory offers a concise and powerful language for analyzing geometric symmetry. For a concise overview of symmetry in mathematics, one may refer to \cite{mostow1996brief}. 
Other symmetries, such as topological symmetries including noninvertible and categorical symmetries in quantum field theory \cite{freed2024topological}, and physical symmetries including supersymmetry and gauge symmetry in theoretical physics \cite{jost2018symmetries}, have attracted significant attention for their rich structures and wide-ranging applications.
In this section, we briefly review geometric symmetry in terms of group language, which provides the mathematical foundation for our study of symmetry in this work.

\subsection{Group action and symmetry}

Let $X$ be a set, and let $G$ be a group. A \textit{group action} of $G$ on $X$ is a map $\phi:G\times X\to X$ satisfying the following conditions:
\begin{itemize}
  \item[$(i)$] for each $x\in X$, $\phi(e,x)=x$, where $e$ is the identity element of $G$;
  \item[$(ii)$] for any $g,h\in G$ and $x\in X$, $\phi(h,\phi(g,x))=\phi(hg,x)$.
\end{itemize}

A classic geometric example is the square $S$. The dihedral group $D_4$ acts on $S$ to capture all of its rigid symmetries. As illustrated in Figure \ref{figure:square_symmetry}, the group $D_4$ consists of eight elements:
\begin{equation*}
e, \tau, \tau^2, \tau^3, \pi, \tau\pi, \tau^2\pi, \tau^3\pi,
\end{equation*}
where $e$ is the identity, $\tau$ is a counterclockwise rotation by $90^\circ$ satisfying $\tau^4 = e$, and $\pi$ is a reflection about the $y$-axis satisfying $\pi^2 = e$.
\begin{figure}[h]
\centering
\begin{tikzpicture}[scale=1.5, >=Stealth]

    \draw[->, color=axisgray, line width=0.4pt] (-1.8,0) -- (1.8,0) node[right, font=\small, text=black] {$x$};
    \draw[->, color=axisgray, line width=0.4pt] (0,-1.8) -- (0,1.8) node[above, font=\small, text=black] {$y$};

    \draw[dash pattern=on 4pt off 2pt, color=symred, line width=0.6pt] (0, -1.4) -- (0, 1.4); 
    \draw[dash pattern=on 4pt off 2pt, color=symred, line width=0.6pt] (-1.4, 0) -- (1.4, 0); 
    \draw[dash pattern=on 4pt off 2pt, color=symred, line width=0.6pt] (-1.3, -1.3) -- (1.3, 1.3); 
    \draw[dash pattern=on 4pt off 2pt, color=symred, line width=0.6pt] (-1.3, 1.3) -- (1.3, -1.3); 

    \draw[line width=1.2pt, color=mainblue] (-1,-1) rectangle (1,1);
    
    \foreach \x/\y in {-1/-1, -1/1, 1/1, 1/-1}
        \fill[color=mainblue] (\x,\y) circle (1.2pt);

    \draw[->, color=arcangle, line width=0.8pt] (0:0.5) arc[start angle=0, end angle=90, radius=0.5];
    \node[color=black, font=\footnotesize] at (30:0.75) {$90^\circ$};

\end{tikzpicture}
\caption{The square $S$ with its four axes of reflectional symmetry and a $90^\circ$ rotational symmetry.}
\label{figure:square_symmetry}
\end{figure}
It is straightforward to observe that $\tau, \tau^2$, and $\tau^3$ correspond to rotations by $90^\circ, 180^\circ$, and $270^\circ$, respectively. Under this setup, the elements $\pi, \tau\pi, \tau^2\pi$, and $\tau^3\pi$ correspond to reflections about the $y$-axis, the line $y = -x$, the $x$-axis, and the line $y = x$, respectively. 

Now, let $X=\mathbb{Z}^{2}$ represent an infinite square network (or grid). Naturally, we have a group action of the local symmetries
\begin{equation*}
D_4 \times X \to X.
\end{equation*}
Beyond rotations and reflections, an infinite square network is also invariant under discrete translations. Thus, we can define the extended group action via the semi-direct product:
\begin{equation*}
(D_4 \ltimes \mathbb{Z}^2) \times X \to X,
\end{equation*}
where the group multiplication in $D_4 \ltimes \mathbb{Z}^2$ is defined by $(g,a)\cdot (h,b)=(gh,a+gb)$ for $g,h\in D_{4}$ and $a,b\in \mathbb{Z}^{2}$. This group action beautifully demonstrates that the infinite square network exhibits translational symmetry alongside the rotational and reflectional symmetries of $D_4$.

Another fundamental example is the action of the symmetric group on a finite set. If $X$ contains $n$ elements, the \textit{symmetric group} $S_n$ acts on $X$ by capturing all possible permutations, where each element $\sigma \in S_n$ defines a bijective map $\sigma: X \to X$.

These examples illustrate how group actions provide an algebraic characterization of the symmetry of a point set. In general, common geometric symmetries in Euclidean space include rotational symmetry, reflection symmetry, and translational symmetry, as the isometries of a finite-dimensional normed space consist of rotations, reflections, and translations.

\subsection{Isometry group and symmetry group}
Let $(M, d)$ be a metric space. An \textit{isometry} is a mapping $\pi: M \to M$ that preserves distances, meaning that for all $x, y \in M$, we have
\begin{equation*}
  d(\pi(x), \pi(y)) = d(x, y).
\end{equation*}
The set of all isometries of $M$, equipped with function composition as the group operation, forms a group known as the \textit{isometry group} of $M$, commonly denoted by $\Iso(M)$.

In two-dimensional Euclidean space, the isometry group is called the \textit{Euclidean group}, denoted by $E(2)$. Recall that the only distance-preserving transformations in $\mathbb{R}^2$ are translations, rotations, and reflections. For example, a rotation by an angle $\theta$ around the origin is represented by the matrix $R(\theta) = \begin{bmatrix} \cos\theta & -\sin\theta \\ \sin\theta & \cos\theta \end{bmatrix}$, while a reflection across a line passing through the origin at an angle $\alpha$ relative to the $x$-axis is given by $M(\alpha) = \begin{bmatrix} \cos 2\alpha & \sin 2\alpha \\ \sin 2\alpha & -\cos 2\alpha \end{bmatrix}$.
Consequently, the Euclidean group can be expressed as the semi-direct product of the translation group $\mathbb{R}^2$ and the orthogonal group:
\begin{equation*}
E(2) \cong \mathbb{R}^2 \rtimes O(2).
\end{equation*}

Let $X$ be a geometric object or a subset within the metric space $(M, d)$. The \textit{symmetry group} of $X$ is defined as the subgroup of isometries of $M$ that map $X$ globally onto itself. Formally, the symmetry group of $X$ is given by
\begin{equation*}
  \Sym(X) = \{\pi \in \Iso(M) \mid \pi(X) = X\},
\end{equation*}
where $\pi(X)$ denotes the image of $X$ under the isometry $\pi$.

An element of $\Sym(X)$ is called a (geometric) \textit{symmetry} of the object $X$. These symmetries describe the transformations of the ambient space that preserve the overall geometric structure of $X$. Depending on the nature of $X$ and the ambient space $M$, these may include rotations, reflections, or translations. In particular, if $X$ is a finite set in a Euclidean space, its symmetry group must fix at least one common point (the centroid of $X$). Consequently, such symmetries exclude non-trivial translations and consist entirely of rotations, reflections, or their compositions.

While the isometry group characterizes all possible distance-preserving transformations of the background space, the symmetry group of a specific geometric object isolates those transformations that leave the object invariant. Consequently, analyzing these symmetry groups provides deep insights into the intrinsic geometric and structural properties of the object itself.

\begin{example}
\begin{enumerate}[label=($\roman*$)]
  \item For an infinite line $\ell$ embedded in $\mathbb{R}^2$, its symmetry group consists of translations along the line and reflections. The structure of this symmetry group $G_{\ell}$ can be expressed as
\begin{equation*}
  G_{\ell} \cong \mathbb{R} \rtimes \mathbb{Z}_2,
\end{equation*}
where $\mathbb{R}$ represents the continuous translation group parallel to the line, and $\mathbb{Z}_2$ represents the reflection group generated by a reflection across a hyperplane perpendicular to $\ell$ (or a reflection across the line $\ell$ itself within the ambient space).
  \item Let $X$ be a regular $n$-gon in the Euclidean plane. The symmetry group of $X$ is the dihedral group $D_n$. The \textit{dihedral group} $D_n$ can be defined via the following group presentation:
\begin{equation*}
D_n = \langle \tau, \pi \mid \tau^n = e, \pi^2 = e, \pi\tau = \tau^{-1}\pi \rangle,
\end{equation*}
where $\tau$ represents a counterclockwise rotation by $\frac{2\pi}{n}$ about the center of the polygon, and $\pi$ represents a reflection across one of its axes of symmetry.
  \item The symmetry group of a cube in $\mathbb{R}^{3}$, denoted by $O_h$, consists of all isometries in $\Iso(\mathbb{R}^3)$ that leave the cube invariant. It can be decomposed as
\begin{equation*}
  O_h \cong S_4 \times \mathbb{Z}_2,
\end{equation*}
where $S_4$ is the symmetric group on four elements, corresponding to the permutations of the four space diagonals of the cube, and $\mathbb{Z}_2$ accounts for the central inversion (or reflections).
\end{enumerate}
\end{example}

\subsection{Persistence object}

The rapid development of persistent homology has spurred the creation of various persistent invariants, such as persistent cohomology, persistent $K$-theory, and persistent Steenrod algebra. Fundamentally, these developments motivate the introduction of \textit{persistence objects} in a categorical sense, furnishing widely used tools to characterize the multi-scale and hierarchical structural information of topological spaces and algebraic structures.

We first recall that a \textit{category} $\mathfrak{C}$ consists of a collection of objects $\mathrm{Ob}(\mathfrak{C})$, a collection of morphisms $\mathrm{Hom}_{\mathfrak{C}}(A, B)$ for each pair of objects $A, B \in \mathrm{Ob}(\mathfrak{C})$, along with an associative composition rule and identity morphisms. Any partially ordered set (poset) $(X, \leq)$ can be naturally viewed as a category, where the objects are the elements of $X$, and a unique morphism $x \to y$ exists if and only if $x \leq y$.

\begin{definition}
Let $(X, \leq)$ be a poset viewed as a category, and let $\mathfrak{C}$ be an arbitrary target category. A \textit{persistence object} in $\mathfrak{C}$ indexed by $X$ is a covariant functor 
\begin{equation*}
\mathcal{F}: (X, \leq) \to \mathfrak{C}.
\end{equation*}
Explicitly, this assignment satisfies:
\begin{enumerate}[label=($\roman*$)]
    \item For each element $x \in X$, there is an associated object $\mathcal{F}(x) \in \mathrm{Ob}(\mathfrak{C})$;
    \item For any pair $x, y \in X$ with $x \leq y$, there is a transition morphism $\mathcal{F}(x \leq y): \mathcal{F}(x) \to \mathcal{F}(y)$ in $\mathfrak{C}$ such that $\mathcal{F}(x \leq x) = \mathrm{id}_{\mathcal{F}(x)}$, and for all $x \leq y \leq z$, the composition satisfies $\mathcal{F}(y \leq z) \circ \mathcal{F}(x \leq y) = \mathcal{F}(x \leq z)$.
\end{enumerate}
\end{definition}

If the indexing poset is chosen as $X = \mathbb{R}^n$ or $X = \mathbb{Z}^n$ (equipped with the standard product order), $\mathcal{F}$ is referred to as a \textit{multiparameter persistence object}. In theoretical developments, the poset $(\mathbb{R}, \leq)$ represents the most classical single-parameter setting. For practical applications, the discrete poset $(\mathbb{Z}, \leq)$ is often preferred, as empirical data collected from real-world observations is inherently discrete.

By reversing the direction of the transition arrows, we obtain the dual notion of a persistence object.

\begin{definition}
Let $(X, \leq)$ be a poset, and let $\mathfrak{C}$ be a category. A \textit{copersistence object} in $\mathfrak{C}$ indexed by $X$ is a contravariant functor 
\begin{equation*}
\mathcal{F}: (X, \leq)^{\mathrm{op}} \to \mathfrak{C},
\end{equation*}
where $(X, \leq)^{\mathrm{op}}$ denotes the opposite category of the poset. Equivalently, a copersistence object can be viewed as a \textit{presheaf} on the poset category $(X, \leq)$, where a unique morphism $\mathcal{F}(x \leq y): \mathcal{F}(y) \to \mathcal{F}(x)$ is assigned whenever $x \leq y$.
\end{definition}

The collection of all persistence (resp. copersistence) objects indexed by $X$ in a fixed category $\mathfrak{C}$ forms a functor category, denoted by $\mathfrak{C}^{(X, \leq)}$ (resp. $\mathfrak{C}^{(X, \leq)^{\mathrm{op}}}$), where the morphisms between two persistence objects are given by \textit{natural transformations}.

\begin{example}
Depending on the choice of the target category $\mathfrak{C}$, the abstract definition of persistence objects specializes into several classical structures foundational to TDA:
\begin{enumerate}[label=($\roman*$)]
    \item When $\mathfrak{C}$ is the category $\mathbf{Top}$ of topological spaces and continuous maps, a persistence object $\mathcal{X}: (X, \leq) \to \mathbf{Top}$ defines a persistent topological space. A classic example is the filtration of a simplicial complex indexed by $\mathbb{R}$ or $\mathbb{Z}$, where the transition maps are inclusion mappings.
    \item When $\mathfrak{C}$ is the category $\mathbf{Ch}_R$ of chain complexes of $R$-modules and chain maps, a functor $\mathcal{C}: (X, \leq) \to \mathbf{Ch}_R$ yields a persistent chain complex. This object typically arises by applying the singular or simplicial chain complex functor to a persistent topological space.
    \item When $\mathfrak{C}$ is the category $\mathbf{Vect}_{\mathbb{K}}$ of vector spaces over a field $\mathbb{K}$ and linear transformations, a functor $\mathcal{V}: (X, \leq) \to \mathbf{Vect}_{\mathbb{K}}$ defines a persistent vector space (often called a persistence module). This is obtained by applying the homology functor $H_n(-; \mathbb{K})$ to a persistent chain complex, serving as the central object of study in persistent homology.
\end{enumerate}
\end{example}

\subsection{Persistent homology and topological signatures}

In this subsection, we provide a brief overview of persistent homology, shwoing how geometric and algebraic filtrations capture multi-scale topological features.

Let $X$ be a topological space, and let $f: X \to \mathbb{R}$ be a real-valued function defined on $X$ (often referred to as a filtering or sublevel-set function). By considering the family of sublevel sets, we obtain a \textit{filtration} of $X$, which is a nested sequence of topological subspaces
\begin{equation*}
X_t = f^{-1}((-\infty, t]), \quad \text{for } t \in \mathbb{R}.
\end{equation*}
Clearly, for any $s \leq t$, we have an inclusion mapping $\iota_{s,t}: X_s \hookrightarrow X_t$. This collection of spaces and inclusions precisely forms a persistent topological space, viewed as a covariant functor $\mathcal{X}: (\mathbb{R}, \leq) \to \mathbf{Top}$.

By applying the $n$-th singular or simplicial homology functor $H_n(-; \mathbb{K})$ with coefficients in a field $\mathbb{K}$ to this persistent topological space, we pass from geometry to algebra. This yield a sequence of vector spaces and induced linear transformations:
\begin{equation*}
\mathcal{V}_t = H_n(X_t; \mathbb{K}), \quad \text{and} \quad v_{s,t} = H_n(\iota_{s,t}): H_n(X_s; \mathbb{K}) \to H_n(X_t; \mathbb{K}).
\end{equation*}
The assignment $\mathcal{V}: (\mathbb{R}, \leq) \to \mathbf{Vect}_{\mathbb{K}}$ is called a persistence module (or a persistent vector space). Under the assumption of pointwise finite-dimensionality (i.e., $\dim \mathcal{V}_t < \infty$ for all $t \in \mathbb{R}$), a structure theorem \cite{zomorodian2004computing} asserts that the persistence module $\mathcal{V}$ can be uniquely decomposed into a direct sum of interval modules
\begin{equation*}
\mathcal{V} \cong \bigoplus_{i \in I} \mathbb{K}_{[b_i, d_i)},
\end{equation*}
where $[b_i, d_i) \subset \mathbb{R}$ represents the persistence interval of a specific topological feature. Here, $b_i$ signifies the \textit{birth time} when a new topological generator (such as a connected component, loop, or void) first appears in the filtration, and $d_i$ signifies its \textit{death time} when it becomes a boundary or merges into an older feature.

This unique algebraic decomposition allows us to summarize the complex multi-scale structural data into intuitive and stable topological signatures:
\begin{itemize}
    \item \textit{Barcode:} A multiset of intervals $\{[b_i, d_i)\}_{i \in I}$ visualized as horizontal bars aligned along the parameter axis. Longer bars denote robust geometric structures that survive across a wide range of scales, whereas shorter bars typically correspond to topological noise or local geometric fluctuations.
    \item \textit{Persistence Diagram:} A multiset of points in the extended plane defined as:
    \begin{equation*}
    \mathrm{Dgm}(\mathcal{V}) = \{(b_i, d_i) \in \mathbb{R} \times (\mathbb{R} \cup \{\infty\}) \mid i \in I \text{ and } b_i < d_i\}.
    \end{equation*}
    Points located far above the diagonal $d=b$ represent significant macroscopic properties of the space $X$, while points close to the diagonal capture minor, short-lived features.
\end{itemize}

A further property of persistent homology that justifies its extensive application is its \textit{stability}. Let $f, g: X \to \mathbb{R}$ be two continuous functions on a triangulable space $X$. The Stability Theorem \cite{cohen2005stability} states that the bottleneck distance $d_{B}$ between their corresponding persistence diagrams is bounded by the supremum norm of the functions
\begin{equation*}
d_{B}(\mathrm{Dgm}(\mathcal{V}_f), \mathrm{Dgm}(\mathcal{V}_g)) \leq \|f - g\|_{\infty}.
\end{equation*}
This structural stability ensures that when the underlying data or continuous objects undergo minor deformations or group perturbations, their topological signatures remain robust. This invariant approach under continuous parameters sets a natural foundation for investigating higher-order continuous symmetries.

%% file: topological_dynamics.tex
\section{The dynamics of symmetries}\label{section:dynamics_symmetries}

In a topological space $M$, the configuration space is defined as the set of all configurations of $n$ distinct points, where each configuration represents a specific spatial arrangement or state of $n$ particles in $M$ \cite{dyer1962homology,fadell1962configuration}. The dynamics of configurations has long been a foundational paradigm in both mathematics and physics, playing a central role in fields such as classical mechanics, robotics, and topological kinematics \cite{farber2009topological,ghrist2007configuration,sussman2015structure}. In this work, we focus on the symmetries of finite configurations and investigate the time-dependent dynamics of such symmetric behaviors.

\subsection{Persistent symmetries}\label{section:persistent_set}

Let $n \geq 2$ be an integer, and let an $n$-point configuration in a Euclidean space $E$ be a finite subset consisting of $n$ distinct points. A \textit{dynamical system} of $n$-point configurations is a parameterized family $\{X_t\}_{t \in \mathbb{R}}$, where each $X_t \subset E$ comprises $n$ distinct points, equipped with a collection of ambient homeomorphisms $\{f_{a,b} : E \to E\}_{a \leq b}$ satisfying the following conditions:
\begin{itemize}
  \item For any $a \leq b$, the restriction of $f_{a,b}$ to $X_a$ induces a bijection
  \[
  {f_{a,b}}|_{X_a} : X_a \to X_b.
  \]
  For brevity, we write $f_{a,b}\big|_{X_a}$ simply as $f_{a,b}$.
  \item For any $a \leq b \leq c$, the classical cocycle relations hold: $f_{a,c} = f_{b,c} \circ f_{a,b}$ and $f_{a,a}=\mathrm{id}_E$.
\end{itemize}

Intuitively, such a dynamical system describes a time-evolving collection of $n$-point configurations whose topological structures are preserved under non-decreasing parameter shifts. This model captures the geometric evolution of configurations over continuous time, giving a natural setting to study changing symmetry profiles.

We denote the configuration at time $t$ by
\[
X_t = \{x_1(t), x_2(t), \dots, x_n(t)\}.
\]
For any real numbers $a \leq b$, the restricted bijection $f_{a,b} : X_a \to X_b$ is defined by setting
\[
f_{a,b}(x_k(a)) = x_k(b)
\]
for each index $k = 1, 2, \dots, n$. This map naturally induces a group isomorphism between the permutation automorphism groups of the point sets:
\begin{equation*}
  f_{a,b}^{\sharp} : \Aut(X_a) \to \Aut(X_b),
\end{equation*}
defined via conjugation:
\begin{equation*}
  f_{a,b}^{\sharp}(g) =  f_{a,b}\circ g\circ f_{a,b}^{-1}
\end{equation*}
for all $g \in \Aut(X_a)$.

Furthermore, this isomorphism endows the map $f_{a,b} : X_a \to X_b$ with a natural $S_n$-equivariant structure. That is, for every $g \in \Aut(X_a) \cong S_n$, the following diagram commutes:
\[
\begin{tikzcd}
X_a \arrow[r, "f_{a,b}"] \arrow[d, "g"'] & X_b \arrow[d, "f_{a,b}^\sharp(g)"] \\
X_a \arrow[r, "f_{a,b}"]                 & X_b
\end{tikzcd}
\]
For convenience, we identify both $\Aut(X_a)$ and $\Aut(X_b)$ with the symmetric group $S_n$ via the natural indexing of the points. Under this identification, the $S_n$-equivariance of $f_{a,b}$ can be expressed explicitly as
\begin{equation*}
  f_{a,b}(g \cdot x_k(a)) = g \cdot f_{a,b}(x_k(a))
\end{equation*}
for any $g \in S_n$ and for all $k = 1, 2, \dots, n$, where the action of $g$ on $x_k(t)$ corresponds to a permutation of the indices.

Crucially, at a given time $t$, an element of the automorphism group $\Aut(X_t)$ does not automatically correspond to a geometric symmetry; that is, permutations in $\Aut(X_t)$ cannot always be extended to global isometries of the ambient space $E$. Let $\Gamma_t = \Sym(X_t) \subset \Iso(E)$ denote the set of geometric symmetries of the configuration $X_t$. Since $X_t$ is finite, $\Gamma_t$ forms a finite subgroup of $\Iso(E)$. When restricted to $X_t$, each symmetry in $\Gamma_t$ induces a distinct permutation belonging to $\Aut(X_t)$. Thus, $\Gamma_t$ can be viewed as a distinguished subgroup of $\Aut(X_t)$ containing those permutations that arise from ambient isometries.

Let $\mathrm{Homeo}(E)$ denote the group of homeomorphisms on $E$. For any $t\in \mathbb{R}$, we define
\[
\mathrm{Homeo}(E; X_t) = \{ u \in \mathrm{Homeo}(E) \mid u(X_t) = X_t \},
\]
which represents the stabilizer subgroup of homeomorphisms fixing $X_t$ globally. Consider the induced map
\[
f_{a,b}^{\sharp} : \mathrm{Homeo}(E; X_a) \to \mathrm{Homeo}(E; X_b), \quad u \mapsto f_{a,b} \circ u \circ f_{a,b}^{-1},
\]
which acts by conjugation. While $f_{a,b}^{\sharp}$ maps $\mathrm{Homeo}(E; X_a)$ to $\mathrm{Homeo}(E; X_b)$ successfully, it does not necessarily map the geometric symmetry subgroup $\Gamma_a$ into $\Gamma_b$.

\begin{example}\label{example:not_persistence}
Consider the point set $X_a = \{x_1(a), x_2(a), x_3(a)\} \subset \mathbb{R}^2$, where
\[
x_1(a) = (0, 0), \quad
x_2(a) = (1, 0), \quad
x_3(a) = \left(\frac{1}{2}, \frac{\sqrt{3}}{2}\right).
\]
These points form an equilateral triangle, which exhibits six geometric symmetries; hence, $\Gamma_a \cong D_3$.

Now consider the configuration at a later time $b$, given by $X_b = \{x_1(b), x_2(b), x_3(b)\} \subset \mathbb{R}^2$, where
\[
x_1(b) = (0, 0), \quad
x_2(b) = (1, 0), \quad
x_3(b) = (0.6, 0.7).
\]
These points form an irregular triangle admitting no non-trivial symmetries, yielding $\Gamma_b = \{\mathrm{id}\}$.

Let $f_{a,b} : \mathbb{R}^2 \to \mathbb{R}^2$ be the unique affine (linear) map determined by the condition $f_{a,b}(x_k(a)) = x_k(b)$ for $k = 1, 2, 3$. The induced conjugation map $f_{a,b}^\sharp$ fails to preserve the geometric symmetries; specifically, $f_{a,b}^\sharp(\Gamma_a) \not\subset \Gamma_b$.
\end{example}

\begin{definition}
Let $\{X_t\}_{t \in \mathbb{R}}$ be a dynamical system of configurations. 
\begin{itemize}
    \item A symmetry $\pi \in \Gamma_{t_0}$ is said to be \textit{born at time $t_0$} if there exists a sufficiently small $\epsilon > 0$ such that for any $t \in (t_0 - \epsilon, t_0)$, there is no symmetry $\pi_t \in \Gamma_t$ satisfying $f_{t,t_0}^\sharp(\pi_t) = \pi$.
    \item A symmetry $\pi \in \Gamma_{t_0}$ is said to \textit{die at time $t_1$} (where $t_0 < t_1$) if $f_{t_0, t}^\sharp(\pi) \in \Gamma_t$ for all $t \in [t_0, t_1)$, but $f_{t_0, t_1}^\sharp(\pi) \notin \Gamma_{t_1}$.
    \item For any real numbers $a \leq b$, a symmetry $\pi \in \Gamma_a$ \textit{persists from time $a$ to time $b$} if, for each $t \in [a, b]$, the tracked configuration symmetry satisfies $f_{a,t}^\sharp(\pi) \in \Gamma_t$.
\end{itemize}
\end{definition}

By definition, a symmetry $\pi \in \Gamma_a$ persists from time $a$ to time $b$ if and only if $f_{a,t}^\sharp(\pi) : E \to E$ remains an ambient isometry for all $t \in [a, b]$. In particular, the identity symmetry trivially persists over any interval $[a, b]$.

\begin{proposition}\label{proposition:persistence}
Let $\{X_t\}_{t \in \mathbb{R}}$ be a dynamical system of finite point sets in a Euclidean space $E$. If for all $t \in [a,b]$, the map $f_{a,t}: E \to E$ is an isometry, then any geometric symmetry $\pi \in \Gamma_a$ persists from time $a$ to time $b$.
\end{proposition}

\begin{proof}
Let $\pi \in \Gamma_a$, meaning that $\pi: E \to E$ is an ambient isometry satisfying $\pi(X_a) = X_a$. Since $f_{a,t}$ is an isometry for each $t \in [a, b]$, its inverse $f_{a,t}^{-1}$ is also an isometry. Consequently, the conjugated mapping
\[
f_{a,t}^{\sharp}(\pi) = f_{a,t} \circ \pi \circ f_{a,t}^{-1}
\]
is a composition of isometries, and thus constitutes an isometry of $E$.

Furthermore, for each $x_k(t) \in X_t$, there exists a unique element $x_k(a) \in X_a$ such that $x_k(t) = f_{a,t}(x_k(a))$. Applying the conjugated operator yields:
\[
f_{a,t}^\sharp(\pi)(x_k(t)) = (f_{a,t} \circ \pi \circ f_{a,t}^{-1})(x_k(t)) = f_{a,t}(\pi(x_k(a))).
\]
Since $\pi \in \Gamma_a$, we have $\pi(X_a) = X_a$, which implies $\pi(x_k(a)) \in X_a$. It follows that $f_{a,t}(\pi(x_k(a))) \in f_{a,t}(X_a) = X_t$. Therefore, the ambient isometry $f_{a,t}^\sharp(\pi)$ maps $X_t$ onto itself, meaning $f_{a,t}^\sharp(\pi) \in \Gamma_t$. We conclude that the symmetry $\pi$ persists throughout the interval $[a,b]$.
\end{proof}

Proposition \ref{proposition:persistence} demonstrates that pure isometric flows preserve the persistence of all structural symmetries. In general dynamical regimes, non-isometric flows distort configurations and break symmetries, as illustrated by Example~\ref{example:not_persistence}. Conversely, the persistence of a non-trivial symmetry $\pi$ from time $a$ to time $b$ implies that the dynamical flow $f_{a,t}$ structurally preserves certain local rigid parameters. Mechanistically, a higher count of persistent symmetries implies that the evolving dynamical system maintains a greater degree of rigidity throughout its temporal evolution.

\subsection{Dynamics of isometric permutations}

In the previous section, we established how geometric symmetries persist within a continuous, ambient topological setting. However, in practical applications and computational topology, combinatorial structures are often significantly easier to describe, track, and compute. To bridge this gap, we develop a purely combinatorial and intrinsic formulation of persistent symmetries.

A \textit{flexible dynamical system} of $n$-point configurations consists of a family of finite configurations $\{X_t\}_{t \in \mathbb{R}}$ indexed by time $t \in \mathbb{R}$, together with a collection of set-theoretic bijections $\{f_{a,b} : X_a \to X_b\}_{a \leq b}$ satisfying the following cocycle conditions: for all $a \leq b \leq c$,
\[
f_{a,c} = f_{b,c} \circ f_{a,b}, \quad \text{and} \quad f_{a,a} = \mathrm{id}_{X_a}.
\]

Unlike the rigid ambient formulation in the preceding section, the structure maps $f_{a,b}$ in a flexible dynamical system are strictly required to be bijections between finite point sets, without any assumption regarding their extendability to homeomorphisms of the ambient Euclidean space $E$. This formulation fundamentally decouples the system from ambient topological constraints, rendering the description and algorithmic construction of configuration trajectories substantially more flexible.

Readily, the following classical geometric result guarantees that if a distance-preserving discrete map can be linearly extended to the entire space, it inherits global rigidity:

\begin{proposition}
Let $f: E \to E$ be an affine transformation of a Euclidean space $E$, and let $X \subseteq E$ be a finite set that spans $E$. If the restricted mapping $f|_X: X \to f(X)$ preserves all pairwise distances in $X$, then the global affine transformation $f$ is an isometry on $E$.
\end{proposition}

For a finite configuration $X$ within a Euclidean space, we define its intrinsic isometry group as:
\[
\Iso(X) = \{ \pi \in \Aut(X) \mid d(\pi(x), \pi(y)) = d(x, y) \text{ for all } x, y \in X \},
\]
which acts as a permutation group on the finite set. The elements of $\Iso(X)$ are referred to as \textit{isometric permutations}. In this combinatorial regime, the dynamic invariants we track are precisely these isometric permutations, which encode the intrinsic structural symmetries of $X$ independently of ambient maps.

For a flexible dynamical system $\{X_t\}_{t \in \mathbb{R}}$, we set $\tilde{\Gamma}_t = \Iso(X_t) \subset \Aut(X_t)$. The induced group isomorphism on total permutations,
\begin{equation*}
  f_{a,b}^{\sharp} : \Aut(X_a) \to \Aut(X_b), \quad g \mapsto f_{a,b} \circ g \circ f_{a,b}^{-1},
\end{equation*}
does not, in general, map the isometric subgroup $\tilde{\Gamma}_a$ into $\tilde{\Gamma}_b$. This discrepancy motivates a combinatorial definition for the persistence of symmetries under flexible flows.

\begin{definition}
Let $\{X_t\}_{t \in \mathbb{R}}$ be a flexible dynamical system with structure maps $\{f_{a,b} : X_a \to X_b\}_{a \leq b}$, and let $\tilde{\Gamma}_t = \Iso(X_t)$.
\begin{itemize}
    \item An isometric permutation $\pi \in \tilde{\Gamma}_{t_0}$ is said to be \textit{born at time $t_0$} if there exists an $\epsilon > 0$ such that for any $t \in (t_0 - \epsilon, t_0)$, there is no isometric permutation $\pi_t \in \tilde{\Gamma}_t$ satisfying $f_{t,t_0}^\sharp(\pi_t) = \pi$.
    \item An isometric permutation $\pi \in \tilde{\Gamma}_{t_0}$ is said to \textit{die at time $t_1$} (where $t_0 < t_1$) if its tracked image satisfies $f_{t_0, t}^\sharp(\pi) \in \tilde{\Gamma}_t$ for all $t \in [t_0, t_1)$, but $f_{t_0, t_1}^\sharp(\pi) \notin \tilde{\Gamma}_{t_1}$.
    \item For any real numbers $a \leq b$, an isometric permutation $\pi \in \tilde{\Gamma}_a$ \textit{persists from time $a$ to time $b$} if, for each $t \in [a, b]$, the dynamically propagated permutation satisfies $f_{a,t}^\sharp(\pi) \in \tilde{\Gamma}_t$.
\end{itemize}
\end{definition}

It is worth noting that the definitions and algebraic persistence mechanisms formulated for flexible dynamical systems closely parallel those established for ambient dynamical systems in Section \ref{section:persistent_set}. The key distinction lies in their methodological emphasis: the results for flexible dynamical systems are combinatorially oriented and computationally tractable via algorithmic matrix representations, whereas the ambient counterparts align more naturally with continuous geometric intuition and classical kinematic flows.

\subsection{Persistence module of symmetries}\label{section:module}

Let $\{X_t\}_{t \in \mathbb{Z}}$ be a discrete-time dynamical system with structure maps $\{f_{a,b}\}_{a \leq b}$, and let $\mathbb{K}$ be a fixed ground field. For each $t \in \mathbb{Z}$, we denote by $M_t$ the $\mathbb{K}$-linear vector space spanned by the geometric symmetries in $\Gamma_t$. Clearly, $M_t$ is a finite-dimensional vector space. Since $\Gamma_t$ invariably contains the identity transformation, $M_t$ is guaranteed to be non-trivial.

To aggregate these local spaces into a unified algebraic object, we define a collection of forward transition linear maps $u_t: M_t \to M_{t+1}$ by specifying their action on the basis elements $\pi \in \Gamma_t$ as follows:
\begin{equation*}
u_t(\pi) = \begin{cases} 
f_{t,t+1}^\sharp(\pi), & \text{if } f_{t,t+1}^\sharp(\pi) \in \Gamma_{t+1}, \\ 
0, & \text{otherwise.} 
\end{cases}
\end{equation*}
Extending these maps linearly yields a global degree-$1$ endomorphism $u: \mathbf{M} \to \mathbf{M}$ on the graded vector space $\mathbf{M} = \bigoplus_{t \in \mathbb{Z}} M_t$. This endomorphism naturally endows $\mathbf{M}$ with a graded $\mathbb{K}[u]$-module structure, where the polynomial action is defined via the module map:
\begin{equation*}
  \mathbb{K}[u] \times \mathbf{M} \to \mathbf{M}, \quad (p(u), x) \mapsto p(u)(x).
\end{equation*}
Explicitly, for a polynomial $p(u) = \sum_{i=0}^{m} a_i u^i$ with $a_i \in \mathbb{K}$, the action on any chain of symmetries $x \in \mathbf{M}$ is given by $p(u)(x) = \sum_{i=0}^{m} a_i u^i(x)$, where $u^i(x)$ denotes the $i$-fold iterated application of the shift operator $u$. 

To formalize a structural decomposition for persistent symmetries and draw a parallel to the barcode representations in persistent homology, we analyze the algebraic evolution of these symmetry generators over time. Much like persistent homology monitors the multi-scale lifespans of topological features, our persistent symmetry module tracks the formation, structural stability, and eventual dissipation of geometric symmetries across temporal scales.

\begin{theorem}\label{theorem:decompostion}
Let $\{X_t\}_{t \in \mathbb{Z}}$ be a discrete dynamical system, and let $\mathbf{M} = \bigoplus_{t \in \mathbb{Z}} M_t$ be its associated persistent symmetry module. Assuming that $\mathbf{M}$ is a finitely generated $\mathbb{K}[u]$-module, there exists a unique canonical decomposition:
\begin{equation}\label{equ:decomposition}
  \mathbf{M} \cong \left( \bigoplus_{i} \mathbb{K}[u](-b_i) \right) \oplus \left( \bigoplus_{j} \mathbb{K}[u]/(u^{d_j})(-c_j) \right),
\end{equation}
where $(-b_i)$ and $(-c_j)$ denote the degree shifts shifting the grading to the respective birth times of the symmetry generators.
\end{theorem}

\begin{proof}
Since the polynomial ring $\mathbb{K}[u]$ over a field is a principal ideal domain (PID), the result follows immediately from the classical structure theorem for finitely generated modules over a PID.
\end{proof}

In the algebraic decomposition \eqref{equ:decomposition}, the free components correspond to symmetry generators that are born at time $t = b_i$. Because their underlying $\mathbb{K}[u]$-submodules are free, these symmetry signatures do not terminate and persist indefinitely into the future as the system evolves. Conversely, the torsion components correspond to transient symmetries born at time $t = c_j$. These elements exhibit a finite lifespan determined by the relation $u^{d_j} \cdot \beta_j = 0$, meaning that the propagated symmetry vanishes entirely at the subsequent time step $t = c_j + d_j$. 

Consequently, Theorem \ref{theorem:decompostion} effectively classifies the structural dynamics of the system by partitioning symmetries into precise mathematical lifetimes. Symmetries associated with free generators capture the permanent geometric rigidity of the dynamical system, while torsion generators reflect localized, short-lived configurations that emerge and dissolve over specific intervals.

This algebraic decompostion gives rise to persistence barcodes that visually trace the lifespan of each individual symmetry. A \textit{persistence barcode} in this context is a multiset of intervals in $\mathbb{Z} \times (\mathbb{Z} \cup \{\infty\})$. Each horizontal bar represents a distinct symmetry line; its start point signifies the birth of the symmetry, and its length corresponds to the duration over which the symmetry remains stable under the dynamical flow. Symmetries represented by significantly longer intervals indicate robust, macroscopically stable geometric properties of the system, whereas short intervals isolate highly sensitive or unstable symmetries.

\begin{remark}
In scenarios where the temporal parameter $t$ is continuously indexed over $\mathbb{R}$, the direct module decomposition in Theorem \ref{theorem:decompostion} cannot be implemented due to the uncountability of the grading. Nevertheless, one can approximate the continuous system $\{X_t\}_{t \in \mathbb{R}}$ by sampling a discrete, countable sequence of time coordinates $\{t_r\}_{r \in \mathbb{Z}}$. Defining the discretized configuration as $Y_r = X_{t_r}$ yields a standard discrete dynamical system $\{Y_r\}_{r \in \mathbb{Z}}$. The decomposition formula \eqref{equ:decomposition} can then be algorithmically deployed to compute the discretized persistent symmetries, providing an arbitrarily precise approximation of the continuous system.
\end{remark}

\begin{example}\label{example:dynamics_symmmetry}
Consider a dynamical system $\{X_t\}_{t \in \mathbb{R}}$ consisting of four points evolving periodically in $\mathbb{R}^2$, defined by the parametric equations:
\begin{equation*}
\begin{aligned}
    & x_{1}(t)=(2-\cos t,\sin t), && x_{2}(t)=(-\sin t,2-\cos t), \\
    & x_{3}(t)=(-2+\cos t,-\sin t), && x_{4}(t)=(\sin t,-2+\cos t).
\end{aligned}
\end{equation*}
This system exhibits a high degree of geometric symmetry. Specifically, the configuration $X_t$ is invariant under a rotation $\tau$ of $90^{\circ}$ counterclockwise about the origin, where $\tau$ acts on $\mathbb{R}^2$ by
\begin{equation*}
  \tau \cdot (x,y) = (-y,x).
\end{equation*}
Consequently, the rotational symmetry group of $X_t$ contains the cyclic subgroup $\{\mathrm{id}, \tau, \tau^{2}, \tau^{3}\} \cong C_4$ for all $t \in \mathbb{R}$.

\begin{figure}[h]
\centering
\definecolor{axisgray}{gray}{0.6}
\definecolor{symred}{RGB}{220, 80, 80}
\definecolor{arcangle}{RGB}{100, 160, 100}
\definecolor{trackcolor}{RGB}{0, 128, 128}

\begin{tikzpicture}[scale=0.8, >=Stealth]

    \begin{scope}
        \draw[->, color=axisgray, line width=0.4pt] (-3.5,0) -- (3.5,0) node[right, font=\small, text=black] {$x$};
        \draw[->, color=axisgray, line width=0.4pt] (0,-3.5) -- (0,3.5) node[right, font=\small, text=black] {$y$};

        \foreach \pos in {(2,0), (-2,0), (0,2), (0,-2)}
            \draw[color=trackcolor, line width=0.7pt] \pos circle(1);

        \draw[dash pattern=on 4pt off 2pt, color=symred, line width=0.6pt] (-2.2, 0) -- (2.2, 0); 
        \draw[dash pattern=on 4pt off 2pt, color=symred, line width=0.6pt] (0, -2.2) -- (0, 2.2); 
        \draw[dash pattern=on 4pt off 2pt, color=symred, line width=0.6pt] (-1.6, -1.6) -- (1.6, 1.6); 
        \draw[dash pattern=on 4pt off 2pt, color=symred, line width=0.6pt] (-1.6, 1.6) -- (1.6, -1.6); 

        \draw[line width=1.2pt, color=mainblue] (1,0) -- (0,1) -- (-1,0) -- (0,-1) -- cycle;
        
        \foreach \pos in {(1,0), (0,1), (-1,0), (0,-1)}
            \fill[color=mainblue] \pos circle (1.2pt);

        \draw[->, color=arcangle, line width=0.8pt] (0:0.7) arc[start angle=0, end angle=90, radius=0.7];
        \node[color=black, font=\footnotesize] at (35:1.0) {$90^\circ$};
    \end{scope}

    \begin{scope}[xshift=9cm]
        \draw[->, color=axisgray, line width=0.4pt] (-3.5,0) -- (3.5,0) node[right, font=\small, text=black] {$x$};
        \draw[->, color=axisgray, line width=0.4pt] (0,-3.5) -- (0,3.5) node[right, font=\small, text=black] {$y$};

        \draw[color=trackcolor, line width=0.7pt] (2,0) circle(0.8);
        \foreach \pos in {(-2,0), (0,2), (0,-2)}
            \draw[color=trackcolor, line width=0.7pt] \pos circle(1);

        \draw[dash pattern=on 4pt off 2pt, color=symred, line width=0.6pt] (-2.2, 0) -- (2.2, 0); 

        \draw[line width=1.2pt, color=mainblue] (1.2,0) -- (0,1) -- (-1,0) -- (0,-1) -- cycle;
        
        \fill[color=mainblue] (1.2,0) circle (1.2pt);
        \foreach \pos in {(0,1), (-1,0), (0,-1)}
            \fill[color=mainblue] \pos circle (1.2pt);
    \end{scope}

\end{tikzpicture}
\caption{The dynamical systems $\{X_t\}$ and $\{Y_t\}$ at $t=0$. The system $X_t$ exhibits full $D_4$ symmetry, whereas the perturbed system $Y_t$ (with $\varepsilon=0.2$) retains only a single axis of reflection and lacks rotational symmetry.}
\label{figure:dynamics_symmetry}
\end{figure}

In addition to rotational symmetry, the system possesses time-dependent reflectional symmetries. The axes of reflection at time $t$ are defined by the linear equations:
\begin{equation*}
  \begin{aligned}
      & (\sin t) x - (2 - \cos t) y  = 0, && (2 - \cos t) x + (\sin t) y  = 0, \\
      & (2-\sin t-\cos t)x + (2+\sin t-\cos t)y=0, && (2+\sin t-\cos t)x - (2-\sin t-\cos t)y=0.
  \end{aligned}
\end{equation*}
Denoting these reflections by $\pi_{1}(t), \dots, \pi_{4}(t)$, the symmetries are periodic with period $2\pi$; that is, $\pi_{i}(t) = \pi_{i}(t + 2\pi)$ for $i = 1,\dots,4$.

To illustrate the sensitivity of symmetry persistence, consider a perturbed dynamical system $\{Y_t\}_{t \in \mathbb{R}}$ given by
\begin{equation*}
\begin{aligned}
    & y_{1}(t) = (2-(1-\varepsilon)\cos t, (1-\varepsilon)\sin t), && y_{2}(t) = (-\sin t, 2-\cos t), \\
    & y_{3}(t) = (-2+\cos t, -\sin t), && y_{4}(t) = (\sin t, -2+\cos t),
\end{aligned}
\end{equation*}
where $0 < \varepsilon < 1$. Unlike $X_t$, the system $Y_t$ possesses no rotational symmetry, and its reflectional symmetry exists only at discrete times $t = m\pi$ ($m \in \mathbb{Z}$). At these instants, the sole axis of reflection is the $x$-axis ($y=0$).

Observe that the trajectory of $Y_t$ lies within an $\varepsilon$-neighborhood of $X_t$, satisfying
\begin{equation*}
  \|x_{k} - y_{k}\| \leq \varepsilon, \quad \text{for all } k = 1, 2, 3, 4.
\end{equation*}
This implies that for an arbitrarily small perturbation magnitude $\varepsilon$, the symmetry persistence modules of $\{X_t\}$ and $\{Y_t\}$ remain fundamentally distinct. Specifically, the persistent symmetries of $X_t$ capture a global rigidity that is structurally fragile under perturbation, demonstrating a significant discrepancy between the two systems despite their close point-wise proximity.
\end{example}

%% file: categorification_symmetries.tex
\section{Categorification of persistent symmetries}\label{section:categorification}

In Section~\ref{section:dynamics_symmetries}, we presented the basic idea of the persistence of symmetries in finite configurations arising from dynamical systems. In this section, we develop a categorical language for persistent symmetry, and study the associated persistent groups, persistent modules, and their barcode representations.

\subsection{Persistence $n$-configuration}\label{section:persistence_configuration}

\begin{definition}
Let $M$ be a topological space. The category $\mathcal{S}_n(M)$ of $n$-configurations is defined as follows.
\begin{itemize}
  \item Objects are subsets $X \subset M$ consisting of $n$ distinct points.
  \item A morphism $f: X \to Y$ is the restriction $f = \tilde{f}|_X$ of a homeomorphism $\tilde{f}: M \to M$ such that $\tilde{f}(X) = Y$ and $f: X \to Y$ is a bijection.
\end{itemize}
Morphisms compose via composition of the underlying homeomorphisms.
\end{definition}

An element in $\mathcal{S}_{n}(M)$ is called an \textit{$n$-point configuration}, or simply an $n$-configuration if there is no ambiguity.

\begin{definition}
A \textit{persistence $n$-configuration} is a functor
\begin{equation*}
\mathcal{F}: (\mathbb{R}, \leq) \to \mathcal{S}_{n}(M).
\end{equation*}
Here, $(\mathbb{R}, \leq)$ is the category of real numbers, where objects are real numbers, and $a \to b$ is a morphism if $a \leq b$ for $a, b \in \mathbb{R}$.
\end{definition}

If $M$ is a Euclidean space, then a persistent $n$-configuration coincides precisely with the notion of a dynamical system of $n$-point configurations introduced in Section~\ref{section:dynamics_symmetries}. Similarly, we may define the flexible category $\widetilde{\mathcal{S}}_n(M)$ of $n$-configurations, where the objects are $n$-configurations in $M$ and the morphisms are bijections between configurations. In this setting, a persistence $n$-configuration valued in $\widetilde{\mathcal{S}}_n(M)$ corresponds to a flexible dynamical system of $n$-point configurations.

In what follows, our focus will be on the study of persistent symmetries of persistence $n$-configurations. To avoid redundancy, we note that the terminology and results developed for persistence $n$-configurations apply parallelly to the flexible case. In certain examples and computations, we will more frequently employ flexible persistence $n$-configurations, as they are often easier to construct. Indeed, by the Whitney extension theorem or the isotopy extension theorem, any bijection between finite subsets of Euclidean space can be extended to a homeomorphism or even a diffeomorphism of the ambient Euclidean space \cite{fefferman2005sharp,hirsch2012differential,whitney1992analytic}. Unless otherwise specified, we will not distinguish between the two settings in the sequel.

From now on, $(M, d)$ is assumed to be a metric space. Given a morphism $f: X \to Y$ in the category $\mathcal{S}_n(M)$, where $X$ and $Y$ are $n$-configurations in the metric space $M$, we obtain a natural group homomorphism between the automorphism groups
\begin{equation*}
  f^{\sharp}: \mathrm{Aut}(X) \to \mathrm{Aut}(Y), \quad \sigma \mapsto f \circ \sigma \circ f^{-1}.
\end{equation*}
However, in general, this does not induce a group homomorphism
\begin{equation*}
  f^\sharp: \Sym(X) \to \Sym(Y).
\end{equation*}
To address this issue, we define the restricted symmetry group
\begin{equation*}
  \Sym_f(X) = \left\{ \sigma \in \Sym(X) \mid f \circ \sigma \circ f^{-1} \in \Sym(Y) \right\}.
\end{equation*}
This is a subgroup of $ \Sym(X) $ because for any $ \sigma, \tau \in \Sym_f(X) $, we have
\begin{equation*}
  (f \circ \sigma \circ f^{-1}) \circ (f \circ \tau \circ f^{-1}) = f \circ (\sigma \circ \tau) \circ f^{-1} \in \Sym(Y),
\end{equation*}
which implies $ \sigma \circ \tau \in \Sym_f(X) $.
The following proposition provides a group homomorphism.

\begin{proposition}\label{proposition:homomorphism}
Let $f: X \to Y$ be a morphism in $\mathcal{S}_n(M)$, where $X$ and $Y$ are $n$-configurations in a metric space $M$. Then the map
\begin{equation*}
f^{\sharp}: \Sym_f(X) \to \Sym(Y), \quad f^\sharp(\sigma)(y) = f(\sigma(f^{-1}(y)))
\end{equation*}
is a well-defined group homomorphism. Here, $\sigma \in \Sym(X)$ and $y \in Y$.
\end{proposition}
\begin{proof}
We first verify that the map $f^\sharp: \Sym_f(X) \to \Sym(Y)$ is well-defined. By the definition of $\Sym_f(X)$, for any $\sigma \in \Sym_f(X)$, the composition $f \circ \sigma \circ f^{-1}$ is an isometry of $Y$. Thus, we have $f^\sharp(\sigma) \in \Sym(Y)$.

Next, we show that $f^\sharp$ is a group homomorphism. For any $\sigma, \tau \in \Sym_f(X)$ and $y \in Y$, we have
\begin{equation*}
f^\sharp(\tau \circ \sigma)(y) = f\left((\tau \circ \sigma)(f^{-1}(y))\right) = f\left(\tau(\sigma(f^{-1}(y)))\right).
\end{equation*}
On the other hand, note that
\begin{equation*}
\left(f^\sharp(\tau) \circ f^\sharp(\sigma)\right)(y) = f^\sharp(\tau)\left(f^\sharp(\sigma)(y)\right) = f\left(\tau(f^{-1}(f(\sigma(f^{-1}(y)))))\right) = f\left(\tau(\sigma(f^{-1}(y)))\right).
\end{equation*}
It follows that
\begin{equation*}
f^\sharp(\tau \circ \sigma)(y) = \left(f^\sharp(\tau) \circ f^\sharp(\sigma)\right)(y),
\end{equation*}
which shows that $f^\sharp(\tau \circ \sigma) = f^\sharp(\tau) \circ f^\sharp(\sigma)$.

It remains to verify that $f^\sharp$ preserves the identity element. Let $e$ be the identity permutation in $\Sym_f(X)$. Then for any $y \in Y$, one has
\begin{equation*}
f^\sharp(e)(y) = f(e(f^{-1}(y))) = f(f^{-1}(y)) = y,
\end{equation*}
so $f^\sharp(e)$ is the identity in $\Sym(Y)$.

Hence, $f^\sharp$ is a well-defined group homomorphism.
\end{proof}

In Proposition~\ref{proposition:homomorphism}, observe that $f$ is a bijection between finite sets and that each $\sigma \in \Sym_f(X)$ is a permutation, hence both $f$ and $\sigma$ are invertible. It follows that the group homomorphism $f^\sharp: \Sym_f(X) \to \Sym(Y)$ is necessarily injective.

\begin{example}\label{example:symmetry_group}
We construct an example where the induced map
\[
f^\sharp: \Sym_f(X) \to \Sym(Y)
\]
is injective but not surjective.

Let $X = \{a, b, c, d\} \subset \mathbb{R}^2$ be the set of vertices of a square, where
\[
a = (0,0), \quad b = (1,0), \quad c = (1,1), \quad d = (0,1).
\]
Then we have $\Sym(X) = D_4$, the dihedral group of the square, which has 8 elements: 4 rotations (including the identity) and 4 reflections.
\begin{figure}[h]
\centering
\begin{tikzpicture}[scale=0.8]
\begin{scope}
\draw[line width=1.2pt, mainblue] (0,0) -- (2,0) -- (2,2) -- (0,2) -- cycle;
\node[below left] at (0,0) {$a$}; \node[below right] at (2,0) {$b$};
\node[above right] at (2,2) {$c$}; \node[above left] at (0,2) {$d$};
\draw[fill=mainblue] (0,0) circle (2pt);
\draw[fill=mainblue] (2,0) circle (2pt);
\draw[fill=mainblue] (2,2) circle (2pt);
\draw[fill=mainblue] (0,2) circle (2pt);
\node at (1,-1) {$X$};
\end{scope}

\begin{scope}[xshift=6cm]
\draw[line width=1.2pt, mainblue] (2,0) -- (0,0) -- (1,1.73) -- (3,1.73) -- cycle;
\node[below] at (2,0) {$a'$}; \node[below] at (0,0) {$b'$};
\node[above] at (1,1.73) {$d'$}; \node[above] at (3,1.73) {$c'$};
\draw[fill=mainblue] (2,0) circle (2pt);
\draw[fill=mainblue] (0,0) circle (2pt);
\draw[fill=mainblue] (1,1.73) circle (2pt);
\draw[fill=mainblue] (3,1.73) circle (2pt);
\node at (1,-1) {$Y$};
\end{scope}
\end{tikzpicture}
\caption{Transformation of a square $X$ into a rhombus $Y$ in Example \ref{example:symmetry_group}.}
\end{figure}
Now define a new set $Y = \{a', b', c', d'\} \subset \mathbb{R}^2$ by
\[
a' = (1,0), \quad b' = (0,0), \quad c' = \left(\frac{3}{2},\frac{\sqrt{3}}{2}\right), \quad d' = \left(\frac{1}{2},\frac{\sqrt{3}}{2}\right),
\]
so that $Y$ forms a rhombus. This rhombus admits three nontrivial symmetries: a reflection across the line $y = \frac{\sqrt{3}}{3}x$, a reflection across the line $y = -\sqrt{3}(x - 1)$, and a $180^\circ$ rotation about its center. Therefore, its full symmetry group is
\[
\Sym(Y) \cong \mathbb{Z}_2 \times \mathbb{Z}_2.
\]
Define a bijection $f: X \to Y$ by
\[
f(a) = a', \quad f(b) = b', \quad f(c) = c', \quad f(d) = d'.
\]
This map $f$ is a bijection of finite sets but is not an isometry. We now consider the induced homomorphism
\[
f^\sharp(\sigma)(y) = f\left( \sigma\left(f^{-1}(y)\right) \right), \quad \sigma \in \Sym_f(X),\ y \in Y.
\]
It follows that
\[
\Sym_f(X) = \{e, \pi_x\},
\]
where $\pi_x$ is the reflection across the horizontal line $y = \frac{1}{2}$, acting as $a \mapsto d$, $b \mapsto c$, $c \mapsto b$, $d \mapsto a$. Under $f^\sharp$, the identity is mapped to the identity in $\Sym(Y)$, and
\[
f^\sharp(\pi_x) = \rho,
\]
where $\rho$ denotes the $180^\circ$ rotation about the center of the rhombus. Thus, the homomorphism $f^\sharp: \Sym_f(X) \to \Sym(Y)$ is injective but not surjective, and hence not an isomorphism.
\end{example}

\subsection{Symmetry groups and the span construction}\label{section:span}

In the previous construction, given an $n$-tuple configuration $X$, we can associate a symmetry group $\Sym(X)$. In Section \ref{section:persistence_configuration}, we observe that the construction $\Sym:\mathcal{S}_{n}(M) \to \mathbf{Grp}$ may not be a functor. Here, $\mathbf{Grp}$ denotes the category of groups. Proposition \ref{proposition:homomorphism} provides an insight that suggests considering the span category of the group category to study this construction. The category of spans has been introduced in various contexts, and a detailed account can be found in \cite{dawson2010span}.

\begin{definition}
Let $\mathfrak{C}$ be a category. A \textit{span} in $\mathfrak{C}$ is a diagram of the form
\begin{equation*}
X \xleftarrow{f} Z \xrightarrow{g} Y,
\end{equation*}
where $X$, $Y$, and $Z$ are objects in $\mathfrak{C}$, and $f: Z \to X$ and $g: Z \to Y$ are morphisms in $\mathfrak{C}$.
\end{definition}
The diagram $X \xleftarrow{f} Z \xrightarrow{g} Y$ represents a span between the objects $X$ and $Y$, with the common source object $Z$, from which both morphisms $f$ and $g$ map.

Let $\mathfrak{C}$ be a category that admits all pullbacks. The \textit{span category} $\text{Span}(\mathfrak{C})$ is defined as follows: the objects of the category $\text{Span}(\mathfrak{C})$ are the objects of $\mathfrak{C}$. Morphisms in $\text{Span}(\mathfrak{C})$ between two objects $X$ and $Y$ are given by spans of the form $X \xleftarrow{f} Z \xrightarrow{g} Y$, where $Z$ is an object in $\mathfrak{C}$, and $f: Z \to X$ and $g: Z \to Y$ are morphisms in $\mathfrak{C}$. In particular, the identity morphism on an object $X$ is given by the span $X \xleftarrow{\mathrm{id}_X} X \xrightarrow{\mathrm{id}_X} X$, where both legs are the identity morphism on $X$ in $\mathfrak{C}$.

Let $A_1 \xleftarrow{f_1} B_1 \xrightarrow{g_1} A_2$ and $A_2 \xleftarrow{f_2} B_2 \xrightarrow{g_2} A_3$ be two spans. The composition in the span category $\text{Span}(\mathfrak{C})$ is given by the pullback $C$ of the two morphisms $g_1: B_1 \to A_2$ and $f_2: B_2 \to A_2$.
\begin{equation*}
  \xymatrix{
  &&C\ar@{->}[ld]_{p_1}\ar@{->}[rd]^{p_2} &&\\
  &B_{1}\ar@{->}[ld]_{f_{1}}\ar@{->}[rd]^{g_{1}} && B_{2}\ar@{->}[ld]_{f_{2}}\ar@{->}[rd]^{g_{2}}&\\
  A_{1}& &A_{2}& &A_{3}
  }
\end{equation*}
More precisely, the composition of spans is the following span
\begin{equation*}
  A_1 \xleftarrow{f_{1}\circ p_{1}} C \xrightarrow{g_{2}\circ p_{2}} A_3,
\end{equation*}
where $C$ is the pullback object, and $p_1$ and $p_2$ are the induced morphisms via the pullback. It is worth noting that the composition defined in this way is not strictly associative. However, due to the universal property of pullbacks, the composition is associative up to canonical isomorphism. As a result, the category $\Span(\mathfrak{C})$ forms a bicategory rather than a strict category. This flexibility in associativity is a hallmark of many categorical constructions involving pullbacks or fiber products.
When regarded as a bicategory, $\Span(\mathfrak{C})$ admits a well-defined notion of 2-morphisms between spans. A detailed treatment of the bicategorical structure of $\Span(\mathfrak{C})$ is beyond the scope of the present discussion, and will not be pursued further here.

\begin{example}
For a morphism $f:X\to Y$ in $\mathcal{S}_{n}(M)$, we can obtain a span
\begin{equation*}
  \Sym(X)\xleftarrow{f^{\flat}} \Sym_{f}(X) \xrightarrow{f^{\sharp}} \Sym(Y),
\end{equation*}
where $f^{\flat}:\Sym_{f}(X)\hookrightarrow \Sym(X)$ and $f^{\sharp}:\Sym_{f}(X)\to \Sym(Y)$ are group homomorphisms.

\end{example}

\begin{definition}
Let $\mathcal{S}_n(M)$ be the category of $n$-point configurations in $M$. We define the \textit{path 2-category} $\mathbb{S}_n(M)$ as follows:
\begin{itemize}
    \item The objects are the configurations $X \in \mathrm{Ob}(\mathcal{S}_n(M))$.

    \item A 1-cell $\gamma: X \to Y$ is a finite path of composable morphisms in $\mathcal{S}_n(M)$. Explicitly, it is a sequence:
    \begin{equation*}
        \gamma = (X = X_0 \xrightarrow{f_1} X_1 \xrightarrow{f_2} \dots \xrightarrow{f_k} X_k = Y).
    \end{equation*}
    The length of the path is $|\gamma| = k$.

    \item A 2-cell is defined as a contraction of a path. 
    For any two composable morphisms $f: X \to Y$ and $g: Y \to Z$, there is a generating 2-cell
    \begin{equation*}
        \alpha_{f,g}: (f, g) \Rightarrow (g \circ f).
    \end{equation*}
    General 2-cells are generated by these elementary contractions under compositions. Formally, a 2-cell $\theta: \gamma \Rightarrow \gamma'$ exists if and only if $\gamma'$ can be obtained from $\gamma$ by replacing a sequence of morphisms $(f_{i}, f_{i+1})$ with their composite $f_{i+1} \circ f_{i}$ one or more times.
\end{itemize}
\end{definition}
The composition of 1-cells is given by concatenation of paths, which is strictly associative. The vertical composition of 2-cells is transitive.

\begin{definition}\label{definition:sym_functor_construction}
We construct a functor $\Sym: \mathbb{S}_n(M) \to \Span(\mathbf{Grp})$ as follows:
\begin{itemize}
    \item \textit{On objects:} $\Sym(X)$ is the group of symmetries of the configuration $X$.

    \item \textit{On 1-cells:} Let $\gamma = (X_0 \xrightarrow{f_1} \dots \xrightarrow{f_k} X_k)$ be a path. We define $\Sym(\gamma)$ to be the span
    \begin{equation*}
        \Sym(\gamma) = \left( \Sym(X_0) \xleftarrow{\gamma^{\flat}} G_\gamma \xrightarrow{\gamma^{\sharp}} \Sym(X_k) \right).
    \end{equation*}
    The apex group $G_\gamma$ is the group of \textit{survival vectors}:
    \begin{equation*}
        G_\gamma = \left\{ (\sigma_0, \dots, \sigma_{k-1}) \in \prod_{i=0}^{k-1} \Sym_{f_{i+1}}(X_i) \ \middle| \ f_{i+1}^\sharp(\sigma_i) = \sigma_{i+1} \right\}.
    \end{equation*}
    The maps are $\gamma^{\flat}(\sigma_0, \dots, \sigma_{k-1}) = \sigma_0$ and $\gamma^{\sharp}(\sigma_0, \dots, \sigma_{k-1}) = f_k^\sharp(\sigma_{k-1})$.

    \item \textit{On 2-cells:} Let $\alpha: (f, g) \Rightarrow (g \circ f)$ be an elementary contraction 2-cell.
    Its image $\Sym(\alpha)$ is a morphism of spans. The domain is $\Sym(f, g)$, with apex $P = \Sym_f(X) \times_{\Sym(Y)} \Sym_g(Y)$.
    The codomain is $\Sym(g \circ f)$, with apex $G = \Sym_{g \circ f}(X)$.
    
    Since any symmetry surviving continuously through $Y$ must survive the composite map, there is a canonical inclusion $\iota: P \hookrightarrow G$. We define:
    \begin{equation*}
        \Sym(\alpha) = \iota.
    \end{equation*}
    This defines a map of spans, preserving the direction of the 2-cell.
\end{itemize}
\end{definition}

\begin{remark}
With the above definitions, $\Sym$ acts as a \textit{lax functor}. 
The 1-cell structure preserves composition up to isomorphism, while the 2-cell structure captures the inclusion relation:
\begin{equation*}
    \text{Continuous Persistence } (G_{(f,g)}) \subseteq \text{Endpoint Persistence } (G_{g \circ f}).
\end{equation*}
\end{remark}

\begin{lemma}\label{lemma:span_injections}
For a path $\gamma = (X_0 \xrightarrow{f_1} \dots \xrightarrow{f_k} X_k)$ in $\mathbb{S}_n(M)$, the leg maps $\gamma^{\flat}:G_{\gamma}\to \Sym(X_0)$ and $\gamma^{\sharp}:G_{\gamma}\to \Sym(X_k)$ are injective group homomorphisms.
\end{lemma}

\begin{proof}
First, note that $G_{\gamma}$ is a subgroup of the direct product $\prod_{i=0}^{k-1} \Sym_{f_{i+1}}(X_i)$. The map $\gamma^{\flat}$ is the projection onto the first component, and $\gamma^{\sharp}$ is the composition of the projection onto the last component with the group homomorphism $f_k^\sharp$. Thus, both maps are group homomorphisms.

To see that $\gamma^{\flat}$ is injective, suppose $\gamma^{\flat}(\sigma_0, \dots, \sigma_{k-1}) = \gamma^{\flat}(\sigma_0', \dots, \sigma_{k-1}')$. This implies $\sigma_0 = \sigma_0'$. By the definition of $G_{\gamma}$, we have $f_{i+1}^\sharp(\sigma_i) = \sigma_{i+1}$ for all $0 \le i \le k-2$. Since $f_1^\sharp$ is a function, $\sigma_1 = f_1^\sharp(\sigma_0) = f_1^\sharp(\sigma_0') = \sigma_1'$. By induction, it follows that $\sigma_i = \sigma_i'$ for all $0 \le i \le k-1$. Thus, the two vectors are identical, proving $\gamma^{\flat}$ is injective.

For $\gamma^{\sharp}$, suppose $\gamma^{\sharp}(\sigma_0, \dots, \sigma_{k-1}) = \gamma^{\sharp}(\sigma_0', \dots, \sigma_{k-1}')$. By definition, this means $f_k^\sharp(\sigma_{k-1}) = f_k^\sharp(\sigma_{k-1}')$. Since the morphisms in $\mathcal{S}_n(M)$ are isometries between $n$-point configurations (and thus bijections), the induced map $f_k^\sharp$ is given by conjugation and is therefore injective. This implies $\sigma_{k-1} = \sigma_{k-1}'$. Using the compatibility condition $f_{k-1}^\sharp(\sigma_{k-2}) = \sigma_{k-1} = \sigma_{k-1}' = f_{k-1}^\sharp(\sigma_{k-2}')$ and the injectivity of $f_{k-1}^\sharp$, we deduce $\sigma_{k-2} = \sigma_{k-2}'$. Repeating this argument inductively backwards yields $\sigma_i = \sigma_i'$ for all $0 \le i \le k-1$. Hence, $\gamma^{\sharp}$ is also injective.
\end{proof}

\begin{lemma}\label{lemma:span_isomorphism}
Let $f: X \to Y$ and $g: Y \to Z$ be morphisms in $\mathcal{S}_n(M)$. Let $g \star f = (f, g)$ denote the path concatenation. Then we have
\begin{equation*}
  \Sym(g) \circ \Sym(f) \cong \Sym(g \star f).
\end{equation*}
\end{lemma}

\begin{proof}
By definition, the spans for $f$ and $g$ are given by
\begin{equation*}
  \xymatrix{
  \Sym(X)  & \Sym_{f}(X) \ar@{->}[l]_-{f^\flat}\ar@{->}[r]^-{f^\sharp}&\Sym(Y) & \Sym_{g}(Y) \ar@{->}[l]_-{g^\flat}\ar@{->}[r]^-{g^\sharp}&\Sym(Z).
  }
\end{equation*}
The composition of these spans in $\Span(\mathbf{Grp})$ is defined by the pullback of the diagram over $\Sym(Y)$. The vertex of the composed span is
\begin{equation*}
   \Sym_{g}(Y) \times_{\Sym(Y)} \Sym_{f}(X) = \left\{ (\tau, \sigma) \in \Sym_{g}(Y) \times \Sym_{f}(X) \mid g^\flat(\tau) = f^\sharp(\sigma) \right\}.
\end{equation*}
Since $g^\flat: \Sym_{g}(Y) \to \Sym(Y)$ is the inclusion map, the compatibility condition reduces to $\tau = f^\sharp(\sigma)$. Thus, we obtain
\begin{equation*}
  P \cong \left\{ \sigma \in \Sym_{f}(X) \mid f^\sharp(\sigma) \in \Sym_{g}(Y) \right\}.
\end{equation*}

On the other hand, the span for the concatenated path $g \star f = (f, g)$ is given
\begin{equation*}
  G_{(f,g)} = \left\{ (\sigma_0, \sigma_1) \in \Sym_{f}(X) \times \Sym_{g}(Y) \mid f^\sharp(\sigma_0) = \sigma_1 \right\}.
\end{equation*}
Therefore, the vertices of the spans are isomorphic
\begin{equation*}
  \Sym_{g}(Y) \times_{\Sym(Y)} \Sym_{f}(X) \cong G_{(f,g)}.
\end{equation*}
It is straightforward to verify that the left and right leg maps of the composed span coincide with $(g \star f)^\flat$ and $(g \star f)^\sharp$ respectively under this isomorphism. Thus, $\Sym(g) \circ \Sym(f) \cong \Sym(g \star f)$.
\end{proof}

\begin{theorem}\label{theorem:lax_functor}
The construction $\Sym: \mathbb{S}_n(M) \to \Span(\mathbf{Grp})$ is a normal lax functor.
\end{theorem}

\begin{proof}
We verify the coherence conditions for a lax functor.

For any object $X$, the identity 1-cell in $\mathbb{S}_n(M)$ is the empty path. By Definition~\ref{definition:sym_functor_construction}, $\Sym(\mathrm{id}_X)$ is the span with apex $\Sym(X)$ and identity legs.
This is exactly the identity morphism in $\Span(\mathbf{Grp})$. Since there are no non-trivial 2-cells involving the identity path, the identity constraints hold strictly (Normality).

Let $\gamma_1$ and $\gamma_2$ be composable paths. By Lemma \ref{lemma:span_isomorphism} and induction on the length of the paths, there exists a canonical isomorphism
\begin{equation*}
    \psi_{\gamma_1, \gamma_2} : \Sym(\gamma_2) \circ \Sym(\gamma_1) \xrightarrow{\cong} \Sym(\gamma_2 \star \gamma_1).
\end{equation*}
This satisfies the associativity pentagon coherence diagram.

The crucial step lies in the interpretation of the contraction 2-cells. Consider an elementary contraction $\alpha: (f, g) \Rightarrow g \circ f$. We analyze the composite comparison map $\phi$ defined by the lax functor structure:
\begin{equation*}
    \phi_{f,g}: \Sym(g) \circ \Sym(f) \Rightarrow \Sym(g \circ f).
\end{equation*}
Using the isomorphism $\psi$ above, this map factors through $\Sym(f, g)$ and we have
\begin{equation*}
    \Sym(g) \circ \Sym(f) \xrightarrow{\psi^{-1}} \Sym(f, g) \xrightarrow{\Sym(\alpha)} \Sym(g \circ f).
\end{equation*}
By Definition~\ref{definition:sym_functor_construction}, $\Sym(\alpha)$ is the inclusion $\iota: P \hookrightarrow G$. Therefore, the comparison map $\phi_{f,g}$ is an inclusion. Since inclusions are not generally invertible, $\Sym$ is not a pseudo-functor but a lax functor.
\end{proof}

\subsection{Persistent symmetry group}\label{section:persistence_group}

Let $\mathcal{F}: (\mathbb{R}, \leq) \to \mathcal{S}_n(M)$ be a persistence $n$-configuration. For real numbers $a\leq b$, let $f_{a,b}:\mathcal{F}_{a}\to \mathcal{F}_{b}$ denote the corresponding morphism. Let $\mathbf{D}_{a,b}$ denote the poset of finite decompositions of the interval $[a,b]$, ordered by refinement.

For any decomposition $\mathcal{D} = (a = t_0 < t_1 < \dots < t_k = b)$, consider the chain of morphisms $\gamma_{\mathcal{D}} = (f_{t_0, t_1}, \dots, f_{t_{k-1}, t_k})$. Let $G_{\gamma_{\mathcal{D}}}$ denote the domain of the composed Span obtained by applying $\Sym$ to $\gamma_{\mathcal{D}}$, as defined in Definition \ref{definition:sym_functor_construction}. Explicitly, $G_{\gamma_{\mathcal{D}}}$ is isomorphic to the iterated pullback of the symmetry groups along the path
\begin{equation*}
G_{\gamma_{\mathcal{D}}} \cong \Sym(f_{t_{k-1}, t_k}) \times_{\Sym(\mathcal{F}_{t_{k-1}})} \dots \times_{\Sym(\mathcal{F}_{t_1})} \Sym(f_{t_0, t_1}).
\end{equation*}
We define the group of symmetries compatible with this fine-grained path as the image of the induced map to the terminal object:
\begin{equation*}
G(\mathcal{D}) = \mathrm{im}\left( G_{\gamma_{\mathcal{D}}} \to \Sym(\mathcal{F}_b) \right).
\end{equation*}
By construction, a refinement $\mathcal{D}' \geq \mathcal{D}$ induces a subgroup inclusion $G(\mathcal{D}') \subseteq G(\mathcal{D})$. Thus, the assignment $\mathcal{D} \mapsto G(\mathcal{D})$ forms a functor from $\mathbf{D}_{a,b}$ to $\mathbf{Grp}$.

\begin{definition}\label{definition:persistence}
For real numbers $a\leq b$, the \textit{$(a,b)$-persistent symmetry group} is defined as the limit
\begin{equation*}
\Sym(\mathcal{F})_{a,b} := \varprojlim_{\mathcal{D} \in \mathbf{D}_{a,b}} G(\mathcal{D}).
\end{equation*}
\end{definition}

Furthermore, let us define the total symmetry space over the interval as the limit of the domains:
\begin{equation*}
   G_{a,b} =  \varprojlim_{\mathcal{D} \in \mathbf{D}_{a,b}} G_{\gamma_{\mathcal{D}}}.
\end{equation*}
Since the morphisms in the diagram are monomorphisms (subgroup inclusions), the limit commutes with the image formation. Consequently, we have
\begin{equation*}
 \Sym(\mathcal{F})_{a,b} = \mathrm{im}\left( G_{a,b} \to \Sym(\mathcal{F}_b) \right).
\end{equation*}

\begin{remark}\label{remark:spans}
By Lemma \ref{lemma:span_injections}, the leg maps of the spans induced by the limit $G_{a,b}$ are injective. Consequently, the spans
\begin{equation*}
  \Sym(\mathcal{F}_a) \leftarrow G_{a,b} \rightarrow \Sym(\mathcal{F}_b),
\end{equation*}
and
\begin{equation*}
  \Sym(\mathcal{F}_a) \leftarrow \Sym(\mathcal{F})_{a,b} \hookrightarrow \Sym(\mathcal{F}_b),
\end{equation*}
are isomorphic, and thus equivalent, in the category $\Span(\mathbf{Grp})$.

Furthermore, let $L(\mathcal{D}) = \mathrm{im}\left( G_{\gamma_{\mathcal{D}}} \to \Sym(\mathcal{F}_a) \right)$ and define $L_{a,b} = \varprojlim_{\mathcal{D} \in \mathbf{D}_{a,b}} L(\mathcal{D})$. Since the left leg maps are injective, $L_{a,b}$ naturally inherits a map to $\Sym(\mathcal{F}_b)$, yielding a third span equivalent to the previous two
\begin{equation*}
  \Sym(\mathcal{F}_a) \hookleftarrow L_{a,b} \rightarrow \Sym(\mathcal{F}_b).
\end{equation*}
The elements of $L_{a,b}$ represent precisely those symmetries in $\Sym(\mathcal{F}_a)$ that can be continuously propagated along the interval $[a,b]$ to yield elements in $\Sym(\mathcal{F}_b)$.
\end{remark}

The persistent symmetry group admits an equivalent computational definition, which offers a more intuitive criterion for identifying valid symmetries.

\begin{proposition}
The $(a,b)$-persistent symmetry group is the intersection of images over all intermediate times. Specifically,
\begin{equation*}
  \Sym(\mathcal{F})_{a,b} = \bigcap_{t \in [a,b]} \mathrm{im}\left( f_{t,b}^\sharp: \Sym_{f_{t,b}}(\mathcal{F}_t) \to \Sym(\mathcal{F}_b) \right).
\end{equation*}
Equivalently, a symmetry $\sigma \in \Sym(\mathcal{F}_b)$ belongs to $\Sym(\mathcal{F})_{a,b}$ if and only if for every intermediate time $t \in [a,b]$, there exists a symmetry $\sigma_t \in \Sym(\mathcal{F}_t)$ that evolves into $\sigma$.
\end{proposition}

\begin{proof}
We prove the equality by establishing mutual subset inclusion.

For any $t \in [a,b]$, define the intermediate group as the image of the span leg
\begin{equation*}
H_t = \mathrm{im}\left( f_{t,b}^{\sharp}: \Sym_{f_{t,b}}(\mathcal{F}_t) \to \Sym(\mathcal{F}_b) \right).
\end{equation*}
For any finite decomposition $\mathcal{D} = (a = t_0 < t_1 < \dots < t_k = b)$, define the decomposition group as the image of the composed span:
\begin{equation*}
G(\mathcal{D}) = \mathrm{im}\left( \Sym(f_{t_{k-1}, t_k}) \circ \dots \circ \Sym(f_{t_0, t_1}) \to \Sym(\mathcal{F}_b) \right).
\end{equation*}
Note that since the morphisms in the diagram are subgroup inclusions, the categorical limit is isomorphic to the set-theoretic intersection:
\begin{equation*}
\Sym(\mathcal{F})_{a,b} = \varprojlim_{\mathcal{D}} G(\mathcal{D}) = \bigcap_{\mathcal{D}} G(\mathcal{D}).
\end{equation*}

Let $\sigma \in \Sym(\mathcal{F})_{a,b} = \bigcap_{\mathcal{D}} G(\mathcal{D})$. For any specific time $t \in [a,b]$, consider the decomposition $\mathcal{D}_t = \{a, t, b\}$. Since $\sigma \in \bigcap_{\mathcal{D}} G(\mathcal{D})$, it follows that $\sigma \in G(\mathcal{D}_t)$. By the functorial composition of spans, we have the inclusion $G(\mathcal{D}_t) \subseteq H_t$. Therefore, $\sigma \in H_t$. Since $t$ was arbitrary, we conclude $\sigma \in \bigcap_{t \in [a,b]} H_t$.

Let $\sigma \in \bigcap_{t \in [a,b]} H_t$. This implies that for every $t \in [a,b]$, there exists a symmetry $\sigma_t \in \Sym_{f_{t,b}}(\mathcal{F}_t)$ such that $\sigma_t \mapsto \sigma$. Consider an arbitrary finite decomposition $\mathcal{D} = (t_0, \dots, t_k)$. Since $\sigma$ survives at every intermediate time $t_i$, there exists a chain of symmetries $\sigma_{t_0} \mapsto \sigma_{t_1} \mapsto \dots \mapsto \sigma_{t_k}=\sigma$. These compatible intermediate symmetries define a valid element in the vertex of the composed span $G_{\gamma_{\mathcal{D}}}$.Consequently, $\sigma$ belongs to the image of the composed span, i.e., $\sigma \in G(\mathcal{D})$. Since $\mathcal{D}$ was arbitrary, $\sigma \in \bigcap_{\mathcal{D}} G(\mathcal{D}) = \Sym(\mathcal{F})_{a,b}$.
\end{proof}

Intuitively, the elements of the $(a,b)$-persistent symmetry group $\Sym(\mathcal{F})_{a,b}$ are those symmetries that persist throughout the entire interval $[a,b]$. Specifically, these elements reside in $\Sym(\mathcal{F}_b)$, originate from $\Sym(\mathcal{F}_a)$, and survive at every intermediate time $t \in [a,b]$.

We denote the direct sum $\Sym(\mathcal{F}) = \bigoplus_{a \in \mathbb{R}} \Sym(\mathcal{F})_a$ as the \textit{persistent symmetry group}. If the filtration parameter is taken over the integers $\mathbb{Z}$, the persistent symmetry group can be regarded as a graded group:
\begin{equation*}
  \Sym(\mathcal{F})= \bigoplus\limits_{a \in \mathbb{Z}} \Sym(\mathcal{F})_a.
\end{equation*}
It is worth noting that the identity element, which corresponds to the identity transformation in metric space, is always a persistent element. Therefore, the persistent symmetry group is never empty. Unless explicitly stated otherwise, we generally do not consider the identity element when analyzing symmetry.

\begin{definition}
Let $\mathcal{F}: (\mathbb{R}, \leq) \to \mathcal{S}_{n}(M)$ be a persistence $n$-configuration.
For a symmetry $\pi \in \Sym(\mathcal{F}_a)$, its \textit{persistent interval} is defined as the maximal interval $[a,b)$ such that for all $t \in [a,b)$, the image of $\pi$ under the map
\[
f_{a,t}^{\sharp}: \Sym_{f_{a,t}}(\mathcal{F}_a) \to \Sym(\mathcal{F}_t)
\]
remains in the image $\mathrm{im}(f_{a,t}^{\sharp})$ and represents a nontrivial symmetry of $\mathcal{F}_t$.

Moreover, if for all $t < a$, the element $\pi$ is not the image of any symmetry in $\Sym(\mathcal{F}_t)$ under the transition map $f_{t,a}^{\sharp}$, then $[a,b)$ is called a \textit{symmetry bar}.
\end{definition}

\begin{definition}
Let $\mathcal{F}: (\mathbb{R}, \leq) \to \mathcal{S}_{n}(M)$ be a persistence $n$-configuration.
The \textit{symmetry barcode} of $\mathcal{F}$ is the multiset of intervals $\{[a_i, b_i)\}_{i \in I}$, where each interval $[a_i, b_i)$ corresponds to a symmetry bar, representing the lifespan of a symmetry that persists under the transition maps of $\mathcal{F}$.
\end{definition}

\begin{proposition}\label{proposition:inverse_bar}
Let $\mathcal{F}: (\mathbb{R}, \leq) \to \mathcal{S}_n(M)$ be a persistence $n$-configuration, and let $\pi \in \Sym(\mathcal{F}_a)$ be a symmetry whose symmetry bar is the interval $[a, b)$. Then the inverse $\pi^{-1}$ also has the same symmetry bar. That is,
\[
\pi \in \Sym_{f_{a,t}}(\mathcal{F}_a) \text{ for all } t \in [a, b) \quad \Longrightarrow \quad \pi^{-1} \in \Sym_{f_{a,t}}(\mathcal{F}_a) \text{ for all } t \in [a, b).
\]
\end{proposition}

\begin{proof}
Assume that for all $t \in [a, b)$, we have $\pi \in \Sym_{f_{a,t}}(\mathcal{F}_a)$. It follows that
\[
f_{a,t} \circ \pi \circ f_{a,t}^{-1} \in \Sym(\mathcal{F}_t).
\]
Since $\Sym(\mathcal{F}_t)$ is a group, it is closed under taking inverses. Thus, the inverse
\[
 \left(f_{a,t} \circ \pi \circ f_{a,t}^{-1}\right)^{-1} = f_{a,t} \circ \pi^{-1} \circ f_{a,t}^{-1}
\]
also belongs to $\Sym(\mathcal{F}_t)$. Hence, $\pi^{-1} \in \Sym_{f_{a,t}}(\mathcal{F}_a)$ for all $t \in [a, b)$.
\end{proof}

\begin{proposition}
Let $\mathcal{F}: (\mathbb{R}, \leq) \to \mathcal{S}_n(M)$ be a persistence $n$-configuration. Suppose $\pi, \varpi \in \Sym(\mathcal{F}_a)$ have persistent intervals $[a, b_\pi)$ and $[a, b_{\varpi})$, respectively. Then the composition $\pi \circ \varpi$ has a symmetry bar that contains the interval $[a, \min(b_\pi, b_{\varpi}))$.
\end{proposition}

\begin{proof}
Fix $t \in [a, \min(b_\pi, b_{\varpi}))$. By the definition of symmetry bars, we have
\[
f_{a,t} \circ \pi \circ f_{a,t}^{-1} \in \Sym(\mathcal{F}_t), \quad
f_{a,t} \circ \varpi \circ f_{a,t}^{-1} \in \Sym(\mathcal{F}_t).
\]
Since $\Sym(\mathcal{F}_t)$ is a group, it is closed under composition. Therefore,
\[
f_{a,t} \circ (\pi \circ \varpi) \circ f_{a,t}^{-1}
= (f_{a,t} \circ \pi \circ f_{a,t}^{-1}) \circ (f_{a,t} \circ \varpi \circ f_{a,t}^{-1}) \in \Sym(\mathcal{F}_t).
\]
This implies that $\pi \circ \varpi \in \Sym_{f_{a,t}}(\mathcal{F}_a)$ for all $t \in [a, \min(b_\pi, b_{\varpi}))$, as claimed.
\end{proof}

\begin{example}\label{example:symmetry_group}
Consider the following persistence configuration
\[
\mathcal{F} = \{X_t\}_{t=0,1,2} \in \mathcal{S}_4(\mathbb{R}^2)^{(\mathbb{Z}, \leq)},
\]
where each $X_t = \{A_t, B_t, C_t, D_t\}$ is a configuration of four points in $\mathbb{R}^2$ at discrete time steps $t = 0, 1, 2$. The coordinates of the points $A_t, B_t, C_t, D_t \in \mathbb{R}^2$ are listed in Table~\ref{table:4configuration}, and the corresponding geometric configurations are illustrated in Figure~\ref{figure:4configuration}. For each pair of time steps $0 \leq i \leq j \leq 2$, the structure map $f_{i,j}: X_i \to X_j$ maps $A_i$, $B_i$, $C_i$, and $D_i$ to $A_j$, $B_j$, $C_j$, and $D_j$, respectively.

\begin{table}[h]
\centering
\begin{tabular}{c|c|c|c|c}
\hline
Time $t$ & Point $A_t$ & Point $B_t$ & Point $C_t$ & Point $D_t$ \\
\hline\hline
$t = 0$ & $(0, -1.2)$ & $(0, 1.2)$ & $(-1, 0)$ & $(1, 0)$ \\
$t = 1$ & $(0, -1)$ & $(0, 1)$ & $(-1, 0)$ & $(1, 0)$ \\
$t = 2$ & $(0, -1)$ & $(0, 1)$ & $(-1, 0)$ & $(1.2, 0)$ \\
\hline
\end{tabular}
\caption{Coordinates of the persistent 4-configuration at different time steps in Example \ref{example:symmetry_group}.}\label{table:4configuration}
\end{table}

\begin{figure}[h]
\centering
\definecolor{axisgray}{gray}{0.6}

\begin{tikzpicture}[scale=0.8, >=Stealth]

\begin{scope}
    \draw[->, color=axisgray, line width=0.4pt] (-1.8,0) -- (1.8,0) node[right, font=\small, text=black] {$x$};
    \draw[->, color=axisgray, line width=0.4pt] (0,-1.8) -- (0,1.8) node[above, font=\small, text=black] {$y$};

    \draw[line width=1.2pt, color=mainblue] (0,-1.2) -- (-1,0) -- (0,1.2) -- (1,0) -- cycle;

    \fill[mainblue] (0,-1.2) circle (2pt); \node[below right, font=\small] at (0,-1.2) {$A_0$};
    \fill[mainblue] (0, 1.2) circle (2pt); \node[above right, font=\small] at (0, 1.2) {$B_0$};
    \fill[mainblue] (-1,0) circle (2pt); \node[above left, font=\small] at (-1,0) {$C_0$};
    \fill[mainblue] ( 1,0) circle (2pt); \node[above right, font=\small] at (1,0) {$D_0$};
    
    \node[below] at (0, -2.2) {$t = 0$};
\end{scope}

\begin{scope}[xshift=6cm]
    \draw[->, color=axisgray, line width=0.4pt] (-1.8,0) -- (1.8,0) node[right, font=\small, text=black] {$x$};
    \draw[->, color=axisgray, line width=0.4pt] (0,-1.8) -- (0,1.8) node[above, font=\small, text=black] {$y$};

    \draw[line width=1.2pt, color=mainblue] (0,-1) -- (-1,0) -- (0,1) -- (1,0) -- cycle;
    
    \fill[mainblue] (0,-1) circle (2pt); \node[below right, font=\small] at (0,-1) {$A_1$};
    \fill[mainblue] (0, 1) circle (2pt); \node[above right, font=\small] at (0, 1) {$B_1$};
    \fill[mainblue] (-1,0) circle (2pt); \node[above left, font=\small] at (-1,0) {$C_1$};
    \fill[mainblue] ( 1,0) circle (2pt); \node[above right, font=\small] at (1,0) {$D_1$};
    
    \node[below] at (0, -2.2) {$t = 1$};
\end{scope}

\begin{scope}[xshift=12cm]
    \draw[->, color=axisgray, line width=0.4pt] (-1.8,0) -- (1.8,0) node[right, font=\small, text=black] {$x$};
    \draw[->, color=axisgray, line width=0.4pt] (0,-1.8) -- (0,1.8) node[above, font=\small, text=black] {$y$};

    \draw[line width=1.2pt, color=mainblue] (0,-1) -- (-1,0) -- (0,1) -- (1.2,0) -- cycle;
    
    \fill[mainblue] (0,-1) circle (2pt); \node[below right, font=\small] at (0,-1) {$A_2$};
    \fill[mainblue] (0, 1) circle (2pt); \node[above right, font=\small] at (0, 1) {$B_2$};
    \fill[mainblue] (-1,0) circle (2pt); \node[above left, font=\small] at (-1,0) {$C_2$};
    \fill[mainblue] ( 1.2,0) circle (2pt); \node[above right, font=\small] at (1.2,0) {$D_2$};
    
    \node[below] at (0, -2.2) {$t = 2$};
\end{scope}

\end{tikzpicture}
\caption{Evolution of a 4-point configuration: from a rhombus, to a square, and then to a quadrilateral.}\label{figure:4configuration}
\end{figure}

We now compute the symmetry barcode of this persistent configuration.
At $t = 0$, the configuration $X_0$ is a rhombus with horizontal and vertical reflection symmetry. Its symmetry group is given by
\[
\Sym(X_0) \cong \mathbb{Z}/2 \times \mathbb{Z}/2.
\]
At $t = 1$, the configuration $X_1$ is a square, whose symmetry group is the dihedral group
\[
\Sym(X_1) \cong D_4.
\]
At $t = 2$, the configuration $X_2$ is a quadrilateral with the symmetry group
\[
\Sym(X_2) \cong \mathbb{Z}/2.
\]

We proceed to compute the persistent symmetry groups $\Sym_{f_{0,1}}(X_0)$ and $\Sym_{f_{1,2}}(X_1)$. The symmetry group $\Sym_{f_{0,1}}(X_0)$ consists of the symmetries of $X_0$ that remain valid under the transition map $f_{0,1}$. Calculation shows that all symmetries of $X_0$ persist to $X_1$, hence
\[
\Sym_{f_{0,1}}(X_0) \cong \mathbb{Z}/2 \times \mathbb{Z}/2.
\]
For the next interval, we find that only the identity and the reflection symmetry $\pi_x$ (reflection with respect to the $x$-axis) persist until $t = 2$. Thus,
\[
\Sym_{f_{1,2}}(X_1) \cong \mathbb{Z}/2.
\]

Excluding the identity element, we obtain the symmetry barcode of the persistence $n$-configuration $\mathcal{F}$ as the multiset:
\[
\mathrm{SymBar}(\mathcal{F}) = \left\{ [0,1), [0,1), [0,2), [1,2), [1,2), [1,2), [1,2) \right\}.
\]
Specifically, the two intervals $[0,1)$ correspond to the reflection symmetry about the $y$-axis and the $180^{\circ}$ rotational symmetry, respectively. The interval $[0,2)$ represents the symmetry bar for the reflection with respect to the $x$-axis. The four intervals $[1,2)$ correspond to the remaining four symmetries in $D_4$ (rotations by $90^{\circ}$ and $270^{\circ}$, and reflections about the diagonals) that emerge at time step $t=1$ but vanish at $t=2$.
\end{example}

The elements of the $(a,b)$-persistent symmetry group introduced above correspond to the notion of $(a,b)$-persistent symmetry defined in Section \ref{section:persistent_set}. In Section~\ref{section:persistent_set}, the collection of persistent symmetries over an interval was treated as a set. In contrast, Definition~\ref{definition:persistence} endows the collection of all persistent symmetries with a group structure. However, this group is fundamentally different from homology groups, as it is generally \textit{non-Abelian}, leading to a significantly more intricate algebraic structure compared to persistent homology. When working with field coefficients, persistent homology yields an Abelian group in the form of a vector space, where the group operation is given by addition, and Betti numbers serve as key invariants capturing topological features. In contrast, the persistent symmetry group encodes richer algebraic information due to its non-Abelian nature, making it inherently more complex and structurally richer than persistent homology. This complexity suggests that persistent symmetry groups may offer deeper insights and broader applications in the study of evolving geometric and topological structures.

\subsection{Persistent symmetry module}

Let $\Span(\mathbf{Set})$ denote the bicategory of spans in the category $\mathbf{Set}$ of sets. Let $\mathbf{Vec}_{\mathbb{K}}$ denote the category of vector spaces over the field $\mathbb{K}$, which we view as a bicategory with only identity 2-cells. Consider the following construction
\begin{equation*}
F\colon \Span(\mathbf{Set}) \to \mathbf{Vec}_{\mathbb{K}},
\end{equation*}
which assigns to each set $S$ the free vector space $F(S) = \mathbb{K}[S]$ generated by the elements of $S$, and to each span $S \xleftarrow{f} S' \xrightarrow{g} T$ a linear map $F(f, g)\colon F(S) \to F(T)$ defined on basis elements $x \in S$ by
\begin{equation*}
F(f, g)(x) = \sum_{\substack{x' \in f^{-1}(x)}} g(x').
\end{equation*}
In particular, if $f^{-1}(x) = \emptyset$, the sum is defined to be $0$.

\begin{proposition}\label{proposition:free}
The construction $F\colon \Span(\mathbf{Set}) \to \mathbf{Vec}_{\mathbb{K}}$ defines a pseudofunctor.
\end{proposition}

\begin{proof}
We verify that $F$ preserves composition and identity up to coherent isomorphism, characteristic of a pseudofunctor. Consider the composition of two spans in $\Span(\mathbf{Set})$ given by the following pullback diagram.
\begin{equation*}
  \xymatrix@=0.6cm{
  &&X'\times_{Y} Y'\ar[ld]_{p_1}\ar[rd]^{p_2} &&\\
  &X'\ar[ld]_{f_{1}}\ar[rd]^{g_{1}} && Y'\ar[ld]_{f_{2}}\ar[rd]^{g_{2}}&\\
  X& &Y& &Z
  }
\end{equation*}
We need to show that the linear map induced by the composed span coincides with the composition of the linear maps induced by each span. Specifically, we check that
\begin{equation*}
F(f_{2},g_{2})\circ F(f_{1},g_{1}) = F(f_{1}\circ p_{1},g_{2}\circ p_{2}),
\end{equation*}
where $p_1\colon X'\times_Y Y' \to X'$, $p_2\colon X'\times_Y Y' \to Y'$ are the canonical projections. By construction, for any basis element $x \in X$, we compute:
\begin{equation*}
\begin{aligned}
\left(F(f_{2},g_{2})\circ F(f_{1},g_{1})\right)(x)
&= F(f_{2},g_{2})\left(\sum_{x' \in f_1^{-1}(x)} g_1(x')\right) \\
&= \sum_{\substack{x' \in f_1^{-1}(x) \\ y' \in f_2^{-1}(g_1(x'))}} g_2(y').
\end{aligned}
\end{equation*}
On the other hand, applying $F$ to the composed span yields:
\begin{equation*}
F(f_1 \circ p_1, g_2 \circ p_2)(x) = \sum_{(x', y') \in (f_1 \circ p_1)^{-1}(x)} g_2(y').
\end{equation*}
Note that the condition $(x', y') \in (f_1 \circ p_1)^{-1}(x)$ is equivalent to $f_1(x') = x$ and $g_1(x') = f_2(y')$. This precisely corresponds to the indexing set of the double summation above. Therefore, the two expressions coincide, verifying the composition coherence condition. 
Moreover, for the identity span $X \xleftarrow{\mathrm{id}_X} X \xrightarrow{\mathrm{id}_X} X$, we have
\begin{equation*}
F(\mathrm{id}_X, \mathrm{id}_X)(x) = \sum_{x' \in \mathrm{id}^{-1}(x)} \mathrm{id}(x') = x,
\end{equation*}
which implies $F(\mathrm{id}_X, \mathrm{id}_X) = \mathrm{id}_{F(X)}$.
Since the associator isomorphisms of span composition arise from the universal property of pullbacks, and we have verified strict equality on the level of elements, $F$ maps these associators to identity maps in $\mathbf{Vec}_{\mathbb{K}}$. Thus, it is a pseudofunctor.
\end{proof}

In conjunction with the previously introduced symmetry construction, we have the following sequence of constructions
\begin{equation*}
\xymatrix{
  \mathbb{S}_{n}(M) \ar[r]^-{\Sym} & \Span(\mathbf{Grp}) \ar[r]^-{U} & \Span(\mathbf{Set}) \ar[r]^-{F} & \mathbf{Vec}_{\mathbb{K}}.
}
\end{equation*}
Here, $U\colon \Span(\mathbf{Grp}) \to \Span(\mathbf{Set})$ is the forgetful pseudofunctor. This is well-defined since the usual forgetful functor $\mathbf{Grp} \to \mathbf{Set}$ preserves limits, including pullbacks.

\begin{proposition}\label{proposition:functor}
The construction $\mathcal{M} = F \circ U \circ \Sym : \mathbb{S}_n(M) \to \mathbf{Vec}_{\mathbb{K}}$ is a strict 2-functor, where $\mathbf{Vec}_{\mathbb{K}}$ is regarded as a 2-category with only identity 2-cells.
\end{proposition}

\begin{proof}
We first analyze the action of $\mathcal{M}$ on 1-cells. Let $\gamma = (X_0 \xrightarrow{f_1} X_1 \xrightarrow{f_2} \dots \xrightarrow{f_k} X_k)$ be a path in $\mathbb{S}_n(M)$. The span $\Sym(\gamma)$ is given by 
\begin{equation*}
  \Sym(X_0) \xleftarrow{\gamma^\flat} G_\gamma \xrightarrow{\gamma^\sharp} \Sym(X_k),
\end{equation*}
where $G_\gamma$ consists of survival vectors $(\sigma_0, \dots, \sigma_{k-1})$ such that $f_{i+1}^\sharp(\sigma_i) = \sigma_{i+1}$, $\gamma^\flat$ is the projection to $\sigma_0$, and $\gamma^\sharp$ is the projection to $\sigma_{k-1}$ followed by $f_k^\sharp$. 
For any basis element $\sigma \in \Sym(X_0)$, applying the free vector space pseudofunctor $F \circ U$ yields:
\begin{equation*}
  \mathcal{M}(\gamma)(\sigma) = (F \circ U)(\Sym(\gamma))(\sigma) = \left\{
                                                              \begin{array}{ll}
                                                                f_k^\sharp \circ \dots \circ f_1^\sharp(\sigma), & \hbox{if $\sigma \in \mathrm{im}(\gamma^\flat)$;} \\
                                                                0, & \hbox{otherwise.}
                                                              \end{array}
                                                            \right.
\end{equation*}
By definition, $\sigma \in \mathrm{im}(\gamma^\flat)$ if and only if there exists a valid survival vector starting with $\sigma_0 = \sigma$. This is equivalent to the conditions that $\sigma \in \Sym_{f_1}(X_0)$, $f_1^\sharp(\sigma) \in \Sym_{f_2}(X_1)$, $\dots$, $f_{k-1}^\sharp \circ \dots \circ f_1^\sharp(\sigma) \in \Sym_{f_k}(X_{k-1})$. 
By iteratively applying the equivalence $\sigma \in \Sym_{g\circ f}(X) \iff (\sigma \in \Sym_f(X) \text{ and } f^\sharp(\sigma) \in \Sym_g(Y))$, we conclude that this survival condition is precisely equivalent to $\sigma \in \Sym_{h}(X_0)$, where $h = f_k \circ \dots \circ f_1$ is the composite morphism of the path $\gamma$. Furthermore, $f_k^\sharp \circ \dots \circ f_1^\sharp(\sigma) = h^\sharp(\sigma)$. Thus, the linear map depends only on the composite morphism:
\begin{equation*}
  \mathcal{M}(\gamma)(\sigma) = \left\{
                                                              \begin{array}{ll}
                                                                h^\sharp(\sigma), & \hbox{if $\sigma \in \Sym_{h}(X_0)$;} \\
                                                                0, & \hbox{otherwise.}
                                                              \end{array}
                                                            \right.
\end{equation*}
This implies $\mathcal{M}(\gamma) = \mathcal{M}(h)$.

Now, let $\gamma_1: X \to Y$ and $\gamma_2: Y \to Z$ be two composable 1-cells in $\mathbb{S}_n(M)$. Let their underlying composite morphisms be $h_1$ and $h_2$ respectively. The composite path $\gamma_2 \star \gamma_1$ has underlying morphism $h_2 \circ h_1$. Using the observation above, we have:
\begin{equation*}
  \mathcal{M}(\gamma_2 \star \gamma_1) = \mathcal{M}(h_2 \circ h_1).
\end{equation*}
For the composite morphism $h_2 \circ h_1$, which can be seen as a path of length 1, a direct computation gives:
\begin{equation*}
  \mathcal{M}(h_2 \circ h_1)(\sigma) = \left\{
                                                              \begin{array}{ll}
                                                                (h_2 \circ h_1)^\sharp(\sigma), & \hbox{if $\sigma \in \Sym_{h_2 \circ h_1}(X)$;} \\
                                                                0, & \hbox{otherwise.}
                                                              \end{array}
                                                            \right.
\end{equation*}
As established, the condition $\sigma \in \Sym_{h_2 \circ h_1}(X)$ is equivalent to $\sigma \in \Sym_{h_1}(X)$ and $h_1^\sharp(\sigma) \in \Sym_{h_2}(Y)$. Therefore, we have
\begin{equation*}
  \mathcal{M}(h_2 \circ h_1)(\sigma) = \mathcal{M}(h_2)\big(\mathcal{M}(h_1)(\sigma)\big) = \mathcal{M}(\gamma_2) \circ \mathcal{M}(\gamma_1)(\sigma).
\end{equation*}
This verifies that $\mathcal{M}$ strictly preserves the composition of 1-cells: $\mathcal{M}(\gamma_2 \star \gamma_1) = \mathcal{M}(\gamma_2) \circ \mathcal{M}(\gamma_1)$. The identity 1-cell is strictly preserved since its underlying morphism is $\mathrm{id}_X$, and $\mathcal{M}(\mathrm{id}_X)(\sigma) = \sigma$.

Finally, we consider the 2-cells. The bicategory $\mathbf{Vec}_{\mathbb{K}}$ is locally discrete, meaning it only possesses identity 2-cells. Since $\mathcal{M}$ is a composite of pseudofunctors, it maps 2-cells in $\mathbb{S}_n(M)$ to 2-cells in $\mathbf{Vec}_{\mathbb{K}}$, which must be the identity 2-cell. Because any 2-cell $\alpha: \gamma \Rightarrow \gamma'$ in $\mathbb{S}_n(M)$ is a contraction of paths, $\gamma$ and $\gamma'$ share the same underlying composite morphism $h$. By our earlier observation, $\mathcal{M}(\gamma) = \mathcal{M}(h) = \mathcal{M}(\gamma')$, so the identity 2-cell is indeed a valid 2-cell between them. 

Thus, $\mathcal{M}$ strictly preserves all 1-cell compositions, identity 1-cells, and maps all 2-cells to identities, satisfying the axioms of a strict 2-functor.
\end{proof}

Let $\mathcal{F}:(\mathbb{R}, \leq) \to \mathcal{S}_{n}(M)$ be a persistence $n$-configuration. Note that $\mathcal{S}_n(M)$ embeds into the path 2-category $\mathbb{S}_n(M)$ as the subcategory of length-1 paths (with only identity 2-cells). By Proposition \ref{proposition:functor}, the strict 2-functor $\mathcal{M}: \mathbb{S}_n(M) \to \mathbf{Vec}_{\mathbb{K}}$ restricts to a strict functor on this subcategory, and thus the composition
\begin{equation*}
  \mathcal{M}\mathcal{F}:(\mathbb{R}, \leq) \to \mathbf{Vec}_{\mathbb{K}}
\end{equation*}
is well-defined, yielding a persistence module. For any real numbers $s\leq t$, the induced $\mathbb{K}$-linear map is
\begin{equation*}
  \mathcal{M}_{s\to t}:\mathcal{M}\mathcal{F}_{s}\to \mathcal{M}\mathcal{F}_{t},
\end{equation*}
where for any basis element $\sigma \in \Sym(\mathcal{F}_s)$,
\begin{equation*}
  \mathcal{M}_{s\to t}(\sigma)=\left\{
                               \begin{array}{ll}
                                 f_{s,t}^{\sharp}(\sigma), & \hbox{$\sigma\in \Sym_{f_{s,t}}(\mathcal{F}_{s})$;} \\
                                 0, & \hbox{otherwise.}
                               \end{array}
                             \right.
\end{equation*}
Here, $f_{s,t}:\mathcal{F}_{s}\to \mathcal{F}_{t}$ is the morphism induced by $s\leq t$. Inspired by the definition of the persistent symmetry group, which requires symmetries to survive at all intermediate times, we define the persistent symmetry module as follows.

\begin{definition}
For real numbers $a\leq b$, the \textit{$(a,b)$-persistent symmetry module} is defined as the intersection of the images of all intermediate linear maps
\begin{equation*}
  \mathcal{M}_{a,b}(\mathcal{F}) := \bigcap_{t \in [a,b]} \mathrm{im}\left( \mathcal{M}_{t \to b}: \mathcal{M}\mathcal{F}_{t}\to \mathcal{M}\mathcal{F}_{b}\right).
\end{equation*}
\end{definition}

For the persistence $n$-configuration $\mathcal{F}: (\mathbb{Z}, \leq) \to \mathcal{S}_{n}(M)$ with an integer filtration, the above definition also applies. In this case, let $\mathbf{M}(\mathcal{F})=\bigoplus\limits_{a\in \mathbb{Z}}\mathcal{M}\mathcal{F}_{a}$. Then $\mathbf{M}(\mathcal{F})$ coincides with the persistence module of symmetries introduced in Section \ref{section:module}.

Let $\mathbb{K}[\Sym(\mathcal{F})_{a,b}]$ denote the $\mathbb{K}$-linear space generated by the elements in $\Sym(\mathcal{F})_{a,b}$. Since $\Sym(\mathcal{F})_{a,b}$ is a finite set for fixed $a \leq b$, the space $\mathbb{K}[\Sym(\mathcal{F})_{a,b}]$ is finite-dimensional. Moreover, as $\Sym(\mathcal{F})_{a,b}$ is a group, $\mathbb{K}[\Sym(\mathcal{F})_{a,b}]$ inherits a natural structure of a finite-dimensional group algebra over $\mathbb{K}$. More precisely, the multiplication in $\mathbb{K}[\Sym(\mathcal{F})_{a,b}]$ is induced by the group operation in $\Sym(\mathcal{F})_{a,b}$. Let $\{\pi_i\}_{i \in I}$ represent the elements of $\Sym(\mathcal{F})_{a,b}$. Then, for any two elements
\begin{equation*}
x = \sum_{i \in I} a_i \pi_i, \quad y = \sum_{j \in I} b_j \pi_j,
\end{equation*}
where $a_i, b_j \in \mathbb{K}$ for all indices $i, j$, the multiplication in $\mathbb{K}[\Sym(\mathcal{F})_{a,b}]$ is given by
\begin{equation*}
x \cdot y = \left(\sum_{i \in I} a_i \pi_i \right) \left(\sum_{j \in I} b_j \pi_j \right) = \sum_{i,j \in I} a_i b_j (\pi_i \cdot \pi_j),
\end{equation*}
where $\pi_i \cdot \pi_j$ is the group multiplication in $\Sym(\mathcal{F})_{a,b}$. Thus, $\mathbb{K}[\Sym(\mathcal{F})_{a,b}]$ is naturally a group algebra, which provides richer algebraic information for studying the persistence of symmetries in $\mathcal{F}$.

\begin{proposition}
For any real numbers $a \leq b$, the $(a,b)$-persistent symmetry module is given by
\begin{equation*}
  \mathcal{M}_{a,b}(\mathcal{F}) = \mathbb{K}[\Sym(\mathcal{F})_{a,b}].
\end{equation*}
Thus, we can endow $\mathcal{M}_{a,b}(\mathcal{F})$ with the group algebra structure inherited from $\mathbb{K}[\Sym(\mathcal{F})_{a,b}]$.
\end{proposition}

\begin{proof}
For any $t \in [a,b]$, the image of the linear map $\mathcal{M}_{t \to b}$ is generated by the elements in the image of $f_{t,b}^{\sharp}$. That is,
\begin{equation*}
  \mathrm{im}(\mathcal{M}_{t \to b}) = \mathbb{K}\left[ \mathrm{im}\left( f_{t,b}^{\sharp}: \Sym_{f_{t,b}}(\mathcal{F}_t) \to \Sym(\mathcal{F}_b) \right) \right].
\end{equation*}
By the equivalent definition of the persistent symmetry group, we have
\begin{equation*}
\Sym(\mathcal{F})_{a,b} = \bigcap_{t \in [a,b]} \mathrm{im}\left( f_{t,b}^{\sharp}: \Sym_{f_{t,b}}(\mathcal{F}_t) \to \Sym(\mathcal{F}_b) \right).
\end{equation*}
Since the elements of $\Sym(\mathcal{F}_b)$ form a basis for $\mathcal{M}\mathcal{F}_b = \mathbb{K}[\Sym(\mathcal{F}_b)]$, the intersection of the linear spans of subsets of this basis is precisely the linear span of their intersection. Therefore, we have
\begin{equation*}
\begin{aligned}
  \mathcal{M}_{a,b}(\mathcal{F}) &= \bigcap_{t \in [a,b]} \mathrm{im}(\mathcal{M}_{t \to b}) \\
  &= \bigcap_{t \in [a,b]} \mathbb{K}\left[ \mathrm{im}\left( f_{t,b}^{\sharp} \right) \right] \\
  &= \mathbb{K}\left[ \bigcap_{t \in [a,b]} \mathrm{im}\left( f_{t,b}^{\sharp} \right) \right] \\
  &= \mathbb{K}[\Sym(\mathcal{F})_{a,b}].
\end{aligned}
\end{equation*}
Moreover, we can endow $\mathcal{M}_{a,b}(\mathcal{F})$ with the group algebra structure of $\mathbb{K}[\Sym(\mathcal{F})_{a,b}]$.
\end{proof}

It is crucial to note that $ \mathcal{M}_{a,b}(\mathcal{F}) $ is not necessarily a commutative algebra. This algebraic tool not only captures the persistence of symmetries but also sheds light on their dynamic interactions, which are inherently noncommutative. As a result, it enriches the study of persistent structures by offering a different perspective on the relationship between symmetry and persistence in topological data analysis.

\subsection{Polybarcode of persistence symmetries}

In Euclidean space, a symmetry corresponds to a linear transformation, which can be expressed in matrix form.
When considering $\mathbf{M}(\mathcal{F}) = \bigoplus_{a} \mathcal{M}\mathcal{F}_{a}$, we can always obtain the symmetries in $\Sym(\mathcal{F})$, where each symmetry has a birth time and a death time. Here, the death time can be infinity.

In this context, it may happen that two symmetries have different birth and death times but correspond to the same linear transformation. In such cases, we regard these two symmetries as the same one (see Example \ref{example:unique}). Hence, the existence time, or lifespan, of a symmetry is not necessarily an interval, but rather a subset of the real numbers $\mathbb{R}$. For a symmetry $\pi$ in $\Sym(\mathcal{F})$, let $\delta_{\pi}(t)=\left\{
                                            \begin{array}{ll}
                                              1, & \hbox{if $\pi$ persists at time $t$;} \\
                                              0, & \hbox{otherwise.}
                                            \end{array}
                                          \right.
$ We denote $I(\pi) = \{t\in \mathbb{R}\mid\delta_{\pi}(t)=1\}$. Then $I(\pi)$ represents the lifespan of the symmetry $\pi$. For any element $\pi \in \Iso(M)$, the set $I(\pi)$ is non-empty if and only if $\pi$ is a symmetry in $\Sym(\mathcal{F})$.

\begin{definition}
Let $\mathcal{F} \colon (\mathbb{R}, \leq) \to \mathcal{S}_{n}(M)$ be a persistence $n$-configuration. For any symmetry $\pi \in \Sym(\mathcal{F})$, the \textit{polybar} of $\pi$ is the set $I(\pi) = \{t\in \mathbb{R}\mid\pi(\mathcal{F}_{t}):=\mathcal{F}_{t}\}$.
\end{definition}

Recall that the $n$-point configuration space of $M$ is the space of all the $n$-configurations, denoted by $\mathrm{Conf}_n(M)$.
A persistence $n$-configuration $\mathcal{F} \colon (\mathbb{R}, \leq) \to \mathcal{S}_{n}(M)$ is \textit{continuous} if the map $f:\mathbb{R}\to \mathrm{Conf}_n(M)$ given by $f(t)=\mathcal{F}_t$ is a continuous map. From now on, we assume that the persistence $n$-configuration under consideration is continuous.

\begin{proposition}\label{proposition:closed}
For each symmetry $\pi$ in $\Sym(\mathcal{F})$, the set $I(\pi)$ is a closed subset of $\mathbb{R}$.
\end{proposition}

\begin{proof}
Let $f:\mathbb{R} \to \mathrm{Conf}_n(M)$ be a continuous map defined by $f(t) = \mathcal{F}_t$. The symmetry $\pi$ survives at time $t$ if and only if $\pi(f(t)) = f(t)$. Thus, the lifespan of $\pi$ is given by
\begin{equation*}
I(\pi) = \{t \in \mathbb{R} \mid \pi(f(t)) = f(t)\}.
\end{equation*}
Additionally, we define a function
\begin{equation*}
    g: \mathbb{R} \to \mathbb{R}, \quad g(t) = d_{H}(f(t), \pi(f(t))).
\end{equation*}
Here, $d_{H}$ is the Hausdorff distance induced by the metric $d$ on $M$. Since the Hausdorff distance $d_H$ is continuous with respect to the topology induced by $d$, and both $f$ and $\pi$ are continuous maps, it follows that $g$ is a continuous function. Now, consider the set
\begin{equation*}
    I(\pi) = g^{-1}(\{0\}) = \{ t \in \mathbb{R} \mid g(t) = 0 \}.
\end{equation*}
Since $g$ is continuous and $\{0\}$ is a closed subset of $\mathbb{R}$, the preimage $I(\pi)$ is also a closed subset of $\mathbb{R}$.
\end{proof}

\begin{remark}
If we consider a persistence $n$-configuration $\mathcal{F} \colon (\mathbb{Z}, \leq) \to \mathcal{S}_n(M)$ with a discrete filtration parameter, then the corresponding polybarcode is a discrete set of points and hence a closed subset. In practical applications, discrete filtration parameters are often used to approximate the general setting. Therefore, in the subsequent discussion, we typically assume that a persistence $n$-configuration is continuous. This assumption does not affect the computation or application of polybarcodes.
\end{remark}

\begin{proposition}
Fix a persistence $n$-configuration $\mathcal{F}:(\mathbb{R},\leq)\to\mathcal{S}_{n}(M)$. Let $\pi$ and $\varpi$ be symmetries in $\Sym(\mathcal{F})$. Then we have
\begin{equation*}
  I(\pi) \cap I(\varpi) \subseteq I(\pi\cdot \varpi).
\end{equation*}
\end{proposition}

\begin{proof}
Observe that at any fixed time $t$, the collection of persistent symmetries forms a group $\Sym(\mathcal{F}_{t})$. This implies that if $\pi, \varpi \in \Sym(\mathcal{F}_{t})$, then the product $\pi \cdot \varpi$ also belongs to $\Sym(\mathcal{F}_{t})$. Consequently, we obtain
\begin{equation*}
  I(\pi) \cap I(\varpi) \subseteq I(\pi\cdot \varpi).
\end{equation*}
which is the desired result.
\end{proof}
It is worth noting that even if the symmetries $\pi$ and $\varpi$ do not exist at a fixed time $t$, their product $\pi \cdot \varpi$ may still exist at $t$.

\begin{proposition}\label{proposition:subset}
Fix a persistence $n$-configuration $\mathcal{F}:(\mathbb{R},\leq)\to\mathcal{S}_{n}(M)$. Let $\pi$ be a symmetry in $\Sym(\mathcal{F})$. Then, for any integer $k\geq 2$, we have
\begin{equation*}
  I(\pi) = I(\pi^{-1}) \subseteq I(\pi^{k}).
\end{equation*}
\end{proposition}

\begin{proof}
At any fixed time $t$, the symmetry group $\Sym(\mathcal{F}_{t})$ always contains at least the identity element $e$, which corresponds to the identity transformation in the underlying metric space. If $\pi \in \Sym(\mathcal{F}_{t})$, then its inverse $\pi^{-1}$ and higher-order powers $\pi^{k}$ (for any integer $k\geq 2$) must also belong to $\Sym(\mathcal{F}_{t})$. Consequently, we obtain the inclusions
\begin{equation*}
  I(\pi) \subseteq I(\pi^{-1}), \quad I(\pi) \subseteq I(\pi^{k}).
\end{equation*}
On the other hand, we observe that
\begin{equation*}
  I(\pi^{-1}) \subseteq I((\pi^{-1})^{-1}) = I(\pi).
\end{equation*}
It follows that $I(\pi) = I(\pi^{-1})$. The desired result follows.
\end{proof}

\begin{definition}
The \textit{polybarcode} of a persistence $n$-configuration $\mathcal{F}: (\mathbb{R}, \leq) \to \mathcal{S}_{n}(M)$ is the collection $\mathcal{B}(\mathcal{F}) := \{I(\pi)\}_{\pi \in \Sym(\mathcal{F})}$ of subsets of $\mathbb{R}$.
\end{definition}

\begin{example}\label{example:unique}
Consider the dynamical system $\{X_{t}\}_{t\geq 0}$ depicted in Figure \ref{figure:poly_bar}\textit{a}. The positions of this system at time $t$ are given by $X_{t}=\{x_{1}(t),x_{2}(t),x_{3}(t),x_{4}(t)\}$, where each $x_i(t)$ is determined by the following equations
\begin{equation}\label{eq:exam}
\begin{split}
    & x_{1}(t)=(1+q(t),2+q(t-1)),\quad x_{2}(t)=(-1-q(t),2+q(t-1)), \\
    & x_{3}(t)=(1+q(t),-2-q(t-1)),\quad x_{4}(t)=(-1-q(t),-2-q(t-1)).
\end{split}
\end{equation}
Here, $q(t)$ is a piecewise function defined as
\begin{equation*}
  q(t)= \left\{
          \begin{array}{ll}
            t, & \hbox{if } t\leq 1; \\
            1, & \hbox{if } 1\leq t\leq 2; \\
            t-1, & \hbox{if } t\geq 2.
          \end{array}
        \right.
\end{equation*}

\begin{figure}[h]
    \centering
    \definecolor{axisgray}{gray}{0.6}

    \begin{tikzpicture}[>=Stealth]

    \draw[->, color=axisgray, line width=0.4pt] (-2, 0) -- (2, 0) node[right, font=\small, text=black] {$x$};
    \draw[->, color=axisgray, line width=0.4pt] (0, -2) -- (0, 2) node[above, font=\small, text=black] {$y$};

    \draw[line width=1.2pt, color=mainblue] (0.4, 0.4) -- (0.8, 0.8) -- (0.8, 1.2) -- (1.2, 1.2) -- (1.8, 1.8);
    \draw[line width=1.2pt, color=mainblue] (-0.4, 0.4) -- (-0.8, 0.8) -- (-0.8, 1.2) -- (-1.2, 1.2) -- (-1.8, 1.8);
    \draw[line width=1.2pt, color=mainblue] (0.4, -0.4) -- (0.8, -0.8) -- (0.8, -1.2) -- (1.2, -1.2) -- (1.8, -1.8);
    \draw[line width=1.2pt, color=mainblue] (-0.4, -0.4) -- (-0.8, -0.8) -- (-0.8, -1.2) -- (-1.2, -1.2) -- (-1.8, -1.8);

    \foreach \pt in {(0.4,0.4), (0.8,0.8), (0.8,1.2), (1.2,1.2), (1.8,1.8), (-0.4,0.4), (-0.8,0.8), (-0.8,1.2), (-1.2,1.2), (-1.8,1.8), (0.4,-0.4), (0.8,-0.8), (0.8,-1.2), (1.2,-1.2), (1.8,-1.8), (-0.4,-0.4), (-0.8,-0.8), (-0.8,-1.2), (-1.2,-1.2), (-1.8,-1.8)}
        \fill[mainblue] \pt circle (1.6pt);

    \node at (-2.4, 1.97) {\textit{a}};

    \draw[->, color=axisgray, line width=0.4pt] (5, -1) -- (11.6, -1) node[right, font=\small, text=black] {$t$};
    \draw[->, color=axisgray, line width=0.4pt] (5, -1.55) -- (5, 1.75) node[above, font=\small, text=black] {Symmetries};
    
    \node at (3.24, 1.97) {\textit{b}};

    \foreach \x/\label in {6.1/1, 8.3/3, 10.5/5} 
        \draw[color=axisgray, line width=0.4pt] (\x, -1) -- (\x, -1.22) node[below, font=\small, text=black] {$\label$};

    \draw[mainblue, very thick, ->] (5, 0.65) -- (11.05, 0.65);
    \draw[mainblue, very thick, ->] (5, 0.98) -- (11.05, 0.98);
    \draw[mainblue, very thick, ->] (5, 1.31) -- (11.05, 1.31);
    \node[mainblue, right, font=\small] at (4.12, 1.31) {$\pi$};
    \node[mainblue, right, font=\small] at (4.12, 0.98) {$\tau^2$};
    \node[mainblue, right, font=\small] at (4.12, 0.65) {$\pi\tau^2$};

    \draw[red!80!black, very thick] (5, -0.67) -- (6.1, -0.67);
    \draw[red!80!black, very thick, ->] (8.3, -0.67) -- (11.05, -0.67);

    \draw[red!80!black, very thick] (5, -0.34) -- (6.1, -0.34);
    \draw[red!80!black, very thick, ->] (8.3, -0.34) -- (11.05, -0.34);

    \draw[red!80!black, very thick] (5, -0.01) -- (6.1, -0.01);
    \draw[red!80!black, very thick, ->] (8.3, -0.01) -- (11.05, -0.01);

    \draw[red!80!black, very thick] (5, 0.32) -- (6.1, 0.32);
    \draw[red!80!black, very thick, ->] (8.3, 0.32) -- (11.05, 0.32);

    \node[red!80!black, right, font=\small] at (4.12, 0.32) {$\tau$};
    \node[red!80!black, right, font=\small] at (4.12, -0.01) {$\tau^3$};
    \node[red!80!black, right, font=\small] at (4.12, -0.34) {$\pi\tau$};
    \node[red!80!black, right, font=\small] at (4.12, -0.67) {$\pi\tau^3$};

    \end{tikzpicture}
    \caption{\textit{a} Illustration of the trajectory of $\{X_{t}\}_{t \geq 0}$ in Example \ref{example:unique}. \textit{b} The polybarcode of persistent symmetries in the dynamical system.}\label{figure:poly_bar}
\end{figure}

In this example, the four points $x_1, x_2, x_3, x_4$ change their positions over time. For $t \in [0,1]$, they form the vertices of a square. In the interval $t \in (1,3)$, their configuration transforms into a rectangle with unequal side lengths. Finally, for $t \in [3,+\infty)$, the points again form the vertices of a square. Consequently, the symmetry group of this dynamical system is the dihedral group
\begin{equation*}
D_4 = \langle \tau, \pi \mid \tau^4 = e, \pi^2 = e, \pi\tau = \tau^{-1}\pi \rangle.
\end{equation*}
for $t \in [0,1]$ and $t \geq 3$. However, for $t \in (1,3)$, the symmetry group is a subgroup of $D_4$ generated by $\tau^2$ and $\pi$, consisting of the elements $e, \tau^{2}, \pi, \pi\tau^{2}$.
This subgroup is abelian and isomorphic to $\mathbb{Z}_{2}\times \mathbb{Z}_{2}$. Figure \ref{figure:poly_bar}\textit{b} illustrates the polybarcode representing the persistent symmetries of this dynamical system. In this polybarcode, the three symmetries $\tau^{2}, \pi, \pi\tau^{2}$ persist over the entire time range $[0,+\infty)$. In contrast, the symmetries $\tau, \tau^{3}, \pi\tau, \pi\tau^{3}$ correspond to polybars given by the union of two disjoint intervals, $[0,1] \cup [3,+\infty)$. This occurs because at $t=3$, the re-emerging symmetry transformations, such as $\tau$, correspond to the same linear transformations $\mathbb{R}^{2} \to \mathbb{R}^{2}$ as in the time interval $t \in [0,1]$.

Furthermore, analyzing the multiplication of persistent symmetries, we observe that the product $\tau \cdot \tau$ exists at $t \in (1,3)$, while $\tau$ itself does not exist in this interval. Notably, we have the inclusion relation
\begin{equation*}
  I(\tau)=[0,1]\cup [3,+\infty)\subseteq I(\tau^{2})=[0,+\infty),
\end{equation*}
which is consistent with Proposition \ref{proposition:subset}, demonstrating a strict inclusion case.

It is worth noting that the dynamical system $\{X_{t}\}_{t\geq 0}$ can be extended to $\{X_{t}\}_{t\in \mathbb{R}}$ by defining $X_{t} = X_{0}$ for all $t < 0$. This extension preserves the initial configuration of the system for negative time values, and $\{X_{t}\}_{t\in \mathbb{R}}$ can be viewed as a persistent 4-configuration, i.e., a map $\mathcal{F}: (\mathbb{R}, \leq) \to \mathcal{S}_{4}(M)$.
\end{example} 

%% file: metrics_symmetries.tex
\section{Metrics on the polybarcodes of symmetries}\label{section:polybarcodes_metrics}

In the previous section, we introduced the symmetry barcode of persistent configurations, which shares a similar structural form with the classical barcode in persistent homology. Standard metrics for comparing such barcodes include the bottleneck distance, Wasserstein distance, matching distance, among others~\cite{bauer2014induced,cohen2005stability,mileyko2011probability}. However, the polybarcode differs fundamentally, as polybars are not easily represented as points in a geometric space. In this section, we introduce several potential distance measures for polybarcodes, serving as a foundation for the stability analysis in later sections.

\subsection{Symmetric difference distance}

\begin{definition}
Let $X$ and $Y$ be Borel subsets of $\mathbb{R}$. The \textit{symmetric difference distance} between $X$ and $Y$ is defined as
\begin{equation*}
    d_{S}(X,Y) := m(X \Delta Y) = m(X \setminus Y) + m(Y \setminus X),
\end{equation*}
where $m(A)$ denotes the Lebesgue measure of a set $A$, and the set difference is given by $X \setminus Y = X - (X \cap Y)$. If either $m(X \setminus Y)$ or $m(Y \setminus X)$ is infinite, then the distance $d_{S}(X,Y)$ can take the value $\infty$.
\end{definition}

If $m(X \cap Y) < \infty$, then we have the equalities
\begin{equation*}
    m(X \setminus Y) = m(X) - m(X \cap Y), \quad m(Y \setminus X) = m(Y) - m(X \cap Y),
\end{equation*}
leading to the alternative expression
\begin{equation*}
    d_{S}(X,Y) = m(X) + m(Y) - 2 m(X \cap Y).
\end{equation*}
In particular, $d_{S}(X,Y) = 0$ if and only if $m(X) = m(Y) = m(X \cap Y)$, which implies that the symmetric difference $X \Delta Y$ is a null set.

\begin{proposition}
The symmetric difference distance defined above is a pseudometric.
\end{proposition}
\begin{proof}
It is straightforward to verify that $d_{S}(X,Y)$ satisfies non-negativity and symmetry. We now prove the triangle inequality. Given Borel subsets $X, Y, Z$, we need to show that
\begin{equation*}
    d_{S}(X,Y) + d_{S}(Y,Z) \geq d_{S}(X,Z).
\end{equation*}
For any $x \in X \setminus Z$, we observe that if $x \in Y$, then $x \in Y \setminus Z$; whereas if $x \notin Y$, then $x \in X \setminus Y$. This implies the inclusion
\begin{equation*}
    X \setminus Z \subseteq (X \setminus Y) \cup (Y \setminus Z).
\end{equation*}
Taking measures on both sides, we obtain
\begin{equation*}
    m(X \setminus Z) \leq m(X \setminus Y) + m(Y \setminus Z).
\end{equation*}
Similarly, we have
\begin{equation*}
    m(Z \setminus X) \leq m(Y \setminus X) + m(Z \setminus Y).
\end{equation*}
Adding these two inequalities yields
\begin{equation*}
    m(X \Delta Z) \leq m(X \Delta Y) + m(Y \Delta Z),
\end{equation*}
which proves the triangle inequality.
\end{proof}

\begin{example}\label{example:symmetry_distance}
Consider the case where $X$ and $Y$ are closed intervals on the real line, given by $[a, b]$ and $[c, d]$, respectively. If $X$ and $Y$ intersect non-trivially, their symmetric difference distance is given by
\begin{equation*}
  d_{S}(X,Y) = |a - c| + |b - d|.
\end{equation*}
If $X$ and $Y$ are disjoint, assuming without loss of generality that $b \leq c$, their symmetric difference distance is
\begin{equation*}
  d_{S}(X,Y) = b + d - a - c.
\end{equation*}
In this disjoint case, we have the inequality
\begin{equation*}
  b + d - a - c \leq |a - c| + |b - d|.
\end{equation*}
If we interpret the intervals $[a, b]$ and $[c, d]$ as points $x = (a, b)$ and $y = (c, d)$ in the plane, their corresponding $\ell_1$ distance is given by
\begin{equation*}
  \|x - y\|_1 = |a - c| + |b - d|.
\end{equation*}
Thus, we obtain the bound
\begin{equation*}
  d_{S}(X,Y) \leq \|x - y\|_1.
\end{equation*}
\end{example}

\begin{definition}
Suppose $\mathcal{F},\mathcal{G}:(\mathbb{R},\leq)\to\mathcal{S}_{n}(M)$ are persistence configurations. The \textit{symmetric difference distance} between their polybarcodes is defined as
\begin{equation*}
  d_{S}(\mathcal{B}(\mathcal{F}),\mathcal{B}(\mathcal{G})) := \sum\limits_{\pi\in \Iso(M)}d_{S}(I^{\mathcal{F}}(\pi),I^{\mathcal{G}}(\pi)).
\end{equation*}
\end{definition}

By Proposition \ref{proposition:closed}, the sets $I^{\mathcal{F}}(\pi)$ and $I^{\mathcal{G}}(\pi)$ are closed, ensuring that the above definition is well-defined. Since the sum is taken over an infinite number of terms, the symmetric difference distance $d_{S}(\mathcal{B}(\mathcal{F}),\mathcal{B}(\mathcal{G}))$ may be infinite, even if each individual term $d_{S}(I^{\mathcal{F}}(\pi), I^{\mathcal{G}}(\pi))$ is finite.

The bottleneck distance between two finite point sets $X$ and $Y$ in a metric space, or equivalently between two persistence diagrams, is defined as
\begin{equation*}
d_B(X, Y) = \inf_{\gamma} \sup_{x \in X} \|x-\gamma(x)\|_{\infty},
\end{equation*}
where $\|x-y\|_{\infty}$ denotes the $\ell_{\infty}$ distance between points $x$ and $y$. Here, $\gamma: X \to Y$ is a matching, which can be a bijection or a partial bijection, allowing points to be matched to the diagonal in the case of persistence diagrams.

A fundamental aspect in the application of the bottleneck distance is the matching between individual points. In persistent homology, each barcode consists of intervals represented by birth-death pairs $(b, d)$, which can be viewed as points in the extended plane. In contrast, persistent symmetries are captured by polybars, which are typically piecewise-defined and represented as unions of multiple intervals. In such cases, identifying them with single points is no longer meaningful. Moreover, for a fixed metric space, each symmetry corresponds to a prescribed isometric transformation. As a result, there is no need to search for an optimal matching between two polybarcodes; the symmetries themselves canonically determine the matching. As the same symmetry may give rise to different polybars across different persistence configurations, comparing polybars instead of individual intervals provides a more coherent measure of variation.

\subsection{Expansion distance}

Let $X$ and $Y$ be two closed sets of $\mathbb{R}$, each consisting of a finite union of disjoint intervals of the form $[a_i, b_i]$, $(-\infty, b_i]$, or $[a_i, \infty)$. Here, we allow the case where $a_i = b_i$. For simplicity, we may also express intervals of the form $[a_i, \infty)$ as $[a_i, b_i]$ by permitting $b_i$ to take infinite values. We can represent $X$ as a union of intervals in the form
\begin{equation*}
  X = \bigcup\limits_{i=1}^{m} A_i,
\end{equation*}
where the intervals $A_1, A_2, \dots, A_m$ are ordered from left to right along the real number line, meaning that the upper bound of $A_i$ is strictly less than the lower bound of $A_{i+1}$ for all $1 \leq i \leq m - 1$. From now on, sets that are unions of disjoint intervals are always assumed to be in the above ordered form.

For an interval $A = [a, b]$, its $\varepsilon$-expansion is defined as $A^{\varepsilon} = [a - \varepsilon, b + \varepsilon]$. In particular, the $\varepsilon$-expansion of $(-\infty, b_i]$ is $(-\infty, b_i + \varepsilon]$, and that of $[a_i, \infty)$ is $[a_i - \varepsilon, \infty)$.

\begin{definition}
The \textit{expansion distance} between $X$ and $Y$ is defined as
\begin{equation*}
    d_E(X, Y) := \max\limits_{1 \leq i \leq m} \inf \{ \varepsilon \geq 0 \mid A_i \subseteq B_i^{\varepsilon}, B_i \subseteq A_i^{\varepsilon} \},
\end{equation*}
where $X = \bigcup\limits_{i=1}^{m} A_i$ and $Y = \bigcup\limits_{i=1}^{n} B_i$, and each $A_i$ and $B_i$ is an interval. If $m \neq n$, we set $d_E(X, Y) = \infty$.
\end{definition}

\begin{example}\label{example:expansion_distance}
Consider the case where $X$ and $Y$ are two closed intervals on the real line, $[a, b]$ and $[c, d]$, respectively. The expansion distance between them is given by
\begin{equation*}
  d_{E}(X, Y) = \max\{|a - c|, |b - d|\}.
\end{equation*}
If we interpret $[a, b]$ and $[c, d]$ as points $(a, b)$ and $(c, d)$ in the plane, then their corresponding $\ell_{\infty}$ distance is exactly $d_{E}(X, Y)$. In this case, the expansion distance coincides with the bottleneck distance, that is,
\begin{equation*}
  d_{E}(X, Y) = d_{B}(X, Y).
\end{equation*}
\end{example}

The expansion distance measures the proximity between the sets $X$ and $Y$. It is straightforward to verify that $d_E(X, Y) = 0$ if and only if $X = Y$. Moreover, one can show that the expansion distance defines a metric.

\begin{proposition}
The expansion distance defined above is a metric.
\end{proposition}
\begin{proof}
We first note that the definition of the expansion distance satisfies symmetry
\begin{equation*}
d_{E}(X, Y) = d_{E}(Y, X),
\end{equation*}
and positive definiteness. We now prove the triangle inequality
\begin{equation*}
d_{E}(X, Z) \leq d_{E}(X, Y) + d_{E}(Y, Z).
\end{equation*}
Clearly, if $d_E(X, Z) = \infty$, then at least one of $d_E(X, Y)$ or $d_E(Y, Z)$ is infinite. Now, suppose that
$X = \bigcup\limits_{i=1}^{n} A_i$, $Y = \bigcup\limits_{i=1}^{n} B_i$, and $Z = \bigcup\limits_{i=1}^{n} C_i$ all have finite expansion distances pairwise.
Let $\epsilon_1 = d_{E}(X,Y)$ and $\epsilon_2 = d_{E}(Y,Z)$. By definition, we have
\begin{equation*}
  A_{i} \subseteq B_{i}^{\epsilon_1} \quad \text{and} \quad B_{i} \subseteq A_{i}^{\epsilon_1},
\end{equation*}
as well as
\begin{equation*}
  B_{i} \subseteq C_{i}^{\epsilon_2} \quad \text{and} \quad C_{i} \subseteq B_{i}^{\epsilon_2},
\end{equation*}
for any $1\leq i\leq n$. Using the additive property of the interval expansion, we have
\begin{equation*}
  A_{i}\subseteq B_{i}^{\epsilon_1}  \subseteq C_{i}^{\epsilon_1+\epsilon_2},\quad 1\leq i\leq n.
\end{equation*}
Similarly, one has
\begin{equation*}
  C_{i}\subseteq B_{i}^{\epsilon_2}  \subseteq A_{i}^{\epsilon_1+\epsilon_2},\quad 1\leq i\leq n.
\end{equation*}
Hence, by the definition of the expansion distance, we have
\begin{equation*}
  d_{E}(X,Z) \leq  \epsilon_1 + \epsilon_2 = d_{E}(X,Y) + d_{E}(Y,Z).
\end{equation*}
This completes the proof of the triangle inequality.
\end{proof}

One motivation for introducing this distance is that, in previous studies of persistent symmetries, polybarcodes often arise as unions of intervals. This metric thus provides a natural tool for quantifying the distance between different persistent symmetries.

\begin{definition}
Let $X= \bigcup\limits_{i=1}^{m} [a_i, b_i]$ and $Y=\bigcup\limits_{i=1}^{n} [c_i, d_i]$ be finite unions of disjoint intervals.
The \textit{left expansion distance} between $X$ and $Y$ is defined as
\begin{equation*}
d_L(X, Y) := \max\limits_{1 \leq i \leq m} |a_i - c_i|, \quad \text{if } m = n,
\end{equation*}
and $d_L(X, Y)=\infty$ otherwise. In particular, if $a_{1}=c_{1}=-\infty$, we make the convention $|a_1 - c_1|=0$.
\end{definition}

\begin{proposition}
The left expansion distance defined above is a pseudometric.
\end{proposition}
\begin{proof}
The left expansion distance satisfies non-negativity and symmetry. However, it does not, in general, satisfy strict definiteness; that is, $d_L(X, Y) = 0$ does not necessarily imply $X = Y$.

Now, we will prove the triangle inequality.
The case where the distance is infinite is straightforward to verify, so we focus on the case where the distance is finite.
Let $X= \bigcup\limits_{i=1}^{n} [a_i, b_i]$, $Y=\bigcup\limits_{i=1}^{n} [c_i, d_i]$, and $Z=\bigcup\limits_{i=1}^{n} [e_i, f_i]$. Our goal is to verify the inequality
\begin{equation*}
  d_L(X, Z) \leq d_L(X, Y) + d_L(Y, Z),
\end{equation*}
which is equivalent to proving
\begin{equation*}
  \max\limits_{1 \leq i \leq n} |a_i - e_i| \leq \max\limits_{1 \leq i \leq n} |a_i - c_i| + \max\limits_{1 \leq i \leq n} |c_i - e_i|.
\end{equation*}
Without loss of generality, assume that
\begin{equation*}
  \max\limits_{1 \leq i \leq n} |a_i - e_i| = |a_k - e_k|
\end{equation*}
for some $1 \leq k \leq n$. Then, by the triangle inequality for absolute values, we obtain
\begin{equation*}
  |a_k - e_k| \leq |a_k - c_k| + |c_k - e_k|.
\end{equation*}
It follows that
\begin{equation*}
  |a_k - e_k| \leq \max\limits_{1 \leq i \leq n} |a_i - c_i| + \max\limits_{1 \leq i \leq n} |c_i - e_i|.
\end{equation*}
Thus, the desired result holds.
\end{proof}

\begin{proposition}
$d_{L}(X,Y)\leq d_{E}(X,Y)$.
\end{proposition}
\begin{proof}
Let $X = \bigcup\limits_{i=1}^{m} [a_i, b_i]$ and $Y = \bigcup\limits_{i=1}^{n} [c_i, d_i]$. If $m \neq n$, both $d_L(X, Y)$ and $d_E(X, Y)$ are infinite. If $m = n$, let $\varepsilon = d_E(X, Y)$. By the definition of the expansion distance, for each $i$, we know that
\begin{equation*}
  [a_i, b_i] \subseteq [c_i - \varepsilon, d_i + \varepsilon], \quad [c_i, d_i] \subseteq [a_i - \varepsilon, b_i + \varepsilon].
\end{equation*}
This implies that $|a_i - c_i| \leq \varepsilon$ for all $1 \leq i \leq n$. Therefore, we have
\begin{equation*}
  \max\limits_{1 \leq i \leq m} |a_i - c_i| \leq \varepsilon.
\end{equation*}
Thus, we conclude that $d_L(X, Y) \leq d_E(X, Y)$, as required.
\end{proof}

\begin{definition}
Let $X= \bigcup\limits_{i=1}^{m} [a_i, b_i]$ and $Y=\bigcup\limits_{i=1}^{n} [c_i, d_i]$ be finite unions of disjoint intervals.
The \textit{matching symmetric difference distance} between $X$ and $Y$ is defined as
\begin{equation*}
  \widetilde{d}_{S}(X,Y) := \frac{1}{m}\sum\limits_{i=1}^{m}d_{S}([a_i, b_i],[c_i, d_i]),\quad \text{if } m = n,
\end{equation*}
and $\widetilde{d}_{S}(X, Y)=\infty$ otherwise.
\end{definition}

\begin{proposition}
$\widetilde{d}_{S}(X,Y)\leq 2d_{E}(X,Y)$.
\end{proposition}
\begin{proof}
By Example \ref{example:symmetry_distance}, we have
\begin{equation*}
  d_{S}([a_i, b_i],[c_i, d_i]) \leq |a_{i}-c_{i}|+|b_{i}-d_{i}|\leq 2 \max\{|a_{i}-c_{i}|,|b_{i}-d_{i}|\}.
\end{equation*}
By the definition of expansion distance, it follows that
\begin{equation*}
  \max\{|a_{i}-c_{i}|,|b_{i}-d_{i}|\} \leq d_{E}(X,Y).
\end{equation*}
Thus, we obtain
\begin{equation*}
  \widetilde{d}_{S}(X,Y)\leq 2\frac{1}{m}\sum\limits_{i=1}^{m}d_{E}(X,Y) = 2d_{E}(X,Y).
\end{equation*}
This completes the proof.
\end{proof}

\subsection{Interleaving distance}
Now, let us recall the interleaving distance. Let $T_x: (\mathbb{R}, \leq) \to (\mathbb{R}, \leq)$ be the translation functor defined by $T_x(a) = a + x$ for all $a, x \in \mathbb{R}$. This induces an endofunctor $\Sigma^x: \mathfrak{C}^{\mathbb{R}} \to \mathfrak{C}^{\mathbb{R}}$ on the category of persistence objects, given by $(\Sigma^x \mathcal{F})(a) = \mathcal{F}(a + x)$. The natural transformation $\Sigma^x|_{\mathcal{F}}: \mathcal{F} \Rightarrow \Sigma^x \mathcal{F}$ defines a morphism in $\mathfrak{C}^{\mathbb{R}}$, which we denote simply by $\Sigma^x$ when no ambiguity arises.

An \textit{$\varepsilon$-interleaving} between two persistence objects $\mathcal{F}, \mathcal{G}: (\mathbb{R}, \leq) \to \mathfrak{C}$ consists of a pair of morphisms $\phi: \mathcal{F} \to \Sigma^\varepsilon \mathcal{G}$ and $\psi: \mathcal{G} \to \Sigma^\varepsilon \mathcal{F}$ in $\mathfrak{C}^{\mathbb{R}}$, satisfying the commutativity conditions
\begin{equation*}
(\Sigma^\varepsilon \psi) \circ \phi = \Sigma^{2\varepsilon}|_{\mathcal{F}}, \quad (\Sigma^\varepsilon \phi) \circ \psi = \Sigma^{2\varepsilon}|_{\mathcal{G}}.
\end{equation*}
These conditions are expressed by the following commutative diagrams:
\begin{equation*}
\xymatrix@=0.6cm{
  & \Sigma^{\varepsilon} \mathcal{G} \ar[rd]^{\Sigma^{\varepsilon} \psi} & \\
  \mathcal{F} \ar[ru]^{\phi} \ar[rr]^{\Sigma^{2\varepsilon}|_{\mathcal{F}}} && \Sigma^{2\varepsilon} \mathcal{F},
}
\qquad \qquad
\xymatrix@=0.6cm{
  & \Sigma^{\varepsilon} \mathcal{F} \ar[rd]^{\Sigma^{\varepsilon} \phi} & \\
  \mathcal{G} \ar[ru]^{\psi} \ar[rr]^{\Sigma^{2\varepsilon}|_{\mathcal{G}}} && \Sigma^{2\varepsilon} \mathcal{G}.
}
\end{equation*}
In the special case $\varepsilon = 0$, the above equations reduce to $\psi \circ \phi = \mathrm{id}_{\mathcal{F}}$ and $\phi \circ \psi = \mathrm{id}_{\mathcal{G}}$, meaning that $\mathcal{F}$ and $\mathcal{G}$ are isomorphic objects in $\mathfrak{C}^{\mathbb{R}}$.

The \textit{interleaving distance} \cite{bubenik2014categorification,chazal2009proximity} between two persistence objects $\mathcal{F}, \mathcal{G}: (\mathbb{R}, \leq) \to \mathfrak{C}$ is defined as
\begin{equation*}
  d_{I}(\mathcal{F},\mathcal{G})=\inf\{\varepsilon\geq 0\mid \text{$\mathcal{F}$ and $\mathcal{G}$ are $\varepsilon$-interleaved}\}.
\end{equation*}
The interleaving distance is an extended pseudometric: ``extended'' means that the distance can be $\infty$, while ``pseudo'' indicates that $d_I(\mathcal{F}, \mathcal{G}) = 0$ does not necessarily imply $\mathcal{F} \cong \mathcal{G}$.

Consider the category $\mathbf{Int}(\mathbb{R})$, where an object $X$ is a finite union of disjoint intervals of the form $[a_i, b_i]$, $(-\infty, b_i]$, or $[a_i, \infty)$ for some real numbers $a_i, b_i$. A morphism $\phi: X \to Y$ in this category is a continuous, monotonically increasing function that maps the intervals of $X$ to those of $Y$ in increasing order. More specifically, for $X=\bigcup\limits_{i=1}^{m}A_{i}$ and $Y=\bigcup\limits_{j=1}^{n}B_{j}$, a morphism $\phi: X \to Y$ in the category $\mathbf{Int}(\mathbb{R})$ satisfies the following conditions:
\begin{itemize}
  \item $\phi$ is a monotonically increasing function;
  \item $m \leq n$, and $\phi$ maps each interval $A_i$ to $B_i$ for $1\leq i\leq m$.
\end{itemize}

Given an object $X$ in the category $\mathbf{Int}(\mathbb{R})$, we define a persistence object $Fr_{X}:(\mathbb{R},\leq)\to \mathbf{Int}(\mathbb{R})$ by
\begin{equation*}
  Fr_{X}(t) = \{x\in X\mid x\leq t\}.
\end{equation*}
Now, consider two persistence objects $Fr_{X},Fr_{Y}:(\mathbb{R},\leq)\to \mathbf{Int}(\mathbb{R})$. We say that $Fr_{X}$ and $Fr_{Y}$ are \textit{$\varepsilon$-interleaved} if there exist morphisms
\begin{equation*}
\phi: Fr_{X} \to \Sigma^\varepsilon Fr_{Y}, \quad \psi: Fr_{Y} \to \Sigma^\varepsilon Fr_{X}
\end{equation*}
in the category $\mathbf{Int}(\mathbb{R})^{\mathbb{R}}$ such that the following commutativity conditions hold:
\begin{equation*}
  (\Sigma^\varepsilon \psi) \circ \phi = \Sigma^{2\varepsilon}|_{Fr_{X}}, \quad (\Sigma^\varepsilon \phi) \circ \psi = \Sigma^{2\varepsilon}|_{Fr_{Y}}.
\end{equation*}

\begin{proposition}\label{proposition:polybar_distance}
Given two persistence objects $Fr_{X},Fr_{Y}:(\mathbb{R},\leq)\to \mathbf{Int}(\mathbb{R})$, we have
\begin{equation*}
  d_{I}(Fr_{X},Fr_{Y}) = d_{L}(X,Y).
\end{equation*}
\end{proposition}

\begin{proof}
If $d_{I}(Fr_{X},Fr_{Y})=\infty$, then there are only two possible cases: either $m \neq n$, or one of $|a_{1}|$ and $|b_{1}|$ is $-\infty$ while the other is finite. In both cases, $d_{L}(X,Y)$ is also infinite.

Now, consider the case where $d_{I}(Fr_{X},Fr_{Y}) < \infty$, and let $\varepsilon = d_{I}(Fr_{X},Fr_{Y})$. If $Fr_{X}$ and $Fr_{Y}$ are $\varepsilon$-interleaved, then for any $t \in \mathbb{R}$, we have a morphism $Fr_{X}(t) \to Fr_{Y}(t+\varepsilon)$ and a morphism $Fr_{Y}(t) \to Fr_{X}(t+\varepsilon)$. Let $X = \bigcup\limits_{i=1}^{m} A_{i}$ and $Y = \bigcup\limits_{j=1}^{n} B_{j}$, where $A_i$ and $B_j$ are intervals for $1\leq i\leq m$ and $1\leq j\leq n$. These intervals are arranged from left to right along the real line, meaning that the upper bound of $A_i$ is always smaller than the lower bound of $A_{i+1}$ for $1 \leq i \leq m-1$. Similarly, we write the intervals as $A_{i} = [a_{i}, a_{i}']$ and $B_{i} = [b_{i}, b_{i}']$. Here, $a_{1}$ and $b_{1}$ may take the value $-\infty$, while $a'_{m}$ and $b'_{n}$ may take the value $\infty$.

First, we can choose $t$ sufficiently large such that $Fr_{X}(t)$ contains $m$ intervals and $Fr_{Y}(t)$ contains $n$ intervals. Since we have the morphisms $Fr_{X}(t) \to Fr_{Y}(t+\varepsilon)$ and $Fr_{Y}(t) \to Fr_{X}(t+\varepsilon)$, it follows that $m = n$.

Next, consider sufficiently small $t$ such that $Fr_{X}(t)$ and $Fr_{Y}(t)$ each contain at most one interval. If $A_1$ is bounded, then since $Fr_{X}$ and $Fr_{Y}$ are $\varepsilon$-interleaved, $B_1$ must also be bounded, and we have $|a_{1} - b_{1}| \leq \varepsilon$. Here, $a_1$ and $b_1$ are the lower bounds of the intervals $A_1$ and $B_1$, respectively. If $A_1$ is unbounded below and $d_{I}(Fr_{X},Fr_{Y}) < \infty$, then $B_1$ must also be unbounded below. In this case, for sufficiently small $t$, both $Fr_{X}(t)$ and $Fr_{Y}(t)$ contain exactly one interval extending to $-\infty$. It follows that $|a_{1} - b_{1}| =0 \leq \varepsilon$.

Now, consider the case where $i > 1$. If $a_{i} = b_{i}$, then clearly $|a_{i} - b_{i}| = 0 < \varepsilon$. If $a_{i} < b_{i}$, we choose a point $t$ such that $a_{i} < t < b_{i}$. Then, $Fr_{X}(t)$ contains $i$ intervals, while $Fr_{Y}(t)$ contains fewer than $i$ intervals. Since $Fr_{X}$ and $Fr_{Y}$ are $\varepsilon$-interleaved, there exists a morphism $Fr_{X}(t) \to Fr_{Y}(t+\varepsilon)$. This implies that $Fr_{Y}(t+\varepsilon)$ must contain at least $i$ intervals, which gives $b_{i} \leq t + \varepsilon$. Letting $t \to a_i$, we obtain $b_{i} \leq a_{i} + \varepsilon$. Similarly, for the case $a_{i} > b_{i}$, we obtain $a_{i} \leq b_{i} + \varepsilon$. Thus, one has $|a_{i}-b_{i}|\leq \varepsilon$.

Combining these results, we conclude that $d_{L}(X,Y)\leq \varepsilon$.

If $d_{L}(X,Y) < \varepsilon$, then there exists $\eta > 0$ such that $d_{L}(X,Y) < \varepsilon - \eta$. Since $Fr_{X}$ and $Fr_{Y}$ are not $(\varepsilon - \eta)$-interleaved, there must exist some $t = t_{0}$ for which at least one of the following morphisms does not exist:
\begin{equation*}
  Fr_{X}(t_{0})\to Fr_{Y}(t_{0}+\varepsilon-\eta),\quad Fr_{Y}(t_{0})\to Fr_{X}(t_{0}+\varepsilon-\eta).
\end{equation*}
Without loss of generality, assume that the morphism $Fr_{X}(t_{0})\to Fr_{Y}(t_{0}+\varepsilon-\eta)$ does not exist. This implies that $Fr_{X}(t_{0})$ contains more intervals than $Fr_{Y}(t_{0}+\varepsilon-\eta)$. Let $Fr_{X}(t_{0})$ contain $k$ intervals and $Fr_{Y}(t_{0}+\varepsilon-\eta)$ contain $l$ intervals, where $k > l$. Then, we have
\begin{equation*}
  a_{k} \leq t_{0} < a_{k+1},\quad b_{l} \leq t_{0}+\varepsilon-\eta < b_{l+1}.
\end{equation*}
It follows that
\begin{equation*}
  a_{k}+\varepsilon-\eta \leq t_{0}+\varepsilon-\eta < b_{l+1} \leq b_{k}.
\end{equation*}
Thus, we obtain
\begin{equation*}
  \varepsilon - \eta < |a_{k} - b_{k}|.
\end{equation*}
Since $\eta$ can be chosen arbitrarily small, it follows that $\varepsilon = |a_{k} - b_{k}| \leq d_{L}(X,Y)$. Combining the previous results, the proposition is proved.
\end{proof}

\subsection{Different distances on polybarcodes}\label{section:distances}
In this section, we investigate several notions of distance between polybarcodes of finite type, including the expansion distance and the interleaving distance. The expansion distance is typically more intuitive and computationally tractable, while the interleaving distance is more suitable for categorical formulations and for establishing relationships among various metrics.

We say that a polybarcode $\{I_{\pi}\}_{\pi \in \Iso(M)}$ is of \textit{finite type} if, for every $\pi \in \Iso(M)$, the set $I_{\pi}$ is either empty or consists of finitely many connected components.

Let $\mathcal{F}:(\mathbb{R},\leq)\to\mathcal{S}_{n}(M)$ be a continuous persistence configuration. We say that $\mathcal{F}$ admits a \textit{finite-type polybarcode} if the associated polybarcode is of finite type. By Proposition~\ref{proposition:closed}, each set $I(\pi)$ is a closed subset of $\mathbb{R}$. Hence, the condition that $\mathcal{F}$ has a finite-type polybarcode ensures that each $I(\pi)$ is an object in the category $\mathbf{Int}(\mathbb{R})$. Building on this idea, we utilize the category $\mathbf{Int}(\mathbb{R})$ to study distances between polybarcodes.

\begin{definition}
Let $B = \{I_{\pi}\}_{\pi \in \Iso(M)}$ and $B' = \{I_{\pi}'\}_{\pi \in \Iso(M)}$ be two finite-type polybarcodes. The \textit{expansion distance} between $B$ and $B'$ is defined by
\begin{equation*}
  d_{E}(B, B') := \sup_{\pi \in \Iso(M)} d_{E}(I_{\pi}, I_{\pi}'),
\end{equation*}
where $d_{E}(I_{\pi}, I_{\pi}')$ denotes the expansion distance between the interval sets $I_{\pi}$ and $I_{\pi}'$.
\end{definition}

Similarly, we define the notion of the \textit{left expansion distance} between polybarcodes.
\begin{definition}
Let $B = \{I_{\pi}\}_{\pi \in \Iso(M)}$ and $B' = \{I_{\pi}'\}_{\pi \in \Iso(M)}$ be two finite-type polybarcodes. The \textit{left expansion distance} between $B$ and $B'$ is defined by
\begin{equation*}
  d_{L}(B, B') := \sup_{\pi \in \Iso(M)} d_{L}(I_{\pi}, I_{\pi}'),
\end{equation*}
where $d_{L}(I_{\pi}, I_{\pi}')$ denotes the left expansion distance between the intervals $I_{\pi}$ and $I_{\pi}'$.
\end{definition}

Recall the construction $Fr_{X} \colon (\mathbb{R}, \leq) \to \mathbf{Int}(\mathbb{R})$ defined by
\begin{equation*}
  Fr_{X}(t) = \{x \in X \mid x \leq t\}.
\end{equation*}
To define the interleaving distance between polybarcodes, we first introduce a category $\mathbf{Polybarc}$. The objects of $\mathbf{Polybarc}$ are polybarcodes of the form $\{I_{\pi}\}_{\pi \in \Iso(M)}$, where for each $\pi \in \Iso(M)$, the set $I_{\pi}$ belongs to $\mathbf{Int}(\mathbb{R})$. A morphism $\{I_{\pi}\}_{\pi \in \Iso(M)} \to \{I_{\pi}'\}_{\pi \in \Iso(M)}$ in $\mathbf{Polybarc}$ is defined to be a collection of morphisms $I_{\pi} \to I_{\pi}'$ in $\mathbf{Int}(\mathbb{R})$, one for each symmetry $\pi$.

Given a polybarcode $B = \{I_{\pi}\}_{\pi \in \Iso(M)}$, we define the associated functor
\begin{equation*}
  Fr_{B} = \{Fr_{I_{\pi}}\}_{\pi \in \Iso(M)},
\end{equation*}
where each $Fr_{I_{\pi}} \colon (\mathbb{R}, \leq) \to \mathbf{Int}(\mathbb{R})$ is functorial. Therefore, $Fr_{B}$ defines a functor
\begin{equation*}
  Fr_{B} \colon (\mathbb{R}, \leq) \to \mathbf{Polybarc}.
\end{equation*}

\begin{definition}
Let $B = \{I_{\pi}\}_{\pi \in \Iso(M)}$ and $B' = \{I_{\pi}'\}_{\pi \in \Iso(M)}$ be two finite-type polybarcodes. The \textit{interleaving distance} between $B$ and $B'$ is defined as
\begin{equation*}
   d_{I}(B, B') := d_{I}(Fr_{B}, Fr_{B'}),
\end{equation*}
where the right-hand side denotes the standard interleaving distance between functors valued in $\mathbf{Polybarc}$.
\end{definition}

\begin{proposition}\label{proposition:polybarcodes_distance}
Let $B = \{I_{\pi}\}_{\pi \in \Iso(M)}$ and $B' = \{I_{\pi}'\}_{\pi \in \Iso(M)}$ be polybarcodes of finite type. Then the interleaving distance satisfies
\begin{equation*}
  d_{I}(B, B') = \sup_{\pi \in \Iso(M)} d_{I}(Fr_{I_{\pi}}, Fr_{I_{\pi}'}).
\end{equation*}
\end{proposition}

\begin{proof}
Let $\varepsilon = d_{I}(B, B')$. By definition, if the functors $Fr_{B}$ and $Fr_{B'}$ are $\varepsilon$-interleaved, then for each $\pi \in \Iso(M)$, the functors $Fr_{I_{\pi}}$ and $Fr_{I_{\pi}'}$ are also $\varepsilon$-interleaved. Therefore,
\begin{equation*}
d_{I}(Fr_{I_{\pi}}, Fr_{I_{\pi}'}) \leq \varepsilon,
\end{equation*}
and taking the supremum over all $\pi$ yields
\begin{equation*}
\sup_{\pi \in \Iso(M)} d_{I}(Fr_{I_{\pi}}, Fr_{I_{\pi}'}) \leq d_{I}(B, B').
\end{equation*}

Conversely, let
\begin{equation*}
\varepsilon_0 = \sup_{\pi \in \Iso(M)} d_{I}(Fr_{I_{\pi}}, Fr_{I_{\pi}'}).
\end{equation*}
Given any $\delta > 0$, there exists $\pi_{\delta} \in \Iso(M)$ such that
\begin{equation*}
d_{I}(Fr_{I_{\pi_{\delta}}}, Fr_{I_{\pi_{\delta}}'}) > \varepsilon_0 - \delta.
\end{equation*}
Since any $(\varepsilon_0 - \delta)$-interleaving between $Fr_{B}$ and $Fr_{B'}$ would induce a componentwise $(\varepsilon_0 - \delta)$-interleaving between $Fr_{I_{\pi_{\delta}}}$ and $Fr_{I_{\pi_{\delta}}'}$, such an interleaving cannot exist. Hence, we have
\begin{equation*}
d_{I}(B, B') \geq \varepsilon_0 - \delta.
\end{equation*}
Since $\delta > 0$ is arbitrary, it follows that $d_{I}(B, B') \geq \varepsilon_0$.

Combining both inequalities gives the desired equality.
\end{proof}

\begin{proposition}\label{proposition:expansion_interleaving}
Let $B = \{I_{\pi}\}_{\pi \in \Iso(M)}$ and $B' = \{I_{\pi}'\}_{\pi \in \Iso(M)}$ be polybarcodes of finite type. Then the interleaving distance and the left expansion distance coincide
\begin{equation*}
  d_{I}(B, B') = d_{L}(B, B').
\end{equation*}
\end{proposition}

\begin{proof}
This follows directly from Propositions~\ref{proposition:polybar_distance} and~\ref{proposition:polybarcodes_distance}, which relate the interleaving distance of polybarcodes to their left expansion distances componentwise.
\end{proof}

\begin{example}\label{example:distance_computation}
Example \ref{example:symmetry_group} continued. Consider another persistent configuration
\[
\mathcal{F}' = \{Y_t\}_{t=0,1,2} \in \mathcal{S}_4(\mathbb{R}^2)^{(\mathbb{Z}, \leq)},
\]
where each $Y_t = \{A_t, B_t, C_t, D_t\}$ is a configuration of four points in $\mathbb{R}^2$ at discrete time steps $t = 0, 1, 2$. The coordinates of the points $A_t, B_t, C_t, D_t \in \mathbb{R}^2$ are listed in Table~\ref{table:4configuration2}.
\begin{table}[h]
\centering
\begin{tabular}{c|c|c|c|c}
\hline
Time $t$ & Point $A_t$ & Point $B_t$ & Point $C_t$ & Point $D_t$ \\
\hline\hline
$t = 0$ & $(0, -1)$ & $(0, 1)$ & $(-1, 0)$ & $(1, 0)$ \\
$t = 1$ & $(0, -1)$ & $(0, 1)$ & $(-1, 0)$ & $(1, 0)$ \\
$t = 2$ & $(0, -1)$ & $(0, 1)$ & $(-1, 0)$ & $(1.2, 0)$ \\
\hline
\end{tabular}
\caption{Coordinates of the persistent 4-configuration at different time steps in Example \ref{example:distance_computation}.}\label{table:4configuration2}
\end{table}
For $0 \leq i \leq j \leq 2$, the structure map $f_{i,j}: Y_i \to Y_j$ maps $A_i$, $B_i$, $C_i$, and $D_i$ to $A_j$, $B_j$, $C_j$, and $D_j$, respectively. The symmetry group of $\mathcal{F}'_0$ is given by
\begin{equation*}
D_4 = \langle \tau, \pi \mid \tau^4 = e,\; \pi^2 = e,\; \pi\tau = \tau^{-1}\pi \rangle,
\end{equation*}
where $\pi$ denotes the reflection across the $x$-axis, and $\tau$ represents the $90^{\circ}$ rotation. The corresponding polybarcode, excluding the identity element $e$, is given by
\begin{align*}
  & \pi: [0,2],\quad   \tau: [0,1],\quad \tau^{2}:[0,1], \quad \tau^{3}:[0,1],\\
  & \pi\tau:[0,1], \quad \pi\tau^{2}:[0,1], \quad \pi\tau^{3}:[0,1].
\end{align*}
Similarly, the polybarcode for the persistent configuration in Example~\ref{example:symmetry_group} is computed as
\begin{align*}
  & \pi: [0,2],\quad   \tau: \{0\},\quad \tau^{2}:[0,1], \quad \tau^{3}:\{0\},\\
  & \pi\tau:\{0\}, \quad \pi\tau^{2}:[0,1], \quad \pi\tau^{3}:\{0\}.
\end{align*}
Therefore, the symmetric difference distance between the polybarcodes of $\mathcal{F}$ and $\mathcal{F}'$ is
\begin{equation*}
  d_{S}(\mathcal{B}(\mathcal{F}),\mathcal{B}(\mathcal{F}')) = 4.
\end{equation*}
The expansion distance is
\begin{equation*}
  d_{E}(\mathcal{B}(\mathcal{F}), \mathcal{B}(\mathcal{F}')) = d_{E}(I^{\mathcal{F}}(\tau), I^{\mathcal{F}'}(\tau)) = 1.
\end{equation*}
The left expansion distance is
\begin{equation*}
  d_{L}(\mathcal{B}(\mathcal{F}), \mathcal{B}(\mathcal{F}')) = 0.
\end{equation*}
A direct verification of the interleaving distance is somewhat involved, but by Proposition~\ref{proposition:expansion_interleaving}, we conclude that
\begin{equation*}
   d_{I}(\mathcal{B}(\mathcal{F}), \mathcal{B}(\mathcal{F}')) = d_{L}(\mathcal{B}(\mathcal{F}), \mathcal{B}(\mathcal{F}')) = 0.
\end{equation*}

In our computations, intervals are used because the configurations in this example can be extended continuously. For instance, between the time steps $t=1$ and $t=2$, the point $D_t$ can be defined as $(0.8 + 0.2t, 0)$ for $t \in [1,2]$, which yields a continuous trajectory from $D_1$ to $D_2$.
\end{example} 

%% file: stability_symmetries.tex
\section{The stability of persistent symmetries}\label{section:stability}

Stability plays a fundamental role in the analysis of symmetries, especially in the context of data subject to noise and/or perturbation. In the study of persistent homology, stability theorems form a theoretical backbone of topological data analysis \cite{bauer2014induced,chazal2009proximity,cohen2005stability,liu2024algebraic}. In the context of persistent symmetries, stability ensures that the extracted symmetry patterns reflect meaningful geometric features. In this section, we study the stability of symmetry barcodes and the stability of polybarcodes associated with symmetries.

\subsection{Persistence of parametrized systems}
Let $(M, d)$ be a metric space. Recall that a persistence $n$-configuration is a functor $\mathcal{F}: (\mathbb{R}, \leq) \to \mathcal{S}_n(M)$, which encodes a time-evolving collection of $n$-point configurations in $M$. Such a functor can be naturally interpreted as a dynamical system, where the evolution of the configuration is governed by the ordering of $\mathbb{R}$.

Consider the $(\mathbb{R}, \leq)$-indexed category $\mathcal{S}_n(M)^{\mathbb{R}}$, which is the functor category consisting of all persistence $n$-configurations. The objects of $\mathcal{S}_n(M)^{\mathbb{R}}$ are functors $\mathcal{F}: (\mathbb{R}, \leq) \to \mathcal{S}_n(M)$, and the morphisms between two such functors $\mathcal{F}$ and $\mathcal{G}$ are natural transformations $T: \mathcal{F} \Rightarrow \mathcal{G}$. A natural transformation $T$ encodes a structured way of deforming one persistence $n$-configuration into another while respecting the parameter order in $\mathbb{R}$.

When studying the stability of symmetry barcodes, the interleaving distance between persistence $n$-configurations serves as a natural notion of perturbation, providing a way to measure the stability of the associated symmetry barcodes. However, in the case of polybarcodes, this approach becomes insufficient. The category of persistence $n$-configurations admits a large class of morphisms, allowing for too much flexibility. As a result, the interleaving distance between persistence $n$-configurations fails to capture enough geometric disparity between $n$-configurations, and thus cannot effectively constrain the perturbation of polybarcodes. To better control the perturbations of polybarcodes, it is natural to consider persistence objects in the functor category $\mathcal{S}_n(M)^{\mathbb{R}}$, as each object in this category encodes the entire trajectory of an $n$-configuration from the past up to the present moment, thereby retaining sufficient dynamical information.

While the category $\mathcal{S}_n(M)^{\mathbb{R}}$ admits objects that model the evolution of $n$-point configurations over time, the category contains too many morphisms, which limits our ability to study interleaving distances in a tractable way. To address this issue, we define a subcategory $\mathbf{Dyn}^{0}_{n}(M)$ of the functor category $\mathcal{S}_n(M)^{\mathbb{R}}$, which we interpret as the \textit{category of truncated dynamical systems}. An object in $\mathbf{Dyn}^{0}_{n}(M)$ is a pair $(\mathcal{F}, r)$, where $\mathcal{F} : (\mathbb{R}, \leq) \to \mathcal{S}_{n}(M)$ is a persistence $n$-configuration and $r \in \mathbb{R}$ is a truncation parameter such that $\mathcal{F}_{t} = \mathcal{F}_{r}$ for all $t \geq r$; that is, the configuration stabilizes at time $r$.
A morphism $\phi : (\mathcal{F}, r) \to (\mathcal{G}, s)$ in the category $\mathbf{Dyn}^{0}_{n}(M)$ is defined whenever $r \leq s$ and $\mathcal{F}_{t} = \mathcal{G}_{t}$ for all $t \leq r$. Then, for each $t \in \mathbb{R}$, the map $\phi_t : \mathcal{F}_t \to \mathcal{G}_t$ is given by
\begin{equation*}
  \phi_t=\left\{
  \begin{array}{ll}
    \mathrm{id}_{\mathcal{F}_t}, & \text{if } t \leq r; \\[0.5em]
    \mathcal{F}_r = \mathcal{G}_r \xrightarrow{\mathcal{G}_{r \leq t}} \mathcal{G}_t, & \text{if } r < t \leq s;\\[0.5em]
    \mathcal{F}_r = \mathcal{G}_r \xrightarrow{\mathcal{G}_{r \leq s}} \mathcal{G}_s, & \text{if } t > s.
  \end{array}
  \right.
\end{equation*}
This construction reflects the fact that the evolution described by $\mathcal{F}$ is entirely contained within that of $\mathcal{G}$: the two dynamics coincide up to time $r$, beyond which $\mathcal{F}$ remains constant while $\mathcal{G}$ may continue to evolve.

\begin{definition}
A \textit{persistence (truncated) dynamical system} is a functor
\begin{equation*}
  \mathcal{P}: (\mathbb{R}, \leq)\to \mathbf{Dyn}^{0}_{n}(M)
\end{equation*}
from the category $(\mathbb{R}, \leq)$ to the category of truncated dynamical systems.
\end{definition}

Given a persistence $n$-configuration $\mathcal{F} : (\mathbb{R}, \leq) \to \mathcal{S}_{n}(M)$, we can always associate to it a persistence dynamical system $\kappa(\mathcal{F}) : (\mathbb{R}, \leq) \to \mathbf{Dyn}^{0}_{n}(M)$ defined by
\begin{equation*}
  \kappa(\mathcal{F})(r)_{t} = \begin{cases}
    \mathcal{F}_{t}, & \text{if } t \leq r, \\
    \mathcal{F}_{r}, & \text{if } t > r.
  \end{cases}
\end{equation*}
Conversely, given a persistence dynamical system $\mathcal{P} : (\mathbb{R}, \leq) \to \mathbf{Dyn}^{0}_{n}(M)$, we can recover a persistence $n$-configuration $\lambda(\mathcal{P}) : (\mathbb{R}, \leq) \to \mathcal{S}_{n}(M)$ by
\begin{equation*}
  \lambda(\mathcal{P})_{t} = \mathcal{P}(t)_{t}.
\end{equation*}

This establishes a one-to-one correspondence of objects between the functor categories $\mathcal{S}_{n}(M)^{\mathbb{R}}$ and $\mathbf{Dyn}^{0}_{n}(M)^{\mathbb{R}}$, indexed by $(\mathbb{R}, \leq)$. However, this correspondence is not necessarily natural in the categorical sense.

Suppose $M$ is a Euclidean space. Even when $\mathcal{F}$ and $\mathcal{G}$ in the functor category $\mathcal{S}_{n}(M)^{\mathbb{R}}$ differ only by an affine transformation, their polybarcodes can exhibit substantial differences. However, interleavings between the objects $\mathcal{F}$ and $\mathcal{G}$ tend to be trivial. In contrast, the corresponding persistence dynamical systems $\kappa(\mathcal{F})$ and $\kappa(\mathcal{G})$ may display significantly richer interleaving behavior. Thus, using the interleaving distance between $\kappa(\mathcal{F})$ and $\kappa(\mathcal{G})$ offers a more meaningful and discriminative way to measure the difference between $\mathcal{F}$ and $\mathcal{G}$. This observation motivates us to define the indexed interleaving distance between persistence $n$-configurations via their associated persistence dynamical systems in what follows.

We can extend the symmetry construction to the category of truncated dynamical systems. Define a pseudofunctor $\Sym: \mathbf{Dyn}^{0}_{n}(M) \to \Span(\mathbf{Grp})$ as follows:
\begin{itemize}
    \item For an object $(\mathcal{F}, r) \in \mathbf{Dyn}^{0}_{n}(M)$, we define $\Sym(\mathcal{F}, r) = \Sym(\mathcal{F}_r)$.
    \item For a morphism $\phi: (\mathcal{F}, r) \to (\mathcal{G}, s)$, we define its image to be the span
    \begin{equation*}
      \Sym(\mathcal{F}_r) \hookleftarrow L_{r,s} \rightarrow \Sym(\mathcal{G}_s),
    \end{equation*}
    where $L_{r,s}$ is the group of survival symmetries along the trajectory of $\mathcal{G}$ from $r$ to $s$, as defined in Remark \ref{remark:spans}.
\end{itemize}

\begin{proposition}\label{proposition:pseudofunctor}
The construction $\Sym:\mathbf{Dyn}^{0}_{n}(M)\to \Span(\mathbf{Grp})$ is a pseudofunctor.
\end{proposition}

\begin{proof}
By Theorem \ref{theorem:lax_functor}, the assignment of survival symmetries to finite paths in $\mathbb{S}_n(M)$ forms a normal lax functor, where refining a path induces a subgroup inclusion. The group $L_{r,s}$ is defined as the limit over all finite decompositions of the interval $[r,s]$, which corresponds to taking the intersection of all such path-dependent survival subgroups.

To verify the pseudofunctoriality, consider composable morphisms $(\mathcal{F}, r) \xrightarrow{\phi} (\mathcal{G}, s) \xrightarrow{\psi} (\mathcal{H}, u)$. The composition $\Sym(\psi) \circ \Sym(\phi)$ in $\Span(\mathbf{Grp})$ is given by the pullback over $\Sym(\mathcal{G}_s)$. Since the leg maps are injective by Lemma \ref{lemma:span_injections}, this pullback is naturally isomorphic to the intersection $L_{r,s} \cap L_{s,u}$ within $\Sym(\mathcal{G}_s)$.

A symmetry belongs to this intersection if and only if it survives from $r$ to $s$ and from $s$ to $u$. This is equivalent to surviving along the entire interval $[r, u]$, which is precisely the definition of $L_{r,u}$. Hence, we have a natural isomorphism:
\begin{equation*}
\Sym(\psi) \circ \Sym(\phi) \cong \Sym(\psi \circ \phi).
\end{equation*}
This shows that the non-invertible inclusions from the lax structure on $\mathbb{S}_n(M)$ become invertible isomorphisms at the limit level. The identity is trivially preserved since $L_{r,r} = \Sym(\mathcal{F}_r)$. Therefore, $\Sym$ is a pseudofunctor.
\end{proof}

Now, consider the following diagram
\begin{equation*}
  \xymatrix{
 (\mathbb{R},\leq )\ar[rr]^{\mathcal{F}}\ar[rd]_{\kappa(\mathcal{F})} && \mathcal{S}_{n}(M) \\
 &\mathbf{Dyn}^{0}_{n}(M)\ar[ru]_{q}. &
 }
\end{equation*}
Here, the functor $q: \mathbf{Dyn}^{0}_{n}(M) \to \mathcal{S}_{n}(M)$ is defined by $q(\mathcal{F}, r) = \mathcal{F}_r$, and assigns to each morphism $(\mathcal{F}, r) \to (\mathcal{G}, s)$ the map $\mathcal{F}_r \to \mathcal{G}_s$. One can readily verify that $q$ is indeed a functor.

\begin{lemma}\label{lemma:functors_equation}
We have the equality $\mathcal{F}=q\circ \kappa(\mathcal{F})$ as functors.
\end{lemma}

\begin{proof}
By the definition of the functor $q: \mathbf{Dyn}^{0}_{n}(M) \to \mathcal{S}_{n}(M)$, we have
\[
(q \circ \kappa(\mathcal{F}))(r) = q(\kappa(\mathcal{F})(r)) = \kappa(\mathcal{F})(r)_r.
\]
From the definition of $\kappa(\mathcal{F})$, we know that
\[
\kappa(\mathcal{F})(r)_r = \mathcal{F}_r.
\]
Therefore, we obtain $(q \circ \kappa(\mathcal{F}))(r) = \mathcal{F}_r$ for all $r \in \mathbb{R}$.

On the other hand, for $r \leq s$ in $\mathbb{R}$, the morphism $\kappa(\mathcal{F})(r \to s)$ in $\mathbf{Dyn}^{0}_{n}(M)$ is a natural transformation, whose components are given by
\[
\kappa(\mathcal{F})(r \to s)_t =
\begin{cases}
\mathrm{id}_{\mathcal{F}_t}, & \text{if } t \leq r; \\
\mathcal{F}(r \to t), & \text{if } r < t \leq s; \\
\mathcal{F}(r \to s), & \text{if } t > s.
\end{cases}
\]
Applying the functor $q$ to this morphism yields the map $\mathcal{F}_r \to \mathcal{F}_s$, which is precisely the component at $t=s$ (or any $t > s$). Hence, it follows that
\[
(q \circ \kappa(\mathcal{F}))(r \to s) = \mathcal{F}(r \to s).
\]

This shows that $\mathcal{F} = q \circ \kappa(\mathcal{F})$ as functors.
\end{proof}

Two persistence dynamical systems $\mathcal{P},\mathcal{Q}: (\mathbb{R}, \leq) \to \mathbf{Dyn}^{0}_{n}(M)$ are said to be $\varepsilon$-interleaved if there are natural transformations $\phi: \mathcal{P} \Rightarrow \Sigma^\varepsilon \mathcal{Q}$ and $\psi: \mathcal{Q} \Rightarrow\Sigma^\varepsilon \mathcal{P}$ satisfying the commutativity conditions
\begin{equation*}
(\Sigma^\varepsilon \psi) \circ \phi = \Sigma^{2\varepsilon}|_{\mathcal{P}}, \quad (\Sigma^\varepsilon \phi) \circ \psi = \Sigma^{2\varepsilon}|_{\mathcal{Q}}.
\end{equation*}
More precisely, for each $r \in \mathbb{R}$, there exist morphisms
\begin{equation*}
\phi_r: \mathcal{P}(r) \to  \mathcal{Q}(r + \varepsilon), \quad \psi_r:  \mathcal{Q}(r) \to  \mathcal{P}(r + \varepsilon)
\end{equation*}
such that the following compositions commute
\begin{equation*}
\psi_{r + \varepsilon} \circ \phi_r = \mathcal{P}(r \xrightarrow{2\varepsilon} r + 2\varepsilon), \quad
\phi_{r + \varepsilon} \circ \psi_r = \mathcal{Q}(r \xrightarrow{2\varepsilon} r + 2\varepsilon).
\end{equation*}

\begin{definition}
Let $\mathcal{F}, \mathcal{G}: (\mathbb{R}, \leq) \to \mathcal{S}_{n}(M)$ be two persistence $n$-configurations. We say that $\mathcal{F}$ and $\mathcal{G}$ are \textit{indexed $\varepsilon$-interleaved} if the corresponding persistence dynamical systems $\kappa(\mathcal{F}), \kappa(\mathcal{G}): (\mathbb{R}, \leq) \to \mathbf{Dyn}^{0}_{n}(M)$ are $\varepsilon$-interleaved. Accordingly, we define the \textit{indexed interleaving distance} between $\mathcal{F}$ and $\mathcal{G}$ as
\begin{equation*}
  d_{II}(\mathcal{F}, \mathcal{G}) := d_{I}(\kappa(\mathcal{F}), \kappa(\mathcal{G})),
\end{equation*}
where $d_I$ denotes the interleaving distance between persistence dynamical systems.
\end{definition}

In the above definition, however, the category $\mathbf{Dyn}^{0}_{n}(M)$ may admit too few morphisms---many pairs of objects are not connected by any morphism at all. As a result, the interleaving distance $d_I$ becomes infinite in many cases. To overcome this limitation, we introduce a refined category of dynamical systems, denoted by $\mathbf{Dyn}^{1}_{n}(M)$, whose objects are also persistence $n$-configurations, but which possesses a richer set of morphisms that better capture the dynamical similarity between configurations.

A morphism $\phi: \mathcal{F} \to \mathcal{G}$ in the category $\mathbf{Dyn}^1_n(M)$ is a natural transformation determined by a continuous function $\theta_{r} = \tau_{r} \circ \theta: \mathbb{R} \to \mathbb{R}$, such that
\begin{equation*}
\mathcal{F} = \mathcal{G} \circ \theta_{r}.
\end{equation*}
Here, $\theta$ is a continuous and strictly increasing function, and $\tau_r(t)= \begin{cases}
t, & \text{if } t\leq r; \\
r, & \text{if } t> r
\end{cases}$ denotes the truncation function at level $r \in \mathbb{R} \cup \{\infty\}$. This implies that for each $t \in \mathbb{R}$, there is a morphism
\begin{equation*}
\phi_t: \mathcal{F}_t =\mathcal{G}_{\theta_{r}(t)} \longrightarrow  \mathcal{G}_t,
\end{equation*}
and the collection $\{\phi_t\}_{t \in \mathbb{R}}$ forms a natural transformation with respect to the parameter $t$. In particular, for any $t\leq t'$, the following square commutes:
\[
\xymatrix{
\mathcal{G}_{\theta_r(t)} \ar[rr]^{\phi_t} \ar[d] & & ~\mathcal{G}_{t}~~\ar[d]\\
\mathcal{G}_{\theta_r(t')}\ar[rr]^{\phi_{t'}} & &  ~\mathcal{G}_{t'},
}
\]
where the vertical maps are induced by $t\to t'$.

In the special case where $\theta(t) = t$, we have $\theta_r(t)= \begin{cases}
t, & \text{if } t\leq r; \\
r, & \text{if } t> r.
\end{cases}$ The morphism $\phi$ corresponds to truncating $\mathcal{G}$ at level $r$, that is, $\mathcal{F}$ is the $r$-truncation of $\mathcal{G}$. Every morphism in the category $\mathbf{Dyn}^0_n(M)$ can be written as a composition of such truncation morphisms. It is worth noting that in the construction of the category $\mathbf{Dyn}^1_n(M)$, we allow $r = \infty$, in which case the truncation function $\theta_r$ reduces to the identity map $\mathrm{id}$.

\begin{definition}
A \textit{persistence dynamical system} is a functor
\begin{equation*}
  \mathcal{P}: (\mathbb{R}, \leq)\to \mathbf{Dyn}^{1}_{n}(M)
\end{equation*}
from the category $(\mathbb{R}, \leq)$ to the category of dynamical systems.
\end{definition}

To compare two persistence dynamical systems, we use the standard interleaving distance between functors valued in the category $\mathbf{Dyn}^{1}_{n}(M)$. This metric quantifies the similarity between the two systems in both structure and dynamics, providing a robust foundation for analyzing their temporal evolution in complex settings.

\subsection{Interleavings between pseudofunctors}

In this section, we extend the theory of persistent homology into the setting of bicategories, building on the categorical foundations of persistence theory. This extension allows for a more flexible treatment of persistent symmetries. The discussion here involves some basic concepts and standard results from the theory of bicategories. For readers interested in the foundational ideas and theoretical background of bicategories, we refer to \cite{benabou1967introduction, lurie2009higher, mac2013categories}.

\begin{definition}
Let $\mathfrak{C}$ be a bicategory. A \textit{persistence object} in $\mathfrak{C}$ is a pseudofunctor
\begin{equation*}
  \mathcal{F}: (\mathbb{R}, \leq) \longrightarrow \mathfrak{C},
\end{equation*}
where $(\mathbb{R}, \leq)$ is the poset of real numbers viewed as a category.
\end{definition}
It is important to note that for each triple $a \leq b \leq c$, the persistence object $\mathcal{F}$ assigns a natural isomorphism
\begin{equation*}
  \Phi_{a,b,c}: \mathcal{F}(b \leq c) \circ \mathcal{F}(a \leq b) \longrightarrow \mathcal{F}(a \leq c).
\end{equation*}
Additionally, for each $a \in \mathbb{R}$, the persistence object $\mathcal{F}$ assigns a natural isomorphism
\begin{equation*}
  \eta_a: \mathrm{id}_{\mathcal{F}(a)} \longrightarrow \mathcal{F}(a \leq a).
\end{equation*}
These isomorphisms must satisfy the usual coherence conditions required of pseudofunctors.

\begin{example}
Recall that the span category in the category of groups, denoted $\mathrm{Span}(\mathbf{Grp})$, forms a bicategory. Proposition~\ref{proposition:pseudofunctor} shows that the symmetry group construction
\begin{equation*}
\Sym: \mathbf{Dyn}^{0}_{n}(M) \to \mathrm{Span}(\mathbf{Grp})
\end{equation*}
is a pseudofunctor. Given a persistence dynamical system
\begin{equation*}
\mathcal{P}: (\mathbb{R}, \leq) \to \mathbf{Dyn}^{0}_{n}(M),
\end{equation*}
we recall that any ordinary category can be regarded as a bicategory in which the hom-categories are discrete (i.e., categories whose only 2-morphisms are identities). Now consider the composition of pseudofunctors:
\begin{equation*}
\xymatrix{
  (\mathbb{R}, \leq)  \ar[r]^{\mathcal{P}}& \mathbf{Dyn}^{0}_{n}(M) \ar[r]^-{\Sym}& \mathrm{Span}(\mathbf{Grp}).
}
\end{equation*}
This composition yields a persistence object $\Sym\circ\mathcal{P}: (\mathbb{R}, \leq) \to \mathrm{Span}(\mathbf{Grp})$.
\end{example}

All pseudofunctors from $(\mathbb{R}, \leq)$ to a bicategory $\mathfrak{C}$ form a category, denoted by $\mathfrak{C}^{\mathbb{R}}$. The objects of this category are persistence objects from $(\mathbb{R}, \leq)$ to the target bicategory $\mathfrak{C}$, and the morphisms are pseudonatural transformations between pseudofunctors.

Recall that the translation functor $T_x: (\mathbb{R}, \leq) \to (\mathbb{R}, \leq)$ is given by $T_x(a) = a + x$ for all $a, x \in \mathbb{R}$. This induces an endopseudofunctor $\Sigma^x: \mathfrak{C}^{\mathbb{R}} \to \mathfrak{C}^{\mathbb{R}}$ on the category of pseudofunctors, where $(\Sigma^x \mathcal{F})(a) = \mathcal{F}(a + x)$ for any pseudofunctor $\mathcal{F}: (\mathbb{R}, \leq) \longrightarrow \mathfrak{C}$. Moreover, there exists a pseudonatural transformation $\Sigma^x|_{\mathcal{F}}: \mathcal{F} \Rightarrow \Sigma^x \mathcal{F}$, which defines a morphism in $\mathfrak{C}^{\mathbb{R}}$. We often denote this transformation simply by $\Sigma^x$ when no ambiguity arises.

\begin{definition}
Two persistence objects $\mathcal{F}, \mathcal{G}: (\mathbb{R}, \leq) \to \mathfrak{C}$ are said to be \textit{$\varepsilon$-interleaved} if there exist pseudonatural transformations $\phi: \mathcal{F} \Rightarrow \Sigma^\varepsilon \mathcal{G}$ and $\psi: \mathcal{G} \Rightarrow \Sigma^\varepsilon \mathcal{F}$ in $\mathfrak{C}^{\mathbb{R}}$ such that the following isomorphisms hold
\begin{equation*}
(\Sigma^\varepsilon \psi) \circ \phi \cong \Sigma^{2\varepsilon}|_{\mathcal{F}}, \quad
(\Sigma^\varepsilon \phi) \circ \psi \cong \Sigma^{2\varepsilon}|_{\mathcal{G}}.
\end{equation*}
\end{definition}
The interleaving distance between two persistence objects $\mathcal{F}, \mathcal{G}: (\mathbb{R}, \leq) \to \mathfrak{C}$ is defined as
\begin{equation*}
d_{I}(\mathcal{F},\mathcal{G}) = \inf \left\{ \varepsilon \geq 0 \ \middle|\ \text{$\mathcal{F}$ and $\mathcal{G}$ are $\varepsilon$-interleaved} \right\}.
\end{equation*}
This defines an extended metric for the collection of persistence objects in the category $\mathfrak{C}^{\mathbb{R}}$. We obtain a result similar to \cite[Proposition 3.6]{chazal2009proximity} as follows.

\begin{proposition}\label{proposition:interleaved}
Let $\mathcal{F}, \mathcal{G}: (\mathbb{R}, \leq) \to \mathfrak{C}$ be persistence objects in the category $\mathfrak{C}^{\mathbb{R}}$, and suppose that they are $\varepsilon$-interleaved.
Then for any pseudofunctor $H: \mathfrak{C} \to \mathcal{D}$ between bicategories, we have that $H \circ \mathcal{F}$ and $H \circ \mathcal{G}$ are also $\varepsilon$-interleaved in $\mathcal{D}^{\mathbb{R}}$.
\end{proposition}

\begin{proof}
Suppose $\mathcal{F}$ and $\mathcal{G}$ are $\varepsilon$-interleaved via pseudonatural transformations
\begin{equation*}
\phi: \mathcal{F} \Rightarrow \Sigma^\varepsilon \mathcal{G}, \quad \psi: \mathcal{G} \Rightarrow \Sigma^\varepsilon \mathcal{F}
\end{equation*}
in $\mathfrak{C}^{\mathbb{R}}$ satisfying the isomorphisms
\begin{equation*}
(\Sigma^\varepsilon \psi) \circ \phi \cong \Sigma^{2\varepsilon}|_{\mathcal{F}}, \quad (\Sigma^\varepsilon \phi) \circ \psi \cong \Sigma^{2\varepsilon}|_{\mathcal{G}}.
\end{equation*}
Since $H$ is a pseudofunctor, it preserves the composition of 1-morphisms and 2-morphisms up to coherent isomorphism.
\begin{equation*}
  \xymatrix@=1.5cm{
\mathcal{F} \ar[r]^{\phi} \ar[d]_{H}\ar@/^1.5pc/[rr]^{\Sigma^{2\varepsilon}|_{\mathcal{F}}} & \Sigma^\varepsilon \mathcal{G} \ar[r]^{\Sigma^{\varepsilon}\psi} \ar[d]^{H} & \Sigma^{2\varepsilon} \mathcal{F} \ar[d]^{H} \\
  H\mathcal{F} \ar[r]^{H\phi}\ar@/_1.5pc/[rr]_{\Sigma^{2\varepsilon}|_{H\mathcal{F}}} & \Sigma^{\varepsilon} H\mathcal{G} \ar[r]^{\Sigma^{\varepsilon}H\psi} & \Sigma^{2\varepsilon} H\mathcal{F}
}
\end{equation*}
Applying $H$ to the pseudonatural transformations $\phi$ and $\psi$ yields pseudonatural transformations
\begin{equation*}
H\phi: H\mathcal{F} \Rightarrow H(\Sigma^\varepsilon \mathcal{G}) \cong \Sigma^\varepsilon (H\mathcal{G}), \quad
H\psi: H\mathcal{G} \Rightarrow \Sigma^\varepsilon (H\mathcal{F}).
\end{equation*}
We now verify the interleaving relations. Because $H$ is a pseudofunctor, it preserves compositions up to coherent isomorphism, yielding
\begin{equation*}
(\Sigma^\varepsilon H\psi) \circ H\phi \cong H((\Sigma^\varepsilon \psi) \circ \phi) \cong H(\Sigma^{2\varepsilon}|_{\mathcal{F}}) \cong \Sigma^{2\varepsilon}|_{H\mathcal{F}},
\end{equation*}
and similarly,
\begin{equation*}
(\Sigma^\varepsilon H\phi) \circ H\psi \cong H((\Sigma^\varepsilon \phi) \circ \psi) \cong H(\Sigma^{2\varepsilon}|_{\mathcal{G}}) \cong \Sigma^{2\varepsilon}|_{H\mathcal{G}}.
\end{equation*}
Therefore, $H\mathcal{F}$ and $H\mathcal{G}$ are $\varepsilon$-interleaved in $\mathcal{D}^{\mathbb{R}}$.
\end{proof}

\begin{remark}
The result follows essentially from the functoriality of bicategorical composition and the naturality of the shift endofunctor $\Sigma^\varepsilon$ on $\mathfrak{C}^{\mathbb{R}}$.
\end{remark}

\subsection{Stability of symmetry barcodes}

Compared to polybarcodes, which heavily depend on the positional evolution of configurations within persistence $n$-configurations, symmetry barcodes rely more fundamentally on the structure of the persistent symmetry groups. This intrinsic difference typically results in smaller distances between symmetry barcodes. Consequently, the interleaving distance between persistence dynamical systems is often effective in controlling the perturbations of symmetry barcodes.

\begin{theorem}[Stability theorem for persistent symmetry groups]\label{theorem:stability_group}
Let $\mathcal{P}, \mathcal{Q}: (\mathbb{R}, \leq) \to \mathbf{Dyn}^{0}_{n}(M)$ be two persistence dynamical systems. Then we have
\begin{equation*}
   d_{I}(\Sym\circ\mathcal{P},\Sym\circ\mathcal{Q})\leq d_{I}(\mathcal{P},\mathcal{Q}) .
\end{equation*}
In particular, for any two persistence $n$-configurations $\mathcal{F}, \mathcal{G}: (\mathbb{R}, \leq) \to \mathcal{S}_{n}(M)$, we have
\begin{equation*}
   d_{I}(\Sym\circ\kappa(\mathcal{F}),\Sym\circ\kappa(\mathcal{G}))\leq d_{II}(\mathcal{F},\mathcal{G}) .
\end{equation*}
\end{theorem}

\begin{proof}
By Proposition \ref{proposition:pseudofunctor}, the construction $\Sym:\mathbf{Dyn}^{0}_{n}(M)\to \mathrm{Span}(\mathbf{Grp})$ is a pseudofunctor. For any two persistence dynamical systems $\mathcal{P}, \mathcal{Q}: (\mathbb{R}, \leq) \to \mathbf{Dyn}^{0}_{n}(M)$, we have the composition of pseudofunctors
\begin{equation*}
\xymatrix{
  (\mathbb{R}, \leq)  \ar@<0.5ex>[r]^{\mathcal{P}} \ar@<-0.5ex>[r]_{\mathcal{Q}}& \mathbf{Dyn}^{0}_{n}(M) \ar@{->}[r]^-{\Sym}& \mathrm{Span}(\mathbf{Grp}).
}
\end{equation*}
Applying Proposition \ref{proposition:interleaved}, if $\mathcal{P}$ and $\mathcal{Q}$ are $\varepsilon$-interleaved, then $\Sym\circ\mathcal{P}$ and $\Sym\circ\mathcal{Q}$ are also $\varepsilon$-interleaved. Therefore, we conclude that
\begin{equation*}
  d_{I}(\Sym\circ\mathcal{P},\Sym\circ\mathcal{Q})\leq d_{I}(\mathcal{P},\mathcal{Q}) .
\end{equation*}

To obtain the inequality for persistence $n$-configurations, let $\mathcal{F}, \mathcal{G}: (\mathbb{R}, \leq) \to \mathcal{S}_{n}(M)$ be two persistence $n$-configurations. By definition, their indexed interleaving distance is given by $d_{II}(\mathcal{F}, \mathcal{G}) = d_{I}(\kappa(\mathcal{F}), \kappa(\mathcal{G}))$. Substituting $\mathcal{P} = \kappa(\mathcal{F})$ and $\mathcal{Q} = \kappa(\mathcal{G})$ into the above inequality yields
\begin{equation*}
   d_{I}(\Sym\circ\kappa(\mathcal{F}),\Sym\circ\kappa(\mathcal{G}))\leq d_{II}(\mathcal{F},\mathcal{G}) .
\end{equation*}
This completes the proof.
\end{proof}

We now proceed to study the stability of the persistence symmetry module. First, we establish that the module construction can be naturally lifted to the category of truncated dynamical systems.

First, we can perform a linearization on the pseudofunctor $\Sym: \mathbf{Dyn}^{0}_{n}(M) \to \Span(\mathbf{Grp})$ to transform it into a strict functor $\mathcal{M}: \mathbf{Dyn}^{0}_{n}(M) \to \mathbf{Vec}_{\mathbb{K}}$ as follows:
\begin{itemize}
    \item For an object $(\mathcal{F}, r) \in \mathbf{Dyn}^{0}_{n}(M)$, we define $\mathcal{M}(\mathcal{F}, r) = \mathbb{K}[\Sym(\mathcal{F}_r)]$, which is the free vector space generated by the symmetry group $\Sym(\mathcal{F}_r)$.
    \item For a morphism $\phi: (\mathcal{F}, r) \to (\mathcal{G}, s)$, we have the span
    \begin{equation*}
      \Sym(\mathcal{F}_r) \xhookleftarrow{\iota} L_{r,s} \xrightarrow{\rho} \Sym(\mathcal{G}_s).
    \end{equation*}
    The linear map $\mathcal{M}(\phi): \mathbb{K}[\Sym(\mathcal{F}_r)] \to \mathbb{K}[\Sym(\mathcal{G}_s)]$ is defined by the composition
    \begin{equation*}
     \mathcal{M} = F \circ U \circ \Sym,
    \end{equation*} 
    where $U: \Span(\mathbf{Grp}) \to \Span(\mathbf{Set})$ is the forgetful pseudofunctor and $F$ is the free vector space construction. Concretely, for a basis element $\sigma \in \Sym(\mathcal{F}_r)$, we have
    \begin{equation*}
      \mathcal{M}(\phi)(\sigma) = \begin{cases} \rho(\sigma), & \text{if } \sigma \in L_{r,s}; \\ 0, & \text{otherwise}. \end{cases}
    \end{equation*}
\end{itemize}

\begin{proposition}\label{proposition:functor}
The construction $\mathcal{M} = F \circ U \circ \Sym : \mathbf{Dyn}^{0}_{n}(M) \to \mathbf{Vec}_{\mathbb{K}}$ is a strict functor.
\end{proposition}

\begin{proof}
The composition of pseudofunctors is again a pseudofunctor, so $\mathcal{M}$ is a pseudofunctor. Since $\mathbf{Dyn}^{0}_{n}(M)$ and $\mathbf{Vec}_{\mathbb{K}}$ are ordinary categories (viewed as locally discrete bicategories with only identity 2-cells), any pseudofunctor between them must strictly preserve compositions and identities. Thus, $\mathcal{M}$ is a strict functor.
\end{proof}

\begin{theorem}[Stability theorem for persistence symmetry modules]\label{theorem:stability_module}
Let $\mathcal{P}, \mathcal{Q}: (\mathbb{R}, \leq) \to \mathbf{Dyn}^{0}_{n}(M)$ be two persistence dynamical systems. Then we have
\begin{equation*}
   d_{I}(\mathcal{M}\circ\mathcal{P},\mathcal{M}\circ\mathcal{Q})\leq d_{I}(\mathcal{P},\mathcal{Q}) .
\end{equation*}
\end{theorem}

\begin{proof}
By Proposition~\ref{proposition:functor}, the construction $\mathcal{M}: \mathbf{Dyn}^{0}_{n}(M) \to \mathbf{Vec}_{\mathbb{K}}$ is a strict functor. Consider the following composition of functors
\begin{equation*}
\xymatrix{
  (\mathbb{R}, \leq)  \ar@<0.5ex>[r]^{\mathcal{P}} \ar@<-0.5ex>[r]_{\mathcal{Q}} & \mathbf{Dyn}^{0}_{n}(M) \ar[r]^-{\mathcal{M}} & \mathbf{Vec}_{\mathbb{K}}.
}
\end{equation*}
Assume that $\mathcal{P}$ and $\mathcal{Q}$ are $\varepsilon$-interleaved. Since $\mathcal{M}$ is a strict functor, it preserves interleavings. Therefore, $\mathcal{M}\circ\mathcal{P}$ and $\mathcal{M}\circ\mathcal{Q}$ are also $\varepsilon$-interleaved. Thus the desired result follows.
\end{proof}

Similar to persistent homology, persistent symmetry modules can also be represented by barcodes. As shown in \cite{bauer2020persistence}, there exists a construction (not necessarily a functor)
\[
\Gamma : (\mathbf{Vec}_{\mathbb{K}})^{\mathbb{R}} \to \mathbf{Barc}
\]
from the category of persistence modules to the category of barcodes, such that for any persistence modules $\mathcal{M}, \mathcal{N}$, we have
\[
d_{I}(\mathcal{M}, \mathcal{N}) = d_{I}(\Gamma(\mathcal{M}), \Gamma(\mathcal{N})) = d_{B}(\Gamma(\mathcal{M}), \Gamma(\mathcal{N})),
\]
where $d_I$ denotes the interleaving distance and $d_B$ the bottleneck distance.

\begin{proposition}
Let $\mathcal{P}, \mathcal{Q} \colon (\mathbb{R}, \leq) \to \mathbf{Dyn}^{0}_{n}(M)$ be two persistence dynamical systems. Then we have
\[
d_B(\mathrm{SymB}(\mathcal{P}), \mathrm{SymB}(\mathcal{Q})) \leq d_I(\mathcal{P}, \mathcal{Q}),
\]
where $\mathrm{SymB}$ denotes the construction that assigns to each persistence dynamical system its symmetry barcode.
\end{proposition}

\begin{proof}
Note that symmetry barcodes are defined as the barcodes of the associated persistent symmetry modules. Thus, we have
\[
\mathrm{SymB}(\mathcal{P}) = \Gamma(\mathcal{M} \circ \mathcal{P}), \quad \mathrm{SymB}(\mathcal{Q}) = \Gamma(\mathcal{M} \circ \mathcal{Q}).
\]
By the stability result in \cite{bauer2020persistence}, we obtain
\[
d_B(\Gamma(\mathcal{M} \circ \mathcal{P}), \Gamma(\mathcal{M} \circ \mathcal{Q})) = d_I(\Gamma(\mathcal{M} \circ \mathcal{P}), \Gamma(\mathcal{M} \circ \mathcal{Q})) = d_I(\mathcal{M} \circ \mathcal{P}, \mathcal{M} \circ \mathcal{Q}).
\]
Combining this with Theorem~\ref{theorem:stability_module}, we conclude that
\[
d_B(\mathrm{SymB}(\mathcal{P}), \mathrm{SymB}(\mathcal{Q})) = d_I(\mathcal{M} \circ \mathcal{P}, \mathcal{M} \circ \mathcal{Q}) \leq d_I(\mathcal{P}, \mathcal{Q}).
\]
This completes the proof.
\end{proof}

Similar to persistent homology, persistent symmetry modules can also be represented by barcodes. As shown in \cite{bauer2020persistence}, there exists a construction (not necessarily a functor)
\[
\Gamma : (\mathbf{Vec}_{\mathbb{K}})^{\mathbb{R}} \to \mathbf{Barc}
\]
from the category of persistence modules to the category of barcodes, such that for any persistence modules $\mathcal{M}, \mathcal{N}$, we have
\[
d_{I}(\mathcal{M}, \mathcal{N}) = d_{I}(\Gamma(\mathcal{M}), \Gamma(\mathcal{N})) = d_{B}(\Gamma(\mathcal{M}), \Gamma(\mathcal{N})),
\]
where $d_I$ denotes the interleaving distance and $d_B$ the bottleneck distance.

\begin{proposition}
Let $\mathcal{F}, \mathcal{G} \colon (\mathbb{R}, \leq) \to \mathcal{S}_n(M)$ be two persistence $n$-configurations. Then we have
\[
d_B(\mathrm{SymB}(\mathcal{F}), \mathrm{SymB}(\mathcal{G})) \leq d_{II}(\mathcal{F}, \mathcal{G}),
\]
where $\mathrm{SymB}$ denotes the construction that assigns to each persistence configuration its symmetry barcode.
\end{proposition}

\begin{proof}
Let $V_{\mathcal{F}} = \mathcal{M} \circ \kappa(\mathcal{F})$ and $V_{\mathcal{G}} = \mathcal{M} \circ \kappa(\mathcal{G})$ be the associated persistent symmetry modules. For $r \leq s$, let $\nu^{\mathcal{F}}_{r,s} : V_{\mathcal{F}}(r) \to V_{\mathcal{F}}(s)$ be the linear transition map. By the definition of $\mathcal{M}$, for a basis element $\sigma \in \Sym(\mathcal{F}_r)$, we have
\[
\nu^{\mathcal{F}}_{r,s}(\sigma) = \begin{cases} \rho_{r,s}(\sigma), & \text{if } \sigma \in L_{r,s}; \\ 0, & \text{otherwise}. \end{cases}
\]
Thus, $\nu^{\mathcal{F}}_{r,s}(\sigma) \neq 0$ if and only if $\sigma$ survives continuously as a non-trivial symmetry along $[r,s]$.

By the structure theorem for persistence modules, $V_{\mathcal{F}}$ decomposes into interval modules. An interval $[a, b)$ appears in the barcode $\Gamma(V_{\mathcal{F}})$ if and only if there exists a symmetry born at $a$ (having no preimage from any $t < a$) that survives continuously over $[a, b)$. This exactly matches the definition of a symmetry bar in Section~\ref{section:persistence_group}. Hence, $\mathrm{SymB}(\mathcal{F}) = \Gamma(V_{\mathcal{F}})$, and similarly $\mathrm{SymB}(\mathcal{G}) = \Gamma(V_{\mathcal{G}})$.

By the stability result in \cite{bauer2020persistence}, we obtain
\[
d_B(\mathrm{SymB}(\mathcal{F}), \mathrm{SymB}(\mathcal{G})) = d_B(\Gamma(V_{\mathcal{F}}), \Gamma(V_{\mathcal{G}})) = d_I(V_{\mathcal{F}}, V_{\mathcal{G}}).
\]
According to Theorem~\ref{theorem:stability_module}, the construction $\mathcal{M}$ preserves interleavings, which implies
\[
d_I(V_{\mathcal{F}}, V_{\mathcal{G}}) = d_I(\mathcal{M} \circ \kappa(\mathcal{F}), \mathcal{M} \circ \kappa(\mathcal{G})) \leq d_I(\kappa(\mathcal{F}), \kappa(\mathcal{G})).
\]
Finally, by the definition of the indexed interleaving distance, $d_I(\kappa(\mathcal{F}), \kappa(\mathcal{G})) = d_{II}(\mathcal{F}, \mathcal{G})$. Combining these yields
\[
d_B(\mathrm{SymB}(\mathcal{F}), \mathrm{SymB}(\mathcal{G})) \leq d_{II}(\mathcal{F}, \mathcal{G}).
\]
This completes the proof.
\end{proof}

\subsection{Stability of polybarcodes}

In this section, we study the stability of polybarcodes. We begin by proving the stability of polybarcodes with respect to persistence dynamical systems. Then, we investigate the stability of polybarcodes arising from persistence $n$-configurations. The interleaving distance is better suited for theoretical analysis and the formulation of general concepts, making it particularly useful for stating and proving the main results. In contrast, the expansion distance is more intuitive and computationally accessible, which facilitates practical applications and numerical implementation. Owing to the close relationship between the expansion and interleaving distances, the stability of polybarcodes can also be studied via the expansion distance.

First, recall that given a persistent $n$-configuration $\mathcal{F}:(\mathbb{R}, \leq) \to \mathcal{S}_{n}(M)$, we have a construction
\begin{equation*}
I^{\mathcal{F}}: \Iso(M) \to \mathbf{Int}(\mathbb{R}), \quad \pi \mapsto I^{\mathcal{F}}(\pi).
\end{equation*}
which assigns to each isometry $\pi$ the set of time points at which $\pi$ preserves the configuration
\begin{equation*}
I^{\mathcal{F}}(\pi) = \{ t \in \mathbb{R} \mid \pi(\mathcal{F}_t) = \mathcal{F}_t \}.
\end{equation*}
Fixing an isometry $\pi \in \Iso(M)$, we obtain a map
\begin{equation*}
I(\pi): \mathcal{S}_n(M)^{(\mathbb{R}, \leq)} \to \mathbf{Int}(\mathbb{R}), \quad \mathcal{F} \mapsto I^{\mathcal{F}}(\pi).
\end{equation*}

Now consider the subcategory $\mathbf{Dyn}^1_n(M)$ of $\mathcal{S}_n(M)^{(\mathbb{R}, \leq)}$, which shares the same objects but restricts the class of morphisms. For any object $\mathcal{F}$ in $\mathbf{Dyn}^1_n(M)$, the construction $I^{\mathcal{F}}(\pi)$ remains valid, and we continue to denote the map as
\begin{equation*}
I(\pi): \mathbf{Dyn}^1_n(M) \to \mathbf{Int}(\mathbb{R}), \quad \mathcal{F} \mapsto I^{\mathcal{F}}(\pi).
\end{equation*}
Take a morphism $\phi: \mathcal{F} = \mathcal{G} \circ \theta_{r}\to \mathcal{G}$ in $\mathbf{Dyn}^1_n(M)$, where $\theta_r : \mathbb{R} \to \mathbb{R}$ is a continuous and monotonic increasing function. Then we observe
\begin{equation}\label{equ:inverse}
I^{\mathcal{F}}(\pi) = \{ t \in \mathbb{R} \mid \pi(\mathcal{G}_{\theta_r(t)}) = \mathcal{G}_{\theta_r(t)} \} = \theta_r^{-1}(I^{\mathcal{G}}(\pi)).
\end{equation}
This implies that the number of connected components (i.e., disjoint intervals) satisfies
\begin{equation}\label{equ:connected}
\#\pi_0(I^{\mathcal{F}}(\pi)) \leq \#\pi_0(I^{\mathcal{G}}(\pi)).
\end{equation}
In particular, if $I^{\mathcal{G}}(\pi)$ contains a left-unbounded interval $(-\infty, a]$ for some $a$, then $I^{\mathcal{F}}(\pi)$ must contain a similar interval $(-\infty, b]$ for some $b$. Likewise, if $I^{\mathcal{G}}(\pi)$ contains $[a, \infty)$, then $I^{\mathcal{F}}(\pi)$ contains $[b, \infty)$ for some $b$.

We now define a map $I^{\phi}(\pi) : I^{\mathcal{F}}(\pi) \to I^{\mathcal{G}}(\pi)$ by
\begin{equation*}
I^{\phi}(\pi)(x) = \theta_r(x).
\end{equation*}
In the following, we shall prove that this map is well-defined and that the assignment $\mathcal{F} \mapsto I^{\mathcal{F}}(\pi)$ extends to a functor.

\begin{proposition}\label{proposition:I_functor}
The construction $I(\pi): \mathbf{Dyn}^{1}_{n}(M)\to \mathbf{Int}(\mathbb{R})$ is functorial.
\end{proposition}

\begin{proof}
$(i)$
First, we will show that for any morphism $\phi: \mathcal{F} = \mathcal{G} \circ \theta_{r}\to \mathcal{G}$ in $\mathbf{Dyn}^{1}_n(M)$, the map $I^{\phi}(\pi): I^{\mathcal{F}}(\pi) \to I^{\mathcal{G}}(\pi)$ is a morphism in the category $\mathbf{Int}(\mathbb{R})$. Specifically, we need to prove that
\begin{itemize}
  \item $I^{\phi}(\pi)$ is a continuous, monotonic increasing function;
  \item $I^{\phi}(\pi)$ maps the intervals in $I^{\mathcal{F}}(\pi)$ to intervals in $I^{\mathcal{G}}(\pi)$ in an ordered manner.
\end{itemize}
Since $\theta_{r}$ is a continuous and monotonic increasing function, the induced map $I^{\phi}(\pi)$ is also continuous and monotonic increasing. Suppose that
\begin{equation*}
  I^{\mathcal{F}}(\pi) = \coprod_{i=1}^{m} A_i, \quad I^{\mathcal{G}}(\pi) = \coprod_{j=1}^{n} B_j
\end{equation*}
are disjoint unions of intervals $A_i= [a_i, a'_i]$ and $B_j= [b_j, b'_j]$, for $1 \leq i \leq m$ and $1 \leq j \leq n$, such that the upper bound of each preceding interval is strictly less than the lower bound of the subsequent one. We denote $A \prec A'$ if interval $A$ lies strictly to the left of $A'$, i.e., the upper bound of $A$ is smaller than the lower bound of $A'$. By Equation~\eqref{equ:connected}, we know that $m \leq n$. We aim to prove the assertion:
\[
\theta_r(A_i) \subseteq B_i \quad \text{for all } 1 \leq i \leq m.
\]

First, suppose for contradiction that $r < a_m$. Then we have
\begin{equation*}
  \theta_r(r) = \theta_r(a_m) \in I^{\mathcal{G}}(\pi),
\end{equation*}
which implies $[r, +\infty) \subseteq I^{\mathcal{F}}(\pi)$. Hence, one has $A_m \subset [r, +\infty)$, contradicting the assumption that $I^{\mathcal{F}}(\pi)$ is a disjoint union of intervals $A_i$. Therefore, we must have
\begin{equation*}
  a_m \leq r.
\end{equation*}

We proceed by induction. For the base case $\ell = 1$, it follows from Equation~\eqref{equ:inverse} and the continuity of $\theta_r$ that
\begin{equation*}
  \theta_r(A_1) \subseteq B_k
\end{equation*}
for some $1 \leq k \leq n$. Thus, we have $b_k \leq \theta_r(a_1) \leq \theta_r(r)$. This implies that the function $\theta_r$ is strictly increasing on $(-\infty, b_k]$.

If $k > 1$, then $B_1 \prec B_k$. Hence, the map $\theta_r$ is strictly increasing and continuous on the interval $\theta_r^{-1}(B_1)$. Let
\begin{equation*}
  \theta_r^{-1}(B_1) = [c, d] \subseteq I^{\mathcal{F}}(\pi).
\end{equation*}
By assumption, we have $a_1 \leq c$. Hence, one has
\begin{equation*}
  b_k \leq \theta_r(a_1) \leq \theta_r(c) \leq b_1'.
\end{equation*}
This contradicts $B_1 \prec B_k$, which implies $b_1' < b_k$. Therefore, we conclude that $k=1$ and
\begin{equation*}
  \theta_r(A_1) \subseteq B_1.
\end{equation*}

Now assume $\theta_r(A_p) \subseteq B_p$ for some $p < m$. Consider $\ell = p+1$. By Equation~\eqref{equ:inverse} and the continuity of $\theta_r$ , we have
\begin{equation*}
  \theta_r(A_{p+1}) \subseteq B_k
\end{equation*}
for some $p+1 \leq k \leq n$. If $k > p+1$, then by a similar argument as in the case $\ell = 1$, we arrive at a contradiction. Therefore,
\begin{equation*}
  \theta_r(A_{p+1}) \subseteq B_{p+1}.
\end{equation*}
By finite induction, we obtain
\[
\theta_r(A_i) \subseteq B_i \quad \text{for all } 1 \leq i \leq m.
\]

Hence, the map $I^{\phi}(\pi): I^{\mathcal{F}}(\pi) \to I^{\mathcal{G}}(\pi)$ defines a morphism in the category $\mathbf{Int}(\mathbb{R})$.

\vspace{1em}

$(ii)$
Secondly, we prove that $I(\pi)$ preserves composition and identity. For maps
\begin{equation*}
  \xymatrix{
   \mathcal{H} = \mathcal{F}  \circ \theta_r \circ \theta'_s\ar@{->}[r]^-{\phi} & \mathcal{F} \circ \theta_r \ar@{->}[r]^-{\psi} & \mathcal{F},
  }
\end{equation*}
we need to show that
\begin{equation*}
  I^{\psi}(\pi) \circ I^{\phi}(\pi) = I^{\psi \circ \phi}(\pi).
\end{equation*}
For any $x \in I^{\mathcal{H}}(\pi)$, note that
\begin{equation*}
  I^{\psi \circ \phi}(\pi)(x) = (\theta_r \circ \theta'_s)(x) = \theta_r \theta'_s(x)=\left(I^{\psi}(\pi) \circ I^{\phi}(\pi)\right)(x).
\end{equation*}
The construction $I(\pi)$ preserves composition.

Finally, when $ \phi: \mathcal{F} \to \mathcal{F} $ is the identity map, we observe that $ I^{\phi}(\pi) $ is also the identity map. Therefore, $ I(\pi) $ preserves identities.

In conclusion, the construction $I(\pi): \mathbf{Dyn}^{1}_{n}(M) \to \mathbf{Int}(\mathbb{R})$ is a functor.
\end{proof}

\begin{proposition}
Let $\mathcal{P},\mathcal{Q}: (\mathbb{R}, \leq) \to \mathbf{Dyn}^{1}_{n}(M)$ be two persistence dynamical systems. Then for any isometry $\pi$, we have
\begin{equation*}
  d_{I}(I(\pi)\mathcal{P},I(\pi)\mathcal{Q}) \leq d_{I}(\mathcal{P},\mathcal{Q}).
\end{equation*}
\end{proposition}
\begin{proof}
By Proposition \ref{proposition:I_functor}, we have the following composition of functors
\begin{equation*}
\xymatrix{
  (\mathbb{R}, \leq) \ar@<0.5ex>[r]^-{\mathcal{P}} \ar@<-0.5ex>[r]_-{\mathcal{Q}}& \mathbf{Dyn}^{1}_{n}(M)\ar@{->}[r]^{I(\pi)}& \mathbf{Int}(\mathbb{R}).
}
\end{equation*}
According to \cite[Proposition 3.6]{chazal2009proximity}, if $\mathcal{P}$ and $\mathcal{Q}$ are $\varepsilon$-interleaved, then $I(\pi)\mathcal{P}$ and $I(\pi)\mathcal{Q}$ are also $\varepsilon$-interleaved. Hence, the desired inequality follows immediately.
\end{proof}

Let $\mathcal{F} : (\mathbb{R}, \leq) \to \mathcal{S}_n(M)$ be a persistence $n$-configuration. For each isometry $\pi \in \Iso(M)$, the induced map $I^{\mathcal{F}}(\pi)$ is functorial in $\mathcal{F}$. This observation allows us to define a functor
\begin{equation*}
  \mathcal{B} : \mathbf{Dyn}^{1}_n(M) \to \mathbf{Polybarc}(M)
\end{equation*}
by setting $\mathcal{B}(\mathcal{F}) = \{ I^{\mathcal{F}}(\pi) \}_{\pi \in \Iso(M)}$.

\begin{theorem}[Stability theorem for polybarcodes I]
Let $\mathcal{P},\mathcal{Q}: (\mathbb{R}, \leq) \to \mathbf{Dyn}^{1}_{n}(M)$ be two persistence dynamical systems. Then we have
\begin{equation*}
  d_{I}(\mathcal{B}\mathcal{P},\mathcal{B}\mathcal{Q}) \leq d_{I}(\mathcal{P},\mathcal{Q}).
\end{equation*}
\end{theorem}
\begin{proof}
Consider the composition of functors
\begin{equation*}
\xymatrix{
  (\mathbb{R}, \leq)  \ar@<0.5ex>[r]^-{\mathcal{P}} \ar@<-0.5ex>[r]_-{\mathcal{Q}}& \mathbf{Dyn}^{1}_{n}(M)\ar@{->}[r]^{\mathcal{B}}& \mathbf{Polybarc}(M).
}
\end{equation*}
By \cite[Proposition 3.6]{chazal2009proximity}, we have the desired result.
\end{proof}

Now, let $\mathcal{F}: (\mathbb{R}, \leq) \to \mathcal{S}_{n}(M)$ be a persistence $n$-configuration. Recall that we have a persistence dynamical system $\kappa(\mathcal{F}):(\mathbb{R}, \leq) \to \mathbf{Dyn}^{0}_{n}(M)$ given by $\kappa(\mathcal{F})(r)=\mathcal{F}\circ \tau_{r}$. Here, $\tau_{r}(x)=\left\{
                                                                   \begin{array}{ll}
                                                                     x, & \hbox{$x\leq r$;} \\
                                                                     r, & \hbox{$x> r$}
                                                                   \end{array}
                                                                 \right.$ is the truncation function. One can embed $\mathbf{Dyn}^{0}_{n}(M)$ into the category $\mathbf{Dyn}^{1}_{n}(M)$ and obtain a sequence of functors
\begin{equation*}
\xymatrix{
  (\mathbb{R}, \leq)  \ar@{->}[r]^-{\kappa(\mathcal{F})}& \mathbf{Dyn}^{0}_{n}(M)\ar@{^{(}->}[r]^{\iota} &\mathbf{Dyn}^{1}_{n}(M)\ar@{->}[r]^{I(\pi)}& \mathbf{Int}(\mathbb{R}).
}
\end{equation*}

\begin{proposition}\label{proposition:inequation1}
Let $\mathcal{F},\mathcal{G}: (\mathbb{R}, \leq) \to \mathcal{S}_{n}(M)$ be two persistence $n$-configurations. Then we have
\begin{equation*}
  d_{I}(I(\pi)\circ \iota\circ \kappa(\mathcal{F}),I(\pi)\circ \iota\circ \kappa(\mathcal{G})) \leq d_{II}(\mathcal{F},\mathcal{G}).
\end{equation*}
\end{proposition}
\begin{proof}
Suppose the functors $\kappa(\mathcal{F})$ and $\kappa(\mathcal{G})$ are $\varepsilon$-interleaved. Since both $I(\pi)$ and $\iota$ are functors, and composition of functors preserves interleavings, it follows that $I(\pi) \circ \iota \circ \kappa(\mathcal{F})$ and $I(\pi) \circ \iota \circ \kappa(\mathcal{G})$ are also $\varepsilon$-interleaved. Therefore, we obtain
\begin{equation*}
   d_I(I(\pi) \circ \iota \circ \kappa(\mathcal{F}), I(\pi) \circ \iota \circ \kappa(\mathcal{G}))\leq  d_I(\kappa(\mathcal{F}), \kappa(\mathcal{G})) =d_{II}(\mathcal{F}, \mathcal{G}) .
\end{equation*}
This concludes the proof.
\end{proof}

\begin{lemma}\label{lemme:natural_injection}
For any isometry $\pi$, there is a natural monomorphism
\begin{equation*}
  Fr_{I^{\mathcal{F}}(\pi)} \Rightarrow I(\pi)\circ \iota \circ \kappa(\mathcal{F})
\end{equation*}
between functors from $(\mathbb{R}, \leq)$ to the category $\mathbf{Int}(\mathbb{R})$.
\end{lemma}

\begin{proof}
Let $r \in \mathbb{R}$ be a fixed real number. We begin by computing the values of the functors at $r$. First, we have
\begin{equation*}
  I(\pi)\circ \iota \circ \kappa(\mathcal{F})(r) =
  I(\pi)(\mathcal{F} \circ \tau_r) =
  \begin{cases}
    \{t \in I^{\mathcal{F}}(\pi) \mid t \leq r\} \cup [r, +\infty), & \text{if } r \in I^{\mathcal{F}}(\pi), \\
    \{t \in I^{\mathcal{F}}(\pi) \mid t \leq r\}, & \text{otherwise}.
  \end{cases}
\end{equation*}
On the other hand, we observe
\begin{equation*}
  Fr_{I^{\mathcal{F}}(\pi)}(r) = \{t \in I^{\mathcal{F}}(\pi) \mid t \leq r\}.
\end{equation*}
Define $\eta_r: Fr_{I^{\mathcal{F}}(\pi)}(r) \hookrightarrow I(\pi) \circ \iota \circ \kappa(\mathcal{F})(r)$ to be the inclusion map.

To show that $\{\eta_r\}_{r \in \mathbb{R}}$ forms a natural transformation, we verify the naturality condition: for any $r \leq s$, the following diagram commutes
\begin{equation*}
 \xymatrix{
Fr_{I^{\mathcal{F}}(\pi)}(r) \ar[d]_{Fr_{I^{\mathcal{F}}(\pi)}(f_{r,s})} \ar[r]^-{\eta_r} & I(\pi) \circ \iota \circ \kappa(\mathcal{F})(r) \ar[d]^{ I(\pi) \circ \iota \circ \kappa(\mathcal{F})(f_{r,s})} \\
Fr_{I^{\mathcal{F}}(\pi)}(s) \ar[r]^-{\eta_s} & I(\pi) \circ \iota \circ \kappa(\mathcal{F})(s)
}
\end{equation*}
Here, the vertical arrows are induced by the functors on $f_{r,s} : r \to s$. Explicitly, the map
\begin{equation*}
I(\pi) \circ \iota \circ \kappa(\mathcal{F})(f_{r,s}) : I(\pi) \circ \iota \circ \kappa(\mathcal{F})(r) \to I(\pi) \circ \iota \circ \kappa(\mathcal{F})(s)
\end{equation*}
is defined by
\begin{equation*}
  I(\pi) \circ \iota \circ \kappa(\mathcal{F})(f_{r,s})(x) = \left\{
                                                                   \begin{array}{ll}
                                                                     x, & \hbox{$x\leq r$;} \\
                                                                     r, & \hbox{$x> r$.}
                                                                   \end{array}
                                                                 \right.
\end{equation*}
By definition, $Fr_{I^{\mathcal{F}}(\pi)}(f_{r,s}):Fr_{I^{\mathcal{F}}(\pi)}(r)\to Fr_{I^{\mathcal{F}}(\pi)}(s)$ is an inclusion.
For any $x \in Fr_{I^{\mathcal{F}}(\pi)}(r)$, we have $x \leq r \leq s$. By evaluating both sides of the naturality condition at $x$, we obtain
\begin{equation*}
\eta_s \circ Fr_{I^{\mathcal{F}}(\pi)}(f_{r,s})(x) = x,
\quad
I(\pi) \circ \iota \circ \kappa(\mathcal{F})(f_{r,s}) \circ \eta_r(x) = x.
\end{equation*}
Therefore, the square commutes and the collection $\{\eta_r\}$ defines a natural monomorphism.
\end{proof}

\begin{lemma}\label{lemme:natural_surjection}
For any isometry $\pi$, there is a natural transformation
\begin{equation*}
  I(\pi)\circ \iota\circ \kappa(\mathcal{F}) \Rightarrow Fr_{I^{\mathcal{F}}(\pi)}
\end{equation*}
between functors from $(\mathbb{R}, \leq)$ to the category $\mathbf{Int}(\mathbb{R})$.
\end{lemma}

\begin{proof}
For a given real number $r \in \mathbb{R}$, recall that
\begin{equation*}
I(\pi) \circ \iota \circ \kappa(\mathcal{F})(r) = I(\pi)(\mathcal{F} \circ \tau_r) = \begin{cases}
Fr_{I^{\mathcal{F}}(\pi)}(r) \cup [r, +\infty), & \text{if } r \in I^{\mathcal{F}}(\pi), \\
Fr_{I^{\mathcal{F}}(\pi)}(r), & \text{otherwise.}
\end{cases}
\end{equation*}
Define the map $\rho_r: I(\pi) \circ \iota \circ \kappa(\mathcal{F})(r) \to Fr_{I^{\mathcal{F}}(\pi)}(r)$ as the retract map, which is given by
\begin{equation*}
\rho_r(x) = \begin{cases}
x, & \text{if } x \leq r, \\
r, & \text{if } x > r.
\end{cases}
\end{equation*}
To show that this defines a natural transformation, we verify the naturality condition.
\begin{equation*}
\xymatrix{
I(\pi) \circ \iota \circ \kappa(\mathcal{F})(r) \ar[d]_{ I(\pi) \circ \iota \circ \kappa(\mathcal{F})(f_{r,s})}  \ar[r]^-{\rho_r} & Fr_{I^{\mathcal{F}}(\pi)}(r) \ar[d]^{Fr_{I^{\mathcal{F}}(\pi)}(f_{r,s})} \\
I(\pi) \circ \iota \circ \kappa(\mathcal{F})(s) \ar[r]^-{\rho_s} & Fr_{I^{\mathcal{F}}(\pi)}(s)
}
\end{equation*}
For any morphism $f_{r,s}$ in $(\mathbb{R}, \leq)$ with $r \leq s$, we need to check the naturality condition
\begin{equation*}
\rho_s \circ  I(\pi) \circ \iota \circ \kappa(\mathcal{F})(f_{r,s}) = Fr_{I^{\mathcal{F}}(\pi)}(f_{r,s}) \circ \rho_r.
\end{equation*}
By definition, $Fr_{I^{\mathcal{F}}(\pi)}(f_{r,s})$ is the inclusion map from $Fr_{I^{\mathcal{F}}(\pi)}(r)$ to $Fr_{I^{\mathcal{F}}(\pi)}(s)$, while $I(\pi) \circ \iota \circ \kappa(\mathcal{F})(f_{r,s})$ is given by
\begin{equation*}
I(\pi) \circ \iota \circ \kappa(\mathcal{F})(f_{r,s})(x) = \begin{cases}
x, & \text{if } x \leq r, \\
r, & \text{if } x > r.
\end{cases}
\end{equation*}
For any $y \in I(\pi) \circ \iota \circ \kappa(\mathcal{F})(r)$, we compute both sides. On one hand, since $\rho_r(y) \leq r \leq s$, we have
\begin{equation*}
Fr_{I^{\mathcal{F}}(\pi)}(f_{r,s}) \circ \rho_{r}(y) = \rho_{r}(y) = \begin{cases}
y, & \text{if } y \leq r, \\
r, & \text{if } y > r.
\end{cases}
\end{equation*}
On the other hand, we have
\begin{equation*}
\rho_{s} \circ I(\pi) \circ \iota \circ \kappa(\mathcal{F})(f_{r,s})(y) =  \begin{cases}
y, & \text{if } y \leq r, \\
r, & \text{if } y > r.
\end{cases}
\end{equation*}
Thus, the naturality condition is satisfied, and we have established the natural transformation.
\end{proof}

\begin{proposition}\label{propostion:equal_distance}
Let $\mathcal{F}, \mathcal{G}: (\mathbb{R}, \leq) \to \mathcal{S}_{n}(M)$ be two persistence $n$-configurations.  Then for any isometry $\pi$, we have
\begin{equation*}
  d_{I}(Fr_{I^{\mathcal{F}}(\pi)},Fr_{I^{\mathcal{G}}(\pi)}) = d_{I}(I(\pi)\circ \iota\circ \kappa(\mathcal{F}),I(\pi)\circ \iota\circ \kappa(\mathcal{G})).
\end{equation*}
\end{proposition}
\begin{proof}
$(i)$
For simplicity, denote
\begin{equation*}
  \mathcal{I}^{\mathcal{F}} = Fr_{I^{\mathcal{F}}(\pi)} , \quad \mathcal{J}^{\mathcal{F}} = I(\pi)\circ \iota\circ \kappa(\mathcal{F}).
\end{equation*}
Let $\varepsilon = d_{I}(\mathcal{J}^{\mathcal{F}},\mathcal{J}^{\mathcal{G}})$. We aim to show that $\mathcal{I}^{\mathcal{F}}$ and $\mathcal{I}^{\mathcal{G}}$ are $\varepsilon$-interleaved.

By definition, there exist natural transformations
\begin{equation*}
  \Phi: \mathcal{J}^{\mathcal{F}} \Rightarrow \Sigma^{\varepsilon}\mathcal{J}^{\mathcal{G}},\quad \Psi: \mathcal{J}^{\mathcal{G}} \Rightarrow \Sigma^{\varepsilon}\mathcal{J}^{\mathcal{F}}
\end{equation*}
such that
\begin{equation*}
  (\Sigma^\varepsilon \Psi) \circ \Phi = \Sigma^{2\varepsilon}|_{\mathcal{J}^{\mathcal{F}}}, \quad (\Sigma^\varepsilon \Phi) \circ \Psi = \Sigma^{2\varepsilon}|_{\mathcal{J}^{\mathcal{G}}}.
\end{equation*}
For any real number $r$ and any point $x\in \mathcal{I}^{\mathcal{F}}(r)\subseteq \mathcal{J}^{\mathcal{F}}(r)$, the naturality of $\Phi$ gives the following commutative diagram
\begin{equation*}
  \xymatrix{
  \mathcal{J}^{\mathcal{F}}(x)\ar@{->}[r]^{\Phi_{x}} \ar@{->}[d]_{\mathcal{J}^{\mathcal{F}}(f_{x\to r})}& \mathcal{J}^{\mathcal{G}}(x+\varepsilon)\ar@{->}[d]^{\mathcal{J}^{\mathcal{G}}(f_{x+\varepsilon\to r+\varepsilon})}\\
  \mathcal{J}^{\mathcal{F}}(r)\ar@{->}[r]^{\Phi_{r}}&  \mathcal{J}^{\mathcal{G}}(r+\varepsilon),\\
  }
\end{equation*}
where $f_{x\to r}: x\to r$ denotes the morphism associated to the ordered pair $x\leq r$. Since $x\in \mathcal{J}^{\mathcal{F}}(x)$, we obtain that
\begin{equation*}
  (\Phi_{r}\circ\mathcal{J}^{\mathcal{F}}(f_{x\to r}))(x) = \Phi_{r}(x)
\end{equation*}
and
\begin{equation*}
  (\mathcal{J}^{\mathcal{G}}(f_{x+\varepsilon\to r+\varepsilon})\circ\Phi_{x})(x) = \left\{
                                                                   \begin{array}{ll}
                                                                     \Phi_{x}(x), & \hbox{$\Phi_{x}(x)\leq x+\varepsilon$;} \\
                                                                     x+\varepsilon, & \hbox{$\Phi_{x}(x)> x+\varepsilon$.}
                                                                   \end{array}
                                                                 \right.
\end{equation*}
The naturality of the above commutative diagram gives that
\begin{equation*}
  \Phi_{r}(x)=(\Phi_{r}\circ\mathcal{J}^{\mathcal{F}}(f_{x\to r}))(x)=(\mathcal{J}^{\mathcal{G}}(f_{x+\varepsilon\to r+\varepsilon})\circ\Phi_{x})(x)\leq r+\varepsilon.
\end{equation*}
This implies
\begin{equation}\label{equation:functor1}
  \Phi_{r}(\mathcal{I}^{\mathcal{F}}(r))\subseteq \mathcal{I}^{\mathcal{G}}(r+\varepsilon).
\end{equation}
 By a similar argument, we obtain
\begin{equation}\label{equation:functor2}
  \Psi_{r}(\mathcal{I}^{\mathcal{G}}(r))\subseteq \mathcal{I}^{\mathcal{F}}(r+\varepsilon).
\end{equation}
By Lemma \ref{lemme:natural_injection}, we have natural injections
\begin{equation}\label{equation:eta}
  \eta^{\mathcal{F}}:\mathcal{I}^{\mathcal{F}}\Rightarrow \mathcal{J}^{\mathcal{F}},\quad \eta^{\mathcal{G}}:\mathcal{I}^{\mathcal{G}}\Rightarrow \mathcal{J}^{\mathcal{G}}
\end{equation}
By Lemma \ref{lemme:natural_surjection}, we also have natural transformations
\begin{equation}\label{equation:rho}
  \rho^{\mathcal{F}}:\mathcal{J}^{\mathcal{F}}\Rightarrow \mathcal{I}^{\mathcal{F}},\quad \rho^{\mathcal{G}}:\mathcal{J}^{\mathcal{G}}\Rightarrow \mathcal{I}^{\mathcal{G}}.
\end{equation}
Here, $\rho^{\mathcal{F}}_{r}$ and $\rho^{\mathcal{G}}_{r}$ are surjection for each $r$.

Now, we define natural transformations
\begin{equation*}
  \phi = (\Sigma^{\varepsilon}\rho^{\mathcal{G}})\circ\Phi\circ\eta^{\mathcal{F}}: \mathcal{I}^{\mathcal{F}} \Rightarrow \Sigma^{\varepsilon}\mathcal{I}^{\mathcal{G}},\quad
  \psi = (\Sigma^{\varepsilon}\rho^{\mathcal{F}})\circ\Psi\circ\eta^{\mathcal{G}}: \mathcal{I}^{\mathcal{G}} \Rightarrow \Sigma^{\varepsilon}\mathcal{I}^{\mathcal{F}}.
\end{equation*}
Given a real number $r$, for $x\in \mathcal{I}^{\mathcal{F}}_{r}$, we have
\begin{align*}
    \left((\Sigma^{\varepsilon}\psi_{r})\circ\phi_{r}\right)(x) =& \left(\rho^{\mathcal{F}}_{r+2\varepsilon}\circ\Psi_{r+\varepsilon}\circ\eta^{\mathcal{G}}_{r+\varepsilon}\circ\rho^{\mathcal{G}}_{r+\varepsilon}\circ\Phi_{r}\circ\eta^{\mathcal{F}}_{r}\right)(x) \\
    =&\Psi_{r+\varepsilon}\Phi_{r}(x)\\
    =& \Sigma^{2\varepsilon}|_{\mathcal{I}^{\mathcal{F}}(r)}(x).
\end{align*}
In the above computation, we used Equations \eqref{equation:functor1} and \eqref{equation:functor2}. Thus, we obtain
\begin{equation*}
  (\Sigma^\varepsilon \psi) \circ \phi = \Sigma^{2\varepsilon}|_{\mathcal{I}^{\mathcal{F}}}
\end{equation*}
Similarly, we have
\begin{equation*}
  (\Sigma^\varepsilon \phi) \circ \psi = \Sigma^{2\varepsilon}|_{\mathcal{I}^{\mathcal{G}}}
\end{equation*}
Therefore, the natural transformations $\phi: \mathcal{I}^{\mathcal{F}} \Rightarrow \Sigma^{\varepsilon}\mathcal{I}^{\mathcal{G}}$ and $\psi: \mathcal{I}^{\mathcal{G}} \Rightarrow \Sigma^{\varepsilon}\mathcal{I}^{\mathcal{F}}$ provide an $\varepsilon$-interleaving between $\mathcal{I}^{\mathcal{F}}$ and $\mathcal{I}^{\mathcal{G}}$.  Hence, we obtain the inequality
\begin{equation}\label{equation:I_to_J}
  d_{I}(\mathcal{I}^{\mathcal{F}},\mathcal{I}^{\mathcal{G}})\leq d_{I}(\mathcal{J}^{\mathcal{F}},\mathcal{J}^{\mathcal{G}}).
\end{equation}

$(ii)$
On the other hand, suppose that $\varepsilon=d_{I}(\mathcal{I}^{\mathcal{F}},\mathcal{I}^{\mathcal{G}})$. Then there exist natural transformations $\phi': \mathcal{I}^{\mathcal{F}} \Rightarrow \Sigma^{\varepsilon}\mathcal{I}^{\mathcal{G}}$ and $\psi': \mathcal{I}^{\mathcal{G}} \Rightarrow \Sigma^{\varepsilon}\mathcal{I}^{\mathcal{F}}$ such that
\begin{equation*}
  (\Sigma^\varepsilon \psi') \circ \phi' = \Sigma^{2\varepsilon}|_{\mathcal{I}^{\mathcal{F}}}, \quad (\Sigma^\varepsilon \phi') \circ \psi' = \Sigma^{2\varepsilon}|_{\mathcal{I}^{\mathcal{G}}}.
\end{equation*}
Recalling Equations \eqref{equation:eta} and \eqref{equation:rho}, we define the natural transformations
\begin{equation*}
  \Phi' = (\Sigma^{\varepsilon}\eta^{\mathcal{G}})\circ\phi'\circ\rho^{\mathcal{F}}: \mathcal{J}^{\mathcal{F}} \Rightarrow \Sigma^{\varepsilon}\mathcal{J}^{\mathcal{G}},\quad
  \Psi' = (\Sigma^{\varepsilon}\eta^{\mathcal{F}})\circ\psi'\circ\rho^{\mathcal{G}}: \mathcal{J}^{\mathcal{G}} \Rightarrow \Sigma^{\varepsilon}\mathcal{J}^{\mathcal{F}}.
\end{equation*}
For any real number $r$ and $y\in \mathcal{J}^{\mathcal{F}}_{r}$, a straightforward computation yields
\begin{align*}
  (\Psi'_{r+\varepsilon}\circ\Phi'_{r})(y) =& (\eta^{\mathcal{F}}_{r+2\varepsilon}\circ\psi'_{r+\varepsilon}\circ\rho^{\mathcal{G}}_{r+\varepsilon}\circ\eta^{\mathcal{G}}_{r+\varepsilon}\circ\phi'_{r}\circ\rho^{\mathcal{F}}_{r})(y)\\
  = & (\eta^{\mathcal{F}}_{r+2\varepsilon}\circ\psi'_{r+\varepsilon}\circ\phi'_{r}\circ\rho^{\mathcal{F}}_{r})(y)\\
  = &  \psi'_{r+\varepsilon} \phi'_{r}\rho^{\mathcal{F}}_{r} (y)\\
  = & \Sigma^{2\varepsilon}|_{\mathcal{I}^{\mathcal{F}}(r)}\rho^{\mathcal{F}}_{r}(y).
\end{align*}
Next, we verify the identity
\begin{equation*}
  \Sigma^{2\varepsilon}|_{\mathcal{I}^{\mathcal{F}}(r)}\rho^{\mathcal{F}}_{r}(y) = \Sigma^{2\varepsilon}|_{\mathcal{J}^{\mathcal{F}}(r)}(y).
\end{equation*}
For $y\in \mathcal{I}^{\mathcal{F}}(r)\subseteq \mathcal{J}^{\mathcal{F}}(r)$, this identity clearly holds. For $y>r$, note that
\begin{equation*}
  \Sigma^{2\varepsilon}|_{\mathcal{J}^{\mathcal{F}}(r)}(y) = \mathcal{J}^{\mathcal{F}}(f_{r\to r+2\varepsilon})(y) =r
\end{equation*}
and
\begin{equation*}
   \Sigma^{2\varepsilon}|_{\mathcal{I}^{\mathcal{F}}(r)}\rho^{\mathcal{F}}_{r}(y) = \Sigma^{2\varepsilon}|_{\mathcal{I}^{\mathcal{F}}(r)}(r) = \mathcal{J}^{\mathcal{F}}(f_{r\to r+2\varepsilon})(r)=r.
\end{equation*}
Hence, we have $ \Sigma^{2\varepsilon}|_{\mathcal{I}^{\mathcal{F}}(r)}\rho^{\mathcal{F}}_{r}(y) = \Sigma^{2\varepsilon}|_{\mathcal{J}^{\mathcal{F}}(r)}(y)$ for any $y\in \mathcal{J}^{\mathcal{F}}(r)$. It follows that
\begin{equation*}
  \Psi'_{r+\varepsilon}\circ\Phi'_{r} = \Sigma^{2\varepsilon}|_{\mathcal{J}^{\mathcal{F}}(r)}.
\end{equation*}
A similar computation shows that
\begin{equation*}
  \Phi'_{r+\varepsilon}\circ\Psi'_{r} = \Sigma^{2\varepsilon}|_{\mathcal{J}^{\mathcal{G}}(r)}.
\end{equation*}
This proves that $\mathcal{J}^{\mathcal{F}}$ and $\mathcal{J}^{\mathcal{G}}$ are $\varepsilon$-interleaved. Hence, we obtain
\begin{equation}\label{equation:J_to_I}
  d_{I}(\mathcal{J}^{\mathcal{F}},\mathcal{J}^{\mathcal{G}})\leq d_{I}(\mathcal{I}^{\mathcal{F}},\mathcal{I}^{\mathcal{G}}).
\end{equation}

Combining Equations \eqref{equation:I_to_J} and \eqref{equation:J_to_I}, we conclude that
\begin{equation*}
  d_{I}(\mathcal{J}^{\mathcal{F}},\mathcal{J}^{\mathcal{G}})= d_{I}(\mathcal{I}^{\mathcal{F}},\mathcal{I}^{\mathcal{G}}),
\end{equation*}
which completes the proof.
\end{proof}

Proposition \ref{propostion:equal_distance} shows that although the constructions $I(\pi)\circ \iota\circ \kappa(\mathcal{F}):(\mathbb{R}, \leq) \to \mathbf{Int}(\mathbb{R})$ and $Fr_{I^{\mathcal{F}}(\pi)}:(\mathbb{R}, \leq) \to \mathbf{Int}(\mathbb{R})$ are different, they induce the same interleaving distance. The construction $I(\pi)\circ \iota\circ \kappa(\mathcal{F})$ is useful for establishing a connection between the interleaving distance on the category of dynamical systems and that on the category of polybarcodes. On the other hand, the construction $Fr_{I^{\mathcal{F}}(\pi)}$ facilitates the relation between the interleaving distance on the category of polybarcodes and the expansion distance.

\begin{theorem}\label{theorem:inequation}
Let $\mathcal{F}, \mathcal{G}: (\mathbb{R}, \leq) \to \mathcal{S}_{n}(M)$ be two persistence $n$-configurations. Then for any isometry $\pi$, we have
\begin{equation*}
   d_{L}(I^{\mathcal{F}}(\pi), I^{\mathcal{G}}(\pi)) \leq d_{II}(\mathcal{F}, \mathcal{G}).
\end{equation*}
\end{theorem}

\begin{proof}
By Proposition~\ref{proposition:inequation1} and Proposition~\ref{propostion:equal_distance}, we have
\begin{equation*}
  d_{I}(Fr_{I^{\mathcal{F}}(\pi)}, Fr_{I^{\mathcal{G}}(\pi)}) \leq d_{II}(\mathcal{F}, \mathcal{G}).
\end{equation*}
Furthermore, by Proposition~\ref{proposition:polybar_distance}, we know that
\begin{equation*}
  d_{L}(I^{\mathcal{F}}(\pi), I^{\mathcal{G}}(\pi)) = d_{I}(Fr_{I^{\mathcal{F}}(\pi)}, Fr_{I^{\mathcal{G}}(\pi)}).
\end{equation*}
Combining the two inequalities above yields the desired result.
\end{proof}

Now, given a persistence $n$-configuration $\mathcal{F}: (\mathbb{R}, \leq) \to \mathcal{S}_{n}(M)$, recall that we have a polybarcode functor
\[
\mathcal{B}: \mathbf{Dyn}^{1}_{n}(M)\to \mathbf{Polybarc}(M).
\]
The following result establishes the stability of polybarcodes associated with persistence $n$-configurations.

\begin{theorem}[Stability theorem for polybarcodes II]\label{theorem:polybarcodes_stability2}
Let $\mathcal{F}, \mathcal{G}: (\mathbb{R}, \leq) \to \mathcal{S}_{n}(M)$ be two persistence $n$-configurations. Then we have
\begin{equation*}
  d_{L}(\mathcal{B}(\mathcal{F}), \mathcal{B}(\mathcal{G}))= d_{I}(\mathcal{B}(\mathcal{F}), \mathcal{B}(\mathcal{G}))\leq d_{II}(\mathcal{F}, \mathcal{G}).
\end{equation*}
\end{theorem}
\begin{proof}
By Theorem~\ref{theorem:inequation}, we have
\begin{equation*}
  d_{L}(\mathcal{B}(\mathcal{F}), \mathcal{B}(\mathcal{G})) \leq d_{II}(\mathcal{F}, \mathcal{G}).
\end{equation*}
Moreover, by Propositions~\ref{proposition:expansion_interleaving} and~\ref{propostion:equal_distance}, we obtain
\begin{equation*}
  d_{I}(\mathcal{B}(\mathcal{F}), \mathcal{B}(\mathcal{G})) = d_{L}(\mathcal{B}(\mathcal{F}), \mathcal{B}(\mathcal{G})).
\end{equation*}
Combining these results completes the proof.
\end{proof} 

%% file: symmetry_types.tex
\section{Symmetry types}\label{section:symmetry_types}

In Euclidean space, two $n$-configurations $X$ and $Y$ may have isomorphic symmetry groups, yet their symmetry structures can differ. However, if $Y = gX$ for some isometry $g$, then $X$ and $Y$ intuitively share the same symmetry, since translations, rotations, or reflections preserve the symmetric structure of a configuration. The only difference lies in how the elements of the symmetry group correspond to specific linear transformations in space. This observation motivates the use of moduli space to study the symmetries of configurations: configurations lying in the same orbit under the isometry group share the same symmetry. To better capture and classify these symmetries, we introduce the notion of symmetry type, which enables a refined classification of configurations according to their symmetric structure.

\subsection{Moduli of symmetric configurations}

Throughout this section, all metric spaces are assumed to be locally compact. Let $(M,d)$ be a metric space. Let $G \leq \Iso(M)$ be a group acting on $M$ by isometries. The unordered configuration space of $n$ points in $M$ is defined by
\[
\Conf_n(M) = \{ X \subset M \mid |X| = n \}.
\]
Note that each element in $\Conf_n(M)$ is an object in $\mathcal{S}_{n}(M)$. The group $G$ acts on $\Conf_n(M)$ by
\[
g \cdot X = \{ g \cdot x \mid x \in X \}, \quad \text{for } g \in G, X\in \Conf_n(M).
\]
We define the associated moduli space of configurations up to isometry as the orbit space
\[
\mathcal{M}_n = \Conf_n(M)/G.
\]

Different configurations may possess different types of symmetry. To capture this, we classify configurations according to their stabilizer subgroups under the $G$-action.

\begin{definition}
For $X \in \Conf_n(M)$, the \textit{stabilizer subgroup} of $X$ in $G$ is defined as
\[
\mathrm{Stab}_G(X) := \{ g \in G \mid g \cdot X = X \}.
\]
Two $n$-configurations $X, Y \in \Conf_n(M)$ are said to have the same \textit{orbit type} if their stabilizer subgroups are conjugate in $G$, that is,
\[
\mathrm{Stab}_G(Y) = g \, \mathrm{Stab}_G(X) \, g^{-1} \quad \text{for some } g \in G.
\]
\end{definition}
In particular, when $G = \Iso(M)$, the stabilizer subgroup $\mathrm{Stab}_G(X)$ coincides with the symmetry group of $X$, that is, $\mathrm{Stab}_G(X) = \Sym(X)$.
\begin{definition}
Two $n$-configurations $X$ and $Y$ are said to have the same \textit{symmetry type} if their symmetry groups are conjugate subgroups of $\Iso(M)$.
\end{definition}
Although the point sets $X$ and $Y$ may differ significantly in their geometric realization, the abstract structure of their symmetry, measured up to conjugacy, remains the same.

Fix $G=\Iso(M)$. Let $\mathsf{SymType}$ denote the set of conjugacy classes of closed subgroups of $G$ that occur as stabilizers of configurations in $\Conf_n(M)$. For each $[H] \in \mathsf{SymType}$, define the stratum
\[
\mathcal{M}_n^{(H)} = \{ G \cdot X \in \mathcal{M}_n \mid \mathrm{Stab}_G(X) \in [H] \}.
\]
Each element of $\mathcal{M}_n^{(H)}$ is a $G$-orbit $G \cdot X$, where the configuration $X$ has a stabilizer subgroup $\mathrm{Stab}_G(X)$ lying in the conjugacy class $[H]$. The stratum $\mathcal{M}_n^{(H)}$ represents the moduli subspace of symmetry type $[H]$. Furthermore, $\mathcal{M}_n$ admits the stratification
\[
\mathcal{M}_n = \bigsqcup_{[H] \in \mathsf{SymType}} \mathcal{M}_n^{(H)},
\]
where each stratum corresponds to a distinct symmetry type under the action of $G$. Since a configuration consists of a finite set of $n$ points, its stabilizer subgroup in $G$ must be a finite group. In many natural settings (such as when $M=\mathbb{R}^d$), the set of conjugacy classes of finite subgroups of $\Iso(M)$ is countably infinite. This implies that, in general, there are infinitely many distinct symmetry types, forming a countable set.

The above stratification classifies $n$-configurations in the metric space according to their symmetry types and provides a different perspective for studying the evolution of symmetry under continuous deformation. This idea does not focus on the evolution of individual symmetries themselves, but rather on the changes in the overall symmetry structure of an $n$-configuration as a whole.

\subsection{The category of symmetry types}

Let $M$ be a metric space. We now fix the isometry group $\Iso(M)$ acting on $M$.

We first define a bicategory $\mathsf{SymType}$ of symmetry types as follows. The objects of $\mathsf{SymType}$ are the conjugacy classes of closed subgroups of the isometry group $\Iso(M)$. A 1-morphism from $[H]$ to $[H']$ in $\mathsf{SymType}$ is defined as an equivalence class of spans
\[
   K\xleftarrow{\phi} A \xrightarrow{\psi} K'
\]
where $K \in [H]$, $K' \in [H']$, $A \leq \Iso(M)$, and $\phi: A \rightarrow K, \psi: A \rightarrow K'$ are group homomorphisms. Two such spans
\[
K \xleftarrow{\phi} A \xrightarrow{\psi} K' \quad \text{and} \quad \tilde{K} \xleftarrow{\tilde{\phi}} \tilde{A} \xrightarrow{\tilde{\psi}} \tilde{K}'
\]
are said to be equivalent if there exists an isometry $g \in \Iso(M)$ such that
\[
\tilde{K} = gKg^{-1}, \quad \tilde{K}' = gK'g^{-1}, \quad \tilde{A} = gAg^{-1},
\]
and the following diagram commutes
\[
\xymatrix@R=1cm@C=1.2cm{
K\ar[d]_{c_g}&A \ar[r]^{\psi}\ar[l]_{\phi} \ar[d]_{c_g} & K' \ar[d]^{c_g} \\
\tilde{K}&\tilde{A} \ar[r]^{\tilde{\psi}}\ar[l]_{\tilde{\phi}} & \tilde{K}'
}
\]
where $c_g$ denotes conjugation by $g$, i.e., $c_g(h) = ghg^{-1}$. The 2-morphisms in $\mathsf{SymType}$ are given by equivalence classes of span morphisms.
The composition of 1-morphisms in $\mathsf{SymType}$ is defined via fibered product over the middle group. This composition is compatible with the equivalence relation, and identity morphisms exist for all objects. Thus, $\mathsf{SymType}$ is a bicategory as composition is associative only up to coherent isomorphism.

Objects in the bicategory $\mathsf{SymType}$ encode classes of symmetry types, while morphisms capture the relationships between them arising from different evolutions of $n$-configurations. A morphism between two distinct symmetry types $[K]$ and $[K']$ is represented by a span of the form
\[
   K \xleftarrow{\phi} A \xrightarrow{\psi} K',
\]
where the intermediate group $A$ plays a role in relating the two types. The homomorphisms $\phi: A \to K$ and $\psi: A \to K'$ provide an algebraic description of how the two symmetry types are connected through $A$.

We consider the construction
\[
\Theta : \mathbb{S}_n(M) \to \mathsf{SymType},
\]
which maps an $n$-configuration $X \subseteq M$ to its symmetry type $[\Sym(X)]$. Here,
\[
\Sym(X) = \{ g \in \Iso(M) \mid g(X) = X \}
\]
is the symmetry group of $X$ under the isometry group of $M$. 
Recall that a morphism $\gamma: X \to Y$ in the category $\mathbb{S}_n(M)$ is a path of morphisms $\gamma = (f_1, \dots, f_k)$. For a single step morphism $f: X \to Y$, we associate to it the equivalence class of the span
\[
   \Sym(X) \xleftarrow{f^{\flat}} \Sym_f(X) \xrightarrow{f^{\sharp}} \Sym(Y),
\]
where
\[
  \Sym_f(X) = \left\{ g \in \Sym(X) \mid f \circ g \circ f^{-1} \in \Sym(Y) \right\},
\]
the map $ f^{\flat} $ is the inclusion of subgroups, and $ f^{\sharp} $ is defined by
\[
f^{\sharp}(g) = f \circ g \circ f^{-1}, \quad \text{for all } g \in \Sym_f(X).
\]
For a general morphism $\gamma = (f_1, \dots, f_k)$, we define $\Theta(\gamma)$ as the equivalence class of the composed span obtained by taking iterated fibered products of the spans associated with each step $f_i$. 

\begin{proposition}\label{proposition:symmetry_type_functor}
The construction $\Theta : \mathbb{S}_n(M) \to \mathsf{SymType}$ is a pseudofunctor.
\end{proposition}
\begin{proof}
Observe that the construction $\Theta$ can be decomposed as
\[
\xymatrix{
  \mathbb{S}_{n}(M) \ar[r]^-{\Sym} & \mathrm{Span}(\mathbf{subGrp}(\Iso(M))) \ar[r] & \mathsf{SymType},
}
\]
where $\mathbf{subGrp}(\Iso(M))$ denotes the category of subgroups of $\Iso(M)$ with morphisms given by group homomorphisms. By Theorem~\ref{theorem:lax_functor}, the assignment $\Sym$, which maps a path of configurations to the composed span of their symmetry groups, defines a pseudofunctor from $\mathbb{S}_n(M)$ to $\mathrm{Span}(\mathbf{subGrp}(\Iso(M)))$. Furthermore, the quotient of the span bicategory by conjugation of subgroups is well-defined, as conjugation preserves group homomorphisms and commutes with the formation of fibered products. This quotient process naturally yields a quotient pseudofunctor. Combining these two facts, we conclude that the overall construction $\Theta$ defines a pseudofunctor.
\end{proof}

We define the construction $\widehat{\Theta}: \mathbf{Dyn}^{0}_{n}(M) \to \mathsf{SymType}$ as follows. 
For an object $(\mathcal{F}, r)$ in $\mathbf{Dyn}^{0}_{n}(M)$, we map it to the symmetry type of its stabilized configuration:
\[
\widehat{\Theta}(\mathcal{F}, r) = [\Sym(\mathcal{F}_r)] \in \mathsf{SymType}.
\]
For a morphism $\phi : (\mathcal{F}, r) \to (\mathcal{G}, s)$ in $\mathbf{Dyn}^{0}_{n}(M)$ (which requires $r \leq s$ and $\mathcal{F}_t = \mathcal{G}_t$ for all $t \leq r$), the underlying evolution from $\mathcal{F}_r$ to $\mathcal{G}_s$ is given by the structure map $\mathcal{G}_{r \leq s}$. We associate to $\phi$ the equivalence class of the span of survival symmetry groups:
\[
\Sym(\mathcal{F}_r) \xleftarrow{\phi^{\flat}} L_{r,s} \xrightarrow{\phi^{\sharp}} \Sym(\mathcal{G}_s),
\]
where $L_{r,s} = \left\{ g \in \Sym(\mathcal{F}_r) \mid \mathcal{G}_{r \leq s}(g) \in \Sym(\mathcal{G}_s) \right\}$ is the group of symmetries that survive the evolution from $r$ to $s$. The map $\phi^{\flat}$ is the inclusion of subgroups, and $\phi^{\sharp}(g) = \mathcal{G}_{r \leq s} \circ g \circ \mathcal{G}_{r \leq s}^{-1}$ is the induced monomorphism.

\begin{proposition}\label{proposition:dyn_to_symtype}
The construction $\widehat{\Theta}: \mathbf{Dyn}^{0}_{n}(M) \to \mathsf{SymType}$ is a pseudofunctor.
\end{proposition}

\begin{proof}
By Proposition \ref{proposition:pseudofunctor}, the assignment $\Sym$, which maps a morphism $\phi: (\mathcal{F}, r) \to (\mathcal{G}, s)$ to the span of its survival symmetry groups $L_{r,s}$, defines a pseudofunctor $\Sym: \mathbf{Dyn}^{0}_{n}(M) \to \mathrm{Span}(\mathbf{subGrp}(\Iso(M)))$. 

Furthermore, as established in the proof of Proposition \ref{proposition:symmetry_type_functor}, taking the quotient by conjugation of subgroups yields a well-defined quotient pseudofunctor $Q: \mathrm{Span}(\mathbf{subGrp}(\Iso(M))) \to \mathsf{SymType}$. 

The construction $\widehat{\Theta}$ is exactly the composition $Q \circ \Sym$. Since the composition of pseudofunctors is again a pseudofunctor, we conclude that $\widehat{\Theta}: \mathbf{Dyn}^{0}_{n}(M) \to \mathsf{SymType}$ is a pseudofunctor.
\end{proof}

It is important to note that for an arbitrary closed subgroup $H$ of $\Iso(M)$, there does not necessarily exist an $n$-point configuration $X$ such that $\Theta(X) = [H]$. In fact, since an $n$-configuration consists of finitely many points, its symmetry group must be a finite subgroup of $\Iso(M)$. Thus, only finite subgroups can be realized as stabilizers of $n$-configurations. In certain settings, such as when $M$ is a Euclidean space, it is known that for every finite subgroup $H \leq \Iso(M)$, there exists a finite configuration whose symmetry type is $[H]$. For more on the realization of topological groups as stabilizers, one can refer to \cite{gao2003classification,melleray2006stabilizers}.

\subsection{Related topological properties}

Let $S$ and $T$ be topological spaces. A \textit{set-valued map} is a function $f : S\to T$ that assigns to each point $x \in S$ a subspace $f(x) \subseteq T$.
A set-valued map $f : S \to T$ is said to be \textit{upper semi-continuous} if for every open set $U \subseteq T$, the set $\{ x \in S \mid f(x) \subseteq U \}$ is open in $S$. Equivalently, for every closed set $F \subseteq T$, the set $\{ x \in S \mid f(x) \cap F \neq \emptyset \}$ is closed in $S$. Similarly, $f$ is said to be \textit{lower semi-continuous} if for every open set $U \subseteq T$, the set $\{ x \in S \mid f(x) \cap U \neq \emptyset \}$ is open in $S$ \cite{aubin2012mutational,aubin1999set}.

In the preceding sections, $\Sym$ was primarily viewed as a pseudofunctor mapping paths of configurations to spans of symmetry groups. From a topological perspective, $\Sym$ can also be naturally regarded as a set-valued map that assigns to each configuration its symmetry group as a subspace of $\Iso(M)$.

\begin{lemma}\label{lemma:upper_continuous}
Suppose $M$ is a compact metric space. Then the set-valued map
\[
\Sym \colon \Conf_n(M) \to \Iso(M) ,\quad X \mapsto \Sym(X)
\]
is upper semi-continuous.
\end{lemma}

\begin{proof}
Let $U \subseteq \Iso(M)$ be an open subset. We aim to show that the set
\[
\{ X \in \Conf_n(M) \mid \Sym(X) \subseteq U \}
\]
is open in $\Conf_n(M)$. Fix a configuration $X_0 \in \Conf_n(M)$ such that $\Sym(X_0) \subseteq U$. It suffices to find an open neighborhood $\mathcal{U}$ of $X_0$ in $\Conf_n(M)$ such that for all $X \in \mathcal{U}$, we have $\Sym(X) \subseteq U$.

Let us denote
\[
F = \{(g, X) \in \Iso(M) \times \Conf_n(M) \mid g \cdot X = X \}.
\]
Since the action of $\Iso(M)$ on $\Conf_n(M)$ is continuous and $\Conf_n(M)$ is Hausdorff, the set $F$ is closed in $\Iso(M) \times \Conf_n(M)$.

Let $K = \Iso(M) \setminus U$. For any $g \in K$, we have $g \notin \Sym(X_0)$. It follows that $(g, X_0) \notin F$. Since $F$ is closed, its complement
\[
(\Iso(M) \times \Conf_n(M)) \setminus F
\]
is open and contains $(g, X_0)$. The space $\Iso(M) \times \Conf_n(M)$ carries the product topology, whose basis consists of sets of the form $V \times W$, where $V \subseteq \Iso(M)$ and $W \subseteq \Conf_n(M)$ are open. Therefore, there exist open neighborhoods
\[
\mathcal{V}_g \subseteq \Iso(M) \quad \text{of } g, \quad \text{and} \quad \mathcal{U}_g \subseteq \Conf_n(M) \quad \text{of } X_0,
\]
such that
\[
\mathcal{V}_g \times \mathcal{U}_g \subseteq (\Iso(M) \times \Conf_n(M)) \setminus F.
\]
This implies that for each $g \in K$, there exist open neighborhoods $\mathcal{V}_g$ of $g$ and $\mathcal{U}_g$ of $X_0$ such that for all $X \in \mathcal{U}_g$, we have $\Sym(X) \cap \mathcal{V}_g = \emptyset$.

Since $M$ is a compact metric space, the Arzel\`{a}-Ascoli theorem implies that the isometry group $\Iso(M)$, equipped with the compact-open topology, is compact. As $K$ is closed in $\Iso(M)$, it follows that $K$ is compact. Therefore, the open cover $\{\mathcal{V}_g\}_{g \in K}$ of $K$ admits a finite subcover, say $\mathcal{V}_{g_1}, \dots, \mathcal{V}_{g_k}$. Correspondingly, we obtain open neighborhoods $\mathcal{U}_{g_1}, \dots, \mathcal{U}_{g_k}$ of $X_0$ such that for all $X \in \mathcal{U}_{g_j}$,
\[
\Sym(X) \cap \mathcal{V}_{g_j} = \emptyset.
\]

Now, we denote
\[
\mathcal{U} = \bigcap_{j=1}^k \mathcal{U}_{g_j}.
\]
Then for any $X \in \mathcal{U}$, we have $\Sym(X) \cap \mathcal{V}_{g_j} = \emptyset$ for all $j$, and hence $\Sym(X) \cap K = \emptyset$. It follows that $\Sym(X) \subseteq U$.

Therefore, $\Sym$ is upper semi-continuous, as desired.
\end{proof}

\begin{example}\label{example:open_dense}
Let $M$ be a compact metric space, and denote $G = \Iso(M)$. For a non-trivial element $g \in G$, the set $\{g\} \subseteq G$ is a closed subset of $G$. By the upper semi-continuity of the map $\mathrm{Stab} : \Conf_n(M) \to \mathrm{Sub}(\Iso(M))$, the subset
\[
A_{g} = \left\{ X \in \Conf_n(M) \mid \mathrm{Stab}(X) \cap \{g\} \neq \emptyset \right\}
\]
is closed. This means that the set of all $n$-configurations in $\Conf_n(M)$ whose stabilizers contain the symmetry $g$ is closed. Let $U_{g}$ be the complement of $A_{g}$. Then, $ U_{g}$ is an open set. Hence, the set
\[
U = \bigcup_{g \neq e} U_{g}
\]
is an open set. Note that
\begin{equation*}
  U = \left\{ X \in \Conf_n(M) \mid \mathrm{Stab}(X) = \{e\} \right\}.
\end{equation*}
Therefore, the set of $n$-configurations with only the trivial symmetry is an open set. Additionally, if $M$ is a connected compact manifold, then $U$ is dense in $\Conf_n(M)$, as any $n$-configuration can be slightly perturbed into one with a trivial symmetry group.
\end{example}

For any point $p \in M$, we define the subgroup of isometries fixing $p$ by
\begin{equation*}
  \Iso_{p}(M) = \{ g \in \Iso(M) \mid g(p) = p \}.
\end{equation*}
Now let $M = \mathbb{R}^k$ be a Euclidean space. Consider the map
\begin{equation*}
  \mu \colon \Conf_n(\mathbb{R}^k) \to \mathbb{R}^k,
\end{equation*}
which sends a configuration $X$ to its barycenter
\[
  \mu(X) = \frac{1}{n} \sum_{x \in X} x.
\]
For each fixed point $p \in \mathbb{R}^k$, define
\begin{equation*}
  H_p(k) = \mu^{-1}(p).
\end{equation*}
Then $H_p(k)$ lies in an affine hyperplane of $\mathbb{R}^{kn}$ and forms a submanifold of $\Conf_n(\mathbb{R}^k)$ of codimension $k$. This gives rise to a set-valued map, also denoted by
\begin{equation*}
  \Sym \colon H_p(k) \to \Iso_p(\mathbb{R}^k),
\end{equation*}
associating to each configuration its symmetry group fixing the barycenter.

In particular, when $p = O$ is the origin, we have $\Iso_p(\mathbb{R}^k) = O(k)$, the orthogonal group. Let us write $H(k) = H_p(k)$ in this case. The above map then reduces to
\begin{equation*}
  \Sym \colon H(k) \to O(k).
\end{equation*}

Since any two configurations that differ by a translation have stabilizer subgroups conjugate via that translation, they share the same symmetry type. Consequently, it suffices to restrict our study of symmetry types to the space $H(k)$ of configurations with barycenter at the origin.

\begin{lemma}\label{lemma:usc_Hk}
The set-valued map
\[
\Sym : H(k) \to O(k), \quad X \mapsto \Sym(X)
\]
is upper semi-continuous.
\end{lemma}

\begin{proof}
Note that $O(k)$ is a compact Lie group. The proof is similar to the proof of Lemma \ref{lemma:upper_continuous}.
\end{proof}

\begin{lemma}\label{lemma:closed_symmetry}
Let $V \subseteq O(k)$ be a subset. Then the set
\[
A_V = \left\{ X \in H(k) \mid \Sym(X) \supseteq V \right\}
\]
is closed.
\end{lemma}
\begin{proof}
For each $g \in O(k)$, let us denote
\[
A_g = \left\{ X \in H(k) \mid g \in \Sym(X) \right\}.
\]
Since the map $\Sym : H(k) \to O(k)$ is upper semi-continuous, the set $A_g$ is closed. Consequently, for any subset $V \subseteq O(k)$, the set
\[
A_V = \bigcap_{g \in V} A_g = \left\{ X \in H(k) \mid V \subseteq \Sym(X) \right\}
\]
is also closed.
\end{proof}

A finite group $V \subseteq O(k)$ is said to be \textit{maximal} if it is not properly contained in the symmetry group of any other configuration.

\begin{proposition}\label{proposition:locally_closed}
Let $V \subseteq O(k)$ be a finite nontrivial subgroup. Then the set
\[
D_V = \left\{ X \in H(k) \mid \Sym(X) = V \right\}
\]
is either empty or locally closed in $H(k)$. In particular, if $V$ is maximal, then $D_V$ is closed in $H(k)$.
\end{proposition}

\begin{proof}
By Lemma~\ref{lemma:closed_symmetry}, for any subgroup $W \subseteq O(k)$, the set
\[
A_W = \left\{ X \in H(k) \mid W \subseteq \Sym(X) \right\}
\]
is closed in $H(k)$.

Note that for any $n$-point configuration $X$, the symmetry group $\Sym(X)$ acts faithfully on the finite set $X$. Thus, $\Sym(X)$ is isomorphic to a subgroup of the symmetric group $S_n$, which implies that the order of $\Sym(X)$ is bounded by $n!$. Therefore, any finite subgroup $W \subseteq O(k)$ containing $V$ and realizable as a symmetry group in $H(k)$ must have order at most $n!$.

Since $O(k)$ is a compact Lie group, it contains only finitely many finite subgroups of any bounded order. Consequently, there are only finitely many such subgroups $W \supsetneq V$ that can possibly be realized as symmetry groups in $H(k)$. Therefore, the union
\[
C_V = \bigcup_{\substack{V < W \\ W \subseteq O(k)}} A_W
\]
is a finite union of closed sets and hence closed. It follows that
\[
D_V = A_V \setminus C_V,
\]
i.e., $D_V$ is the difference of two closed sets in $H(k)$ and hence locally closed.
\end{proof}

\begin{example}\label{example:perturbation}
Let $M = \mathbb{R}^2$, and let $G = \Iso(\mathbb{R}^2)$ be the group of Euclidean isometries. Consider the configuration space $\Conf_3(\mathbb{R}^2)$ of unlabeled triples of distinct points in the plane. Let $H(2)$ be the subspace of $\Conf_3(\mathbb{R}^2)$ consisting of configurations whose centroid is at the origin.

Let $X_0 = \{A, B, C\} \subset \mathbb{R}^2$ be an equilateral triangle with centroid at the origin given by
\begin{equation*}
  A=\left( 2,  0 \right),\quad B = \left( -1,  \sqrt{3} \right), C = \left( -1,   -\sqrt{3} \right).
\end{equation*}
Then the symmetry group $\Sym(X_0) < G$ is isomorphic to the dihedral group $D_3$ of order six, consisting of three rotations and three reflections that preserve the triangle.

Now, for each small $\varepsilon > 0$, let $X_\varepsilon=\{A_\varepsilon,B_\varepsilon,C_\varepsilon\}$ be a configuration obtained by slightly perturbing the points to break the symmetry of the triangle, so that $X_\varepsilon$ is an asymmetric triangle with centroid at the origin. The coordinates of the perturbed points $A_\varepsilon, B_\varepsilon, C_\varepsilon$ are given by
\begin{equation*}
  A_\varepsilon = \left( 2+\varepsilon, \varepsilon \right), \quad
  B_\varepsilon = \left( -1-\varepsilon,  \sqrt{3} \right), \quad
  C_\varepsilon = \left( -1,   -\sqrt{3}-\varepsilon \right).
\end{equation*}
Then the symmetry group $\Sym(X_\varepsilon)$ is trivial (i.e., the identity isometry is the only one preserving the configuration).

As $\varepsilon \to 0$, we have $X_\varepsilon \to X_0$ in $H(2)$, and hence
\[
\lim_{\varepsilon \to 0} \Sym(X_\varepsilon) = \{e\} \neq D_3 = \Sym(X_0).
\]
The set
\begin{equation*}
  D = \left\{ X \in H(2) \mid \Sym(X) = D_3 \right\}
\end{equation*}
is not open, since arbitrarily small perturbations of $X_0$ leave the set $D$; however, it is closed. On the other hand, the set
\begin{equation*}
  U = \left\{ X \in H(2) \mid \Sym(X) = \{e\} \right\}
\end{equation*}
is not closed, as the limit of configurations with trivial symmetry may lie in $D$, but it is open by the upper semi-continuity of the map $\Sym \colon H(2) \to O(2)$.
\end{example}

\begin{definition}[Chabauty topology \cite{benedetti1992lectures,chabauty1950limite}]
Let $ G $ be a locally compact Hausdorff topological group. Denote by $\mathrm{Sub}(G)$ the set of all closed subgroups of $G$. The \textit{Chabauty topology} on $\mathrm{Sub}(G)$ is the topology generated by the subbasis consisting of all sets of the form
\[
\mathcal{O}_K = \{ H \in \mathrm{Sub}(G) \mid H \cap K = \emptyset \},
\]
where $K \subseteq G$ is compact, and all sets of the form
\[
\mathcal{O}^U = \{ H \in \mathrm{Sub}(G) \mid H \cap U \neq \emptyset \},
\]
where $U \subseteq G$ is open.

With this topology, the space $\mathrm{Sub}(G)$ is compact, Hausdorff, and metrizable whenever $G$ is second-countable.
\end{definition}

\begin{definition}[Quotient Chabauty topology on conjugacy classes]
Let $ G $ be a locally compact Hausdorff topological group. Denote by $ \mathrm{Conj}(G) $ the set of all conjugacy classes of closed subgroups of $ G $. The \textit{quotient Chabauty topology} on $ \mathrm{Conj}(G) $ is defined as the quotient topology induced by the Chabauty topology on $ \mathrm{Sub}(G) $ under the conjugation action of $ G $. That is, a set $ \mathcal{U} \subseteq \mathrm{Conj}(G) $ is open if and only if
\[
p^{-1}(\mathcal{U}) = \{ H \in \mathrm{Sub}(G) \mid [H] \in \mathcal{U} \}
\]
is open in $ \mathrm{Sub}(G) $, where $ p: \mathrm{Sub}(G) \to \mathrm{Conj}(G) $ is the canonical projection mapping a subgroup to its conjugacy class.
\end{definition}

\begin{remark}
It is worth noting that $\mathrm{Conj}(G)$ is not a Hausdorff space in general. Even when $G$ is a Euclidean space, the space $\mathrm{Conj}(G)$ may still fail to be Hausdorff. As a consequence, a conjugacy class $[H]$, viewed as a singleton in $\mathrm{Conj}(G)$, is not necessarily a closed point.
\end{remark}

From now on, we denote $|\mathsf{SymType}| \subseteq \mathrm{Conj}(\Iso(M))$ as the topological space underlying the category of symmetry types. We also define the map
\begin{equation*}
  \widetilde{\Theta} : \Conf_n(M) \to |\mathsf{SymType}|,\quad X \mapsto [\Sym(X)],
\end{equation*}
which assigns to each configuration its conjugacy class of symmetry groups.

Furthermore, the space $|\mathsf{SymType}|$ carries a natural partial order defined by $[H]\leq [H']$ if there exists an element $g\in G$ such that $gHg^{-1}\leq H'$. This partial order reflects the relative richness of symmetries between different symmetry types.

\begin{theorem}\label{theorem:symmetry_type}
Let $M$ be a compact metric space.
Let $\widetilde{\Theta} : \Conf_n(M) \to |\mathsf{SymType}|$ be the map given by $\widetilde{\Theta}(X) = [\Sym(X)]$. For each $\omega \in |\mathsf{SymType}|$, denote
\[
  \Conf_n^{\omega}(M) = \{ X \in \Conf_n(M) \mid \widetilde{\Theta}(X) = \omega \}.
\]
Then we have
\begin{enumerate}[label=\textit{(\roman*)}]
  \item The space $\Conf_n(M)$ admits a disjoint decomposition into these subsets
  \[
  \Conf_n(M) = \bigsqcup_{\omega \in |\mathsf{SymType}|} \Conf_n^{\omega}(M).
  \]
  \item For any $\omega, \omega' \in |\mathsf{SymType}|$, if
  \[
  \overline{\Conf_n^{\omega}(M)} \cap \Conf_n^{\omega'}(M) \neq \emptyset,
  \]
  then $ \omega \leq \omega'$ in the natural partial order on $|\mathsf{SymType}|$.
\end{enumerate}
\end{theorem}

\begin{proof}
$(i)$ The map $\widetilde{\Theta}$ assigns to each configuration $X \in \Conf_n(M)$ a unique symmetry type $\omega = [\Sym(X)]$, where $\Sym(X)$ is the symmetry group of $X$. Consequently, the fibers $\Conf_n^{\omega}(M)$ form a partition of $\Conf_n(M)$. In other words, we have a disjoint decomposition
\[
\Conf_n(M) = \bigsqcup_{\omega \in |\mathsf{SymType}|} \Conf_n^{\omega}(M).
\]

$(ii)$ Suppose that
\[
\overline{\Conf_n^{\omega}(M)} \cap \Conf_n^{\omega'}(M) \neq \emptyset.
\]
Then there exists a configuration $X \in \Conf_n(M)$ such that
\[
X \in \overline{\Conf_n^{\omega}(M)} \quad \text{and} \quad X \in \Conf_n^{\omega'}(M).
\]
Since $X \in \overline{\Conf_n^{\omega}(M)}$, it follows that $X$ is either in $\Conf_n^{\omega}(M)$ or is a limit of configurations from $\Conf_n^{\omega}(M)$. Thus, there exists a sequence $\{X_k\} \subseteq \Conf_n^{\omega}(M)$ such that $X_k \to X$ as $k \to \infty$. For each $k$, we have
\[
\widetilde{\Theta}(X_k) = [\Sym(X_k)] = \omega.
\]
On the other hand, since $X \in \Conf_n^{\omega'}(M)$, we know
\[
\widetilde{\Theta}(X) = [\Sym(X)] = \omega'.
\]

Let $H$ be a fixed representative of the conjugacy class $\omega$. Since $\widetilde{\Theta}(X_k) = \omega$, for each $k$ there exists an isometry $g_k \in \Iso(M)$ such that $\Sym(X_k) = g_k H g_k^{-1}$. 

Because $M$ is a compact metric space, $\Iso(M)$ is compact. Therefore, the sequence $\{g_k\}$ in $\Iso(M)$ has a convergent subsequence. Without loss of generality (by passing to this subsequence), we may assume $g_k \to g \in \Iso(M)$. 

By the continuity of the conjugation action, the sequence of subgroups $g_k H g_k^{-1}$ converges to $g H g^{-1}$ in the Chabauty topology. Since $\Sym(X_k) = g_k H g_k^{-1}$, this implies that $\limsup_{k \to \infty} \Sym(X_k) = g H g^{-1}$.

Recall that the set-valued map $\Sym : \Conf_n(M) \to \Iso(M)$ is upper semi-continuous. This implies that
\[
\limsup_{k \to \infty} \Sym(X_k) \subseteq \Sym(X).
\]
Therefore, we obtain
\[
g H g^{-1} \subseteq \Sym(X).
\]
By conjugating both sides with $g^{-1}$, we get $H \subseteq g^{-1} \Sym(X) g$. By the definition of the partial order on $|\mathsf{SymType}|$, this means $\omega = [H] \leq [\Sym(X)] = \omega'$.
\end{proof}

Theorem \ref{theorem:symmetry_type} does not imply that $\omega = \omega'$, because configurations in the closure $\overline{\Conf_n^{\omega}(M)}$ that are not in $\Conf_n^{\omega}(M)$ may acquire larger symmetry groups in the limit. A special case is when $\omega = \{e\}$ is the trivial group. For any $\omega' \in |\mathsf{SymType}|$, Example \ref{example:open_dense} shows that
\[
\overline{\Conf_n^{[e]}(M)} \cap \Conf_n^{\omega'}(M) \neq \emptyset.
\]
Hence, we have $[\{e\}] \leq [\omega']$, but equality $[\{e\}] = [\omega']$ does not necessarily hold. This illustrates the phenomenon described in Theorem~\ref{theorem:symmetry_type}, where limits of configurations may acquire larger symmetry groups, leading to a strict inequality in the partial order.

Let $M = \mathbb{R}^k$ be Euclidean space. Consider the map $\Sym : H(k) \to O(k)$, which sends an $n$-configuration with centroid at the origin to its symmetry group. This gives a map
\[
\widetilde{\Theta} : H(k) \to |\mathsf{SymType}|,
\]
where $\mathsf{SymType}$ denotes the set of conjugacy classes of finite subgroups of $O(k)$. Analogous to Theorem~\ref{theorem:symmetry_type}, we obtain the following result.

\begin{theorem}\label{theorem:symmetry_type2}
For each $\omega \in |\mathsf{SymType}|$, denote
\[
  H^{\omega}(k) = \{ X \in H(k) \mid \widetilde{\Theta}(X) = \omega \}.
\]
Then we have
\begin{enumerate}[label=\textit{(\roman*)}]
  \item The space $H(k)$ admits a disjoint decomposition into these subsets
  \[
   H(k) = \bigsqcup_{\omega \in |\mathsf{SymType}|} H^{\omega}(k).
  \]
  \item For any $\omega, \omega' \in |\mathsf{SymType}|$, if
  \[
  \overline{H^{\omega}(k)} \cap H^{\omega'}(k) \neq \emptyset,
  \]
  then $ \omega \leq \omega'$ in the natural partial order on $|\mathsf{SymType}|$.
\end{enumerate}
\end{theorem}

\subsection{Persistent symmetry types}

We now consider the variation of symmetry types under evolution of configurations. This leads to the notion of persistent symmetry types, which captures how symmetry evolves in families of configurations.

\begin{definition}
A \textit{persistent symmetry type} is a pseudofunctor $\mathcal{P}: (\mathbb{R}, \leq) \to \mathsf{SymType}$, where $\mathsf{SymType}$ is the bicategory of symmetry types.

We say that the persistent symmetry type $\mathcal{P}: (\mathbb{R}, \leq) \to \mathsf{SymType}$ is \textit{constant on an interval} $ I \subseteq \mathbb{R} $ if $\mathcal{P}_t = \mathcal{P}_{t'}$ for all $ t, t' \in I $.

A parameter $ t_0 \in \mathbb{R}$ is called a \textit{transition point} of $\mathcal{P}$ if there is no open interval containing $t_0$ on which $\mathcal{P}$ is constant. Equivalently, for every $\varepsilon > 0$, there exist $t_1, t_2 \in (t_0 - \varepsilon, t_0 + \varepsilon)$ such that $\mathcal{P}_{t_1} \neq \mathcal{P}_{t_2}$.
\end{definition}

Let $\mathcal{F}: (\mathbb{R}, \leq) \to \mathcal{S}_{n}(M)$ be a persistence $n$-configuration. By Proposition \ref{proposition:dyn_to_symtype}, the composition
\begin{equation*}
  \xymatrix{
  (\mathbb{R}, \leq) \ar@{->}[r]^-{\kappa(\mathcal{F})}&\mathbf{Dyn}^{0}_{n}(M) \ar@{->}[r]^{\widetilde{\Theta}}& \mathsf{SymType}
  }
\end{equation*}
is a persistent symmetry type. If $\widetilde{\Theta}\circ \kappa(\mathcal{F})$ admits only finitely many transition points $t_0 < t_1 < \cdots < t_k$, then it determines a sequence of symmetry types
\[
[\Sym(\mathcal{F}_{t_0})] \to [\Sym(\mathcal{F}_{t_1})] \to \cdots \to [\Sym(\mathcal{F}_{t_k})].
\]
The data
\[
\{ (t_i, [H_i]) \}_{i=0}^k, \quad \text{where } [H_i] = [\Sym(\mathcal{F}_{t_i})].
\]
encodes the symmetry type information at these transition points. We refer to this data as the \textit{symmetry type trajectory}.

\begin{proposition}
Let $M$ be a Euclidean space.
Let $\mathcal{F}: (\mathbb{R}, \leq) \to \mathcal{S}_n(M)$ be a persistence $n$-configuration such that the image $\mathcal{F}(\mathbb{R})$ is a compact subset of $\Conf_n(M)$. Then the associated symmetry type trajectory
\[
\{ (t, [\Sym(\mathcal{F}_{t})]) \}_{t\in \mathbb{R}}
\]
admits only finitely many distinct symmetry types.
\end{proposition}

\begin{proof}
Since symmetry types are invariant under translations, we may project the configurations to the subspace $H(k)$ consisting of configurations with barycenter at the origin. The projected image of $\mathcal{F}(\mathbb{R})$ remains compact in $H(k)$. The orthogonal group $O(k)$ is a compact Lie group acting locally smoothly on the manifold $H(k)$. By \cite[IV, Proposition 1.2]{bredon1972introduction}, any compact subset of $H(k)$ intersects only finitely many orbit types. Therefore, the trajectory admits only finitely many distinct symmetry types.
\end{proof}

\begin{definition}
Let $\mathcal{F}: (\mathbb{R}, \leq) \to \mathcal{S}_n(M)$ be a persistence $n$-configuration.
For any real numbers $a \leq b$, let $\phi_{a,b}: (\mathcal{F}_a, a) \to (\mathcal{F}_b, b)$ be the morphism in $\mathbf{Dyn}^{0}_{n}(M)$ induced by the persistence. The \textit{$(a,b)$-persistent symmetry type} is defined as the equivalence class of the span
\[
\Sym(\mathcal{F}_a) \xleftarrow{\phi_{a,b}^\flat} L_{a,b} \xrightarrow{\phi_{a,b}^\sharp} \Sym(\mathcal{F}_b),
\]
where $L_{a,b}$ is the group of symmetries that survive the evolution from $a$ to $b$. This span describes the relationship between the symmetry types $[\Sym(\mathcal{F}_a)]$ and $[\Sym(\mathcal{F}_b)]$.
\end{definition}

In particular, when $a=b$, the $(a,b)$-persistent symmetry type coincides with the symmetry type $[\Sym(\mathcal{F}_a)]$.

We now use the notion of a digraph of groups to describe the evolution of a persistent symmetry type. Recall that a \textit{digraph of groups} is a directed graph $ D = (V, E) $, equipped with a group $ G_v $ for each vertex $ v \in V $, and a group $ G_e $ for each edge $ e \in E $, together with a family of monomorphisms
\[
\phi_{e,v} : G_e \hookrightarrow G_v
\]
defined whenever the vertex $ v $ is either the source or the target of the edge $ e $.

Let $\mathcal{F}: (\mathbb{R}, \leq) \to \mathcal{S}_n(M)$ be a persistence $n$-configuration. Suppose that the associated persistent symmetry type
\[
\widetilde{\Theta} \circ \kappa(\mathcal{F}): (\mathbb{R}, \leq) \to \mathsf{SymType}
\]
has finitely many transition points on $\mathbb{R}$, say $ t_0 < t_1 < \cdots < t_k $. Then we obtain the data
\[
\left\{ (t_i,\; [\Sym(\mathcal{F}_{t_i})]) \right\}_{i=0}^k.
\]
We define a digraph of groups $ D(\mathcal{F}) $ whose vertex set is $ V = \{ v_i \mid 0 \leq i \leq k \} $ and whose edge set is $ E = \{ e_{ij} \mid 0 \leq i < j \leq k \} $. Each vertex $ v_i $ is assigned the group $ G_{v_i} = \Sym(\mathcal{F}_{t_i}) $, and each edge $ e_{ij} $ is assigned the group $ G_{e_{ij}} = L_{t_i, t_j} $, where $ \phi_{ij}: (\mathcal{F}_{t_i}, t_i) \to (\mathcal{F}_{t_j}, t_j) $ is the structure morphism in $\mathbf{Dyn}^{0}_{n}(M)$ induced by the persistence. These data naturally determine a digraph of groups via the canonical monomorphisms $\phi_{ij}^\flat: G_{e_{ij}} \hookrightarrow G_{v_i}$ and $\phi_{ij}^\sharp: G_{e_{ij}} \hookrightarrow G_{v_j}$.

\begin{example}\label{example:digraph_of_groups}
We consider a persistence $n$-configuration $\mathcal{F}$ evolving through three distinct stages at times $t = 0, 1, 2$.
At time $t_0 = 0$, the configuration $\mathcal{F}_{0} = \{A_0, B_0, C_0, D_0\}$ forms a rhombus, whose symmetry group is the dihedral group $G_{v_0} = D_2$ of order 4, generated by two reflections and a $180^\circ$ rotation.
At time $t_1 = 1$, the points $\mathcal{F}_{1} = \{A_1, B_1, C_1, D_1\}$ form a square, with full dihedral symmetry group $G_{v_1} = D_4$.
At time $t_2 = 2$, the configuration $\mathcal{F}_{2} = \{A_2, B_2, C_2, D_2\}$ consists of an isosceles right triangle $\{A_2, B_2, C_2\}$ together with a central point $D_2$. Its symmetry group is $G_{v_2} = \mathbb{Z}/2$, generated by reflection about the triangle's altitude.

\begin{figure}[h]
\centering
\definecolor{axisgray}{gray}{0.6}

\begin{tikzpicture}[scale=0.8, >=Stealth]

\begin{scope}
    \draw[->, color=axisgray, line width=0.4pt] (-1.8,0) -- (1.8,0) node[right, font=\small, text=black] {$x$};
    \draw[->, color=axisgray, line width=0.4pt] (0,-1.8) -- (0,1.8) node[above, font=\small, text=black] {$y$};

    \draw[line width=1.2pt, color=mainblue] (0,-1.2) -- (-1,0) -- (0,1.2) -- (1,0) -- cycle;

    \fill[mainblue] (0,-1.2) circle (2pt); \node[below right, font=\small] at (0,-1.2) {$A_0$};
    \fill[mainblue] (0, 1.2) circle (2pt); \node[above right, font=\small] at (0, 1.2) {$B_0$};
    \fill[mainblue] (-1,0) circle (2pt); \node[above left, font=\small] at (-1,0) {$C_0$};
    \fill[mainblue] ( 1,0) circle (2pt); \node[above right, font=\small] at (1,0) {$D_0$};
    
    \node[below] at (0, -2.2) {$t = 0$};
\end{scope}

\begin{scope}[xshift=6cm]
    \draw[->, color=axisgray, line width=0.4pt] (-1.8,0) -- (1.8,0) node[right, font=\small, text=black] {$x$};
    \draw[->, color=axisgray, line width=0.4pt] (0,-1.8) -- (0,1.8) node[above, font=\small, text=black] {$y$};

    \draw[line width=1.2pt, color=mainblue] (0,-1) -- (-1,0) -- (0,1) -- (1,0) -- cycle;
    
    \fill[mainblue] (0,-1) circle (2pt); \node[below right, font=\small] at (0,-1) {$A_1$};
    \fill[mainblue] (0, 1) circle (2pt); \node[above right, font=\small] at (0, 1) {$B_1$};
    \fill[mainblue] (-1,0) circle (2pt); \node[above left, font=\small] at (-1,0) {$C_1$};
    \fill[mainblue] ( 1,0) circle (2pt); \node[above right, font=\small] at (1,0) {$D_1$};
    
    \node[below] at (0, -2.2) {$t = 1$};
\end{scope}

\begin{scope}[xshift=12cm]
    \draw[->, color=axisgray, line width=0.4pt] (-1.8,0) -- (1.8,0) node[right, font=\small, text=black] {$x$};
    \draw[->, color=axisgray, line width=0.4pt] (0,-1.8) -- (0,1.8) node[above, font=\small, text=black] {$y$};

    \draw[line width=1.2pt, color=mainblue] (0,-1) -- (-1,0) -- (0,1) -- cycle;
    
    \fill[mainblue] (0,-1) circle (2pt); \node[below right, font=\small] at (0,-1) {$A_2$};
    \fill[mainblue] (0, 1) circle (2pt); \node[above right, font=\small] at (0, 1) {$B_2$};
    \fill[mainblue] (-1,0) circle (2pt); \node[above left, font=\small] at (-1,0) {$C_2$};
    \fill[mainblue] (0,0) circle (2pt); \node[above right, font=\small] at (0,0) {$D_2$};
    
    \node[below] at (0, -2.2) {$t = 2$};
\end{scope}

\end{tikzpicture}
\caption{Evolution of a 4-configuration from rhombus, to square, and then to isosceles right triangle.}
\end{figure}

The induced structure maps $\phi_{ij}: (\mathcal{F}_{i}, i) \to (\mathcal{F}_{j}, j)$ for $i < j$ send $A_i$ to $A_j$, $B_i$ to $B_j$, $C_i$ to $C_j$, and $D_i$ to $D_j$. According to the definition, we compute the edge groups as
\[
G_{e_{01}} = D_2, \quad G_{e_{12}} = \mathbb{Z}/2, \quad G_{e_{02}} = \mathbb{Z}/2.
\]
The above computations yield a directed graph of groups with three vertices and three edges.

\begin{figure}[h]
\centering

\begin{tikzpicture}[>=Stealth, font=\small, every node/.style={inner sep=2pt}]

\node (v0) at (0,0) [circle, draw=mainblue, line width=1.2pt, fill=mainblue!10, minimum size=1.4cm] {$G_{v_0} = D_2$};
\node (v1) at (4,0) [circle, draw=mainblue, line width=1.2pt, fill=mainblue!10, minimum size=1.4cm] {$G_{v_1} = D_4$};
\node (v2) at (8,0) [circle, draw=mainblue, line width=1.2pt, fill=mainblue!10, minimum size=1.4cm] {$G_{v_2} = \mathbb{Z}/2$};

\draw[->, line width=1.2pt] (v0) -- (v1) node[midway, above] {$G_{e_{01}} = D_2$};
\draw[->, line width=1.2pt] (v1) -- (v2) node[midway, above] {$G_{e_{12}} = \mathbb{Z}/2$};
\draw[->, line width=1.2pt, dashed] (v0) .. controls (4,-2) .. (v2) node[midway, below] {$G_{e_{02}} = \mathbb{Z}/2$};

\end{tikzpicture}

\caption{The digraph of groups associated with the persistent configuration $\mathcal{F}$ from Example~\ref{example:digraph_of_groups}.}
\end{figure}

\end{example} 

%% file: symmetry_degree_analysis.tex
\section{Symmetry degree analysis}\label{section:symmetry_degree}

Symmetry groups and their invariants offer a fundamental approach for quantifying and understanding the intrinsic structure of mathematical and physical systems.
A central theme in the study of symmetry is how to characterize the degree of symmetry of a given system, as well as the distribution and evolution of its symmetrical features. This motivates the search for suitable isomorphism invariants of symmetry groups, which can be used as effective descriptors of the overall symmetry, its localization, and dynamic behavior.
In this section, we introduce several quantitative and visual tools including the degree of symmetry and asymmetry, symmetry entropy, weighted paths, and Cayley graphs to analyze and visualize the symmetry properties of both finite configurations and persistence configurations. These tools help to formalize intuitive notions of symmetry and provide a structured way to compare symmetry across systems, over time, and parameter varying.

\subsection{Assignments and persistent invariants}

In many practical settings, although the group itself contains rich structural information, it is often too complex to be used directly as a feature. In contrast, features that are commutative, linear, and numerical are typically preferred, as they are easier to store, manipulate, and apply in data-driven pipelines.  As a result, it is desirable to transform the group structure into more tractable numerical or algebraic forms.

Let $\mathfrak{C}$ and $\mathfrak{D}$ be categories.
An \textit{object-wise assignment} $F: \mathfrak{C}\to \mathfrak{D}$ is a map that assigns to each object $X \in \mathfrak{C}$ an object $F(X) \in \mathfrak{D}$. Unlike a functor, an object-wise assignment does not act on morphisms and therefore does not preserve the categorical structure.

\begin{definition}
Let $\mathfrak{C}$ and $\mathfrak{D}$ be categories.
An object-wise assignment $F : \mathfrak{C} \to \mathfrak{D}$ satisfies \textit{isomorphism invariance} if, for any objects $X, Y \in \mathfrak{C}$ such that $X \cong Y$, it follows that $F(X) \cong F(Y)$ in $\mathfrak{D}$.
\end{definition}

It is worth noting that a functor is an object-wise assignment that satisfies isomorphism invariance, since a functor maps isomorphisms to isomorphisms.

\begin{definition}
Let $\mathcal{G} : (\mathbb{R}, \leq) \to \mathbf{Grp}$ be a persistence group, and let $F : \mathbf{Grp} \to \mathfrak{C}$ be an object-wise assignment into a category $\mathfrak{C}$ that satisfies isomorphism invariance.
The \textit{$F$-invariant} of the persistence group $\mathcal{G}$ is defined to be the object-wise composition
\[
F \circ \mathcal{G} : (\mathbb{R}, \leq) \to \mathfrak{C},
\]
which is naturally understood as a parameterized family of objects $F(\mathcal{G}_t)$ in $\mathfrak{C}$.
\end{definition}

\begin{example}
Consider the set $\mathbb{N}$ as a poset category. The group order function
\[
\mathrm{ord} : \mathbf{Grp}_{\mathrm{fin}} \to (\mathbb{N}, \leq), \quad G \mapsto |G|
\]
is an object-wise assignment that satisfies isomorphism invariance. Therefore, for a persistence group $\mathcal{G} : (\mathbb{R}, \leq) \to \mathbf{Grp}_{\mathrm{fin}}$, we obtain an $\mathrm{ord}$-invariant of $\mathcal{G}$ defined by the composition
\[
\mathrm{ord} \circ \mathcal{G} : (\mathbb{R}, \leq) \to (\mathbb{N}, \leq).
\]
\end{example}

\begin{example}
The center of a group, defined by
\[
Z(G) = \{ g \in G \mid gh = hg \text{ for all } h \in G \},
\]
assigns to each group an abelian subgroup. The center construction $Z : \mathbf{Grp} \to \mathbf{Ab}$ has the isomorphism invariance, though it is not a functor in general. For a persistence group $\mathcal{G} : (\mathbb{R}, \leq) \to \mathbf{Grp}_{\mathrm{fin}}$, the composition $Z \circ \mathcal{G}: (\mathbb{R}, \leq) \to \mathbf{Ab}$ is a $Z$-invariant.
\end{example}

\begin{definition}
Let $\mathcal{G} : (\mathbb{R}, \leq) \to \mathbf{Grp}$ be a persistence group, and let $F : \mathbf{Grp} \to \mathfrak{C}$ be a (covariant or contravariant) functor, where $\mathfrak{C}$ is a category. The strong \textit{$F$-invariant} of the persistence group $\mathcal{G}$ is defined to be the functor composition
\[
F \circ \mathcal{G} : (\mathbb{R}, \leq) \to \mathfrak{C},
\]
if $F$ is covariant, or $F \circ \mathcal{G}^{\mathrm{op}} : (\mathbb{R}, \geq) \to \mathfrak{C}$, if $F$ is contravariant.
\end{definition}

\begin{example}
Let $F_{\mathrm{ab}} : \mathbf{Grp} \to \mathbf{Ab}$ be the canonical abelianization functor $G \mapsto G / [G, G]$.
Then the composition
\[
F_{\mathrm{ab}} \circ \mathcal{G} \colon (\mathbb{R},\leq) \to \mathbf{Ab}
\]
is a strong $F_{\mathrm{ab}}$-invariant that captures the evolution of the abelian part of the persistence group $\mathcal{G}$.
\end{example}

\begin{example}
Let $A$ be a commutative ring with identity, and let $H^p(-; A) : \mathbf{Grp}^{\mathrm{op}} \to \mathbf{Mod}_A$ be the $p$-th group cohomology functor. For a group $G$, the group cohomology $H^p(G; A)$ can be computed as $\mathrm{Ext}^p_{A[G]}(A, A)$, where $A[G]$ denotes the group algebra of $G$ over $A$ and $A$ is equipped with a trivial $G$-module structure. Then the composition
\[
H^p(-; A) \circ \mathcal{G}^{\mathrm{op}} \colon (\mathbb{R}, \geq) \to \mathbf{Mod}_A
\]
defines a strong invariant of $\mathcal{G}$ taking values in the category of $A$-modules.

In particular, the zeroth cohomology group $H^0(G, A)$ is isomorphic to $A$, and the first cohomology group $H^1(G, A)$ is naturally isomorphic to $\mathrm{Hom}(G^{\mathrm{ab}}, A)$, where $G^{\mathrm{ab}} = G / [G, G]$ denotes the abelianization of $G$.
\end{example}

While a strong invariant of a persistence group is stable in a certain sense, a general invariant of a persistence group is not necessarily stable.

\begin{proposition}\label{proposition:invariant_stablity}
Let $\mathcal{G}, \mathcal{G}' : (\mathbb{R}, \leq) \to \mathbf{Grp}$ be persistence groups, and let $F : \mathbf{Grp} \to \mathfrak{C}$ be a functor into a category $\mathfrak{C}$ of persistence modules. Then we have
\[
  d_I(F \circ \mathcal{G}, F \circ \mathcal{G}') \leq d_I(\mathcal{G}, \mathcal{G}'),
\]
where $d_I$ denotes the interleaving distance.
\end{proposition}

\begin{proof}
It is directly obtained from \cite[Proposition 3.6]{chazal2009proximity}.
\end{proof}

In the subsequent sections, we introduce the concepts of degree of symmetries, symmetry entropy, and the symmetry degree polynomial, which are the object-wise assignments of primary interest in the characterization of symmetries.

\subsection{Degree of symmetries}

Recall that in Euclidean space, the isometries of a finite configuration $X$ consist of rotations and reflections. Intuitively, in two-dimensional space, a configuration with more axes of reflectional symmetry is considered more symmetric; similarly, a configuration admitting a higher order of rotational symmetry is viewed as more symmetric. As an extreme example, a circle possesses infinitely many reflection axes and is invariant under all rotations about its center, highlighting its maximal symmetry, in alignment with our geometric intuition. Motivated by this perspective, we propose the following general definition of the degree of symmetry in the context of arbitrary metric spaces.

From now on, unless otherwise specified, $M$ is the fixed metric space.
\begin{definition}
Let $X$ be an $n$-configuration in $\mathcal{S}_n(M)$. The \textit{degree of symmetry} of $X$ is defined as
\[
\mathrm{DegSym}(X) := \sum_{\pi \in \Sym(X)} \mathrm{ord}(\pi).
\]
Here, $\Sym(X)$ denotes the symmetry group of $X$, and $\mathrm{ord}(\pi)$ denotes the order of the group element $\pi$, that is, the smallest positive integer $k$ such that $\pi^k = e$.
\end{definition}

\begin{proposition}
Let $X$ be an $n$-configuration in $\mathcal{S}_n(M)$. Then $\mathrm{DegSym}(X) = 1$ if and only if $\Sym(X) = \{e\}$.
\end{proposition}

\begin{proof}
If $\Sym(X) = \{e\}$, then the only group element is the identity $e$, which has order $1$, so $\mathrm{DegSym}(X) = 1$. Conversely, if $\mathrm{DegSym}(X) = 1$, then only one group element of order $1$ exists, namely $e$, so the group must be trivial.
\end{proof}

\begin{proposition}
Let $X \in \mathcal{S}_n(M)$ be a configuration with symmetry group $G = \Sym(X)$. Then the degree of symmetry satisfies
\[
2|G| - 1 \leq \mathrm{DegSym}(X) \leq |G| \cdot \mathrm{ord}_{\max}(G),
\]
where $\mathrm{ord}_{\max}(G) = \max\{\mathrm{ord}(\pi) \mid \pi \in G\}$ denotes the maximal order of the elements in $G$.
\end{proposition}

\begin{proof}
By definition, the degree of symmetry of $X$ is
\[
\mathrm{DegSym}(X) = \sum_{\pi \in G} \mathrm{ord}(\pi).
\]
Observe that the identity element $e \in G$ has order $\mathrm{ord}(e) = 1$, while every non-identity element $\pi \neq e$ satisfies $\mathrm{ord}(\pi) \geq 2$. Thus, we obtain
\[
\mathrm{DegSym}(X) = \mathrm{ord}(e) + \sum_{\substack{\pi \in G \\ \pi \neq e}} \mathrm{ord}(\pi) \geq 1 + (|G| - 1) \cdot 2 = 2|G| - 1,
\]
establishing the lower bound.

For the upper bound, since by definition $\mathrm{ord}(\pi) \leq \mathrm{ord}_{\max}(G)$ for all $\pi \in G$, it follows that
\[
\mathrm{DegSym}(X) = \sum_{\pi \in G} \mathrm{ord}(\pi) \leq \sum_{\pi \in G} \mathrm{ord}_{\max}(G) = |G| \cdot \mathrm{ord}_{\max}(G).
\]
This completes the proof.
\end{proof}

\begin{remark}
For any configuration $X$ consisting of $n$ points, its symmetry group $\Sym(X)$ is a subgroup of the symmetric group $S_n$, and hence the order of any element in $\Sym(X)$ is bounded above by $\mathrm{ord}_{\max}(S_n)$, which is Landau's function $g(n)$.
\end{remark}

\begin{proposition}
Let $X, Y \in \mathcal{S}_n(M)$ be two configurations such that there exists an isometry $\varphi : M \to M$ satisfying $\varphi(X) = Y$. Then
\[
\mathrm{DegSym}(X) = \mathrm{DegSym}(Y).
\]
\end{proposition}

\begin{proof}
Any isometry $\varphi : M \to M$ induces an isomorphism between the symmetry groups $\Sym(X)$ and $\Sym(Y)$ via conjugation. Since group isomorphisms preserve the order of elements, the multiset of orders in $\Sym(X)$ and $\Sym(Y)$ coincide. Thus, their degree of symmetry is equal.
\end{proof}

\begin{remark}
The degree of symmetry is in general not stable under perturbations. That is, for arbitrarily small perturbations of a highly symmetric configuration $X$, the resulting configuration $X'$ may satisfy $\mathrm{DegSym}(X') = 1$.
\end{remark}

\begin{definition}
Let $X \in \mathcal{S}_n(M)$ be a configuration with symmetry group $G = \Sym(X)$. For each positive integer $k$, let
\[
p_k = \frac{|\{\sigma \in G \mid \mathrm{ord}(\sigma) = k\}|}{|G|}.
\]
The \textit{symmetry entropy} of $X$ is defined as
\[
\mathrm{SymEntropy}(X) := - \sum_{k \geq 1} p_k \log p_k.
\]
\end{definition}

Larger entropy indicates a richer variety of symmetry types, and can be viewed as a measure of structural complexity within the symmetry group.

\begin{definition}
Let $X$ be an $n$-configuration in $\mathcal{S}_n(M)$. The \textit{symmetry degree polynomial} of $X$ is defined to be the polynomial
\[
\mathrm{DegSym}_X(t) := \sum_{\pi \in \Sym(X)} t^{\mathrm{ord}(\pi)},
\]
where $\mathrm{ord}(\pi)$ denotes the order of the group element $\pi \in \Sym(X)$.
\end{definition}

Note that the symmetry degree polynomial recovers the order of the symmetry group when evaluated at $t = 1$, that is,
\[
\mathrm{DegSym}_X(1) = |\Sym(X)|.
\]
Moreover, its derivative evaluated at $t = 1$ recovers the degree of symmetry, that is,
\[
\mathrm{DegSym}_X'(1) = \sum_{\pi \in \Sym(X)} \mathrm{ord}(\pi) = \mathrm{DegSym}(X).
\]
Furthermore, the symmetry degree polynomial can be expressed in terms of the order distribution of group elements as
\[
\mathrm{DegSym}_X(t) = \sum_{k} a_k t^k,
\]
where $a_k = |G| \cdot p_k$, with $G = \Sym(X)$ and $p_k$ denoting the proportion of elements in $G$ of order $k$. Thus, the symmetry entropy is exactly the Shannon entropy associated with the normalized coefficient distribution of the symmetry degree polynomial.

\begin{example}
Let $X = \{x_1, x_2, x_3\} \subset \mathbb{R}^2$ be a configuration consisting of the vertices of an equilateral triangle centered at the origin. The coordinates of the three points are
\[
x_1 = (1, 0), \quad x_2 = \left(-\frac{1}{2}, \frac{\sqrt{3}}{2} \right), \quad x_3 = \left(-\frac{1}{2}, -\frac{\sqrt{3}}{2} \right).
\]
The symmetry group $\Sym(X)$ is the dihedral group $D_3$, which has six elements. These include the identity transformation, two nontrivial rotations by $120^\circ$ and $240^\circ$, and three reflections across the axes of symmetry of the triangle.

The orders of the elements in $\Sym(X)$ are as follows: the identity has order $1$, each of the three reflections has order $2$, and each of the two rotations has order $3$. Therefore, the degree of symmetry of the configuration is
\[
\mathrm{DegSym}(X) = \sum_{\pi \in \Sym(X)} \mathrm{ord}(\pi) = 1 + 3 \cdot 2 + 2 \cdot 3 = 13.
\]
The corresponding symmetry degree polynomial is given by
\[
\mathrm{DegSym}_X(t) = t + 3t^2 + 2t^3.
\]

Now, we perturb the point $x_1$ slightly to obtain $x_1' = (0.99, 0)$, and consider the new configuration $X' = \{x_1', x_2, x_3\}$. This represents a small geometric perturbation of the original configuration $X$. However, the symmetry group of $X'$ is dramatically reduced. We have $\Sym(X') \cong \mathbb{Z}/2\mathbb{Z}$, since the only nontrivial symmetry remaining is the reflection across the $x$-axis passing through the origin. As a result, the degree of symmetry drops significantly:
\[
\mathrm{DegSym}(X') = 1 + 2 = 3.
\]
This example demonstrates that although $X'$ is geometrically very close to $X$, their symmetry degrees differ substantially.
\end{example}

\subsection{Evolution of symmetry degree}

Let $\mathcal{F} : (\mathbb{R}, \leq) \to \mathcal{S}_{n}(M)$ be a persistence configuration. Recall that $\mathcal{F}$ induces a pseudofunctor $\kappa(\mathcal{F}) : (\mathbb{R}, \leq) \to \mathbf{Dyn}^{0}_{n}(M)$ mapping $t \in \mathbb{R}$ to the truncated dynamical system $(\mathcal{F}, t)$. By Proposition~\ref{proposition:pseudofunctor}, the assignment $\Sym : \mathbf{Dyn}^{0}_{n}(M) \to \mathrm{Span}(\mathbf{Grp})$ is a pseudofunctor.
Recall that the degree of symmetry construction
\[
\mathrm{DegSym} : \mathrm{Span}(\mathbf{Grp}) \to \mathbb{N}
\]
is an object-wise assignment that satisfies isomorphism invariance. By composing these maps, we obtain a function
\[
\xymatrix{
  (\mathbb{R}, \leq) \ar[r]^-{\kappa(\mathcal{F})} & \mathbf{Dyn}^{0}_{n}(M) \ar[r]^-{\Sym} & \mathrm{Span}(\mathbf{Grp}) \ar[rr]^-{\mathrm{DegSym}} && \mathbb{N},
}
\]
which we denote by $\delta : \mathbb{R} \to \mathbb{N}$, given by $\delta(t) = \mathrm{DegSym}(\mathcal{F}_t)$. This function can be interpreted as the $\mathrm{DegSym}$-invariant of the persistence configuration $\mathcal{F}$, and it captures the evolution of symmetry degree along the filtration.

Recall that for any real numbers $a \leq b$, we have a morphism $\phi_{a,b} : (\mathcal{F}_a, a) \to (\mathcal{F}_b, b)$ in $\mathbf{Dyn}^{0}_{n}(M)$ induced by the persistence. This morphism yields a span of symmetry groups
\[
\Sym(\mathcal{F}_a) \leftarrow \Sym(\mathcal{F})_{a,b} \rightarrow \Sym(\mathcal{F}_b),
\]
where $\Sym(\mathcal{F})_{a,b}$ is the persistence symmetry group introduced in Section~\ref{section:persistence_group}, consisting of symmetries that survive the evolution of $\mathcal{F}$ from $a$ to $b$. This leads to the following definition of a persistent degree of symmetry.

\begin{definition}
Let $\mathcal{F} : (\mathbb{R}, \leq) \to \mathcal{S}_n(M)$ be a persistence configuration. For any real numbers $a \leq b$, the \textit{$(a,b)$-persistent degree of symmetry} of $\mathcal{F}$ is defined by
\[
\mathrm{DegSym}_{a,b}(\mathcal{F}) := \sum_{\pi \in \Sym(\mathcal{F})_{a,b}} \mathrm{ord}(\pi),
\]
where $\Sym(\mathcal{F})_{a,b}$ is the group of symmetries that survive the evolution from $a$ to $b$.
\end{definition}

In practical computations, we often work with a discrete filtration parameter. Consider a persistence configuration $\mathcal{F} : (\mathbb{Z}, \leq) \to \mathcal{S}_n(M)$. This gives rise to a function
\[
\delta : \mathbb{Z} \to \mathbb{N}, \quad \delta(k) = \mathrm{DegSym}(\mathcal{F}_k).
\]
Furthermore, for each $k \in \mathbb{Z}$, we define a function
\[
\bar{\delta} : \mathbb{Z} \to \mathbb{N}, \quad \bar{\delta}(k) = \mathrm{DegSym}_{k,k+1}(\mathcal{F}),
\]
which measures the $(k,k+1)$-persistent degree of symmetry. Thus, we obtain a weighted path, where each vertex is assigned the symmetry degree $\delta(k)$, and each edge is weighted by the corresponding persistent symmetry degree $\bar{\delta}(k)$. This structure is represented schematically as follows.
\begin{figure}[H]
\centering
\begin{tikzpicture}[scale=1.1, >=Stealth, thick,
    vertex/.style={circle, fill=mainblue, inner sep=2.5pt},
    edge label/.style={midway, fill=white, font=\small, inner sep=1pt}
  ]

\node[vertex] (A) at (0,0) {};
\node[vertex] (B) at (3,0) {};
\node[vertex] (C) at (6,0) {};

\node at (0,0.5) {$\delta(k)$};
\node at (3,0.5) {$\delta(k+1)$};
\node at (6,0.5) {$\delta(k+2)$};

\draw[->, mainblue] (A) -- (B) node[edge label, above] {$\bar{\delta}(k)$};
\draw[->, mainblue] (B) -- (C) node[edge label, above] {$\bar{\delta}(k+1)$};
\draw[->, mainblue] (-1.5,0) -- (A); \node at (-2,0) {$\cdots$};
\draw[->, mainblue] (C) -- (7.5, 0); \node at (8,0) {$\cdots$};
\end{tikzpicture}
\end{figure}

\begin{example}\label{persistence_symmetry_degree}
Consider the three vertices $A, B, C$ of an equilateral triangle, given by the coordinates
\[
A = (0, 1), \quad B = \left(-\frac{\sqrt{3}}{2}, -\frac{1}{2}\right), \quad C = \left(\frac{\sqrt{3}}{2}, -\frac{1}{2}\right).
\]
In addition, we set the following points
\[
P_0 = \left(-\sqrt{3}, 1\right), \quad P_1 = (0,0), \quad P_2 = \left(\sqrt{3}, 1\right).
\]
We now consider a persistence configuration given by
\[
\mathcal{F}_0 = \{A, B, C, P_0\}, \quad \mathcal{F}_1 = \{A, B, C, P_1\}, \quad \mathcal{F}_2 = \{A, B, C, P_2\}.
\]
And the structure maps $f_{ij} : \mathcal{F}_i \to \mathcal{F}_j$ for $i < j$ are defined as follows: each map fixes the points $A, B, C$, and maps $P_i$ to $P_j$. That is, for all $0\leq i < j\leq 2$, we have
\[
f_{ij}(A) = A, \quad f_{ij}(B) = B, \quad f_{ij}(C) = C, \quad f_{ij}(P_i) = P_j.
\]

\begin{figure}[h]
\centering
\definecolor{axisgray}{gray}{0.6}
\begin{tikzpicture}[scale=1.2, >=Stealth]

\coordinate (A) at (0, 1);
\coordinate (B) at (-0.866, -0.5);  
\coordinate (C) at (0.866, -0.5);   
\coordinate (P0) at (-1.732, 1);
\coordinate (P1) at (0, 0);
\coordinate (P2) at (1.732, 1);

\draw[->, color=axisgray, line width=0.4pt] (-2.2, 0) -- (2.2, 0) node[right, font=\small, text=black] {$x$};
\draw[->, color=axisgray, line width=0.4pt] (0, -1.2) -- (0, 1.8) node[above, font=\small, text=black] {$y$};

\draw[line width=1.2pt, color=mainblue] (A) -- (B) -- (C) -- cycle;

\draw[->, dashed, color=symred, line width=0.8pt] (P0) to[bend right=20] (P1);
\draw[->, dashed, color=symred, line width=0.8pt] (P1) to[bend right=20] (P2);

\fill[mainblue] (A) circle (2pt); \node[above right, font=\small] at (A) {$A$};
\fill[mainblue] (B) circle (2pt); \node[below left, font=\small] at (B) {$B$};
\fill[mainblue] (C) circle (2pt); \node[below right, font=\small] at (C) {$C$};
\fill[mainblue] (P0) circle (2pt); \node[above left, font=\small] at (P0) {$P_0$};
\fill[mainblue] (P1) circle (2pt); \node[below right, font=\small] at (P1) {$P_1$};
\fill[mainblue] (P2) circle (2pt); \node[above right, font=\small] at (P2) {$P_2$};

\end{tikzpicture}
\caption{Geometric illustration of point positions in the persistent configuration in Example \ref{persistence_symmetry_degree}.}
\end{figure}

The configuration of $\mathcal{F}_{0}$ forms a rhombus, whose symmetry group is the dihedral group $D_2$, consisting of four isometries: the identity, one $180^\circ$ rotation about the center of the rhombus, and two reflections passing through its center. Hence, we have
\[
\Sym(\mathcal{F}_0) \cong D_2 \cong \mathbb{Z}/2 \times \mathbb{Z}/2.
\]
Moreover, the degree of symmetry of $\mathcal{F}_0$ is given by
\[
\mathrm{DegSym}(\mathcal{F}_0) = \sum_{\pi \in D_2} \mathrm{ord}(\pi) = 1 + 2 + 2 + 2 = 7.
\]
The configuration of $\mathcal{F}_{1}$ form an equilateral triangle centered at the origin. The symmetry group of an equilateral triangle is the dihedral group $D_3$, which consists of six isometries: three rotations (including the identity) and three reflections. By a straightforward calculation, we have the degree of symmetry
\[
\mathrm{DegSym}(\mathcal{F}_1)
= \sum_{\pi \in D_3} \mathrm{ord}(\pi)
= 1 + 3 + 3 + 2 + 2 + 2 = 13.
\]
On the other hand, by a step-by-step calculation, we obtain
\begin{equation*}
  \Sym(\mathcal{F})_{0,1} \cong \mathbb{Z}/2.
\end{equation*}
Therefore, the $(0,1)$-persistent degree of symmetry of the configuration $\mathcal{F}$ is given by
\begin{equation*}
  \mathrm{DegSym}_{0,1}(\mathcal{F}) = \sum_{\pi \in \Sym(\mathcal{F})_{0,1}} \mathrm{ord}(\pi) = 1 + 2 = 3.
\end{equation*}
Furthermore, the configuration $\mathcal{F}_2$ has the same geometric shape as $\mathcal{F}_0$, differing only in position. As a result, its symmetry group is also isomorphic to the dihedral group $D_2$, and hence its degree of symmetry is
\[
\mathrm{DegSym}(\mathcal{F}_2) = 7.
\]
A similar calculation shows that
\[
\Sym(\mathcal{F})_{1,2} \cong \mathbb{Z}/2,
\quad
\mathrm{DegSym}_{1,2}(\mathcal{F}) = 1 + 2 = 3.
\]
Hence, we have the following concept diagram.
\begin{figure}[H]
\centering
\begin{tikzpicture}[scale=1.1, >=Stealth, thick,
    vertex/.style={circle, fill=mainblue, inner sep=2.5pt},
    edge label/.style={midway, fill=white, font=\small, inner sep=1pt}
  ]

\node[vertex] (A) at (0,0) {};
\node[vertex] (B) at (3,0) {};
\node[vertex] (C) at (6,0) {};

\node at (0,0.5) {$\delta(0) = 7$};
\node at (3,0.5) {$\delta(1) = 13$};
\node at (6,0.5) {$\delta(2) = 7$};

\draw[->, mainblue] (A) -- (B) node[edge label, above] {$\bar{\delta}(0) = 3$};
\draw[->, mainblue] (B) -- (C) node[edge label, above] {$\bar{\delta}(1) = 3$};

\end{tikzpicture}
\caption{Illustration of the weighted path of the degree of symmetry in Example \ref{persistence_symmetry_degree}.}
\end{figure}
\end{example}

\subsection{Cayley graph of symmetries}

Cayley graphs are highly symmetric graphs that encode the abstract structure of a group \cite{cayley1878desiderata,magnus2004combinatorial}. We use Cayley graphs as combinatorial representations to explore the symmetries arising in our setting.

\begin{definition}
Let $G$ be a group and let $S \subseteq G$ be a subset such that $e \notin S$ and $S = S^{-1}$. The \textit{Cayley graph} $\mathrm{Cay}(G, S)$ is the graph with vertex set $G$, where an edge connects $g$ and $h$ if and only if $g^{-1}h \in S$.
\end{definition}
By definition, the edge set of the Cayley graph $\mathrm{Cay}(G, S)$ is given by
\[
  E = \big\{ \{g, gs\} \mid g \in G,\, s \in S \big\}.
\]

According to \cite[Lemma 3.7.3]{royle2001algebraic}, Cayley graphs are isomorphic under group isomorphisms: if $\phi: G \to G'$ is a group isomorphism, then the Cayley graphs $\mathrm{Cay}(G, S)$ and $\mathrm{Cay}(G', \phi(S))$ are isomorphic. Consequently, the construction of Cayley graphs exhibits isomorphism invariance.

\begin{example}
Let $G = \mathbb{Z}/n\mathbb{Z} = \{0,1,\dots,n-1\}$ be the cyclic group of order $n$ under addition modulo $n$, and consider the generating set
\[
  S = \{1, -1\} = \{1, n-1\}.
\]
The Cayley graph $\mathrm{Cay}(G, S)$ has vertex set $G$, where two vertices $g$ and $h$ are connected by an edge if and only if $g - h \equiv \pm 1 \pmod{n}$. This graph is isomorphic to the cycle graph $C_n$.

Let $S'= G \setminus \{e\}$. Then the Cayley graph $\mathrm{Cay}(G, S')$ is the complete graph $K_{|G|}$, in which every pair of distinct vertices is connected by an edge.
\end{example}

We now construct Cayley graphs from spans of groups. Specifically, for any real numbers $a \leq b$, we consider the span of symmetry groups
\[
\Sym(\mathcal{F}_a) \leftarrow \Sym(\mathcal{F})_{a,b} \rightarrow \Sym(\mathcal{F}_b).
\]
Let $S_a = \Sym(\mathcal{F})_{a,b} \setminus \{e\}$. Then we obtain the Cayley graph
\[
\mathrm{Cay}(\Sym(\mathcal{F}_a), S_a).
\]
In particular, if $S_a = \emptyset$, then the Cayley graph $\mathrm{Cay}(\Sym(\mathcal{F}_a), S_a)$ is a discrete graph, that is, a graph with no edges.

Note that the Cayley graph $\mathrm{Cay}(\Sym(\mathcal{F}_a), S_a)$ is the disjoint union of $[\Sym(\mathcal{F}_a):\Sym(\mathcal{F})_{a,b}]$ copies of the complete graph $K_{m}$. Here, $m$ is the order of the group $\Sym(\mathcal{F})_{a,b}$. Hence, the structure of the Cayley graph $\mathrm{Cay}(\Sym(\mathcal{F}_a), S_a)$ is largely determined by the subgroup $\Sym(\mathcal{F})_{a,b}$.

The order function $\mathrm{ord}: \Sym(\mathcal{F})_{a,b} \to \mathbb{N}$ maps each group element to its order. Define the subset
\[
S_a^k = \{ \pi \in \Sym(\mathcal{F})_{a,b} \mid \pi \neq e, \, \mathrm{ord}(\pi) \leq k \}.
\]
The set $S_a^k$ is a subset of $\Sym(\mathcal{F})_{a,b}$, and for any $\pi \in S_a^k$, its inverse $\pi^{-1} \in S_a^k$. Consequently, we can construct the Cayley graph
\[
\mathrm{Cay}(\Sym(\mathcal{F})_{a,b}, S_a^k).
\]
Thus, for distinct values of the order $k$, we obtain different Cayley graphs, enabling us to analyze their properties to characterize the symmetries of $\Sym(\mathcal{F})_{a,b}$. For instance, we may consider the spectral properties of these Cayley graphs, which are frequently employed as a key feature in data analysis.

\begin{example}\label{example:cayley_symmetry}

Consider the configuration $\mathcal{F}_0 = \{A_0, B_0, C_0, D_0\}$, which forms a square, and the configuration $\mathcal{F}_1 = \{A_1, B_1, C_1, D_1\}$, which forms a rhombus. The symmetry groups of the configurations $\mathcal{F}_0$ and $\mathcal{F}_1$ are $D_4$ and $D_2$, respectively. Let
\[
\phi_{01}: \mathcal{F}_0 \to \mathcal{F}_1
\]
be a mapping such that $A_0 \mapsto A_1$, $B_0 \mapsto B_1$, $C_0 \mapsto C_1$, and $D_0 \mapsto D_1$.
\begin{figure}[h]
\centering
\definecolor{axisgray}{gray}{0.6}
\definecolor{mainblue}{RGB}{26, 89, 142}
\begin{tikzpicture}[scale=0.9, >=Stealth]

\begin{scope}
\draw[->, color=axisgray, line width=0.4pt] (-1.5, 0) -- (1.5, 0) node[right, font=\small, text=black] {$x$};
\draw[->, color=axisgray, line width=0.4pt] (0, -1.5) -- (0, 1.5) node[above, font=\small, text=black] {$y$};

\draw[line width=1.2pt, color=mainblue] (0,-1) -- (-1,0) -- (0,1) -- (1,0) -- cycle;

\fill[mainblue] (0,-1) circle (2pt); \node[below right, font=\small] at (0,-1) {$A_0$};
\fill[mainblue] (0,1) circle (2pt); \node[above right, font=\small] at (0,1) {$B_0$};
\fill[mainblue] (-1,0) circle (2pt); \node[above left, font=\small] at (-1,0) {$C_0$};
\fill[mainblue] (1,0) circle (2pt); \node[above right, font=\small] at (1,0) {$D_0$};

\node at (0,-2) {$t = 0$};
\end{scope}

\begin{scope}[xshift=5cm]
\draw[->, color=axisgray, line width=0.4pt] (-1.7, 0) -- (1.7, 0) node[right, font=\small, text=black] {$x$};
\draw[->, color=axisgray, line width=0.4pt] (0, -1.5) -- (0, 1.5) node[above, font=\small, text=black] {$y$};

\draw[line width=1.2pt, color=mainblue] (0,-1) -- (-1.2,0) -- (0,1) -- (1.2,0) -- cycle;

\fill[mainblue] (0,-1) circle (2pt); \node[below right, font=\small] at (0,-1) {$A_1$};
\fill[mainblue] (0,1) circle (2pt); \node[above right, font=\small] at (0,1) {$B_1$};
\fill[mainblue] (-1.2,0) circle (2pt); \node[above left, font=\small] at (-1.2,0) {$C_1$};
\fill[mainblue] (1.2,0) circle (2pt); \node[above right, font=\small] at (1.2,0) {$D_1$};

\node at (0,-2) {$t = 1$};
\end{scope}
\end{tikzpicture}
\caption{Geometric illustration of point positions in the configurations in Example \ref{example:cayley_symmetry}.}
\end{figure}
It follows that $\Sym(\mathcal{F})_{0,1} = D_2$. Moreover, we have the span of groups
\[
D_4 \leftarrow D_2 \rightarrow D_2.
\]
Consequently, we obtain the Cayley graph $\mathrm{Cay}(D_4, S)$, where $S = D_2 \setminus \{e\}$. Let
\begin{equation*}
  D_4 = \{e, r, r^2, r^3, s, sr, sr^2, sr^3\}
\end{equation*}
and let $S = \{r^2, s, sr^2\}$. The adjacency matrix of the Cayley graph $\mathrm{Cay}(D_4, S)$ represented as follows.
\[
\begin{array}{c|cccccccc}
       & e & r & r^2 & r^3 & s & sr & sr^2 & sr^3 \\
\hline
e      & 0 & 0 & 1 & 0 & 1 & 0 & 1 & 0 \\
r      & 0 & 0 & 0 & 1 & 0 & 1 & 0 & 1 \\
r^2    & 1 & 0 & 0 & 0 & 1 & 0 & 1 & 0 \\
r^3    & 0 & 1 & 0 & 0 & 0 & 1 & 0 & 1 \\
s      & 1 & 0 & 1 & 0 & 0 & 0 & 1 & 0 \\
sr     & 0 & 1 & 0 & 1 & 0 & 0 & 0 & 1 \\
sr^2   & 1 & 0 & 1 & 0 & 1 & 0 & 0 & 0 \\
sr^3   & 0 & 1 & 0 & 1 & 0 & 1 & 0 & 0
\end{array}
\]
Thus, the Cayley graph $\mathrm{Cay}(D_4, S)$ consists of two connected components, as shown in Figure~\ref{figure:cayley_graph}.

\begin{figure}[h]
\centering
\begin{tikzpicture}[>=Stealth, font=\small, every node/.style={circle, draw=mainblue, line width=1.2pt, fill=mainblue!10, minimum size=7mm, inner sep=0pt}]
\node (e) at (0, 1.5) {$e$};
\node (s) at (-1.5, 1.5) {$s$};
\node (r2) at (-1.5, 0) {$r^2$};
\node (sr2) at (0, 0) {$sr^2$};

\node (r) at (2.5, 1.5) {$r$};
\node (sr) at (4, 1.5) {$sr$};
\node (r3) at (4, 0) {$r^3$};
\node (sr3) at (2.5, 0) {$sr^3$};

\draw[line width=1.2pt, mainblue] (e) -- (s);
\draw[line width=1.2pt, mainblue] (s) -- (r2);
\draw[line width=1.2pt, mainblue] (r2) -- (sr2);
\draw[line width=1.2pt, mainblue] (sr2) -- (e);
\draw[line width=1.2pt, mainblue] (e) -- (r2);
\draw[line width=1.2pt, mainblue] (s) -- (sr2);

\draw[line width=1.2pt, mainblue] (r) -- (sr);
\draw[line width=1.2pt, mainblue] (sr) -- (r3);
\draw[line width=1.2pt, mainblue] (r3) -- (sr3);
\draw[line width=1.2pt, mainblue] (sr3) -- (r);
\draw[line width=1.2pt, mainblue] (r) -- (r3);
\draw[line width=1.2pt, mainblue] (sr) -- (sr3);
\end{tikzpicture}
\caption{Illustration of Cayley graph $\mathrm{Cay}(D_4, S)$ with $S = \{r^2, s, sr^2\}$.}\label{figure:cayley_graph}
\end{figure}

The Laplacian spectrum of the Cayley graph $\mathrm{Cay}(D_4, S)$ is $\{0,0,4,4,4,4,4,4\}$.
\end{example} 

%% file: symmetry_defect_analysis.tex
\section{Symmetry defect analysis}\label{section:symmetry_defect}

Symmetry appears everywhere, from patterns in buildings and artwork to the laws of nature and the structure of molecules. But in the real world, perfect symmetry is rare. Objects and systems often show only approximate symmetry because of imperfections, noise, or external disturbances. To better understand and describe these situations, we need ways to measure how far something is from being perfectly symmetric. Concepts like symmetry defect, degree of asymmetry (asymmetricity), and symmetry/asymmetry measure help us do this by giving a precise way to evaluate and compare the degree of symmetry (asymmetry) in different settings. These tools let us study symmetry even when it is only partial or approximate, which is often the case in real-world systems.

\subsection{Symmetry defect}

We now define a symmetry defect on the objects of $\mathcal{S}_n(M)$, which measures the degree of asymmetry of a finite configuration.

Let $(M,d)$ be a metric space. Recall that $\mathcal{S}_n(M)$ is the category whose objects are $n$-point subsets of $M$. The set of all objects in $\mathcal{S}_n(M)$ constitutes the (unordered) configuration space of $M$, denoted by $\Conf_n(M)$. Let $D : \Conf_n(M) \times \Conf_n(M) \to \mathbb{R}$ be a suitable metric defined on this space.

To more accurately capture the discrepancy between $n$-configurations $X$ and $Y$, the metric $D$ is typically chosen to reflect the geometric dissimilarity or matching cost between two point sets. The $p$-Wasserstein distance is
\[
W_p(X,Y) = \left(\inf_{\gamma} \sum_{x \in X} d(x, \gamma(x))^p \right)^{1/p},
\]
where the infimum is taken over all bijections $\gamma : X \to Y$.

\begin{definition}
Let $X$ be an $n$-configuration in $\mathcal{S}_n(M)$. For a given isometry $\pi$ of $M$, the \textit{standard symmetry defect of $X$ with respect to $\pi$} is defined by
\[
  \widetilde{\mu}(X,\pi) := W_p(X,\pi(X)).
\]
\end{definition}

According to the above definition, we have $\widetilde{\mu}(X,e) = 0$, where $e$ is the identity isometry. Although the identity isometry theoretically represents a trivial symmetry, it provides no useful information in computations. This motivates the introduction of a refined definition of symmetry defect later, which better captures meaningful symmetries in practical applications.

For any isometry $\pi \in \mathrm{Iso}(M)$ and any $n$-configuration $X \in \mathcal{S}_n(M)$, the restriction $\pi|_X : X \to \pi(X)$ is a bijection between finite subsets of $M$.

\begin{definition}
Let $X$ be an $n$-configuration in $\mathcal{S}_n(M)$. For a given isometry $\pi$ of $M$, the \textit{symmetry defect of $X$ with respect to $\pi$} is defined by
\[
  \mu(X,\pi) := \left( \inf_{\gamma \neq \pi|_X}   \sum_{x \in X} d(x, \gamma(x))^p  \right)^{1/p},
\]
where the infimum is taken over all bijections $\gamma : X \to \pi(X)$ other than $\pi|_{X}$.
\end{definition}

\begin{remark}
In the above definition of symmetry defect, the infimum is taken over all bijections $\gamma : X \to \pi(X)$ such that $\gamma \neq \pi|_X$ in order to exclude contributions from infinitesimal perturbations. For instance, in Euclidean space, a small rotation of a configuration $X$ may yield a set $\pi(X)$ that is arbitrarily close to $X$ in the Wasserstein distance. To avoid such negligible perturbations from being counted as genuine symmetries, we exclude the canonical bijection $\pi|_X : X \to \pi(X)$ when computing the cost of matching points in $X$ with those in $\pi(X)$.
\end{remark}

In particular, when $M$ is a Euclidean space, let $X = \{x_1, \dots, x_n\}$ and $Y = \{y_1, \dots, y_n\}$, and let $S_n$ denote the symmetric group on $n$ elements. Then the bijection $\pi|_X : X \to \pi(X)$ corresponds to a unique permutation $\sigma_\pi \in S_n$. In this case, the symmetry defect reduces to
\[
  \mu(X, \pi) = \left(\inf_{\substack{\sigma \in S_n \\ \sigma \neq \sigma_\pi}}  \sum_{i=1}^n \| x_i - y_{\sigma(i)} \|^p \right)^{1/p},
\]
where $\|\cdot\|$ denotes the Euclidean norm. From now on, whenever $M$ is a Euclidean space, we equip it with the standard Euclidean metric.

\begin{lemma}\label{lemme:continuous}
Let $M = \mathbb{R}^k$ be a $k$-dimensional Euclidean space, and let $X \in \mathcal{S}_n(M)$ be an $n$-configuration. Fix $p \geq 1$. Then the map
\[
\mu(X, -) : \Iso(\mathbb{R}^k) \to \mathbb{R}, \quad \pi \mapsto \mu(X,\pi)
\]
is continuous.
\end{lemma}

\begin{proof}
Let $\pi \in \Iso(M)$ be any isometry, and let $Y = \pi(X) = \{y_1, \dots, y_n\}$. Then there exists a unique permutation $\sigma_\pi \in S_n$ such that $\pi(x_i) = y_{\sigma_\pi(i)}$ for all $i$. The symmetry defect is given by
\[
f(\pi) = \mu(X, \pi)^p = \inf_{\substack{\sigma \in S_n \\ \sigma \neq \sigma_\pi}} \sum_{i=1}^n \| x_i - y_{\sigma(i)} \|^p.
\]
For each $\sigma \in S_n$, define
\[
f_\sigma(\pi) = \sum_{i=1}^n \| x_i - \pi(x_{\sigma(i)}) \|^p.
\]
Then $f(\pi) = \inf_{\sigma \neq \text{id}} f_{\sigma \circ \sigma_\pi^{-1}}(\pi)$, so $f$ is the pointwise minimum of finitely many continuous functions (since each $f_\sigma$ is continuous in $\pi$). Therefore, $f$ is upper semi-continuous.

To prove lower semi-continuity, fix $\pi_0 \in \Iso(M)$ and $\varepsilon > 0$. Choose $\sigma_0 \in S_n \setminus \{\sigma_{\pi_0}\}$ such that
\[
f(\pi_0) \geq \sum_{i=1}^n \| x_i - \pi_0(x_{\sigma_0(i)}) \|^p - \frac{\varepsilon}{2}.
\]
Since $\pi \mapsto \sum_{i=1}^n \| x_i - \pi(x_{\sigma_0(i)}) \|^p$ is continuous, there exists a neighborhood $U$ of $\pi_0$ such that for all $\pi \in U$,
\[
\sum_{i=1}^n \| x_i - \pi(x_{\sigma_0(i)}) \|^p < \sum_{i=1}^n \| x_i - \pi_0(x_{\sigma_0(i)}) \|^p + \frac{\varepsilon}{2}.
\]
Hence, we have
\[
f(\pi) \leq \sum_{i=1}^n \| x_i - \pi(x_{\sigma_0(i)}) \|^p < f(\pi_0) + \varepsilon,
\]
which shows that $f$ is also lower semi-continuous at $\pi_0$.

Therefore, $f$ is continuous on $\Iso(M)$. Finally, since $\mu(X, \pi) = f(\pi)^{1/p}$ and the map $x \mapsto x^{1/p}$ is continuous on $[0, \infty)$, we conclude that $\mu(X, -)$ is continuous.
\end{proof}

\begin{definition}
Let $X$ be an $n$-configuration in $\mathcal{S}_n(M)$. The \textit{symmetry defect} of $X$ is defined by
\[
  \mu(X) := \inf_{\pi \in \mathrm{Iso}(M) } \mu(X,\pi),
\]
where the infimum is taken over all isometries $\pi$ of $M$.
\end{definition}

Clearly, if there exists a isometry $\pi$ of $M$ such that $\pi(X) = X$, then $\mu(X) = 0$, indicating that the configuration $X$ is symmetric. Conversely, if $X$ is asymmetric, then $\mu(X) > 0$. The larger the value of $\mu(X)$, the more asymmetric the configuration is.

\begin{proposition}
The symmetry defect is invariant under isometries of the metric space $M$. Specifically, for any isometry $\sigma \in \mathrm{Iso}(M)$ and any configuration $X \in \mathcal{S}_n(M)$, we have
\[
\mu(\sigma(X)) = \mu(X).
\]
\end{proposition}

\begin{proof}
By definition, we have
\[
\mu(\sigma(X)) = \inf_{\pi \in \mathrm{Iso}(M)} \mu(\sigma(X), \pi).
\]
For any isometry $\pi \in \mathrm{Iso}(M)$, we can write
\[
\mu(\sigma(X), \pi) = \left( \inf_{\gamma' \neq \pi|_{\sigma(X)}} \sum_{x' \in \sigma(X)} d(x', \gamma'(x'))^p \right)^{1/p}.
\]
Now consider the change of variable $\tau = \sigma^{-1} \circ \pi \circ \sigma \in \mathrm{Iso}(M)$, so that $\pi = \sigma \circ \tau \circ \sigma^{-1}$. For each bijection $\gamma' : \sigma(X) \to \pi(\sigma(X))$, let
\[
\gamma = \sigma^{-1} \circ \gamma' \circ \sigma : X \to \tau(X).
\]
Then for any $x \in X$, setting $x' = \sigma(x)$, we have
\[
d(x, \gamma(x)) = d(x, \sigma^{-1}(\gamma'(\sigma(x)))) = d(\sigma(x), \gamma'(\sigma(x))) = d(x', \gamma'(x')).
\]
It follows that
\begin{align*}
  \mu(\sigma(X), \pi)^p
  &= \inf_{\gamma' \neq \pi|_{\sigma(X)}} \sum_{x' \in \sigma(X)} d(x', \gamma'(x'))^p \\
  &= \inf_{\gamma \neq \tau|_X} \sum_{x \in X} d(x, \gamma(x))^p \\
  &= \mu(X, \tau)^p.
\end{align*}
Hence, we obtain
\[
\mu(\sigma(X), \pi) = \mu(X, \tau).
\]
Taking the infimum over all $\pi \in \mathrm{Iso}(M)$ is equivalent to taking the infimum over all $\tau \in \mathrm{Iso}(M)$ via conjugation by $\sigma$. Therefore,
\[
\mu(\sigma(X)) = \inf_{\pi} \mu(\sigma(X), \pi) = \inf_{\tau} \mu(X, \tau) = \mu(X),
\]
which completes the proof.
\end{proof}

\begin{proposition}\label{proposition:continuous_defect}
Let $M = \mathbb{R}^k$ be a $k$-dimensional Euclidean space.
The symmetry defect function
\[
\mu : \Conf_n(\mathbb{R}^k) \to \mathbb{R}
\]
is continuous with respect to the configuration $X$.
\end{proposition}

\begin{proof}
The configuration space $\Conf_n(\mathbb{R}^k)$ is naturally embedded in the product space $(\mathbb{R}^k)^n$, equipped with the product topology. For any fixed isometry $\pi \in \mathrm{Iso}(\mathbb{R}^k)$, consider the function
\[
\mu(X, \pi) = \left(\inf_{\substack{\gamma : X \to \pi(X) \\ \gamma \neq \pi|_X}}  \sum_{x \in X} \| x - \gamma(x) \|^p \right)^{1/p}.
\]
Since both $X$ and $\pi(X)$ are finite sets of size $n$, there are only finitely many bijections $\gamma: X \to \pi(X)$. For each such $\gamma$, the function
\[
X \mapsto \left( \sum_{x \in X} \| x - \gamma(x) \|^p \right)^{1/p}
\]
is continuous with respect to the product topology, because the norm $\|\cdot\|$ is continuous on $\mathbb{R}^k$ and the sum is finite. Taking the minimum over finitely many continuous functions, we conclude that $X \mapsto \mu(X, \pi)$ is continuous for each fixed $\pi \in \mathrm{Iso}(\mathbb{R}^k)$.

\vspace{1ex}
Now consider the symmetry defect function
\[
\mu(X) = \inf_{ \pi \in \mathrm{Iso}(\mathbb{R}^k) } \mu(X, \pi).
\]
This is the infimum over a (noncompact) family of continuous functions indexed by the isometry group $\mathrm{Iso}(\mathbb{R}^k) \cong O(k) \ltimes \mathbb{R}^k$. Since this group is not compact (due to unbounded translations), the pointwise infimum may not a priori be continuous. However, we can restrict the domain of minimization to a compact subset as follows.

Fix a configuration $X \in \Conf_n(\mathbb{R}^k)$. By Lemma~\ref{lemme:continuous}, the map $\pi \mapsto \mu(X, \pi)$ is continuous on $\mathrm{Iso}(\mathbb{R}^k)$. For any sufficiently large constant $K > 0$, the sublevel set
\[
S_K = \left\{ \pi \in \mathrm{Iso}(\mathbb{R}^k) \mid \mu(X, \pi) \leq K \right\}
\]
is closed by continuity and bounded due to the finiteness of $X$ and the growth of $\mu(X, \pi)$ under large translations. Since $\mathrm{Iso}(\mathbb{R}^k)$ is a locally compact and complete metric space, it follows that $S_K$ is compact.

Since $\mu(X,\pi)$ is continuous in both $X$ and $\pi$, and the group action is continuous, the map
\[
(X, \pi) \mapsto \mu(X, \pi)
\]
is jointly continuous on $\Conf_n(\mathbb{R}^k) \times \mathrm{Iso}(\mathbb{R}^k)$. Thus, for fixed compact $S_K$, the infimum
\[
\mu(X) = \inf_{\pi \in S_K} \mu(X,\pi)
\]
is continuous in $X$ as the infimum of a compact family of continuous functions.

\smallskip
Therefore, the symmetry defect function $\mu$ is continuous on $\Conf_n(\mathbb{R}^k)$.
\end{proof}

In practical applications, the ambient metric space is often taken to be a Euclidean space. Given that it is sometimes desirable to restrict the class of isometries under consideration, we may focus on a subset $\Gamma \subseteq \Iso(M)$ of isometries. In this case, we define the \textit{symmetry defect} of a configuration $X$ relative to $\Gamma$ by
\[
  \mu_\Gamma(X) = \inf_{\pi \in \Gamma} \mu (X,\pi).
\]
In Euclidean space, the only isometries that preserve a finite point set are reflections and rotations. Thus, one often considers $\Gamma$ to be the set of all reflections or all rotations in $\mathbb{R}^k$. Even more restricted choices of $\Gamma$ may be useful in practice. For example, one may consider only those rotations centered at the centroid of $X$, or reflections whose axes pass through the centroid. Such restrictions can significantly reduce computational complexity while still capturing meaningful symmetry behavior, with little compromise in theoretical accuracy.

\subsection{Symmetries and approximate groups}

In this section, we introduce the notion of approximate symmetry and study the properties of the associated set of approximate symmetries for a finite configuration. Notably, this set exhibits structural similarities with the concept of approximate groups, as formulated by T.~Tao et al \cite{breuillard2012structure}.

Let $M$ be a $k$-dimensional Euclidean space $\mathbb{R}^k$, and let $X$ be an $n$-configuration in $\mathcal{S}_n(M)$. For convenience, we assume that the finite configurations under consideration are centered at the origin. Under this assumption, all candidate symmetry transformations preserve the origin, and the corresponding symmetry groups are subgroups of the orthogonal group $O(k)$.
For any $\varepsilon \geq 0$, we define the set of $\varepsilon$-approximate symmetries of $X$ by
\[
  \Sym_\varepsilon(X) = \left\{ \pi \in O(k) \mid \widetilde{\mu}(X, \pi) \leq \varepsilon \right\}.
\]
In particular, when $\varepsilon = 0$, the set $\Sym_0(X)$ coincides with the exact symmetry group $\Sym(X)$. The motivation for introducing the set $\Sym_\varepsilon(X)$ is to capture approximate symmetries of a configuration, recognizing that multiple transformations may contribute significantly to the near-symmetry structure of $X$.

Note that $\Sym_\varepsilon(X)$ is generally not a subgroup of the Lie group $O(k)$, as the composition of two $\varepsilon$-approximate symmetries may not remain within the same error bound $\varepsilon$.

\begin{proposition}\label{proposition:compact}
The set $\Sym_\varepsilon(X)$ is a compact subset of the orthogonal group $O(k)$.
\end{proposition}

\begin{proof}
By a similar proof of Lemma~\ref{lemme:continuous}, the map $\pi \mapsto \widetilde{\mu}(X, \pi)$ is continuous on $O(k)$. Therefore, the set of sublevels $\Sym_\varepsilon(X) = \{\pi \in O(k) \mid \widetilde{\mu}(X, \pi) \leq \varepsilon\}$ is closed in $O(k)$. Since $O(k)$ is a compact Lie group, any closed subset of $O(k)$ is also compact. Hence, $\Sym_\varepsilon(X)$ is compact.
\end{proof}

\begin{example}
Let $X = \{(1, 0), (-1, 0)\} \subset \mathbb{R}^2$, and for each $\theta \in [0, 2\pi)$, let $R_\theta \in O(2)$ denote the counterclockwise rotation by angle $\theta$ about the origin. Then we have
\[
R_\theta(X) = \{(\cos \theta, \sin \theta), (-\cos \theta, -\sin \theta)\}.
\]
Note that $\widetilde{\mu}(X, R_0) = 0$, since $R_0$ is the identity transformation. Since the map $\theta \mapsto \widetilde{\mu}(X, R_\theta)$ is continuous and satisfies $\widetilde{\mu}(X, R_0) = 0$, it follows that for any fixed $\varepsilon > 0$, there exists an open interval $I \subset \mathbb{R}$ containing $0$ such that
\[
R_\theta \in \Sym_\varepsilon(X) \quad \text{for all } \theta \in I.
\]

This example shows that even for a simple configuration, the set of $\varepsilon$-approximate symmetries may be uncountably infinite.
\end{example}

\begin{lemma}\label{lemma:composition_approximate}
For any $\varepsilon \geq 0$, the composition of two $\varepsilon$-approximate symmetries yields a $2\varepsilon$-approximate symmetry. That is,
\[
\Sym_\varepsilon(X) \cdot \Sym_\varepsilon(X) \subseteq \Sym_{2\varepsilon}(X).
\]
\end{lemma}

\begin{proof}
For $\pi_1, \pi_2 \in \Sym_\varepsilon(X)$, we have
\[
\widetilde{\mu}(X, \pi_1) \leq \varepsilon \quad \text{and} \quad \widetilde{\mu}(X, \pi_2) \leq \varepsilon.
\]
We aim to estimate the defect of the composition $\pi_1 \circ \pi_2$.

Using the triangle inequality for the distance function $\widetilde{\mu}$, we obtain
\begin{align*}
\widetilde{\mu}(X, \pi_1 \circ \pi_2)
&= W_p\big((\pi_1 \circ \pi_2)(X), X\big) \\
&\leq W_p\big((\pi_1 \circ \pi_2)(X), \pi_1(X)\big) + W_p\big(\pi_1(X), X\big) \\
&= W_p\big(\pi_1(\pi_2(X)), \pi_1(X)\big) + \widetilde{\mu}(X, \pi_1).
\end{align*}
Since $\pi_1$ is an isometry, we have
\[
W_p\big(\pi_1(\pi_2(X)), \pi_1(X)\big) = W_p\big(\pi_2(X), X\big) = \widetilde{\mu}(X, \pi_2).
\]
It follows that
\[
\widetilde{\mu}(X, \pi_1 \circ \pi_2) \leq \widetilde{\mu}(X, \pi_2) + \widetilde{\mu}(X, \pi_1) \leq \varepsilon + \varepsilon = 2\varepsilon.
\]
This shows that $\pi_1 \circ \pi_2 \in \Sym_{2\varepsilon}(X)$, completing the proof.
\end{proof}

Although the set $\Sym_\varepsilon(X)$ is not a subgroup of the orthogonal group, its compactness and partial closure under composition suggest that it possesses certain group-like properties. This behavior can be partially understood within the notion of \textit{approximate groups}, as developed by Tao and his coauthors~\cite{breuillard2012structure,tao2008product,tao2006additive}.

\begin{definition}[Approximate subgroup in a topological group]
Let $G$ be a locally compact topological group. A non-empty subset $A \subseteq G$ is said to be a $K$-approximate subgroup if it satisfies the following properties:
\begin{itemize}
  \item[(i)] $A$ is symmetric, that is, $A = A^{-1}$;
  \item[(ii)] $A$ contains the identity element $e \in G$;
  \item[(iii)] There exists a finite subset $X \subseteq G$ with cardinality at most $K$ such that
  \[
    A \cdot A \subseteq X \cdot A.
  \]
\end{itemize}
If $A$ is also compact, we refer to $A$ as a compact $K$-approximate subgroup.
\end{definition}

\begin{remark}
The classical definition of an approximate subgroup typically assumes that the subset $A$ is finite, focusing on combinatorial and algebraic properties within discrete groups. However, as noted by T.~Tao, the concept of approximate groups can be extended to the setting of topological groups. The definition given above is adapted here primarily to serve the setting of this work.
\end{remark}

\begin{definition}[$K$-approximate group]
Let $G$ be a group. A finite subset $A \subseteq G$ is called a \textit{$K$-approximate group} if it satisfies the following conditions:
\begin{itemize}
  \item[(i)] $A$ is symmetric: $A = A^{-1}$, i.e., for all $a \in A$, we have $a^{-1} \in A$;
  \item[(ii)] $A$ contains the identity element: $e \in A$;
  \item[(iii)] There exists a finite subset $X \subseteq G$ with $|X| \leq K$ such that
  \[
    A \cdot A \subseteq X \cdot A.
  \]
\end{itemize}
\end{definition}

\begin{theorem}
Let $X$ be an $n$-configuration in $\mathcal{S}_n(M)$, where $M = \mathbb{R}^{k}$ is a Euclidean space. Then the set $\Sym_\varepsilon(X)$ is a compact $K$-approximate subgroup of the orthogonal group $O(k)$ for some finite $K$.
\end{theorem}

\begin{proof}
We verify the three defining properties of a compact approximate subgroup.

First, since
\[
W_p(e(X), X) = 0 \leq \varepsilon,
\]
the identity isometry $e$ belongs to $\Sym_\varepsilon(X)$.

Next, let $\pi \in \Sym_\varepsilon(X)$. As $\pi$ is an isometry and $W_p$ is a symmetric metric, we have
\[
W_p(\pi^{-1}(X), X) = W_p(X, \pi(X)) \leq \varepsilon,
\]
so $\pi^{-1} \in \Sym_\varepsilon(X)$. Hence, $\Sym_\varepsilon(X)$ is symmetric.

By Proposition~\ref{proposition:compact}, the set $\Sym_{2\varepsilon}(X)$ is a compact subset of the Lie group $O(k)$. Therefore, there exists a finite subset $F \subseteq \Sym_{2\varepsilon}(X)$ such that
\[
\Sym_{2\varepsilon}(X) \subseteq \bigcup_{f \in F} f \cdot \Sym_\varepsilon(X).
\]
By Lemma~\ref{lemma:composition_approximate}, for any $\pi_1, \pi_2 \in \Sym_\varepsilon(X)$, we have
\[
\pi_1 \pi_2 \in \Sym_{2\varepsilon}(X),
\]
and thus there exists $f \in F$ such that
\[
\pi_1 \pi_2 \in f \cdot \Sym_\varepsilon(X).
\]
Therefore, we obtain
\[
\Sym_\varepsilon(X) \cdot \Sym_\varepsilon(X) \subseteq F \cdot \Sym_\varepsilon(X),
\]
where $F$ is finite.

We conclude that $\Sym_\varepsilon(X)$ is a compact $K$-approximate subgroup of $O(k)$ with $K = |F|$.
\end{proof}

The approximate symmetry set $\Sym_\varepsilon(X)$ provides not only a qualitative description of the near-symmetries of the configuration $X$, but also serves as a quantitative indicator of its global symmetry structure. In particular, the Haar measure of $\Sym_\varepsilon(X)$ reflects the volume of symmetries that approximately preserve $X$ within a prescribed tolerance $\varepsilon$. A larger Haar measure indicates that $X$ admits a greater abundance of approximate symmetries, which suggests a higher degree of geometric regularity.

This perspective is consistent with standard practices in geometric group theory and topological data analysis, where invariant measures of transformation groups are often used to characterize the intrinsic structure of geometric and combinatorial objects.

\subsection{Symmetry measure}

While symmetry defect captures how far a configuration deviates from being symmetric, it is also useful to consider how symmetric a configuration actually is. To this end, we introduce the notion of a symmetry measure, which provides a way to quantify the degree of symmetry present in a configuration. In this subsection, we introduce the concept of the symmetry measure and establish its relationship with the symmetry defect.

Let $X$ be an $n$-configuration in $\mathcal{S}_n(M)$, where $M$ is $k$-dimensional Euclidean space $\mathbb{R}^k$. Let $\pi\in \Iso(\mathbb{R}^k)$ be an isometry on $\mathbb{R}^k$. The \textit{standard symmetry measure} of $X$ with respect to $\pi$ is defined by
\begin{equation*}
  \widetilde{\varpi}(X, \pi) = \sup_{\gamma} \left( \frac{\sum_{x \in X} \gamma(x) \cdot \pi(x)}{\sqrt{\sum_{x \in X} \|x\|^2} \sqrt{\sum_{x \in X} \|\pi(x)\|^2}} \right),
\end{equation*}
where the supremum is taken over all bijections $\gamma: X \to X$.

Since $X$ is finite, we may write $X = \{x_1, \dots, x_n\}$. Then the symmetry measure can be equivalently expressed as
\begin{equation*}
  \widetilde{\varpi}(X, \pi) = \sup_{\sigma \in S_n} \left( \frac{\sum_{i=1}^{n} x_{\sigma(i)} \cdot \pi(x_i)}{\sqrt{\sum_{i=1}^{n} \|x_i\|^2} \sqrt{\sum_{i=1}^{n} \|\pi(x_i)\|^2}} \right),
\end{equation*}
where $S_n$ denotes the symmetric group on $n$ elements.

Similar to the symmetry defect, we exclude perturbations and define the symmetry measure as follows.
\begin{definition}
Let $X$ be an $n$-configuration in $\mathcal{S}_n(\mathbb{R}^{k})$. For a given isometry $\pi$ of $\mathbb{R}^{k}$, the \textit{symmetry measure of $X$ with respect to $\pi$} is defined by
\[
  \varpi(X,\pi) := \sup_{\gamma \neq \mathrm{id}_X} \left( \frac{\sum_{x \in X} \gamma(x) \cdot \pi(x)}{\sqrt{\sum_{x \in X} \|x\|^2} \sqrt{\sum_{x \in X} \|\pi(x)\|^2}} \right),
\]
where the supremum is taken over all bijections $\gamma: X \to X$ other than the identity.
\end{definition}

Let $\Gamma\subseteq \Iso(\mathbb{R}^{k})$ be a collection of symmetries of $\mathbb{R}^{k}$. The \textit{symmetry measure} of $X$ with respect to $\Gamma$ is then defined by
\begin{equation*}
  \varpi(X, \Gamma) = \sup_{\pi \in \Gamma} \varpi(X, \pi),
\end{equation*}
where the supremum is taken over all $\pi \in \Gamma$.

\begin{example}
Now, consider the point set $X = \{(0,0), (0,2), (1,0)\}$. Let $\pi_{xy}$ be the reflection symmetry across the line $y = x$. Then, the set $\pi_{xy}(X) = \{(0,0), (2,0), (0,1)\}$ is obtained as the reflection of the points in $X$ under $\pi_{xy}$. We can now compute the symmetry measure of $X$ with respect to the reflection symmetry $\pi_{xy}$ as
\begin{equation*}
   \varpi(X, \pi_{xy}) = \frac{(0,2) \cdot (0,1) + (1,0) \cdot (2,0)}{\sqrt{5} \times \sqrt{5}} = 0.8.
\end{equation*}
Now, let $\Gamma_{1}$ be the set of reflection symmetries across all lines passing through the origin. Suppose the reflection occurs across the line $\ell_{\theta}$, which forms an angle $\theta$ with the x-axis. Through elementary calculations, we obtain the reflected point set as
\begin{equation*}
 \sigma_{\theta}(X) = \{(0,0), (2\sin 2\theta, -2\cos 2\theta), (\cos 2\theta, \sin 2\theta)\}.
\end{equation*}
Noting that $(4\sin 2\theta)/5 \leq 0.8$ and $(3\cos 2\theta)/5 \leq 0.6 \leq 0.8$, we conclude that the symmetry measure of the point set $X$ with respect to $\Gamma_{1}$ is
\begin{equation*}
 \varpi(X, \Gamma_{1}) = 0.8.
\end{equation*}
Now, let $\Gamma_{2}$ be the set of reflection symmetries across all lines through the centroid $(1/3, 2/3)$ of the set $X$. We use computer code to compute this process. Consider the reflection axis passing through $(1/3, 2/3)$ and forming an angle $\theta$ with the x-axis. We divide the range of $\theta$ from 0 to $\pi$ radians into 180 intervals, i.e., $\theta = k\pi/180$ for $k = 0, 1, \dots, 179$. For each value of $\theta$, we compute the symmetry measure $\varpi(X, \pi_{\theta})$. The maximum symmetry measure obtained is $\varpi(X, \Gamma_{2})=0.9899$. The corresponding reflection axis forms an angle of $168^{\circ}$ with the $x$-axis, or equivalently, an angle of $14\pi/15$ radians. The corresponding reflected point set is
\begin{equation*}
  (0.1932, 0.9090), (2.0203, 0.0955), (-0.2135, -0.0045).
\end{equation*}
Let $\Gamma$ be the set of reflection symmetries across all lines in the plane. To the best of our knowledge, it is not clear whether there is a result that can prove the following equality
\begin{equation*}
  \varpi(X, \Gamma) = \varpi(X, \Gamma_{2}).
\end{equation*}
The incenter of a triangle is the point of intersection of the angle bisectors, and it plays an important role in balancing the triangle's angles and side lengths. The circumcenter, on the other hand, is the center of the circumscribed circle, and its symmetry focuses more on the external geometry of the triangle. However, intuitively, the centroid of the triangle, as the center of symmetry, can be seen as the ``balance point'' of all the points. In contrast, the incenter and circumcenter's reflection point sets and the original point set's similarity might be influenced by the asymmetry of the geometric shape.

For simple geometric shapes like triangles, the centroid's contribution to all edges and angles is balanced, making the reflected point set, when using the centroid as the axis of symmetry, typically align or overlap with the original point set to a larger extent. For the examples mentioned above, we also computed the symmetry measure of the reflection symmetries across all lines through the incenter and circumcenter, yielding values of 0.9897 and 0.9507, respectively.
\end{example}

\begin{definition}
Let $X$ be an $n$-configuration in $\mathcal{S}_n(M)$. For a given isometry $\pi$ of $\mathbb{R}^{k}$, the \textit{symmetry measure} of $X$ is defined by
\[
  \varpi(X) := \sup_{\pi\in \Iso(\mathbb{R}^{k})} \varpi(X, \pi).
\]
\end{definition}

We further investigate the relationship between the symmetry defect and the symmetry measure. From now on, we assume that the metric space $M$ is Euclidean and that the symmetry defect is measured using the 2-Wasserstein distance. In this setting, the symmetry defect with respect to an isometry $\pi$ of $M$ can be written explicitly as
\[
\mu(X, \pi) = \left( \inf_{\substack{\gamma : X \to \pi(X) \\ \gamma \neq \pi|_X}} \sum_{x \in X} \| x - \gamma(x) \|^2 \right)^{1/2} =\left( \inf_{\gamma'\neq \mathrm{id}_{X}} \sum_{x \in X} \|\pi(x') - {\gamma'}(x')\|^2 \right)^{1/2}.
\]
Here, we make the change of variables $\gamma=\pi|_{X}\circ{\gamma'}^{-1}$ and $x=\gamma'(x')$. Hence, the symmetry defect can be written as the form
\[
  \mu(X,\pi) = \left( \inf_{\gamma\neq \mathrm{id}_{X}} \sum_{x \in X} \|\pi(x) - \gamma(x)\|^2 \right)^{1/2},
\]
where the infimum is taken over all bijections $\gamma: X \to X$ except the identity map.

\begin{proposition}\label{proposition:symmetry_inequality}
Let $X$ be an $n$-configuration in $\mathcal{S}_n(M)$. Let $\pi$ be an isometry on $M$. Then we have 
\[
  \varpi(X, \pi) \geq 1 - \frac{1}{2} \cdot \frac{\mu(X, \pi)^2}{\|X\| \cdot \|\pi(X)\|}.
\]
Here, $\|X\| = \sqrt{\sum\limits_{x\in X}\|x\|^2}$. In particular, if $\|X\| = \|\pi(X)\| = 1$, the inequality becomes an equality
\[
  \mu(X, \pi)^2 = 2(1 - \varpi(X, \pi)).
\]
\end{proposition}

\begin{proof}
By definition, we have
\[
\mu(X, \pi)^2 = \inf_{\gamma\neq \mathrm{id}_{X}} \sum_{x \in X} \|\pi(x) - \gamma(x)\|^2.
\]
For any bijection $\gamma: X \to X$, we use the identity
\[
\|\pi(x) - \gamma(x)\|^2 = \|\pi(x)\|^2 + \|\gamma(x)\|^2 - 2\pi(x) \cdot \gamma(x).
\]
Summing over $x \in X$ gives
\[
\sum_{x \in X} \|\pi(x) - \gamma(x)\|^2 = \sum_{x \in X} \|\pi(x)\|^2 + \sum_{x \in X} \|\gamma(x)\|^2 - 2 \sum_{x \in X} \pi(x) \cdot \gamma(x).
\]
Since $\gamma$ is a bijection, we have $\sum_{x \in X} \|\gamma(x)\|^2 = \sum_{x \in X} \|x\|^2 = \|X\|^2$. Similarly, let $\|\pi(X)\|^2 = \sum_{x \in X} \|\pi(x)\|^2$. Then we obtain
\[
\sum_{x \in X} \|\pi(x) - \gamma(x)\|^2 = \|\pi(X)\|^2 + \|X\|^2 - 2 \sum_{x \in X} \pi(x) \cdot \gamma(x).
\]
Taking the infimum over $\gamma \neq \mathrm{id}_X$ (which corresponds to taking the supremum for the dot product term), we get
\[
\mu(X, \pi)^2 = \|\pi(X)\|^2 + \|X\|^2 - 2 \cdot \sup_{\gamma\neq \mathrm{id}_X} \sum_{x \in X} \pi(x) \cdot \gamma(x).
\]
Dividing both sides by $\|X\| \cdot \|\pi(X)\|$ yields
\[
\frac{\mu(X, \pi)^2}{\|X\| \cdot \|\pi(X)\|} = \frac{\|\pi(X)\|}{\|X\|} + \frac{\|X\|}{\|\pi(X)\|} - 2 \varpi(X, \pi).
\]
Since the arithmetic-geometric mean inequality implies
\[
\frac{\|X\|}{\|\pi(X)\|} + \frac{\|\pi(X)\|}{\|X\|} \geq 2,
\]
we obtain
\[
\frac{\mu(X, \pi)^2}{\|X\| \cdot \|\pi(X)\|} \geq 2 - 2\varpi(X, \pi),
\]
which leads to the inequality
\[
\varpi(X, \pi) \geq 1 - \frac{1}{2} \cdot \frac{\mu(X, \pi)^2}{\|X\| \cdot \|\pi(X)\|}.
\]

In the special case where $\|\pi(X)\| = \|X\| = 1$, equality is achieved throughout, yielding the identity $\mu(X, \pi)^2 = 2(1 - \varpi(X, \pi))$.
\end{proof}

\begin{definition}
Let $X$ be a finite configuration in Euclidean space. Let $\pi$ be an isometry on $M$. The \textit{normalized symmetry defect} of $X$ with respect to $\pi$ is defined by
\[
  \overline{\mu}(X, \pi) := \frac{\mu(X, \pi)}{\sqrt{\|X\| \cdot \|\pi(X)\|}}.
\]
\end{definition}

By the proof of Proposition~\ref{proposition:symmetry_inequality}, we obtain the identity
\[
  \overline{\mu}(X, \pi)^2 = \frac{\|\pi(X)\|}{\|X\|} + \frac{\|X\|}{\|\pi(X)\|} - 2 \varpi(X, \pi),
\]
which provides an explicit quantitative relation between the normalized symmetry defect and the symmetry measure. In particular, when $\pi$ is a rotation about the origin or a reflection across a line (or hyperplane) passing through the origin, we have $\|X\| = \|\pi(X)\|$, and the identity simplifies to
\[
  \overline{\mu}(X, \pi)^2 = 2(1 - \varpi(X, \pi)).
\]

In practical computations, it is common to translate the configuration $X$ so that its centroid lies at the origin, and to consider symmetries whose centers of rotation or axes of reflection pass through the origin. This normalization not only reduces computational complexity, but also ensures that the considered symmetries are maximally representative of the object's intrinsic (approximate) symmetry.

\subsection{The stability of symmetry defect}

Consider the case where the underlying space is the Euclidean space $\mathbb{R}^k$. As previously discussed, the symmetry defect function $\mu(X, \pi)$ is continuous with respect to both the configuration $X \in \Conf_n(\mathbb{R}^k)$ and the isometry $\pi \in \Iso(\mathbb{R}^k)$. More generally, the symmetry defect also satisfies a Lipschitz continuity property, which is the notion of stability we aim to investigate.

\begin{theorem}\label{theorem:lipnitz_space}
Let $X, Y \in \Conf_n(\mathbb{R}^k)$ be two $n$-configurations in $\mathbb{R}^k$, and let $\pi \in \Iso(\mathbb{R}^k)$ be a fixed isometry. Then
\[
|\mu(X, \pi) - \mu(Y, \pi)| \leq 2 W_p(X, Y).
\]
\end{theorem}

\begin{proof}
Recall that the $p$-Wasserstein distance between $X$ and $Y$ is
\[
W_p(X,Y) = \left( \inf_{\gamma} \sum_{x \in X} d(x, \gamma(x))^p \right)^{1/p},
\]
where the infimum is taken over all bijections $\gamma: X \to Y$. Fix such a bijection $\delta: X \to Y$ realizing the infimum, then we have
\[
W_p(X,Y) = \left( \sum_{x \in X} d(x, \delta(x))^p \right)^{1/p}.
\]
Define the induced bijection $\pi(\delta): \pi(X) \to \pi(Y)$ by
\[
\pi(\delta)(\pi(x)) = \pi(\delta(x)).
\]
By definition, we have
\[
\mu(X,\pi) = \left( \inf_{\gamma \neq \pi|_X} \sum_{x \in X} d(x, \gamma(x))^p \right)^{1/p}.
\]
For each such $\gamma: X \to \pi(X)$, define $\gamma' : Y \to \pi(Y)$ by requiring the following diagram to commute:
\[
\xymatrix{
X \ar[r]^{\gamma} \ar[d]_{\delta} & \pi(X) \ar[d]^{\pi(\delta)} \\
Y \ar[r]_{\gamma'} & \pi(Y).
}
\]
It follows that $\gamma' = \pi(\delta) \circ \gamma \circ \delta^{-1}$. Then $\gamma = \pi|_X$ if and only if $\gamma' = \pi|_Y$. Therefore, one has
\[
\mu(Y, \pi) = \left( \inf_{\gamma \neq \pi|_X} \sum_{y \in Y} d\big(y, \pi(\delta)(\gamma(\delta^{-1}(y)))\big)^p \right)^{1/p}.
\]

We now estimate $\mu(X, \pi)$. For any $\gamma \neq \pi|_X$, by the triangle inequality
\[
d(x, \gamma(x)) \leq d(x, \delta(x)) + d(\delta(x), \gamma(x)),
\]
and applying Minkowski's inequality, we obtain
\begin{align*}
  \mu(X,\pi)  =& \inf_{\gamma \neq \pi|_X}    \left( \sum_{x \in X} d(x, \gamma(x))^p  \right)^{1/p}\\
   \leq & \inf_{\gamma\neq \pi|_{X}}\left[\left( \sum_{x \in X} d(x, \delta(x))^p \right)^{1/p}+\left( \sum_{x \in X} d(\delta(x), \gamma(x))^p \right)^{1/p}\right]\\
   = & W_p(X,Y) + A.
\end{align*}
Here, we denote
\[
A = \inf_{\gamma \neq \pi|_X} \left( \sum_{x \in X} d(\delta(x), \gamma(x))^p \right)^{1/p}.
\]

Similarly, we estimate $\mu(Y, \pi)$ as follows:
\begin{align*}
\mu(Y,\pi)
= &\inf_{\gamma \neq \pi|_X}\left(  \sum_{y \in Y} d\left( y, \pi(\delta)(\gamma(\delta^{-1}(y))) \right)^p \right)^{1/p} \\
   \leq & \inf_{\gamma\neq \pi|_{X}}\left[\left(\sum_{y\in Y}d(y, \gamma(\delta^{-1}(y)))^p \right)^{1/p}+\left(\sum_{y\in Y}d(\gamma(\delta^{-1}(y)), \pi(\delta)(\gamma(\delta^{-1}(y))))^p \right)^{1/p}\right]\\
   =&\inf_{\gamma \neq \pi|_{X}}\left[ \left( \sum_{x \in X} d(\delta(x), \gamma(x))^p \right)^{1/p}
+ \left( \sum_{x \in X} d\left( \gamma(x), \pi(\delta)(\gamma(x)) \right)^p \right)^{1/p}\right]\\
    = &  \inf_{\gamma\neq \pi|_{X}}\left[\left(\sum_{x\in X}d(\delta(x), \gamma(x))^p \right)^{1/p}+\left(\sum_{x\in X}d(x', \delta(x'))^p \right)^{1/p}\right]\\
    = & A +W_p(X,Y).
\end{align*}
In the above computation, we apply the change of variable $x' = \gamma(x)$, so that
\[
\pi(\delta)(\gamma(x)) = \pi(\delta(x')).
\]
Therefore, we have
\[
d\left(\gamma(x), \pi(\delta)(\gamma(x))\right)
= d\left( \pi(x'), \pi(\delta(x')) \right)
= d(x', \delta(x')),
\]
where the last equality follows from the fact that $\pi$ is an isometry.

Combining both estimates, we conclude that
\[
|\mu(X, \pi) - \mu(Y, \pi)| \leq |\mu(X, \pi) - A| + |\mu(Y, \pi) - A| \leq 2 W_p(X,Y),
\]
which completes the proof.
\end{proof}

\begin{example}
Let $\pi : \mathbb{R} \to \mathbb{R}$ be the reflection about the origin, defined by $\pi(x) = -x$. Consider the two configurations
\[
X = \{-1, 1\}, \quad Y = \{0, 2\} \in \Conf_2(\mathbb{R}).
\]
We use the symmetry defect function with $p = 1$, given by
\[
  \mu(X, \pi) = \inf_{\gamma \neq \pi|_X} \sum_{x \in X} d(x, \gamma(x)).
\]

We begin by computing $\mu(X, \pi)$. Since $\pi (X) = \{1, -1\} = X$, it follows that
\[
\mu(X, \pi)  = 0.
\]
For the configuration $Y = \{0, 2\}$, its image under $\pi$ is $\pi (Y) = \{0, -2\}$. The optimal bijection between $Y$ and $\pi (Y)$ matches $0$ with $-2$ and $2$ with $0$, giving a total cost of
\[
\mu(Y, \pi) =  |0 - (-2)| + |2 - 0| = 4.
\]
To compute the Wasserstein distance between $X$ and $Y$, observe that the optimal matching pairs $-1$ with $0$ and $1$ with $2$, resulting in a total cost of
\[
W_1(X, Y) = |-1 - 0| + |1 - 2|  = 2.
\]
Putting everything together, we find
\[
|\mu(X, \pi) - \mu(Y, \pi)| = |0 - 4| = 4 = 2 \cdot W_1(X, Y),
\]
so the Lipschitz inequality becomes an equality in this case. This example demonstrates that the constant $2$ in the Lipschitz bound for the symmetry defect function $\mu(X, \pi)$ with respect to the $1$-Wasserstein distance is indeed sharp.
\end{example}

Now consider two isometries $\pi, \sigma \in \Iso(\mathbb{R}^k)$. Let $X \subset \mathbb{R}^k$ be a fixed $n$-configuration.
A natural idea is to use
\[
W_p(\pi(X), \sigma(X))
\]
as a distance between $\pi$ and $\sigma$. However, this definition does not satisfy the separation property of a metric. In fact, when the symmetry group $\Sym(X)$ is nontrivial, for any $\pi, \sigma \in \Sym(X)$, we have
\[
W_p(\pi(X), \sigma(X)) = W_p(X, X) = 0.
\]
Taking $X$ into account, we instead define
\[
d_{X,p}(\pi, \sigma) = \left( \sum_{x\in X} d(\pi(x) ,\sigma(x))^p \right)^{1/p}.
\]
It is clear that $d_{X,p}(\pi, \sigma) = 0$ if and only if $\pi = \sigma$.

\begin{theorem}\label{theorem:symmetry_inequality}
Let $X \in \Conf_n(\mathbb{R}^k)$ be a fixed configuration, and let $\pi, \sigma \in \Iso(\mathbb{R}^k)$ be two isometries. Then we have
\[
\left| \mu(X, \pi) - \mu(X, \sigma) \right| \leq d_{X,p}(\pi, \sigma).
\]
\end{theorem}

\begin{proof}
For each bijection $\gamma: X \to \pi(X)$, define a corresponding bijection $\gamma' : X \to \sigma(X)$ by
\[
\gamma' = \sigma \circ \pi^{-1} \circ \gamma.
\]
It is clear that $\gamma \neq \pi|_X$ if and only if $\gamma' \neq \sigma|_X$.

By the triangle inequality, for each $x \in X$, we have
\[
d(x, \gamma(x)) \leq d(x, \gamma'(x)) + d(\gamma'(x), \gamma(x)).
\]
Applying Minkowski's inequality, we obtain
\begin{align*}
\mu(X,\pi) &= \inf_{\gamma \neq \pi|_X} \left( \sum_{x \in X} d(x, \gamma(x))^p \right)^{1/p} \\
&\leq \inf_{\gamma \neq \pi|_X} \left[ \left( \sum_{x \in X} d(x, \gamma'(x))^p \right)^{1/p} + \left( \sum_{x \in X} d(\gamma'(x), \gamma(x))^p \right)^{1/p} \right].
\end{align*}
We now simplify the second term. Note that
\begin{align*}
\sum_{x \in X} d(\gamma'(x), \gamma(x))^p
&= \sum_{x \in X} d\left( \sigma(\pi^{-1}(\gamma(x))), \gamma(x) \right)^p \\
&= \sum_{x' \in \pi(X)} d\left(\sigma(x'), \pi(x')\right)^p \\
&=  d_{X,p}(\pi, \sigma)^p,
\end{align*}
where we used the substitution $x=\gamma^{-1}(\pi(x'))$.

Hence, we conclude that
\[
\mu(X, \pi) \leq \mu(X, \sigma) + d_{X,p}(\pi, \sigma).
\]
Similarly, we also obtain
\[
\mu(X, \sigma) \leq \mu(X, \pi) + d_{X,p}(\pi, \sigma).
\]
Combining the two inequalities yields
\[
\left| \mu(X, \pi) - \mu(X, \sigma) \right| \leq d_{X,p}(\pi, \sigma),
\]
as desired.
\end{proof}

On the other hand, to define a distance between $\pi, \sigma \in \Iso(\mathbb{R}^k)$ that is independent of the configuration $X$, one may use the operator norm defined by
\[
\|\pi - \sigma\|_{\mathrm{op}}  = \sup_{\|x\| = 1} \| \pi x - \sigma x \|,
\]
where $x \in \mathbb{R}^k$.

\begin{theorem}
Let $X \in \Conf_n(\mathbb{R}^k)$ be a fixed configuration, and let $\pi, \sigma \in \Iso(\mathbb{R}^k)$ be two isometries. Then we have
\[
\left| \mu(X, \pi) - \mu(X, \sigma) \right| \leq C \|\pi - \sigma\|_{\mathrm{op}},
\]
where $C= \left( \sum_{x\in X} \|x\|^p \right)^{1/p}$.
\end{theorem}

\begin{proof}
By Theorem \ref{theorem:symmetry_inequality}, we have
\[
|\mu(X, \pi) - \mu(X, \sigma)| \leq \left( \sum_{x\in X} d(\pi(x) ,\sigma(x))^p \right)^{1/p}.
\]
Each term can be estimated using the operator norm
\[
d(\pi(x), \sigma(x)) = \| (\pi - \sigma)(x) \| \leq \|\pi - \sigma\|_{\mathrm{op}} \cdot \|x\|.
\]
Hence, we obtain
\[
|\mu(X, \pi) - \mu(X, \sigma)| \leq \|\pi - \sigma\|_{\mathrm{op}} \cdot \left( \sum_{x\in X} \|x\|^p \right)^{1/p} = C \|\pi - \sigma\|_{\mathrm{op}}.
\]
This proves the inequality.
\end{proof}

%% file: persistent_representation.tex
\section{Persistent representation theory}\label{section:persistent_representation}

In this section, we study persistent groups that encode the evolution of group structures over time or parametrization. By considering their linear representations, we gain access to powerful algebraic methods for analyzing time-varying or parametrized symmetries. In our analysis, `time' can be regarded as a parameter. Our focus lies in the structure and decomposition of persistent representations, the behavior of persistent irreducible subrepresentations, their module-theoretic interpretation, and the regular representations of persistence groups.

Computational tools such as SageMath, GAP, and MAGMA provide effective methods to calculate irreducible representations and character tables of finite groups \cite{bosma1997magma,gray2008sage,linton2007gap}. These software packages enable practical exploration and verification of representation-theoretic properties that are often difficult to handle manually.

For the representation theory of finite groups, classical references include \cite{curtis1966representation,fulton2013representation, serre1977linear}. Throughout this section, all groups considered are assumed to be finite, and the ground field $\mathbb{K}$ is assumed to have characteristic zero.

\subsection{Persistent group representations}

Let $G$ be a finite group. A \textit{representation} of $G$ over a field $\mathbb{K}$ is a homomorphism of groups
\[
\rho: G \to \mathrm{GL}(V).
\]
Here, $V$ is a finite-dimensional vector space over $\mathbb{K}$ and $\mathrm{GL}(V)$ denotes the group of invertible linear transformations on $V$. The space $V$ is called the \textit{representation space} and the dimension of $V$ is called the \textit{dimension} of the representation. From now on, all representation spaces considered henceforth are assumed to be finite-dimensional.

An important example of group representation is the \textit{regular representation}. Let $\mathbb{K}[G]$ be the group algebra of $G$ over $\mathbb{K}$. Then $\mathbb{K}[G]$ is a $\mathbb{K}$-vector space with basis elements indexed by elements of $G$, and group multiplication extends linearly to an action on the space
\[
g \cdot \left( \sum_{h \in G} a_h h \right) = \sum_{h \in G} a_h (g h), \quad \text{for } g \in G, a_{h}\in \mathbb{K}.
\]
This defines a representation $\rho_{G}:G\to \mathrm{GL}(\mathbb{K}[G])$ of $G$ on $\mathbb{K}[G]$, called the \textit{left regular representation}. The regular representation has dimension $|G|$ and contains every irreducible representation of $G$ as a subrepresentation, with multiplicities equal to their dimensions.

Let $\mathbf{Rep}_{\mathbb{K}}$ denote the category of finite group representations over a field $\mathbb{K}$. An object in $\mathbf{Rep}_{\mathbb{K}}$ is a triple $(G, V, \rho)$, where $G$ is a finite group, $V$ is a finite-dimensional $\mathbb{K}$-vector space, and $\rho: G \to \mathrm{GL}(V)$ is a group homomorphism. A morphism between two such objects $(G, V, \rho)$ and $(G', V', \rho')$ is a pair $(\phi, f)$, where $\phi: G \to G'$ is a group homomorphism and $f: V \to V'$ is a $\mathbb{K}$-linear map, such that the following condition
\[
f(\rho(g)v) = \rho'(\phi(g))f(v)
\]
holds for all $g \in G$ and $v \in V$.

\begin{definition}
Let $(T, \leq)$ be a poset, regarded as a small category. A \textit{persistence representation} over a field $\mathbb{K}$ is a covariant functor
\[
\mathcal{R}: T \to \mathbf{Rep}_{\mathbb{K}}.
\]
\end{definition}

Explicitly, the functor $\mathcal{R}$ assigns to each $t \in T$ a representation $(G_t, V_t, \rho_t)$, where $G_t$ is a finite group, $V_t$ is a finite-dimensional $\mathbb{K}$-vector space, and $\rho_t: G_t \to \mathrm{GL}(V_t)$ is a group homomorphism. For every pair $s \leq t$ in $T$, the functor assigns a morphism of representations
\[
(\phi_{s,t}, f_{s,t}) : (G_s, V_s, \rho_s) \to (G_t, V_t, \rho_t),
\]
where $\phi_{s,t}: G_s \to G_t$ is a group homomorphism and $f_{s,t}: V_s \to V_t$ is a $\mathbb{K}$-linear map satisfying the following commutative diagram
\[
\xymatrix{
V_s \ar[d]_{f_{s,t}} \ar[rr]^{\rho_s(g)} && V_s \ar[d]^{f_{s,t}} \\
V_t \ar[rr]^{\rho_t(\phi_{s,t}(g))} && V_t
}
\]
for any $g\in G_{s}$.

For $s \leq t$, the morphism of representations $(\phi_{s,t}, f_{s,t})$ induces a new group representation, which we call the \textit{$(s,t)$-persistent group representation}, denoted by $\mathcal{R}_{s,t}$. The associated representation space is defined as
\[
V_{s,t} = \mathrm{im}(f_{s,t}: V_s \to V_t),
\]
and is referred to as the \textit{$(s,t)$-persistent representation space}.

The representation $\mathcal{R}_{s,t}$ is a group homomorphism
\[
\rho_{s,t} : G_s \to \mathrm{GL}(V_{s,t}),
\]
where the group action is given by
\[
\rho_{s,t}(g)(v) = \rho_t(\phi_{s,t}(g))(v), \quad \text{for all } g \in G_s,\ v \in V_{s,t}.
\]
This action is well-defined. Indeed, since $v \in V_{s,t}$, there exists $v_{s} \in V_s$ such that $v = f_{s,t}(v_{s})$. Then
\[
\rho_{s,t}(g)(v) = \rho_t(\phi_{s,t}(g))(f_{s,t}(v_{s})) = f_{s,t}(\rho_s(g)(v_{s})) \in V_{s,t}.
\]
Hence, $V_{s,t}$ is stable under the action of $G_s$, and $\mathcal{R}_{s,t}$ defines a well-defined $G_s$-representation on the subspace $V_{s,t}$.

\begin{definition}
Let $\mathcal{R} : (\mathbb{R}, \leq) \to \mathbf{Rep}_{\mathbb{K}}$ be a persistence representation. A \textit{persistent subrepresentation} of $\mathcal{R}$ is a subfunctor $\mathcal{W} : (\mathbb{R}, \leq) \to \mathbf{Rep}_{\mathbb{K}}$ such that, for each $t \in \mathbb{R}$, $\mathcal{W}_t = (G_t, W_t, \rho_t|_{W_t})$ with $W_t \subseteq V_t$ a $G_t$-invariant subspace of $\mathcal{R}_t = (G_t, V_t, \rho_t)$,  and the structure maps $f_{s,t}^{\mathcal{W}} : W_s \to W_t$ are the restrictions of the maps $f_{s,t} : V_s \to V_t$ in $\mathcal{R}$, meaning $f_{s,t}^{\mathcal{W}} = f_{s,t}|_{W_s}$ with $f_{s,t}(W_s) \subseteq W_t$ for all $s \leq t$.
\end{definition}

Now, let $\mathbb{K} = \mathbb{C}$ be the field of complex numbers. Given a representation $\rho : G \to \mathrm{GL}(V)$ of a finite group $G$, the \textit{character} associated to $\rho$, denoted by $\chi_\rho$, is the function
\[
\chi_\rho : G \to \mathbb{C}, \quad \chi_\rho(g) = \tr(\rho(g)),
\]
where $\tr(\rho(g))$ denotes the trace of the linear operator $\rho(g)$ acting on the vector space $V$.

\begin{definition}
Let $\mathcal{R} : T \to \mathbf{Rep}_{\mathbb{K}}$ be a persistence representation.
For $s\leq t$, the \textit{$(s,t)$-persistent character} associated to $\mathcal{R}$ is the class function
\[
\chi_{\mathcal{R}_{s,t}} : G_s \to \mathbb{K}, \quad \chi_{\mathcal{R}_{s,t}}(g) = \tr \left( \rho_t(\phi_{s,t}(g))|_{V_{s,t}} \right),
\]
where $\tr$ denotes the trace of the linear operator restricted to $V_{s,t}$.
\end{definition}

\begin{example}
Consider the two-point poset $T = \{0 \le 1\}$ and work over the field $\mathbb{C}$. At time $t = 0$, let $\mathcal{R}_0 = (G_0, V_0, \rho_0)$ be a representation given by
\[
G_0 = \mathbb{Z}/2 = \{e, a\}, \qquad
V_0 = \mathbb{C}^2, \qquad
\rho_0(a) =
\begin{bmatrix}
1 & 0 \\
0 & -1
\end{bmatrix}.
\]
At time $t = 1$, define the representation $\mathcal{R}_1 = (G_1, V_1, \rho_1)$ as follows
\[
G_1 = \mathbb{Z}/2 \times \mathbb{Z}/2 = \{(e,e), (a,e), (e,b), (a,b)\}, \qquad
V_1 = \mathbb{C}^3,
\]
with representation defined on generators by
\[
\rho_1(a,e) =
\begin{bmatrix}
1 & 0 & 0 \\
0 & -1 & 0 \\
0 & 0 & 1
\end{bmatrix}, \qquad
\rho_1(e,b) =
\begin{bmatrix}
1 & 0 & 0 \\
0 & 1 & 0 \\
0 & 0 & -1
\end{bmatrix}.
\]
The morphism $(\phi_{0,1}, f_{0,1}) : \mathcal{R}_0 \to \mathcal{R}_1$ is given by the group homomorphism
\[
\phi_{0,1} : G_0 \to G_1, \qquad \phi_{0,1}(e) = (e,e), \quad \phi_{0,1}(a) = (a,e),
\]
and the $\mathbb{C}$-linear map
\[
f_{0,1} : V_0 \to V_1, \qquad
f_{0,1}\left(\begin{bmatrix}x \\ y\end{bmatrix}\right) =
\begin{bmatrix}
x \\
y \\
0
\end{bmatrix}.
\]
On can check that $(\phi_{0,1}, f_{0,1})$ is a morphism of representations, i.e.,
\[
f_{0,1}(\rho_0(g)v) = \rho_1(\phi_{0,1}(g)) f_{0,1}(v), \quad \forall g \in G_0,\ v \in V_0.
\]
The $(0,1)$-persistent representation space is given by
\[
V_{0,1} = \operatorname{im}(f_{0,1}) = \left\{
\begin{bmatrix}
x \\
y \\
0
\end{bmatrix} : x, y \in \mathbb{C}
\right\} \subset V_1.
\]
This subspace is preserved by the action of $G_0$ via the map $\rho_{0,1} : G_0 \to \mathrm{GL}(V_{0,1})$, defined by
\[
\rho_{0,1}(g) = \rho_1(\phi_{0,1}(g))|_{V_{0,1}}, \quad \forall g \in G_0.
\]
Restricting to $V_{0,1}$ yields
\[
\rho_{0,1}(e) =
\begin{bmatrix}
1 & 0 \\
0 & 1
\end{bmatrix}, \qquad
\rho_{0,1}(a) =
\begin{bmatrix}
1 & 0 \\
0 & -1
\end{bmatrix}.
\]
Therefore, the associated $(0,1)$-persistent character is the class function
\[
\chi_{\mathcal{R}_{0,1}} : G_0 \to \mathbb{C}, \qquad
\chi_{\mathcal{R}_{0,1}}(e) = \tr(\rho_{0,1}(e)) = 2, \quad
\chi_{\mathcal{R}_{0,1}}(a) = \tr(\rho_{0,1}(a)) = 0.
\]
\end{example}

In addition to the persistent character, we define the \textit{$(s,t)$-persistent dimension} of a persistence representation $\mathcal{R}$ as $d_{s,t}(\mathcal{R}) = \dim V_{s,t}$ the dimension of the $(s,t)$-persistent representation space. A fundamental observation is that the persistent dimension equals the value of the persistent character at the identity element
\[
d_{s,t}(\mathcal{R}) = \chi_{\mathcal{R}_{s,t}}(e).
\]
This follows from the fact that $\rho_{s,t}(e)$ acts as the identity on $V_{s,t}$, so its trace equals $\dim V_{s,t}$.

\subsection{Persistent irreducible subrepresentations}

Let $\mathcal{R} : T \to \mathbf{Rep}_{\mathbb{K}}$ denote a persistence representation. For any pair of parameters $s \leq t$, the representation $\mathcal{R}_{s} = (G_s, V_s, \rho_s)$ admits a decomposition into a direct sum of irreducible subrepresentations, and consequently, the representation space $V_s$ decomposes into a direct sum of irreducible $G_{s}$-invariant subspaces.

Let $W_s \subseteq V_s$ be an irreducible $G_{s}$-invariant subspace. We define its persistent image at time $t$ as
\[
W_{s,t} = f_{s,t}(W_s),
\]
where $f_{s,t}$ is the structure map of the persistence representation $\mathcal{R}$. The image $W_{s,t}$ is a subspace of $V_{s,t} = f_{s,t}(V_s)$ that is invariant under the action of $G_{s}$. Specifically, for any $w \in W_{s,t}$, there exists $w_s \in W_s$ such that $w = f_{s,t}(w_s)$, and for any $g \in G_{s}$,
\[
\rho_{s,t}(g)(w) = \rho_t(\phi_{s,t}(g))(f_{s,t}(w_{s})) = f_{s,t}(\rho_s(g)(w_{s})) \in f_{s,t}(W_s) = W_{s,t}.
\]
Thus, $W_{s,t}$ is indeed a $G_{s}$-invariant subspace of $V_{s,t}$.

\begin{proposition}\label{proposition:irreducible}
Let $W_s \subseteq V_s$ be an irreducible $G_{s}$-invariant subspace. Then, the $G_{s}$-invariant subspace $W_{s,t}$ of $V_{s,t}$ is either irreducible or zero space.
\end{proposition}

\begin{proof}
Suppose $W_{s,t} \neq \{0\}$ and is reducible. Then there exist nontrivial $G_{s}$-invariant subspaces $U_1, U_2 \subseteq W_{s,t}$ such that
\[
W_{s,t} = U_1 \oplus U_2.
\]
Consider the preimages under the linear map $f_{s,t}: V_s \to V_t$ given by
\[
U_1' = f_{s,t}^{-1}(U_1) \cap W_s, \quad U_2' = f_{s,t}^{-1}(U_2) \cap W_s.
\]
We claim that $U_1', U_2'$ are $G_{s}$-invariant subspaces of $W_s$. Indeed, for any $g \in G_{s}$ and $u_1' \in U_1'$, we have
\[
f_{s,t}(\rho_{s}(g)  u_1') = \rho_{t}(\phi_{s,t}(g)) f_{s,t}(u_1') \in \rho_{t}(\phi_{s,t}(g)) U_1 \subseteq U_1,
\]
so $\rho_{s}(g)  u_1' \in f_{s,t}^{-1}(U_1) \cap W_s = U_1'$. Thus $U_1'$ is $G_{s}$-invariant; similarly, so is $U_2'$.

Since $f_{s,t}|_{W_s}: W_s \to W_{s,t}$ is a surjective $G_{s}$-module homomorphism and $U_1, U_2$ are nontrivial, it follows that $U_1', U_2'$ are nontrivial subspaces of $W_s$.

Next, we show that $U_1' \cap U_2' = \{0\}$. Suppose there is a nonzero element $x \in U_1' \cap U_2'$. Then
\[
f_{s,t}(x) \in U_1 \cap U_2 = \{0\}.
\]
So $x \in \ker f_{s,t} \cap W_s$. Thus, $\ker f_{s,t} \cap W_s$ is a nonzero $G_{s}$-invariant subspace of $W_{s}$. Since $W_s$ is irreducible, it must be that $\ker f_{s,t} \cap W_s = W_s$, i.e., $f_{s,t}(W_s) = \{0\}$, contradicting the assumption that $W_{s,t} \neq \{0\}$.

Therefore, $U_1' \cap U_2' = \{0\}$. Now, since $W_s$ is irreducible, any nonzero $G_s$-invariant subspace of $W_s$ must be equal to $W_s$ itself. We have shown that $U_1'$ is a nonzero $G_s$-invariant subspace, so $U_1' = W_s$.
This implies $f_{s,t}(W_s) \subseteq U_1$. However, $f_{s,t}(W_s) = W_{s,t} = U_1 \oplus U_2$. This forces $U_2$ to be the zero subspace, which contradicts the assumption that $W_{s,t}$ is reducible (requiring $U_2$ to be nontrivial).

We conclude that $W_{s,t}$ must be either irreducible or the zero subspace.
\end{proof}

\begin{corollary}\label{corollary:irreducible}
Let $W_s \subseteq V_s$ be an irreducible $G_{s}$-invariant subspace. If the restriction $f_{s,t}|_{W_s}$ is injective, then $W_{s,t}$ is an irreducible $G_{s}$-invariant subspace of $V_{s,t}$.
\end{corollary}

\begin{proof}
Since the restriction $f_{s,t}|_{W_s}$ is injective, the image $W_{s,t} = f_{s,t}(W_s)$ is nonzero. By Proposition~\ref{proposition:irreducible}, any nonzero $G_s$-invariant image of an irreducible $G_s$-representation under a $G_s$-equivariant linear map remains irreducible. Hence, $W_{s,t}$ is an irreducible $G_s$-invariant subspace of $V_{s,t}$.
\end{proof}

\begin{remark}
Let $W_s \subseteq V_s$ be an irreducible $G_s$-invariant subspace. It is worth noting that although the image $W_{s,t} = f_{s,t}(W_s)$ is a subspace of $V_t$, it is not necessarily $G_t$-invariant. However, when the group homomorphism $\phi_{s,t} : G_s \to G_t$ is surjective, $W_{s,t}$ can be regarded as a $G_t$-invariant subspace. 
For any $h \in G_t$, since $\phi_{s,t}$ is surjective, there exists $g \in G_s$ such that $\phi_{s,t}(g) = h$. The group action is given by
\[
\rho_t(h)(w) = \rho_t(\phi_{s,t}(g))(f_{s,t}(w_s)) = f_{s,t}(\rho_s(g)(w_s)) \in W_{s,t},
\]
where $w = f_{s,t}(w_s) \in W_{s,t}$. This action is well-defined: if $\phi_{s,t}(g') = h$ as well, then $\phi_{s,t}(g^{-1}g') = e$, and one can verify that $f_{s,t}(\rho_s(g)w_s - \rho_s(g')w_s) = 0$, ensuring consistency.
\end{remark}

\begin{example}
Let $\mathcal{R} : (\mathbb{R}, \leq) \to \mathbf{Rep}_{\mathbb{R}}$ be a persistence representation defined as follows. At parameter $s$, let $\mathcal{R}_s = (G, V_s, \rho_s)$, where $G = \mathbb{Z}/2 = \{e, g\}$ is the cyclic group of order $2$ and $V_s = \mathbb{K}[G]$ is the regular representation. At parameter $t > s$, let $\mathcal{R}_t = (H, V_t, \rho_t)$, where $H = \{e\}$ is the trivial group and $V_t = \mathbb{K}[H]$ is its regular representation.

The group homomorphism $\phi_{s,t} : G \to H$ is defined by $\phi_{s,t}(g) = e$ for all $g \in G$. This induces a linear map $f_{s,t} : V_s \to V_t$ given by
\[
f_{s,t}(e) = f_{s,t}(g) = e \in V_t = \mathbb{K}[H].
\]
Let $W_s = \langle e - g \rangle \subset V_s$ be the one-dimensional irreducible subrepresentation corresponding to the sign representation of $G$. Then,
\[
f_{s,t}(W_s) = f_{s,t}(\langle e - g \rangle) = \langle f_{s,t}(e - g) \rangle = \langle e - e \rangle = \{0\}.
\]
Thus, the persistent image at time $t$ is the zero subspace
\[
W_{s,t} = f_{s,t}(W_s) = \{0\}.
\]
Therefore, $W_{s,t}$ is not an irreducible subspace of $V_t$. This example illustrates that if the restriction $f_{s,t}|_{W_s}$ is not injective, then the irreducibility may not be preserved.
\end{example}

\begin{example}
Let $\mathcal{R} : (\mathbb{R}, \leq) \to \mathbf{Rep}_{\mathbb{R}}$ be a persistence representation defined as follows. At parameter $s$, let $\mathcal{R}_s = (G, V_s, \rho_s)$, where $G = \mathbb{Z}/2 \times \mathbb{Z}/2 = \langle a, b \rangle$ is the Klein four-group, and let $V_s = \mathbb{K}[G]$ be the regular representation of $G$. Then $V_s \cong \mathbb{K}^4$, with basis elements $e, a, b, ab$.

At parameter $t > s$, let $\mathcal{R}_t = (H, V_t, \rho_t)$, where $H = \mathbb{Z}/2 = \langle g \rangle$ is the cyclic group of order 2, and $V_t = \mathbb{K}[H]$ is its regular representation, with basis $e', g$.

Define the group homomorphism $\phi_{s,t} : G \to H$ by $\phi_{s,t}(a) = g$ and $\phi_{s,t}(b) = e'$. This induces a $\mathbb{K}$-linear map $f_{s,t} : V_s \to V_t$, defined by linearly extending the map on group basis elements
\[
f_{s,t}(x) = \phi_{s,t}(x) \in H \subset V_t.
\]
Explicitly, we have
\[
f_{s,t}(e) = e', \quad f_{s,t}(a) = g, \quad f_{s,t}(b) = e', \quad f_{s,t}(ab) = g.
\]

The regular representation $V_s$ of $G = \mathbb{Z}/2 \times \mathbb{Z}/2$ decomposes into four one-dimensional irreducible subrepresentations. Let us consider two such subspaces
\[
W_s = \langle e - a - b + ab \rangle, \quad W_s' = \langle e + a - b - ab \rangle.
\]
These are one-dimensional subrepresentations corresponding to distinct nontrivial characters of $G$.

Now we consider their images under $f_{s,t}$. A straightforward calculation shows that
\begin{align*}
    & f_{s,t}(e - a - b + ab) = e' - g - e' + g = 0,  \\
    & f_{s,t}(e + a - b - ab) = e' + g - e' - g = 0.
\end{align*}

Therefore, we have
\[
f_{s,t}(W_s) = \{0\}, \quad f_{s,t}(W_s') = \{0\}.
\]
In this example, although $W_s$ and $W_s'$ are distinct irreducible subrepresentations of $V_s$, their persistent images in $V_t$ are both trivial.
\end{example}

\begin{definition}
Let $\mathcal{R} : (\mathbb{R}, \leq) \to \mathbf{Rep}_{\mathbb{K}}$ be a persistence representation.
Fix $s \in \mathbb{R}$, and let $B_s \subseteq \mathbb{R}$ be an interval of the form $[s, t)$, $[s, t]$, or $[s, \infty)$ for some $t$.
An irreducible $G_s$-invariant subspace $W_s \subseteq V_s$ is said to be \textit{stably persistently irreducible} on $B_s$ if, for every $r \in B_s$, the persistent image
\[
W_{s,r} = f_{s,r}(W_s)
\]
remains an irreducible $G_s$-invariant subspace of $V_{s,r} \subseteq V_r$ under the induced $G_s$-action via $\rho_{s,r}(g) = \rho_r(\phi_{s,r}(g))$ for all $g \in G_s$.
\end{definition}

Given a persistence representation $\mathcal{R} : (\mathbb{R}, \leq) \to \mathbf{Rep}_{\mathbb{K}}$, an irreducible $G_s$-invariant subspace $W_s \subseteq V_s$ at time $s$ may either persist across subsequent time parameters or vanish at some later time. The behavior of such subspaces over time is entirely determined by the structure of the persistence representation $\mathcal{R}$.

\begin{definition}
Let $\mathcal{R} : (\mathbb{R}, \leq) \to \mathbf{Rep}_{\mathbb{K}}$ be a persistence representation. For each parameter $s \in \mathbb{R}$, let $\mathcal{R}_s = (G_s, V_s, \rho_s)$.
Fix an irreducible $G_s$-invariant subspace $W_s \subseteq V_s$. The \textit{persistence interval} of $W_s$ is the maximal interval $B_s = [s, b)$, $[s,b]$, or $[s, \infty)$ with $b \geq s$, such that for every $r \in B_s$, the persistent image
\[
W_{s,r} = f_{s,r}(W_s) \subseteq V_r
\]
is nonzero.

In addition, if for all $t < s$, there does not exist an irreducible $G_t$-invariant subspace $W_t \subseteq V_t$ such that
\[
f_{t,s}(W_t) = W_s,
\]
then $B_s$ is called a \textit{persistence bar} (or simply a \textit{bar}) of the irreducible subspace $W_s$.
\end{definition}

\begin{definition}
The \textit{irreducible barcode} of a persistence representation $\mathcal{R}$ is the multiset of persistence bars
\[
\{ B_{s_i} \}_{i \in I},
\]
where each $B_{s_i}$ corresponds to a persistence bar birth at $s_i$.
\end{definition}

Each bar can be further equipped with the dimension $\dim(W_{s_i})$ as a weight, encoding additional information beyond the classical barcode.

In practical computations, the set of filtration parameters is typically discrete. However, this does not prevent us from approximating the continuous case. Therefore, it is also natural to represent barcodes using continuous intervals.

\begin{example}\label{example:representation_barcode}
Consider the group $G = \mathbb{Z}/2 = \{e, g\}$ and the persistence representation
$\mathcal{V} : \{0,1,2\} \to \mathbf{Rep}_{\mathbb{C}}$, where each $V_t$ is a finite-dimensional $G$-representation.

At time $t=0$, set $V_0 = T_0 \oplus S_0$ with $T_0 = \langle v_1 \rangle$ the trivial representation defined by $g \cdot v_1 = v_1$, and $S_0 = \langle v_2 \rangle$ the sign representation with $g \cdot v_2 = -v_2$.

At time $t=1$, define $V_1 = T_1 \oplus S_1$, where $T_1 = \langle w_1 \rangle$ and $S_1 = \langle w_2 \rangle$ are the trivial and the sign representations, respectively.

At $t=2$, let $V_2 = T_2 = \langle u \rangle$ be a one-dimensional trivial representation.

The $G$-equivariant linear maps are given by
\[
f_{0,1} : V_0 \to V_1, \quad f_{0,1}(v_1) = w_1, \quad f_{0,1}(v_2) = w_2,
\]
\[
f_{1,2} : V_1 \to V_2, \quad f_{1,2}(w_1) = u, \quad f_{1,2}(w_2) = 0.
\]
By composition, we have
\[
f_{0,2} := f_{1,2} \circ f_{0,1}, \quad f_{0,2}(v_1) = u, \quad f_{0,2}(v_2) = 0.
\]

Analyzing the persistence of irreducible $G$-invariant subspaces, the trivial subspace $T_0 = \langle v_1 \rangle$ persists through $T_1 = \langle w_1 \rangle$ to $T_2 = \langle u \rangle$, giving a bar $[0,2]$. The sign subspace $S_0 = \langle v_2 \rangle$ persists to $S_1 = \langle w_2 \rangle$ but vanishes at time 2, giving a bar $[0,1]$.

Hence, the barcode associated with $\mathcal{V}$ is $\{ [0,2], \; [0,1] \}$.
\end{example}

\subsection{Decomposition of persistence representations}

In this section, we consider the decomposition of persistence representations as direct sums of irreducible persistence subrepresentations. To ensure that the decomposition has finitely many components, we assume that the persistence representations under consideration change only at finitely many time points. Under this assumption, the filtration parameter can be taken to be the poset $(\mathbb{Z}, \leq)$. Consequently, the collection of persistent irreducible subrepresentations gives rise to only finitely many bars. Each bar corresponds to an interval of one of the following forms: $[a,b]$, $[a,\infty)$, or $(-\infty,b]$, where $a, b \in \mathbb{Z}$. Furthermore, we fix the group actions over time to ensure that each irreducible subrepresentation remains irreducible at every time step.

Let $G$ be a finite group, and let $\mathbf{Rep}_{\mathbb{K}}^{G}$ denote the category of finite-dimensional $\mathbb{K}$-representations of the fixed group $G$. The persistence representation
\[
\mathcal{R} : (\mathbb{Z}, \leq) \to \mathbf{Rep}_{\mathbb{K}}^{G}
\]
is a functor in which the underlying $G$-representation structure is preserved as the parameter varies.

\begin{definition}[Irreducible persistence representation]
A persistence representation $\mathcal{R} : (\mathbb{Z}, \leq) \to \mathbf{Rep}_{\mathbb{K}}^{G}$ is said to be \textit{irreducible} on an interval $B$ if
\begin{itemize}
    \item For all $r \in B$, the representation $\mathcal{R}_r = (G, V_r, \rho_r)$ is irreducible; that is, $V_r$ contains no nontrivial $G$-invariant proper subspaces;
    \item For all $r \notin B$, we have $V_r = 0$;
    \item For all $r, s \in B$ with $r \leq s$, the structure maps $f_{r,s} \colon V_r \to V_s$ are isomorphisms.
\end{itemize}
We refer to such a representation as an \textit{irreducible persistence representation supported on $B$}.
\end{definition}

\begin{example}[Irreducible persistence subrepresentation]\label{example:irreducible}
Let $\mathcal{R} : (\mathbb{Z}, \leq) \to \mathbf{Rep}_{\mathbb{K}}^{G}$ be a persistence representation, and let $\mathcal{R}_t = (G, V_t, \rho_t)$ for some $t \in \mathbb{Z}$. Given a parameter $t$ and an irreducible $G$-subrepresentation $W_t \subseteq V_t$, we can always find a minimal index $s \leq t$ (possibly $s = -\infty$) and an irreducible $G$-subrepresentation $W_s \subseteq V_s$ such that
\[
f_{s,t}(W_s) = W_t,
\]
where $f_{s,t}$ denotes the structure map of the persistence representation $\mathcal{R}$.

Such a pair $(s, W_s)$ always exists since at least $f_{t,t}(W_t) = W_t$ holds. However, the choice of $W_s$ may not be unique: even when the minimal $s$ is fixed, there may exist distinct irreducible $G$-subrepresentations $W_s, W_s' \subseteq V_s$ with $f_{s,t}(W_s) = f_{s,t}(W_s') = W_t$. In such a case, we may arbitrarily select one, say $W_s$. Thus, we always obtain a valid pair $(W_s, W_t)$. Moreover, for any $r \in [s, t]$, the composition of the structure maps satisfies
\[
f_{r,t}(f_{s,r}(W_s)) = f_{s,t}(W_s) = W_t.
\]

We now define an irreducible persistence subrepresentation
\[
\mathcal{R}^{(W_s,W_t)} : (\mathbb{Z}, \leq) \to \mathbf{Rep}_{\mathbb{K}}^{G}
\]
supported on the interval $[s,t]$ as follows:
\begin{itemize}
    \item For each $r \in [s,t]$, let
    \[
    \mathcal{R}^{(W_s,W_t)}_{r} = \left(G, W_{s,r}, \rho_{r}|_{W_{s,r}}\right),
    \]
    where $W_{s,r} = f_{s,r}(W_s) \subseteq V_r$ is the image of $W_s$ under the structure map, and $\rho_{r}|_{W_{s,r}}$ denotes the restriction of $\rho_r$ to the $G$-invariant subspace $W_{s,r}$.

    \item For $r \notin [s,t]$, define $\mathcal{R}^{(W_s,W_t)}_{r} = (G, 0, 0)$.

    \item For all $r_1 \leq r_2$ in $[s,t]$, the structure map
    \[
    f^{\mathcal{I}}_{r_1, r_2} := f_{r_1, r_2}|_{W_{s,r_1}} : W_{s,r_1} \to W_{s,r_2}
    \]
    is defined as the restriction of $f_{r_1, r_2}$ to $W_{s,r_1}$. This map is $G$-equivariant and an isomorphism since
    \[
    f^{\mathcal{I}}_{r_1, r_2}(W_{s,r_1}) = f_{r_1, r_2}(f_{s,r_1}(W_s)) = f_{s,r_2}(W_s) = W_{s,r_2}.
    \]
\end{itemize}

To verify that $\mathcal{R}^{(W_s,W_t)}$ defines a functor, it suffices to check the commutativity of the following diagram for all $g \in G$.
\[
\xymatrix{
W_{s,r_1} \ar[d]_{f^{\mathcal{I}}_{r_1, r_2}} \ar[rrr]^{\rho_{r_1}(g)} &&& W_{s,r_1} \ar[d]^{f^{\mathcal{I}}_{r_1, r_2}} \\
W_{s,r_2} \ar[rrr]^{\rho_{r_2}(g)} &&& W_{s,r_2}
}
\]
Indeed, for any $w \in W_{s,r_1}$, we compute
\begin{align*}
\rho_{r_2}(g)(f^{\mathcal{I}}_{r_1, r_2}(w))
&= \rho_{r_2}(g)(f_{r_1,r_2}(w)) \\
&= f_{r_1,r_2}(\rho_{r_1}(g)(w)) \\
&= f^{\mathcal{I}}_{r_1, r_2}(\rho_{r_1}(g)(w)),
\end{align*}
so the diagram commutes.

By construction, the functor $\mathcal{R}^{(W_s,W_t)}$ is irreducible and is supported exactly on the interval $[s,t]$. Moreover, $\mathcal{R}^{(W_s,W_t)}$ can be regarded as a subrepresentation of the original persistence representation $\mathcal{R}$, since each $(G, W_{s,r}, \rho_{r}|_{W_{s,r}})$ is a $G$-subrepresentation of $(G, V_r, \rho_r)$.
\end{example}

Notations as in Example \ref{example:irreducible}. Let $\mathcal{R}^{(W_s, W_t)} : (\mathbb{Z}, \leq) \to \mathbf{Rep}_{\mathbb{K}}^{G}$ denote the irreducible persistence subrepresentation of $\mathcal{R}$ constructed from the pair $(W_s, W_t)$. We say that an irreducible subrepresentation $W_r$ \textit{lies in} $\mathcal{R}^{(W_s, W_t)}$ if $W_r = f_{s,r}(W_s)$.

\begin{definition}
A \textit{decomposition of a persistence representation} is a decomposition of the functor $\mathcal{R}: (\mathbb{Z}, \leq) \to \mathbf{Rep}_{\mathbb{K}}$ into a direct sum of two persistence representations, $\mathcal{R}'$ and $\mathcal{R}''$, i.e., for each $t \in \mathbb{Z}$, we have
\[
\mathcal{R}_t \cong \mathcal{R}'_t \oplus \mathcal{R}''_t,
\]
where $\mathcal{R}'_t = (G_t, V'_t, \rho'_t)$ and $\mathcal{R}''_t = (G_t, V''_t, \rho''_t)$, with $V_t = V'_t \oplus V''_t$ and $\rho_t = \rho'_t \oplus \rho''_t$.

For each pair $s \leq t$, the morphism
\[
(\phi_{s,t}, f_{s,t}) : \mathcal{R}_s \to \mathcal{R}_t
\]
is the direct sum of the morphisms of representations
\[
(\phi_{s,t}, f'_{s,t}) \oplus (\phi_{s,t}, f''_{s,t}),
\]
where the group homomorphism component $\phi_{s,t}$ is the same in both summands, and the linear map satisfies
\[
f_{s,t} = f'_{s,t} \oplus f''_{s,t}.
\]
Each $f'_{s,t}$ and $f''_{s,t}$ satisfies the following commutative diagram for all $g \in G_s$.
\[
\xymatrix{
V'_s \ar[d]_{f'_{s,t}} \ar[rr]^{\rho'_s(g)} && V'_s \ar[d]^{f'_{s,t}} \\
V'_t \ar[rr]^{\rho'_t(\phi_{s,t}(g))} && V'_t
}
\qquad
\xymatrix{
V''_s \ar[d]_{f''_{s,t}} \ar[rr]^{\rho''_s(g)} && V''_s \ar[d]^{f''_{s,t}} \\
V''_t \ar[rr]^{\rho''_t(\phi_{s,t}(g))} && V''_t
}
\]
\end{definition}

To study the decomposition of persistence representations, a finiteness condition is required. Specifically, a persistence representation $\mathcal{R} : (\mathbb{Z}, \leq) \to \mathbf{Rep}_{\mathbb{K}}^{G}$ is said to be \textit{finite} if there exist only finitely many pairs $r \leq s$ such that the structure map $\mathcal{R}_{r\leq s}:\mathcal{R}_{r}\to \mathcal{R}_{s}$ is not an isomorphism.
Equivalently, $\mathcal{R}$ is finite if there exist integers $p$ and $q$ such that for any $p_1 \leq p_2 \leq p$, the structure map $\mathcal{R}_{p_1\leq p_2}: \mathcal{R}_{p_1}\to \mathcal{R}_{p_2}$ is an isomorphism, and for any $q \leq q_1 \leq q_2$, the structure map $\mathcal{R}_{q_1\leq q_2}:\mathcal{R}_{q_1}\to \mathcal{R}_{q_2}$ is also an isomorphism.

\begin{theorem}[Decomposition of persistence representations]\label{theorem:decomposition_representation}
Let $\mathcal{R} : (\mathbb{Z}, \leq) \to \mathbf{Rep}_{\mathbb{K}}^{G}$ be a finite persistence representation. Then there exists a finite collection $\{(U_i, U_{i}')\}_{i \in I}$ of irreducible $G$-invariant subspaces $U_i \subseteq V_{s_i}$ and $U_{i}' \subseteq V_{t_i}$ such that
\[
\mathcal{R} \cong \bigoplus_{i \in I} \mathcal{R}^{(U_i, U_{i}')},
\]
where each $\mathcal{R}^{(U_i, U_{i}')}$ is an irreducible persistence representation supported on the interval $[s_i, t_i]$.
\end{theorem}

\begin{proof}
Since $\mathcal{R} : (\mathbb{Z}, \leq) \to \mathbf{Rep}_{\mathbb{K}}^{G}$ is finite, there exists an integer $q \in \mathbb{Z}$ such that for all $t \geq q$, the morphism of representations
\[
\mathcal{R}_{q\leq t}: \mathcal{R}_q \xrightarrow{\cong} \mathcal{R}_t
\]
are isomorphisms. In particular, we define the stabilized representation at infinity by
\[
\mathcal{R}_{\infty} = \mathcal{R}_{t} \quad \text{for any } t \geq q,
\]
which is well-defined up to canonical isomorphism. The representation $V_{\infty} = \mathcal{R}_{\infty}$ is completely reducible and can be written as a finite direct sum of irreducible subspaces
\[
V_{\infty} = \bigoplus_j W_{\infty}^{(j)}.
\]
By the construction in Example~\ref{example:irreducible}, each $W_{\infty}^{(j)}$ determines an irreducible persistence subrepresentation $\mathcal{R}^{(U^{(j)}, W_{\infty}^{(j)})}$ supported on an interval $[a, \infty)$ for some $a \leq q$. This gives a collection $\mathcal{Z}_{\infty} = \{(U^{(j)}, W_{\infty}^{(j)})\}_j$ of irreducible persistence pairs. By the finiteness assumption, every irreducible subrepresentation $W_q \subseteq V_q$ must lie in one such subrepresentation, i.e., $f_{s,\infty}(W_q) = W_{\infty}^{(j)}$ for some $j$.

We now construct the collections $\mathcal{Z}_t$ inductively so that for each $r \geq t$, any irreducible subrepresentation $W_r \subseteq V_r$ lies in some $\mathcal{R}^{(U_i, U_i')}$ indexed by the pair $(U_i, U_i')$ in $\mathcal{Z}_t$.

Start with $\mathcal{Z}_q = \mathcal{Z}_{\infty}$. Suppose $\mathcal{Z}_t$ has been constructed for some $t \leq q$. We construct $\mathcal{Z}_{t-1}$ as follows. Decompose $V_{t-1}$ into irreducible subspaces
\[
V_{t-1} = \bigoplus_j W_{t-1}^{(j)}.
\]
For any $W_{t-1}^{(j)}$ not already lying in a subrepresentation indexed by a pair in $\mathcal{Z}_t$, construct a new irreducible subrepresentation $\mathcal{R}^{(U^{(j)}, W_{t-1}^{(j)})}$ supported on $[a, t-1]$ for some $a \leq t-1$, again using Example~\ref{example:irreducible}. Add each such pair $(U^{(j)}, W_{t-1}^{(j)})$ to $\mathcal{Z}_t$ to form $\mathcal{Z}_{t-1}$.

Inductively, we eventually reach an integer $p \in \mathbb{Z}$ such that for all $t \leq p$, the maps
\[
\mathcal{R}_{t\leq p} : \mathcal{R}_{t} \xrightarrow{\cong} \mathcal{R}_{p}
\]
are isomorphisms. Thus, every irreducible subspace $W_t \subseteq V_t$ for $t \leq p$ lies in a persistence subrepresentation indexed by some pair in $\mathcal{Z}_p$. We obtain a nested family of sets
\[
\mathcal{Z}_{\infty} = \mathcal{Z}_v = \mathcal{Z}_q \subseteq \mathcal{Z}_s \subseteq \mathcal{Z}_r \subseteq \mathcal{Z}_p = \mathcal{Z}_u = \mathcal{Z}_{-\infty}
\]
for all $u \leq p \leq r \leq s \leq q \leq v$. For simplicity, denote the final collection by $\mathcal{Z} = \mathcal{Z}_p$.

We now show that each irreducible subspace $W_r \subseteq V_r$ lies in a unique irreducible persistence subrepresentation indexed by a pair in $\mathcal{Z}$.

Suppose $W_r \subseteq V_r$ lies in both $\mathcal{R}^{(U_1, U_1')}$ and $\mathcal{R}^{(U_2, U_2')}$ with $U_1 \subseteq V_{s_1}$, $U_1' \subseteq V_{t_1}$, $U_2 \subseteq V_{s_2}$, $ U_2' \subseteq V_{t_2}$, and $s_1, s_2 \leq r\leq t_1, t_2 $. Then we have
\[
W_r = f_{s_1,r}(U_1) = f_{s_2,r}(U_2), \quad f_{r,t_1}(W_r) = U_1', \quad f_{r,t_2}(W_r) = U_2'.
\]
If $t_1 \neq t_2$, we may assume without loss of generality that $t_1 < t_2$, then one has
\[
f_{t_1,t_2}(U_1') = f_{t_1,t_2} f_{r,t_1}(W_r) = f_{r,t_2}(W_r) = U_2'.
\]
This implies that the component $\mathcal{R}^{(U_1, U_1')}$ does not terminate at $t_1$ but continues to $t_2$, contradicting the construction of the intervals as maximal lifespans. If $t_1 = t_2$, then $U_1' = U_2' = f_{r,t_1}(W_r)$, which also contradicts the uniqueness of the decomposition at time $t_1$. Hence, each $W_r \subseteq V_r$ lies in a unique irreducible subrepresentation indexed by $\mathcal{Z}$.

Now, we prove the decomposition
\[
\mathcal{R} \cong \bigoplus_{(U, U') \in \mathcal{Z}} \mathcal{R}^{(U, U')}.
\]
For any fixed $r \in \mathbb{Z}$, consider the direct sum of the subrepresentations $\mathcal{R}^{(U, U')}_{r}$ over all $(U, U') \in \mathcal{Z}$. The underlying representation space of $\bigoplus_{(U, U') \in \mathcal{Z}} \mathcal{R}^{(U, U')}_{r}$ is given by
\[
\bigoplus_{\substack{
U \subseteq V_s,\ s \leq r \\
(U, U') \in \mathcal{Z}
}} f_{s,r}(U).
\]
Each summand $f_{s,r}(U)$ is an irreducible $G$-invariant subspace of $V_r$. Moreover, these subspaces are pairwise distinct. Indeed, by construction, every irreducible $G$-subrepresentation $W_r \subseteq V_r$ lies in a unique irreducible persistence subrepresentation indexed by a pair $(U, U') \in \mathcal{Z}$, which guarantees the mutual disjointness of the summands.

On the other hand, again by construction, every irreducible $G$-subrepresentation $W_r \subseteq V_r$ arises as $W_r = f_{s,r}(U)$ for some $s \leq r$ and some irreducible $G$-invariant subspace $U \subseteq V_s$, where $(U, U') \in \mathcal{Z}$. Therefore, the union of all such $f_{s,r}(U)$ exhausts the entire space $V_r$. That is,
\[
V_r = \bigoplus_{\substack{
U \subseteq V_s,\ s \leq r \\
(U, U') \in \mathcal{Z}
}} f_{s,r}(U).
\]
Moreover, for any $r_1 \leq r_2$, we have
\[
f_{r_1, r_2} = \bigoplus_{\substack{U \subseteq V_s,\ s \leq r_1 \\ (U, U') \in \mathcal{Z}}} f_{r_1, r_2}|_{f_{s, r_1}(U)}.
\]
This shows that $\mathcal{R}$ admits a decomposition as a direct sum of irreducible persistence subrepresentations indexed by the finite collection $\mathcal{Z}$. Hence,
\[
\mathcal{R} \cong \bigoplus_{(U, U') \in \mathcal{Z}} \mathcal{R}^{(U, U')}.
\]
This completes the proof.
\end{proof}

The above decomposition is interpreted, in the context of persistence representations, as a direct sum decomposition of $\mathcal{R}$ in the functor category $\mathrm{Fun}(\mathbb{Z}, \mathbf{Rep}_{\mathbb{K}}^{G})$, where each summand is an irreducible persistence representation that persists over a time interval. This decomposition gives rise to the barcode associated with $\mathcal{R}$, in which each pair $(U_i, U'_i)$ corresponds to a support interval $[s_i, t_i]$, representing the lifespan of an irreducible summand, or \textit{bar}. The finiteness condition ensures that there are only finitely many such bars. These bars collectively form the barcode, which fully encodes the essential features of the original persistence representation.

\begin{example}
Let us recall Example \ref{example:representation_barcode}. Consider again the persistence representation
\[
\mathcal{V} : \{0,1,2\} \to \mathbf{Rep}_{\mathbb{C}}^{G},
\]
where $G = \mathbb{Z}/2 = \{e,g\}$ and each $V_t$ decomposes as
\[
V_0 = T_0 \oplus S_0, \quad V_1 = T_1 \oplus S_1, \quad V_2 = T_2,
\]
with $T_t$ the trivial and $S_t$ the sign representations of $G$.

By the structure maps
\[
f_{0,1}(v_1) = w_1, \quad f_{0,1}(v_2) = w_2, \quad f_{1,2}(w_1) = u, \quad f_{1,2}(w_2) = 0,
\]
the trivial summand persists from time $0$ to time $2$, whereas the sign summand vanishes at time $2$.

According to Theorem~\ref{theorem:decomposition_representation}, the persistence representation $\mathcal{V}$ admits a decomposition into irreducible persistence subrepresentations supported on intervals corresponding to their lifespans
\[
\mathcal{V} \cong \mathcal{R}^{(T_0, T_2)} \oplus \mathcal{R}^{(S_0, S_1)}.
\]
Here, the persistence representation
\[
\mathcal{R}^{(T_0, T_2)} : \{0,1,2\} \to \mathbf{Rep}_{\mathbb{C}}^{G}
\]
is an irreducible persistence representation supported on the interval $[0,2]$, corresponding to the trivial subrepresentation, and the persistence representation
\[
\mathcal{R}^{(S_0, S_1)} : \{0,1\} \to \mathbf{Rep}_{\mathbb{C}}^{G}
\]
is an irreducible persistence representation supported on $[0,1]$, corresponding to the sign subrepresentation.

This decomposition matches precisely the barcode
\[
\{ [0,2], \; [0,1] \},
\]
where each bar represents the time interval over which an irreducible subrepresentation persists.

Thus, the example illustrates concretely how the general decomposition theorem produces a barcode decomposition of a persistence representation into irreducible persistence subrepresentations, each supported on a finite interval.
\end{example}

\subsection{The module structure of persistence representations}\label{section:module_representations}

In this section, we interpret persistence representations as graded modules over a polynomial ring generated by the shift operator induced by persistence morphisms. This algebraic structure encodes both the vector spaces and group actions, allowing the application of module decomposition theory. As a result, classical persistence module decompositions correspond to decompositions of persistence representations into irreducible components supported on intervals.

Let $\mathcal{R} \colon (\mathbb{Z}, \leq) \to \mathbf{Rep}_{\mathbb{K}}$ be a persistence representation. We define the associated graded representation
\begin{equation*}
  \mathbf{R} = \bigoplus_{k \in \mathbb{Z}} \mathcal{R}_{k},
\end{equation*}
which is a $\mathbb{Z}$-graded $\mathbb{K}$-vector space. Let $\tau \colon \mathbf{R} \to \mathbf{R}$ be a $\mathbb{K}$-linear map of degree $+1$ defined componentwise by the morphisms in the persistence representation. That is, for each $k \in \mathbb{Z}$, the restriction
\[
\tau|_{\mathcal{R}_{k}} = (\phi_{k,k+1}, f_{k,k+1}) \colon \mathcal{R}_{k} \to \mathcal{R}_{k+1},
\]
where
\[
(\phi_{s,t}, f_{s,t}) \colon (G_s, V_s, \rho_s) \to (G_t, V_t, \rho_t)
\]
is a morphism of representations in the category $\mathbf{Rep}_{\mathbb{K}}$ induced by the persistence structure.

We now consider the polynomial ring $\mathbb{K}[\tau]$, where $\tau$ is treated as an indeterminate of degree $+1$. The graded representation $\mathbf{R}$ naturally acquires the structure of a graded $\mathbb{K}[\tau]$-module via the action
\begin{equation*}
  \mathbb{K}[\tau] \times \mathbf{R} \to \mathbf{R}, \quad (p(\tau), v) \mapsto p(\tau)(v),
\end{equation*}
where $p(\tau)$ acts on $v \in \mathcal{R}_{k}$ by iterated composition of the morphisms $(\phi_{\ast,\ast+1}, f_{\ast,\ast+1})$ along the grading. This module structure encodes the entire persistence representation in algebraic form.

\begin{remark}\label{remark:persistence_module}
The construction of the $\mathbb{K}[\tau]$-module structure on the graded representation $\mathbf{R}$ is closely related to the standard construction of a persistence module. Consider the graded vector space
\begin{equation*}
  \mathbf{V} = \bigoplus_{k \in \mathbb{Z}} V_k,
\end{equation*}
where each $V_k$ is the representation space corresponding to $\mathcal{R}_{k}$. Define a linear map
\begin{equation*}
  \theta \colon \mathbf{V} \to \mathbf{V}
\end{equation*}
of degree $+1$, given by $\theta|_{V_k} = f_{k,k+1} \colon V_k \to V_{k+1}$. Then $\mathbf{V}$ becomes a $\mathbb{K}[\theta]$-module, where $\mathbb{K}[\theta]$ denotes the polynomial ring in the variable $\theta$.

Observe that $\mathbf{V}$ carries the usual structure of a persistence module. In comparison, the $\mathbb{K}[\tau]$-module structure on the graded representation $\mathbf{R}$ encodes additional information, as it also preserves the group actions. Specifically, for $v \in V_k$ and $g \in G_k$, we have
\begin{equation*}
  \tau (\rho(g)(v)) = \rho(\phi_{k,k+1}(g))(\tau(v)).
\end{equation*}
Under the setting of the persistence representation $\mathcal{R}$, even when we consider only the $\mathbb{K}[\theta]$-module structure on the graded space $\mathbf{V}$, it still retains the structural information coming from the group actions. Suppose that the underlying graded vector space $\mathbf{V}$ is a finitely generated $\mathbb{K}[\theta]$-module. Then $\mathbf{V}$ admits the following decomposition
\begin{equation*}
  \mathbf{V} \cong \left( \bigoplus_{i} \mathbb{K}[\theta] \cdot e_i \right) \oplus \left( \bigoplus_{j} \frac{\mathbb{K}[\theta]}{\theta^{t_j} \mathbb{K}[\theta]} \cdot \tilde{e}_j \right),
\end{equation*}
where each $e_i \in V_{s_i}$ and each $\tilde{e}_j \in V_{r_j}$ for some indices $s_i, r_j \in \mathbb{Z}$.
\end{remark}

\begin{remark}
When we fix the group action, that is, consider a persistence representation $\mathcal{R} \colon (\mathbb{Z}, \leq) \to \mathbf{Rep}_{\mathbb{K}}^{G}$, the persistence representation $\mathcal{R}$ can be viewed as a persistence $\mathbb{K}[G]$-module. More precisely, we define a functor
\[
\mathcal{M}_{\mathcal{R}} \colon (\mathbb{Z}, \leq) \to \mathbf{Mod}_{\mathbb{K}[G]}
\]
by setting $\mathcal{M}_{\mathcal{R}}(t) = V_t$ for each $t \in \mathbb{Z}$. For any $s \leq t$, the structure map $f_{s,t} \colon V_s \to V_t$ is a $\mathbb{K}[G]$-module homomorphism.

In this setting, the linear map
\[
\theta \colon \mathbf{V} \to \mathbf{V}
\]
can also be regarded as a $\mathbb{K}[G]$-module homomorphism of degree $+1$. Consider the polynomial ring $\mathbb{K}[G][\theta]$, which is generally non-commutative. Then the graded vector space $\mathbf{V}$ naturally acquires the structure of a graded $\mathbb{K}[G][\theta]$-module, with the action defined as
\[
\mathbb{K}[G][\theta] \times \mathbf{V} \to \mathbf{V},\quad (g\theta,v)\mapsto g\theta(v).
\]
\end{remark}

An \textit{interval module} $\mathcal{I}_{[b,d]} \colon (\mathbb{Z}, \leq) \to \mathbf{Vec}_{\mathbb{K}}$ supported on the interval $[b,d] \subseteq \mathbb{Z}$ is a persistence module defined by
\[
\mathcal{I}_{[b,d]}(k) =
\begin{cases}
\mathbb{K} & \text{if } b \leq k \leq d, \\
0 & \text{otherwise},
\end{cases}
\quad \text{for each } k \in \mathbb{Z},
\]
and for $k \leq l$,
\[
\mathcal{I}_{[b,d]}(k \leq l) =
\begin{cases}
\mathrm{id}_{\mathbb{K}} & \text{if } b \leq k \leq l \leq d, \\
0 & \text{otherwise}.
\end{cases}
\]

Let $\mathcal{R} \colon (\mathbb{Z}, \leq) \to \mathbf{Rep}_{\mathbb{K}}$ be a persistence representation. Consider the associated persistence module $\mathcal{V} \colon (\mathbb{Z}, \leq) \to \mathbf{Vec}_{\mathbb{K}}$ given by taking the underlying representation spaces of $\mathcal{R}$. Specifically, for each $s \in \mathbb{Z}$, we set $\mathcal{V}_{s} = V_s$, and for each $s \leq t$, the structure map $\mathcal{V}_{s\leq t} \colon \mathcal{V}_{s} \to \mathcal{V}_{t}$ is defined by $\mathcal{V}_{s\leq t} = f_{s,t}$, where $f_{s,t}$ is the structure map from the morphism $(\phi_{s,t}, f_{s,t})$ in $\mathcal{R}$.

\begin{theorem}[Decomposition theorem of persistence module]\label{theorem:decomposition_module}
Let $\mathcal{R} \colon (\mathbb{Z}, \leq) \to \mathbf{Rep}_{\mathbb{K}}$ be a finite persistence representation. Then there exists a finite multiset of intervals $\{[b_i, d_i]\}_{i \in I}$ such that
\[
\mathcal{V} \cong \bigoplus_{i \in I} \mathcal{I}_{[b_i, d_i]},
\]
where $\mathcal{I}_{[b_i, d_i]}$ denotes the interval module supported on $[b_i, d_i]$.
\end{theorem}
The above decomposition, viewed from the perspective of persistence modules, is essentially the same as the decomposition of $\mathbb{K}[\theta]$-module $\mathbf{V}$ presented in Remark~\ref{remark:persistence_module}.

In the above, we studied the persistence representations from the viewpoint of persistence modules. We now reverse the viewpoint and interpret persistence modules in terms of persistence representations.

Let $\mathcal{V} \colon (\mathbb{Z}, \leq) \to \mathbf{Vec}_{\mathbb{K}}$ be a persistence module. One can associate to it a persistence representation
\begin{equation*}
  \mathcal{R}_{\mathcal{V}} : (\mathbb{Z}, \leq) \to \mathbf{Rep}_{\mathbb{K}}^{G},
\end{equation*}
where each $\mathcal{R}_{\mathcal{V},t} = (G, \mathcal{V}_{t}, \rho)$, with $G = \{e\}$ the trivial group and $\rho \colon G \to \mathrm{GL}(\mathcal{V}_{t})$ the trivial representation. In this case, every $G$-invariant irreducible subspace of $\mathcal{V}_{t}$ is simply a one-dimensional subspace, and any decomposition of $\mathcal{V}_{t}$ into $G$-invariant subspaces corresponds to a choice of basis of $\mathcal{V}_{t}$.

Theorem~\ref{theorem:decomposition_representation} asserts that there exists a finite collection $\{(U_i, U_i')\}_{i \in I}$, where $U_i \subseteq V_{s_i}$ and $U_i' \subseteq V_{t_i}$ are irreducible $G$-invariant subspaces, such that
\[
\mathcal{R}_{\mathcal{V}} \cong \bigoplus_{i \in I} \mathcal{R}_{\mathcal{V}}^{(U_i, U_i')},
\]
where each $\mathcal{R}_{\mathcal{V}}^{(U_i, U_i')}$ is an irreducible persistence representation supported on the interval $[s_i, t_i]$. Note that each $U_i$ can be written as $\mathbb{K} \cdot \alpha_i$ for some generator $\alpha_i \in V_{s_i}$. The pair $(U_i, U_i')$ then corresponds to a generator $\alpha_i$ that persists over the interval $[s_i, t_i]$. Hence, the above decomposition of $\mathcal{R}_{\mathcal{V}}$ reduces to the decomposition of $\mathcal{V}$ described in Theorem~\ref{theorem:decomposition_module}.

\subsection{Regular representations of persistence groups}

The structure of a representation is determined by its irreducible constituents. For a finite group, the regular representation contains all irreducible representations as direct summands. In a persistent setting, understanding how these irreducible constituents evolve with respect to the filtration parameter provides valuable insight into the dynamics of the group system and its symmetries.

A \textit{persistence group} is a functor $\mathcal{G} : (T, \leq) \to \mathbf{Grp}$ from a poset $(T, \leq)$, viewed as a category, to the category of groups. A \textit{representation of the persistence group} $\mathcal{G}$ over a field $\mathbb{K}$ is a functor
\[
\mathcal{R} : (T, \leq) \to \mathbf{Rep}_{\mathbb{K}},
\]
such that for each $t \in T$, the object $\mathcal{R}_{t} = (G_t, V_t, \rho_t)$ is a representation of the group $G_t = \mathcal{G}_{t}$ on a $\mathbb{K}$-vector space $V_t$, and for each $s \leq t$, the morphism
\[
\mathcal{R}_{s \leq t} = (\phi_{s,t}, f_{s,t}) : (G_s, V_s, \rho_s) \to (G_t, V_t, \rho_t)
\]
is a morphism of group representations, where $\phi_{s,t} = \mathcal{G}_{s \leq t}$ is the group homomorphism assigned by $\mathcal{G}$. In other words, the representation $\mathcal{R}$ is compatible with the group structure provided by $\mathcal{G}$.

Let $\mathcal{G} : (T, \leq) \to \mathbf{Grp}$ be a persistence group. For each $t \in T$, define $G_t = \mathcal{G}_{t}$ and consider the regular representation of $G_t$. Let $V_t = \mathbb{K}[G_t]$ be the group algebra, viewed as a $\mathbb{K}$-vector space with basis corresponding to the elements of $G_t$. Let $\rho_t : G_t \to \mathrm{GL}(V_t)$ be the regular representation. This gives a representation $(G_t, V_t, \rho_t)$ at each time point $t \in T$.
For each pair $s \leq t$, the persistence group $\mathcal{G}$ induces a group homomorphism $\phi_{s,t} = \mathcal{G}_{s \leq t} : G_s \to G_t$. We define a $\mathbb{K}$-linear map $f_{s,t} : V_s \to V_t$ as follows
\[
f_{s,t}(g) = \phi_{s,t}(g) \quad \text{for } g \in G_s.
\]
It is straightforward to verify that this pair $(\phi_{s,t}, f_{s,t})$ is a morphism of representations, meaning that for all $h \in G_s$, we have
\[
f_{s,t} \circ \rho_s(h) = \rho_t(\phi_{s,t}(h)) \circ f_{s,t}.
\]

Thus, we obtain a functor $\mathcal{R}_{\mathcal{G}} : (T, \leq) \to \mathbf{Rep}_{\mathbb{K}}$ defined as follows. For each $t \in T$, the object $\mathcal{R}_{\mathcal{G},t}$ is the regular representation $(G_t, V_t, \rho_t)$ of the group $G_t$. For each morphism $s \leq t$ in $(T, \leq)$, the corresponding morphism $\mathcal{R}_{\mathcal{G},s\leq t}$ is given by the pair $(\phi_{s,t}, f_{s,t})$ constructed above. The functor $\mathcal{R}_{\mathcal{G}}$ is referred to as the \textit{regular persistence representation} of the persistence group $\mathcal{G}$.

A \textit{filtered group} $\mathcal{G} : (\mathbb{Z}, \leq) \to \mathbf{Grp}$ is a persistence group such that $\mathcal{G}_s\leq \mathcal{G}_t$ for any $s\leq t$.  We say that a filtered group $\mathcal{G}$ is left constant if there exists an integer $r \in \mathbb{Z}$ such that for all $r_1 \leq r_2 \leq r$, the structure map
\[
\mathcal{G}_{r_1 \leq r_2} : \mathcal{G}_{r_1} \xrightarrow{\cong} \mathcal{G}_{r_2}
\]
is an isomorphism.

We now recall the construction of induced representations \cite{fulton2013representation,serre1977linear}. Given a subgroup $H \leq G$ and a representation $(\rho, V)$ of $H$, the induced representation $\mathrm{Ind}_H^G \rho$ is a representation of $G$ defined on the vector space
\[
\mathrm{Ind}_H^G V = \left\{ f : G \to V \;\middle|\; f(hg) = \rho(h)f(g),\ \forall h \in H,\ g \in G \right\}.
\]
The group $G$ acts on $\mathrm{Ind}_H^G V$ by right translation: for $g \in G$ and $f \in \mathrm{Ind}_H^G V$, define
\[
(\tilde{\rho}(g)f)(x) = f(xg), \quad \forall x \in G.
\]
This action gives a well-defined representation $(G,\mathrm{Ind}_H^G V, \tilde{\rho})$. In particular, if $\rho : H \to \mathrm{GL}(\mathbb{C}[H])$ is the regular representation of $H$, then its induced representation is the regular representation of $G$, i.e., $\tilde{\rho} : G \to \mathrm{GL}(\mathbb{C}[G])$.

Let $\mathcal{G} : (\mathbb{Z}, \leq) \to \mathbf{Grp}$ be a filtered group that is constant below some fixed parameter $r \in \mathbb{Z}$. This induces a persistence representation $\mathcal{R}_{\mathcal{G}} : (\mathbb{Z}, \leq) \to \mathbf{Rep}_{\mathbb{K}}$. For $r \leq s$, and for any subrepresentation $W_r$ of the regular representation $(\mathcal{G}_{r}, \mathbb{K}[\mathcal{G}_{r}], \rho_r)$ at parameter $r$, we define the induced representation at parameter $s$ as
\begin{equation*}
   \left( \mathcal{G}_{s},\ \mathrm{Ind}_{\mathcal{G}_{r}}^{\mathcal{G}_{s}} W_r,\ \tilde{\rho}_s(W_r) \right),
\end{equation*}
where the underlying vector space is given by
\[
\mathrm{Ind}_{\mathcal{G}_{r}}^{\mathcal{G}_{s}} W_r := \mathbb{K}[\mathcal{G}_{s}] \otimes_{\mathbb{K}[\mathcal{G}_{r}]} W_r.
\]
Moreover, for $r \leq s \leq t$, the induction functor is transitive in the sense that
\begin{equation*}
   \mathrm{Ind}_{\mathcal{G}_{r}}^{\mathcal{G}_{t}} W_r \cong
   \mathrm{Ind}_{\mathcal{G}_{s}}^{\mathcal{G}_{t}} \left( \mathrm{Ind}_{\mathcal{G}_{r}}^{\mathcal{G}_{s}} W_r \right).
\end{equation*}
Thus, one can obtain a persistence representation
\begin{equation*}
  \mathcal{R}^{W_r}:(\mathbb{Z}, \leq) \to \mathbf{Rep}_{\mathbb{K}}
\end{equation*}
given by
\begin{equation*}
  \mathcal{R}^{W_r}_{t} =\left\{
                                      \begin{array}{ll}
                                        (\mathcal{G}_{r},\ W_r,\ \rho_{r}|_{W_r}),\ & \text{if } t\leq r; \\
                                        (\mathcal{G}_{t},\ \mathrm{Ind}_{\mathcal{G}_{r}}^{\mathcal{G}_{t}} W_r,\ \tilde{\rho}_t(W_r)) & \text{if } t>r.
                                      \end{array}
                                    \right.
\end{equation*}

The structure $\mathcal{R}^{W_r}$ encodes the evolution of $\mathcal{G}_{t}$-invariant subspaces as the parameter $t$ varies. The induced representation $\mathrm{Ind}_{\mathcal{G}_{t}}^{\mathcal{G}_{t+1}} U$ is always nontrivial for any nontrivial representation $U$ of $\mathcal{G}_{t}$, but it is not necessarily irreducible. In general, the induction of an irreducible representation decomposes into a direct sum of several irreducible components. This implies that the $\mathcal{G}_{t}$-invariant subspaces of $\mathcal{R}^{W_r}$ undergo a branching process over time, where irreducible components at a given stage $t$ may split into multiple irreducible components at stage $t+1$. This branching behavior can be visualized as a graded digraph. In this graph, the grading corresponds to the filtration parameter $t$, each vertex at grading $t$ represents an irreducible representation of the group $\mathcal{G}_{t}$, and a directed edge $(U, V)$ indicates that $V$ is an irreducible component of the induced representation $\mathrm{Ind}_{\mathcal{G}_{t}}^{\mathcal{G}_{t+1}} U$ for some irreducible $U$ of $\mathcal{G}_{t}$. It is worth noting that two distinct irreducible representations may have isomorphic induced representations, or their induced representations may share one or more irreducible components. This non-injectivity of the induction process reflects the complexity of how symmetry and invariant structure evolve under group actions in persistent settings.

\begin{example}
Let $\mathcal{G} : (\mathbb{Z}, \leq) \to \mathbf{Grp}$ be a persistence group defined by
\[
\mathcal{G}_{0} = \{e\}, \quad
\mathcal{G}_{1} = \langle (12) \rangle \cong C_2, \quad
\mathcal{G}_{2} = \langle (12), (123) \rangle \cong S_3, \quad
\mathcal{G}_{3} = S_4.
\]
All representations are taken over the field $\mathbb{C}$. The group $\mathcal{G}_{0}$ admits a unique irreducible representation $V^{(1)}$, which is trivial and $1$-dimensional. Inducing $V^{(1)}$ to $\mathcal{G}_{1}$ yields
\[
U = \mathrm{Ind}_{\mathcal{G}_{0}}^{\mathcal{G}_{1}}(V^{(1)}) \cong \mathbb{C}[\mathcal{G}_{1}] \cong V^{(2)} \oplus V^{(1^2)},
\]
where $V^{(2)}$ and $V^{(1^2)}$ denote the trivial and sign representations of $C_2$ respectively.

Each irreducible summand of $U$ is then induced to $\mathcal{G}_{2} \cong S_3$. One obtains
\[
  \mathrm{Ind}_{\mathcal{G}_{1}}^{\mathcal{G}_{2}}(V^{(2)}) \cong V^{(3)} \oplus V^{(2,1)}, \quad
  \mathrm{Ind}_{\mathcal{G}_{1}}^{\mathcal{G}_{2}}(V^{(1^2)}) \cong V^{(1^3)} \oplus V^{(2,1)},
\]
where the irreducible representations of $S_3$ are indexed by partitions of $3$.

Each of the summands $V^{(3)}$, $V^{(1^3)}$, and $V^{(2,1)}$ is now induced to $\mathcal{G}_{3} = S_4$, yielding
\[
\mathrm{Ind}_{\mathcal{G}_{2}}^{\mathcal{G}_{3}}(V^{(3)}) \cong V^{(4)} \oplus V^{(3,1)}, \quad
\mathrm{Ind}_{\mathcal{G}_{2}}^{\mathcal{G}_{3}}(V^{(1^3)}) \cong V^{(1^4)} \oplus V^{(2,1,1)},
\]
\[
\mathrm{Ind}_{\mathcal{G}_{2}}^{\mathcal{G}_{3}}(V^{(2,1)}) \cong V^{(3,1)} \oplus V^{(2,2)} \oplus V^{(2,1,1)}.
\]

The following figure illustrates the branching behavior of irreducible representations under induction as the filtration parameter increases. The evolution of irreducible components can be naturally encoded by a graded directed graph, where each vertex corresponds to an irreducible representation at a given filtration level, and edges represent the inclusion of an irreducible component in an induced representation at the next level.
\begin{figure}[h]
\centering
\begin{tikzpicture}[
    >=Stealth, 
    font=\small, 
    every node/.style={
        circle, 
        draw=mainblue, 
        line width=1.2pt, 
        fill=mainblue!10, 
        minimum size=1.2cm, 
        inner sep=0pt
    },
    shorten >= 3pt, 
    shorten <= 3pt
]

\begin{scope}[scale=0.8, transform shape]
\node (L0) at (0, 4.8) {$V^{(1)}$};

\node (L1L) at (-2.2, 3.2) {$V^{(2)}$};
\node (L1R) at (2.2, 3.2) {$V^{(1^2)}$};

\node (L2L) at (-4.4, 1.6) {$V^{(3)}$};
\node (L2M) at (0, 1.6) {$V^{(2,1)}$}; 
\node (L2R) at (4.4, 1.6) {$V^{(1^3)}$};

\node (L3L) at (-6.6, 0) {$V^{(4)}$};
\node (L3ML) at (-2.2, 0) {$V^{(3,1)}$}; 
\node (L3MM) at (0, 0) {$V^{(2,2)}$};
\node (L3MR) at (2.2, 0) {$V^{(2,1,1)}$}; 
\node (L3R) at (6.6, 0) {$V^{(1^4)}$};

\draw[->, line width=1.2pt, mainblue] (L0) -- (L1L);
\draw[->, line width=1.2pt, mainblue] (L0) -- (L1R);

\draw[->, line width=1.2pt, mainblue] (L1L) -- (L2L);
\draw[->, line width=1.2pt, mainblue] (L1L) -- (L2M);
\draw[->, line width=1.2pt, mainblue] (L1R) -- (L2M);
\draw[->, line width=1.2pt, mainblue] (L1R) -- (L2R);

\draw[->, line width=1.2pt, mainblue] (L2L) -- (L3L);
\draw[->, line width=1.2pt, mainblue] (L2L) -- (L3ML);
\draw[->, line width=1.2pt, mainblue] (L2M) -- (L3ML);
\draw[->, line width=1.2pt, mainblue] (L2M) -- (L3MM);
\draw[->, line width=1.2pt, mainblue] (L2M) -- (L3MR);
\draw[->, line width=1.2pt, mainblue] (L2R) -- (L3MR);
\draw[->, line width=1.2pt, mainblue] (L2R) -- (L3R);
\end{scope}
\end{tikzpicture}
\caption{Branching behavior of irreducible representations under induction.}\label{figure:induction_branching}
\end{figure}

It is worth noting that this directed graph is not a tree; multiple irreducible representations at one level may induce to the same irreducible representation at the next level, and different paths may merge. Nonetheless, the overall branching structure of irreducible components remains clearly discernible and faithfully captures the decomposition patterns induced by the filtration.
\end{example} 

%% file: persistent_Fourier_analysis.tex
\section{Fourier analysis on persistence groups}\label{section:Fourier_analysis}

In many settings involving the evolution or parametrization of symmetric structures in data, a fundamental inquiry concerns the evolution of group structures over time or with respect to a varying parameter. The concept of the persistent Fourier transform offers a robust methodology for investigating the temporal evolution of frequency components within functions defined on persistence groups.

Compared to standard Fourier analysis, the persistent version captures long-term patterns in the spectral energy and changes in dominant frequencies. These spectral patterns may reflect important processes such as symmetry breaking, transition, or stabilization. Quantities like spectral energy curves, persistent entropy, and dominant frequency ratios provide stable and interpretable features for tasks such as classification, clustering, and model analysis. Such spectral methods are expected to provide useful insights into structural dynamics in settings with evolving symmetries, including dynamic point clouds, time-dependent or parametrization-dependent graphs, and algebraic signals.

\subsection{Persistent Fourier transform}

For the Fourier transform on finite groups, one can refer to \cite{folland2016course,terras1999fourier} for foundational and comprehensive treatments. Let $G$ be a finite group, and let $\mathbb{C}[G]$ denote its group algebra. Denote by $\widehat{G}$ the set of equivalence classes of irreducible complex representations of $G$. For each $\rho \in \widehat{G}$, let $d_\rho$ be the dimension of the representation $\rho$. Given a function $\theta : G \to \mathbb{C}$, its Fourier transform at $\rho$ is defined by
\begin{equation*}
\widehat{\theta}(\rho) = \sum_{g \in G} \theta(g)\, \rho(g^{-1}) \in \mathrm{End}(V_\rho) \cong \mathbb{C}^{d_\rho \times d_\rho}.
\end{equation*}
This operator captures the frequency content of $\theta$ along the irreducible representation $\rho$.

\begin{definition}
The \textit{spectral energy} of $\theta$ at $\rho \in \widehat{G}$ is defined as
\begin{equation*}
E_\theta(\rho) := \frac{1}{|G|} \| \widehat{\theta}(\rho) \|_{\mathrm{HS}}^2,
\end{equation*}
where $\| A \|_{\mathrm{HS}}^2 := \tr(A^\dagger A)$ denotes the Hilbert-Schmidt norm.
\end{definition}

The total energy of the function $\theta$ is given by
\[
E_\theta := \sum_{g \in G} |\theta(g)|^2.
\]

\begin{theorem}[Parseval identity]
For any function $\theta : G \to \mathbb{C}$, the total energy is equal to the sum of the spectral energies weighted by representation dimensions
\begin{equation*}
\sum_{g \in G} |\theta(g)|^2 = \frac{1}{|G|} \sum_{\rho \in \widehat{G}} d_\rho \cdot \| \widehat{\theta}(\rho) \|_{\mathrm{HS}}^2 = \sum_{\rho \in \widehat{G}} d_\rho \cdot E_\theta(\rho).
\end{equation*}
\end{theorem}

\begin{theorem}[Fourier inversion]
Any function $\theta : G \to \mathbb{C}$ can be recovered from its Fourier transform via
\begin{equation*}
\theta(g) = \frac{1}{|G|} \sum_{\rho \in \widehat{G}} d_\rho \cdot \tr \left( \widehat{\theta}(\rho)\, \rho(g) \right).
\end{equation*}
\end{theorem}

\begin{example}
Let $\theta: G \to \mathbb{C}$ be a function on a finite group $G$. If $\theta$ satisfies
\[
\theta(hgh^{-1}) = \theta(g) \quad \text{for all } g, h \in G,
\]
then $\theta$ is called a \textit{class function}. Class functions are constant on conjugacy classes and form a subspace of $\mathbb{C}[G]$ invariant under conjugation.

A fundamental property of class functions is that their Fourier transforms are scalar multiples of identity matrices. That is, for each irreducible representation $\rho \in \widehat{G}$, there exists a scalar $\lambda_\rho \in \mathbb{C}$ such that
\[
\widehat{\theta}(\rho) = \lambda_\rho \cdot I_{d_\rho}.
\]
This implies that the spectral content of a class function is aligned uniformly along each isotypic component of the representation space.

The spectral energy of $\theta$ at $\rho$ then simplifies to
\[
E_\theta(\rho) = \frac{1}{|G|} \| \widehat{\theta}(\rho) \|_{\mathrm{HS}}^2
= \frac{1}{|G|} \cdot \| \lambda_\rho I_{d_\rho} \|_{\mathrm{HS}}^2
= \frac{d_\rho}{|G|} \cdot |\lambda_\rho|^2.
\]
Hence, the contribution of each irreducible representation to the total energy is determined by the scalar Fourier coefficient $\lambda_\rho$ and the dimension $d_\rho$.

A canonical example of a class function is the character $\chi_\rho$ associated with an irreducible representation $\rho$, defined as
\[
\chi_\rho(g) := \tr(\rho(g)).
\]
The spectral energy of the character function is given by
\begin{equation*}
E_{\chi_\rho}(\sigma) := \frac{1}{|G|} \left\| \widehat{\chi}_\rho(\sigma) \right\|_{\mathrm{HS}}^2 =
\begin{cases}
\frac{|G|}{d_\rho}, & \text{if } \rho \cong \sigma, \\
0, & \text{otherwise}.
\end{cases}
\end{equation*}
\end{example}

Let $\mathcal{G} : (\mathbb{Z}, \leq) \to \mathbf{Grp}$ be a persistence group. Let $\mathrm{Hom}(\mathcal{G}_t,\mathbb{C})$ denote the linear space of complex-valued functions on $\mathcal{G}_t$. Consider the inverse system of sets
\[
\{ \mathrm{Hom}(\mathcal{G}_t,\mathbb{C}), \phi_{s,t}^* \}_{t \in \mathbb{Z}},
\]
where $\phi_{s,t}:\mathcal{G}_s\to \mathcal{G}_t$ are group homomorphisms, and the transition maps
\[
\phi_{s,t}^* : \mathrm{Hom}(\mathcal{G}_t,\mathbb{C}) \to \mathrm{Hom}(\mathcal{G}_s,\mathbb{C}), \quad h \mapsto h \circ \phi_{s,t}
\]
are given by precomposition.

\begin{definition}
A \textit{persistent function on $\mathcal{G}$} is an element in the inverse limit
\[
\theta=\{ \theta_t \}_{t \in \mathbb{Z}} \in \varprojlim_{t \in \mathbb{Z}} \mathrm{Hom}(\mathcal{G}_t,\mathbb{C}).
\]
\end{definition}
More precisely, a persistent function is a family of functions $\{ \theta_t\}_{t \in \mathbb{Z}}$ with $\theta_t\in \mathrm{Hom}(\mathcal{G}_t,\mathbb{C})$ such that for all $s \leq t$, we have $\phi_{s,t}^*(\theta_t) = \theta_s$, or equivalently, $\theta_s = \theta_t \circ \phi_{s,t}$.
\begin{equation*}
  \xymatrix{
 \mathcal{G}_{s}\ar[dr]_{\theta_s} \ar[rr]^{\phi_{s,t}}&&\mathcal{G}_{t}\ar[ld]^{\theta_t} \\
 &\mathbb{C}&
 }
\end{equation*}
\begin{definition}
Let $\mathcal{G} : (\mathbb{Z}, \leq) \to \mathbf{Grp}$ be a persistence group, and let $\{ \theta_t \}_{t \in \mathbb{Z}}$ be a persistent function on $\mathcal{G}$. The \textit{$(s,t)$-persistent Fourier transform} of $\{\theta_t\}_{t \in \mathbb{Z}}$ is given by
\[
\widehat{\theta}_{s,t}(\varrho_{s}) := \sum_{g_s \in \mathcal{G}_{s}} \theta_s(g_s) \cdot \varrho_{s,t}(g_s^{-1}),
\]
where $\varrho_{s,t}$ is the persistent representation of the irreducible representation $\varrho_{s} \in \widehat{\mathcal{G}}_{s}$.
\end{definition}

It is straightforward to verify that when $s = t$, the $(s,t)$-persistent Fourier transform coincides with the classical Fourier transform on the group $\mathcal{G}_s$. Moreover, if the persistent representation $\varrho_{s,t}$ degenerates to the zero representation at parameter $t$, then its representation space is the zero space, i.e., $\dim W_{s,t} = 0$. In this case, the persistent Fourier transform $\widehat{\theta}_{s,t}(\varrho_s)$ is the zero endomorphism in $\mathrm{End}(\{0\})$.

\begin{remark}\label{remark:persistent_function}
We may also define a persistent function indexed in reverse filtration parameter, a collection $\{\theta_t \}_{t \leq r}$ on the persistence group $\mathcal{G} : (\mathbb{Z}, \leq) \to \mathbf{Grp}$. Such a family can naturally be interpreted as an element of the inverse limit $\varprojlim_{t \leq r} \mathrm{Hom}(\mathcal{G}_t, \mathbb{C})$. Moreover, any function $\theta_r : \mathcal{G}_r \to \mathbb{C}$ canonically extends to a persistent function $\theta = \{ \theta_t \}_{t \leq r}$ via pullback along the structure maps
\[
\theta_t = \theta_r \circ \phi_{t,r}, \quad \forall t \leq r,
\]
where $\phi_{t,r} : \mathcal{G}_t \to \mathcal{G}_r$ is the morphism in the persistence group system. In addition, a compatible system $\{ \theta_t \}_{t \geq r}$ on $\mathcal{G}$ can always extend to a persistent function $\{ \theta_t \}_{t \in \mathbb{Z}}$ on $\mathcal{G}$.
\end{remark}

Recall that the regular persistence representation of a persistence group $\mathcal{G} : (\mathbb{Z}, \leq) \to \mathbf{Grp}$ has representation spaces given by the group algebras $\mathbb{C}[\mathcal{G}_s]$, and the structure maps $f_{s,t} : \mathbb{C}[\mathcal{G}_s] \to \mathbb{C}[\mathcal{G}_t]$ are the linear extensions of the group homomorphisms $\phi_{s,t} : \mathcal{G}_s \to \mathcal{G}_t$. Let $\widehat{\theta}_s(\varrho_s)$ denote the Fourier transform of a function $\theta_s \in \mathrm{Hom}(\mathcal{G}_s, \mathbb{C})$ at the irreducible representation $\varrho_s \in \widehat{\mathcal{G}}_s$. Then the $(s,t)$-persistent Fourier transform satisfies the intertwining relation
\[
f_{s,t} \circ \widehat{\theta}_s(\varrho_s) = \widehat{\theta}_{s,t}(\varrho_s) \circ f_{s,t},
\]
which gives rise to the following commutative diagram.
\[
\xymatrix{
W_s \ar[rr]^{\widehat{\theta}_s(\varrho_s)} \ar[d]_{f_{s,t}|_{W_s}} && W_s \ar[d]^{f_{s,t}|_{W_s}} \\
W_{s,t} \ar[rr]^{\widehat{\theta}_{s,t}(\varrho_s)} && W_{s,t}
}
\]
Here, $W_s$ denotes the representation space of $\varrho_s$, and $W_{s,t} = f_{s,t}(W_s)$ is the corresponding $(s,t)$-persistent representation space.

\begin{example}\label{example:persistent_fourier}
Let $\mathcal{G}_t = \mathbb{Z}/2^t\mathbb{Z}$ be the cyclic group of order $2^t$, and let $\mathcal{G}: (\mathbb{Z}_{\geq 0}, \leq) \to \mathbf{Grp}$ be the persistence group assigning $t \mapsto \mathcal{G}_t$, with structure maps $\phi_{s,t}: \mathbb{Z}/2^s\mathbb{Z} \rightarrow \mathbb{Z}/2^t\mathbb{Z}$ given by $\phi_{s,t}([x]_s) = [2^{t-s}x]_t$ for $s \leq t$. Define a family of functions $\{ \theta_t \}_{t \in \mathbb{Z}_{\geq 0}}$ with
\[
\theta_t([x]_t) = \cos\left( \frac{2\pi x}{2^t} \right).
\]
Then for each $s \leq t$, we have
\[
\theta_s([x]_s) = \cos\left( \frac{2\pi x}{2^s} \right) = \cos\left( \frac{2\pi \cdot 2^{t-s}x}{2^t} \right) = \theta_t\left( \phi_{s,t}([x]_s) \right).
\]
Hence, $\theta_s = \theta_t \circ \phi_{s,t}$, and the family $\{ \theta_t \}$ defines a persistent function on $\mathcal{G}$.

The irreducible representations of $\mathcal{G}_t$ are the 1-dimensional characters $\chi_k^{(t)} : \mathcal{G}_t \to \mathbb{C}^\times$ defined by
\[
\chi_k^{(t)}([x]_t) = e^{2\pi i k x / 2^t}, \quad k = 0, 1, \dots, 2^t - 1.
\]
Let $\rho_t = \chi_k^{(t)}$ be a representation of $\mathcal{G}_t$. Then the \textit{persistent representation} $\rho_{s,t} : \mathcal{G}_s \to \mathbb{C}$ is defined as
\[
\rho_{s,t}(g) := \rho_t(\phi_{s,t}(g)) = \chi_k^{(t)}([2^{t-s}x]_t) = e^{2\pi i k \cdot 2^{t-s} x / 2^t} = e^{2\pi i k x / 2^s}.
\]
Therefore, $\rho_{s,t} = \chi_k^{(s)}$, i.e., the persistent representation simply pulls back $\chi_k^{(t)}$ along $\phi_{s,t}$ to get the corresponding character of $\mathcal{G}_s$.

Now the $(s,t)$-persistent Fourier transform is given by
\[
\widehat{\theta}_{s,t}(k) = \sum_{x = 0}^{2^s - 1} \theta_s([x]_s) \cdot \overline{\rho_{s,t}([x]_s)} = \sum_{x=0}^{2^s - 1} \cos\left( \frac{2\pi x}{2^s} \right) \cdot e^{-2\pi i k x / 2^s}.
\]
Expanding the cosine term, we have
\begin{align*}
   \widehat{\theta}_{s,t}(k) = &  \frac{1}{2} \sum_{x=0}^{2^s - 1} \left( e^{2\pi i x / 2^s} + e^{-2\pi i x / 2^s} \right) \cdot e^{-2\pi i k x / 2^s}, \\
    = & \frac{1}{2} \left( \sum_{x=0}^{2^s - 1} e^{-2\pi i (k-1)x / 2^s} + \sum_{x=0}^{2^s - 1} e^{-2\pi i (k+1)x / 2^s} \right).
\end{align*}
Note that
\[
\sum_{x=0}^{2^s - 1} e^{-2\pi i m x / 2^s} =
\begin{cases}
2^s, & \text{if } m \equiv 0 \mod 2^s, \\
0, & \text{otherwise}.
\end{cases}
\]
Thus, we have
\[
\widehat{\theta}_{s,t}(k) =
\begin{cases}
2^{s-1}, & \text{if } k \equiv \pm 1 \mod 2^s, \\
0, & \text{otherwise}.
\end{cases}
\]
This shows that $\theta_s$ has persistent frequency support at indices $k \equiv \pm 1 \mod 2^s$ under the lifted representation from $\mathcal{G}_s$ to $\mathcal{G}_t$.
\end{example}

\subsection{Persistent Fourier inversion formula}

Let $H \trianglelefteq G$ be a normal subgroup of a finite group $G$, and let $\pi: G \to G/H$ denote the canonical projection onto the quotient group. Given a representation $\bar{\rho}: G/H \to \mathrm{GL}(V)$ of the quotient group $G/H$, we obtain a representation of $G$ by pullback via the composition
\[
\rho = \bar{\rho} \circ \pi: G \to \mathrm{GL}(V),
\]
explicitly defined by $\rho(g)(v) = \bar{\rho}(gH)(v)$ for all $g \in G$ and $v \in V$.

Conversely, suppose that $\rho: G \to \mathrm{GL}(V)$ is a representation of $G$ satisfying $\rho(h) = \mathrm{id}_V$ for all $h \in H$. Then there exists a unique representation
\[
\bar{\rho}: G/H \to \mathrm{GL}(V)
\]
such that the following diagram commutes
\[
\xymatrix@C=1.2cm@R=1.2cm{
G/H \ar@{-->}[dr]^-{\bar{\rho}} &\\
G \ar[r]^-{\rho} \ar[u]^-{\pi} & \mathrm{GL}(V) .
}
\]
That is, $\rho = \bar{\rho} \circ \pi$.

This construction gives a bijection between the set of representations of the quotient group $G/H$ and the set of representations of $G$ that are trivial on $H$. In particular, there is a bijection between the sets of irreducible representations
\[
\widehat{G/H} \cong \widehat{G}_H,
\]
where $\widehat{G}_H = \left\{ \rho \in \widehat{G} \mid \rho(h) = \mathrm{id}_V \text{ for all } h \in H \right\}$ denotes the set of irreducible representations of $G$ that are trivial on $H$.

\begin{proposition}\label{proposition:fourier_note}
Let $(\phi_{s,t}, f_{s,t}) : (\mathcal{G}_s, V_s, \rho_s) \to (\mathcal{G}_t, V_t, \rho_t)$ be a morphism of regular representations, which means that for all $g \in \mathcal{G}_s$, the following diagram commutes.
\[
\xymatrix{
V_s \ar[d]_{f_{s,t}} \ar[rr]^{\rho_s(g)} && V_s \ar[d]^{f_{s,t}} \\
V_t \ar[rr]^{\rho_t(\phi_{s,t}(g))} && V_t
}
\]
Let $\varrho_s: \mathcal{G}_s \to \mathrm{GL}(W_s)$ be an irreducible subrepresentation of $\rho_s$ that is trivial on $\ker \phi_{s,t}$. Then the image
\[
W_{s,t} = f_{s,t}(W_s) \subseteq V_t
\]
is a nonzero $\mathcal{G}_s$-invariant subspace.
\end{proposition}

\begin{proof}
Let $K = \ker \phi_{s,t}$. Since $\varrho_s$ is a subrepresentation of the regular representation $\rho_s$, we can view $W_s$ as a subspace of the group algebra $\mathbb{C}[\mathcal{G}_s]$ spanned by basis vectors $\{g\}_{g \in \mathcal{G}_s}$.
Let $w \in W_s$ be any nonzero vector, expressed as $w = \sum_{g \in \mathcal{G}_s} c_g g$ with $c_g \in \mathbb{C}$.
Because $\varrho_s$ is trivial on $K$, we have $\rho_s(k)w = w$ for all $k \in K$. The action of $\rho_s(k)$ on the basis is permutation, so $c_{k^{-1}g} = c_g$ for all $k \in K$ and $g \in \mathcal{G}_s$. This implies that the coefficients $c_g$ are constant on the cosets of $K$ in $\mathcal{G}_s$.

Now, consider the image $f_{s,t}(w) = \sum_{g \in \mathcal{G}_s} c_g \phi_{s,t}(g) \in \mathbb{C}[\mathcal{G}_t]$.
Suppose, for contradiction, that $f_{s,t}(w) = 0$. Then for each $y \in \phi_{s,t}(\mathcal{G}_s) \subseteq \mathcal{G}_t$, the sum of coefficients of the basis elements corresponding to $y$ must be zero.
For a fixed coset $\bar{g} \in \mathcal{G}_s/K$, let $g$ be a representative. The elements in this coset map to $\phi_{s,t}(g)$ in $\mathcal{G}_t$.
The coefficient sum for this image is $\sum_{k \in K} c_{gk}$.
Since $c_{gk} = c_g$ for all $k \in K$, this sum equals $\sum_{k \in K} c_g = |K| \cdot c_g$.
Thus, $|K| \cdot c_g = 0$, which implies $c_g = 0$ (as we are working over $\mathbb{C}$).
This holds for every coset, so $c_g = 0$ for all $g \in \mathcal{G}_s$, which means $w = 0$. This contradicts the assumption that $w$ is nonzero.
Therefore, $f_{s,t}(w) \neq 0$ for any nonzero $w \in W_s$, implying that $W_{s,t} = f_{s,t}(W_s)$ is nonzero.
\end{proof}

\begin{theorem}[Persistent Fourier inversion formula]
Let $\{\theta_t\}_{t \in \mathbb{Z}}$ be a persistent function on a persistence group $\mathcal{G} : (\mathbb{Z}, \leq) \to \mathbf{Grp}$. For any $s \leq t$, let $K = \ker(\phi_{s,t})$ be the kernel of the group homomorphism $\phi_{s,t} : \mathcal{G}_s \to \mathcal{G}_t$. Then for any $g_s \in \mathcal{G}_s$, the value $\theta_s(g_s)$ can be recovered via the formula
\[
\theta_s(g_s) = \frac{1}{|\mathcal{G}_s|} \sum_{\varrho_s \in \widehat{(\mathcal{G}_s)}_K} d_{\varrho_s} \cdot \tr\left( \widehat{\theta}_{s,t}(\varrho_s) \cdot \varrho_{s,t}(g_s) \right),
\]
where $\widehat{(\mathcal{G}_s)}_K$ denotes the set of inequivalent irreducible representations of $\mathcal{G}_s$ that are trivial on $K$.
\end{theorem}

\begin{proof}
For $s \leq t$, the homomorphism $\phi_{s,t} : \mathcal{G}_s \to \mathcal{G}_t$ induces a quotient map
\[
\bar{\phi}_{s,t} : \mathcal{G}_s / K \to \mathcal{G}_t,
\]
where $K = \ker(\phi_{s,t})$. Let us denote $\mathcal{G}_s^\flat := \mathcal{G}_s / K$. The persistence relation $\theta_s = \theta_t \circ \phi_{s,t}$ naturally induces a well-defined function
\[
\theta_s^\flat : \mathcal{G}_s^\flat \to \mathbb{C}, \quad \theta_s^\flat(\bar{g}) = \theta_s(g),
\]
where $g$ is any representative of the coset $\bar{g} \in \mathcal{G}_s^\flat$.

Applying the classical Fourier inversion formula to the function $\theta_s^\flat$ on the finite group $\mathcal{G}_s^\flat$, we obtain
\[
\theta_s^\flat(\bar{g}_s) = \frac{1}{|\mathcal{G}_s^\flat|} \sum_{\bar{\varrho}_s \in \widehat{\mathcal{G}_s^\flat}} d_{\bar{\varrho}_s} \cdot \tr\left( \widehat{\theta_s^\flat}(\bar{\varrho}_s) \cdot \bar{\varrho}_s(\bar{g}_s) \right),
\]
where $\bar{g}_s \in \mathcal{G}_s^\flat$ and the Fourier transform is given by
\[
\widehat{\theta_s^\flat}(\bar{\varrho}_s) = \sum_{\bar{h} \in \mathcal{G}_s^\flat} \theta_s^\flat(\bar{h}) \cdot \bar{\varrho}_s(\bar{h}^{-1}).
\]

Let $(\mathcal{G}_s, W_s, \varrho_s)$ be an irreducible subrepresentation of the regular representation $(\mathcal{G}_s, V_s, \rho_s)$ such that $\varrho_s$ is trivial on $K$. By Propositions~\ref{proposition:irreducible} and~\ref{proposition:fourier_note}, the corresponding persistent representation $(\mathcal{G}_s, W_{s,t}, \varrho_{s,t})$ is irreducible. By Schur's lemma, the restriction $f_{s,t}|_{W_s} : W_s \to W_{s,t}$ is an isomorphism of irreducible representations.
This induces a one-to-one correspondence between irreducible representations of $\mathcal{G}_s$ that are trivial on $K$ and the persistent irreducible representations arising in the persistent regular representation. Moreover, there exists a canonical correspondence between irreducible representations of $\mathcal{G}_s^\flat$ and persistent irreducible representations of $\mathcal{G}_s$. Specifically, for each $\bar{\varrho}_s \in \widehat{\mathcal{G}_s^\flat}$, we have
\[
\varrho_{s,t} \cong \bar{\varrho}_s \circ \pi,
\]
where $\pi : \mathcal{G}_s \to \mathcal{G}_s^\flat$ is the natural projection.

Consequently, we analyze the term $\tr\left( \widehat{\theta}_{s,t}(\varrho_s) \cdot \varrho_{s,t}(g_s) \right)$. Using the definition of the persistent Fourier transform:
\begin{align*}
\tr\left( \widehat{\theta}_{s,t}(\varrho_s) \cdot \varrho_{s,t}(g_s) \right)
&= \tr\left( \sum_{h_s \in \mathcal{G}_s} \theta_s(h_s) \cdot \varrho_{s,t}(h_s^{-1}) \cdot \varrho_{s,t}(g_s) \right) \\
&= \tr\left( \sum_{h_s \in \mathcal{G}_s} \theta_s(h_s) \cdot \varrho_{s,t}(g_s h_s^{-1}) \right).
\end{align*}
Since $\varrho_{s,t}$ is trivial on $K$, it factors through $\pi$, i.e., $\varrho_{s,t}(g) = \bar{\varrho}_s(\pi(g))$. The sum over $\mathcal{G}_s$ decomposes into a sum over cosets $\bar{h} \in \mathcal{G}_s^\flat$ and elements $k \in K$. It follows that
\begin{align*}
\tr\left( \sum_{h_s \in \mathcal{G}_s} \theta_s(h_s) \cdot \varrho_{s,t}(g_s h_s^{-1}) \right)&= \tr\left( \sum_{\bar{h}_s \in \mathcal{G}_s^\flat} \sum_{k \in K} \theta_s^\flat(\bar{h}_s) \cdot \bar{\varrho}_s(\pi(g_s h_s^{-1} k)) \right) \\
&= \tr\left( \sum_{\bar{h}_s \in \mathcal{G}_s^\flat} \theta_s^\flat(\bar{h}_s) \bar{\varrho}_s(\bar{g}_s \bar{h}_s^{-1}) \sum_{k \in K} \mathbf{1} \right) \\
&= |K| \cdot \tr\left( \sum_{\bar{h}_s \in \mathcal{G}_s^\flat} \theta_s^\flat(\bar{h}_s) \cdot \bar{\varrho}_s(\bar{h}_s^{-1}) \cdot \bar{\varrho}_s(\bar{g}_s) \right) \\
&= |K| \cdot \tr\left( \widehat{\theta_s^\flat}(\bar{\varrho}_s) \cdot \bar{\varrho}_s(\bar{g}_s) \right).
\end{align*}

Now substitute this back into the formula from the theorem statement. Note that $d_{\varrho_s} = d_{\bar{\varrho}_s}$ and $\sum_{\varrho_s \in \widehat{(\mathcal{G}_s)}_K}$ is equivalent to $\sum_{\bar{\varrho}_s \in \widehat{\mathcal{G}_s^\flat}}$. Also $|\mathcal{G}_s| = |K| \cdot |\mathcal{G}_s^\flat|$. Thus, we have
\begin{align*}
\frac{1}{|\mathcal{G}_s|} \sum_{\varrho_s} d_{\varrho_s} \tr\left( \widehat{\theta}_{s,t}(\varrho_s) \varrho_{s,t}(g_s) \right)
&= \frac{1}{|K| \cdot |\mathcal{G}_s^\flat|} \sum_{\bar{\varrho}_s} d_{\bar{\varrho}_s} \cdot |K| \cdot \tr\left( \widehat{\theta_s^\flat}(\bar{\varrho}_s) \cdot \bar{\varrho}_s(\bar{g}_s) \right) \\
&= \frac{1}{|\mathcal{G}_s^\flat|} \sum_{\bar{\varrho}_s} d_{\bar{\varrho}_s} \cdot \tr\left( \widehat{\theta_s^\flat}(\bar{\varrho}_s) \cdot \bar{\varrho}_s(\bar{g}_s) \right) \\
&= \theta_s^\flat(\bar{g}_s) = \theta_s(g_s).
\end{align*}
This completes the proof.
\end{proof}

\subsection{Persistent convolution and Fourier transform}

Persistent convolution provides a theoretical foundation for subsequent investigations into spectral feature computation and the analysis of persistent Laplacian operators. In this section, we formally define the notion of persistent convolution and establish the convolution-multiplication theorem, which generalizes the classical convolution theorem to the persistent setting.

\begin{definition}
Let $\mathcal{G} : (\mathbb{Z}, \leq) \to \mathbf{Grp}$ be a persistence group. Let $\{\theta_t\}$ and $\{\eta_t\}$ be persistent functions on $\mathcal{G}$. The \textit{$(s,t)$-persistent convolution} of $\{\theta_t\}$ and $\{\eta_t\}$ is defined by

\[
(\theta \ast_{s} \eta)_t(g) := \sum_{h \in \mathcal{G}_s} \theta_s(h) \cdot \eta_t(\phi_{s,t}(h^{-1})g), \quad \forall g \in \mathcal{G}_t.
\]
\end{definition}

\begin{proposition}
Let $\{\theta_t\}$ and $\{\eta_t\}$ be persistent functions on $\mathcal{G}$. The convolution $\{(\theta \ast \eta)_{t}\}_{t\geq s}$ is also a persistent function on $\mathcal{G}$.
\end{proposition}

\begin{proof}
We want to show that the convolution family $\{(\theta \ast_s \eta)_t\}_{t \geq s}$ defines a persistent function on $\mathcal{G}$, that is, for any $s \leq r \leq t$,
\[
(\theta \ast_s \eta)_r = (\theta \ast_s \eta)_t \circ \phi_{r,t}.
\]
Since $\{\theta_t\}$ and $\{\eta_t\}$ are persistent functions, we have
\[
\theta_s = \theta_t \circ \phi_{s,t}, \quad \eta_r = \eta_t \circ \phi_{r,t}.
\]
It follows that
\[
\eta_r(\phi_{s,r}(h^{-1}) g) = \eta_t\big(\phi_{r,t}(\phi_{s,r}(h^{-1}) g)\big) = \eta_t\big(\phi_{s,t}(h^{-1}) \phi_{r,t}(g)\big).
\]
A straightforward calculation shows
\begin{align*}
  (\theta \ast_s \eta)_r(g)= & \sum_{h \in \mathcal{G}_s} \theta_s(h) \cdot \eta_r(\phi_{s,r}(h^{-1}) g) \\
  = &\sum_{h \in \mathcal{G}_s} \theta_s(h) \cdot \eta_t(\phi_{s,t}(h^{-1})\phi_{r,t}(g))\\
  = &(\theta \ast_s \eta)_t(\phi_{r,t}(g)).
\end{align*}
Thus, $\{(\theta \ast_s \eta)_t\}_{t \geq s}$ is a persistent function on $\mathcal{G}$.
\end{proof}

The persistent convolution defined above is generally non-commutative. However, it becomes commutative either when the two persistent functions coincide, i.e., ${\theta_t} = {\eta_t}$, or when the time indices coincide, i.e., $s = t$.

\begin{theorem}
Let $\mathcal{G} : (\mathbb{Z}, \leq) \to \mathbf{Grp}$ be a persistence group with structure maps $\phi_{s,t} : \mathcal{G}_s \to \mathcal{G}_t$, and let $\theta = \{ \theta_t \}$ and $\eta = \{ \eta_t \}$ be two persistent functions on $\mathcal{G}$. Then we have
\[
\widehat{(\theta \ast_{s} \eta)}_{s,t}(\varrho_s) = \widehat{\theta}_{s,t}(\varrho_s) \cdot \widehat{\eta}_{s,t}(\varrho_s),
\]
for all $\varrho_s \in \widehat{\mathcal{G}}_s$.
\end{theorem}

\begin{proof}
By definition, the $(s,t)$-persistent Fourier transform of the convolution is given by
\begin{align*}
  \widehat{(\theta \ast_{s} \eta)}_{s,t}(\varrho_s)
  &= \sum_{g \in \mathcal{G}_s} (\theta \ast_{s} \eta)_s(g) \cdot \varrho_{s,t}(g^{-1}) \\
  &= \sum_{g \in \mathcal{G}_s} \left( \sum_{h \in \mathcal{G}_s} \theta_s(h) \cdot \eta_s(h^{-1}g) \right) \cdot \varrho_{s,t}(g^{-1}) \\
  &= \sum_{h \in \mathcal{G}_s} \theta_s(h) \sum_{g \in \mathcal{G}_s} \eta_s(h^{-1}g) \cdot \varrho_{s,t}(g^{-1}) .
\end{align*}
Here, we use the change of variable $u = h^{-1}g$ in the inner sum, which implies $g = hu$ and $g^{-1} = u^{-1}h^{-1}$. Since $g$ varies over $\mathcal{G}_s$, $u$ also varies over $\mathcal{G}_s$. The expression becomes
\begin{align*}
  & \sum_{h \in \mathcal{G}_s} \theta_s(h) \sum_{g \in \mathcal{G}_s} \eta_s(h^{-1}g) \cdot \varrho_{s,t}(g^{-1})  \\
  =& \sum_{h \in \mathcal{G}_s} \theta_s(h) \sum_{u \in \mathcal{G}_s} \eta_s(u) \cdot \varrho_{s,t}(u^{-1}h^{-1}) \\
  =& \sum_{h \in \mathcal{G}_s} \theta_s(h) \sum_{u \in \mathcal{G}_s} \eta_s(u) \cdot \varrho_{s,t}(u^{-1}) \varrho_{s,t}(h^{-1}).
\end{align*}
Using the representation property and linearity, we can separate the sums:
\[
  \sum_{h \in \mathcal{G}_s} \theta_s(h) \varrho_{s,t}(h^{-1}) \cdot \left( \sum_{u \in \mathcal{G}_s} \eta_s(u) \cdot \varrho_{s,t}(u^{-1}) \right).
\]
Noting that $\widehat{\eta}_{s,t}(\varrho_s) = \sum_{g \in \mathcal{G}_s} \eta_s(g) \cdot \varrho_{s,t}(g^{-1})$, we obtain
\[
  \widehat{(\theta \ast_{s} \eta)}_{s,t}(\varrho_s)
  = \left( \sum_{h \in \mathcal{G}_s} \theta_s(h) \cdot \varrho_{s,t}(h^{-1}) \right) \cdot \widehat{\eta}_{s,t}(\varrho_s)
  = \widehat{\theta}_{s,t}(\varrho_s) \cdot \widehat{\eta}_{s,t}(\varrho_s).
\]
This completes the proof.
\end{proof}

\subsection{Persistent Fourier spectral analysis}

Let $\mathcal{G} : (\mathbb{Z}, \leq) \to \mathbf{Grp}$ be a persistence group and $\{ \theta_t \}_{t \in \mathbb{Z}}$ a persistent function on $\mathcal{G}$, where $\theta_t \in \mathrm{Hom}(\mathcal{G}_t, \mathbb{C})$. Denote the $(s,t)$-persistent Fourier transform by
\[
\widehat{\theta}_{s,t}(\varrho_s) = \sum_{g_s \in \mathcal{G}_s} \theta_s(g_s) \cdot \varrho_{s,t}(g_s^{-1}), \quad \varrho_s \in \widehat{\mathcal{G}}_s,
\]
where $\varrho_{s,t}$ is the $(s,t)$-persistent lift of the irreducible representation $\varrho_s$ along the structure map $\phi_{s,t}$.

\begin{definition}[Persistent Fourier spectrum]
The \textit{$(s,t)$-persistent Fourier spectrum} of $\theta$ is the function
\[
\mathrm{Spec}_{s,t}(\theta) : \widehat{\mathcal{G}}_s \to \mathbb{R}_{\geq 0}, \quad
\varrho_s \mapsto \left\| \widehat{\theta}_{s,t}(\varrho_s) \right\|,
\]
where $\| \cdot \|$ denotes a chosen matrix norm, such as the Hilbert-Schmidt norm.
\end{definition}

For each irreducible representation $\varrho_s \in \widehat{\mathcal{G}}_s$, the value $\mathrm{Spec}_{s,t}(\theta)(\varrho_s)$ quantifies the persistent spectral energy contributed by the frequency component $\varrho_s$ at time $t$. The $(s,t)$-persistent Fourier energy of the function $\theta$ is defined as the total energy across all frequency components
\[
\mathcal{E}_{s,t}(\theta) = \sum_{\varrho_s \in \widehat{\mathcal{G}}_s} \left\| \widehat{\theta}_{s,t}(\varrho_s) \right\|^2.
\]
This quantity measures the total preserved energy of the signal $\theta_s$ after being propagated to time $t$ through the persistent group structure.

\begin{theorem}\label{theorem:persistent_spectral}
Let $\varrho_s \in \widehat{\mathcal{G}}_s$ be an irreducible representation. If the persistent representation $\varrho_{s,t}$ is nonzero, then the persistent Fourier transform $\widehat{\theta}_{s,t}(\varrho_s)$ has the same spectrum as the classical Fourier transform $\widehat{\theta}_s(\varrho_s)$. That is,
\[
  \mathrm{Spec}_{s,t}(\theta)(\varrho_s) = \mathrm{Spec}_{s}(\theta)(\varrho_s),
\]
where $\mathrm{Spec}_{s}(\theta)(\varrho_s)= \left\| \widehat{\theta}_s(\varrho_s) \right\|$.
\end{theorem}

\begin{proof}
Let $\varrho_s : \mathcal{G}_s \to \mathrm{GL}(W_s)$ be an irreducible representation of $\mathcal{G}_s$, and let $\phi_{s,t} : \mathcal{G}_s \to \mathcal{G}_t$ be the structure morphism of the persistence group. The persistent representation $\varrho_{s,t}$ acts on a space $W_{s,t}=f_{s,t}(W_s)$ isomorphic to $W_s$. Here, $f_{s,t} : W_s \to W_{s,t}$ is the linear extension of $\phi_{s,t}$ restricted to $W_s$. Note that
\[
f_{s,t} \circ \varrho_s(g_s) = \varrho_{s,t}(g_s) \circ f_{s,t},\quad \forall g_s \in \mathcal{G}_s.
\]
This implies that for the inverse element, we also have
\[
f_{s,t} \circ \varrho_s(g_s^{-1}) = \varrho_{s,t}(g_s^{-1}) \circ f_{s,t}.
\]
By definition of the Fourier transforms, we have
\begin{align*}
\widehat{\theta}_{s,t}(\varrho_s)
&= \sum_{g_s \in \mathcal{G}_s} \theta_s(g_s) \cdot \varrho_{s,t}(g_s^{-1}) \\
&= \sum_{g_s \in \mathcal{G}_s} \theta_s(g_s) \cdot f_{s,t} \circ \varrho_s(g_s^{-1}) \circ f_{s,t}^{-1} \\
&= f_{s,t} \circ \left( \sum_{g_s \in \mathcal{G}_s} \theta_s(g_s) \cdot \varrho_s(g_s^{-1}) \right) \circ f_{s,t}^{-1} \\
&= f_{s,t} \circ \widehat{\theta}_s(\varrho_s) \circ f_{s,t}^{-1}.
\end{align*}
Therefore, $\widehat{\theta}_{s,t}(\varrho_s)$ is conjugate to $\widehat{\theta}_s(\varrho_s)$, and in particular, they have the same spectrum. Since matrix norms like the Hilbert--Schmidt norm are unitarily invariant (i.e., invariant under conjugation), it follows that
\[
\left\| \widehat{\theta}_{s,t}(\varrho_s) \right\| = \left\| \widehat{\theta}_s(\varrho_s) \right\|.
\]
Hence, the equation
\[
\mathrm{Spec}_{s,t}(\theta)(\varrho_s) = \mathrm{Spec}_s(\theta)(\varrho_s),
\]
holds as claimed.
\end{proof}

Based on the persistent Fourier spectrum, we can further introduce tools such as the persistent spectral entropy, persistent Fourier spectrum trajectory, and persistent Fourier correlation, to analyze the structure and evolution of signals in the persistence setting.

\begin{definition}[Persistent spectral entropy]
Let
\[
  D_{s,t} = \sum_{\psi \in \widehat{\mathcal{G}}_s} \| \widehat{\theta}_{s,t}(\psi) \|^2.
\]
If $D_{s,t} = 0$, we define the \textit{$(s,t)$-persistent spectral entropy} to be $H_{s,t}(\theta) := 0$. 
Otherwise, let
\[
  p_{s,t}(\varrho_s) = \frac{ \| \widehat{\theta}_{s,t}(\varrho_s) \|^2 }{ D_{s,t} }
\]
be the normalized spectral weight. The \textit{$(s,t)$-persistent spectral entropy} is defined by
\[
  H_{s,t}(\theta) := - \sum_{\varrho_s \in \widehat{\mathcal{G}}_s} p_{s,t}(\varrho_s) \log p_{s,t}(\varrho_s).
\]
\end{definition}

\begin{definition}[Persistent Fourier spectrum trajectory]
Fix $\varrho_s \in \widehat{\mathcal{G}}_s$. The \textit{persistent spectrum trajectory} of $\varrho_s$ is the function
\[
\mathsf{Traj}_{\varrho_s}(t) := \left\| \widehat{\theta}_{s,t}(\varrho_s) \right\|, \quad \text{for } t \geq s.
\]
\end{definition}

\begin{definition}[Persistent Fourier correlation]
Let $\{ \theta_t \}$ and $\{ \eta_t \}$ be two persistent functions on $\mathcal{G}$. Their \textit{$(s,t)$-persistent Fourier correlation} is given by
\[
\mathrm{Corr}_{s,t}(\theta,\eta) := \sum_{\varrho_s \in \widehat{\mathcal{G}}_s} \operatorname{Tr} \left( \widehat{\theta}_{s,t}(\varrho_s)^* \cdot \widehat{\eta}_{s,t}(\varrho_s) \right).
\]
\end{definition}

\begin{example}

Continuing from Example~\ref{example:persistent_fourier}, consider the cosine signal $\theta_s([x]_s) = \cos\left(  2\pi x/2^s  \right)$ defined on the group $\mathcal{G}_s = \mathbb{Z}/2^s\mathbb{Z}$. Its $(s,t)$-persistent Fourier transform is supported only at frequencies $k \equiv \pm 1 \mod 2^s$. The squared spectral norms are given by
\[
\| \widehat{\theta}_{s,t}(k) \|^2 =
\begin{cases}
2^{2(s-1)}, & \text{if } k \equiv \pm 1 \mod 2^s, \\
0, & \text{otherwise}.
\end{cases}
\]
Accordingly, the normalized spectral weights become
\[
p_{s,t}(k) =
\begin{cases}
\frac{1}{2}, & \text{if } k \equiv \pm 1 \mod 2^s, \\
0, & \text{otherwise}.
\end{cases}
\]
This yields a persistent spectral entropy of
\[
H_{s,t}(\theta) = -\left( \frac{1}{2} \log \frac{1}{2} + \frac{1}{2} \log \frac{1}{2} \right) = \log 2,
\]
indicating that the spectral energy is evenly distributed between two frequency modes.

Furthermore, the persistent spectrum trajectory for $k \equiv \pm 1 \mod 2^s$ is the constant function
\[
\mathsf{Traj}_k(t) = \left\| \widehat{\theta}_{s,t}(k) \right\| = 2^{s-1},
\]
while for all other frequencies $k$, the trajectory is identically zero.

Let $\eta_t([x]_t) = \sin\left(2\pi x/2^t \right)$. Then the family $\{\eta_{t}\}_{t\in \mathbb{Z}_{\geq 0}}$ also defines a persistent function in $\varprojlim_{t \in \mathbb{Z}} \mathrm{Hom}(\mathcal{G}_t,\mathbb{C})$ on the persistence group $\mathcal{G}$. A straightforward calculation shows that
\[
\widehat{\eta}_{s,t}(k) =
\begin{cases}
-2^{s-1}i, & \text{if } k \equiv 1 \mod 2^s, \\
2^{s-1}i, & \text{if } k \equiv -1 \mod 2^s, \\
0, & \text{otherwise}.
\end{cases}
\]
Then the $(s,t)$-persistent Fourier correlation is given by
\[
\mathrm{Corr}_{s,t}(\theta,\eta) = (2^{s-1})\cdot (-2^{s-1}i) + (2^{s-1})\cdot (2^{s-1}i) = 0.
\]
This reflects the orthogonality of $\cos$ and $\sin$ at the same frequency under the persistent Fourier transform.
\end{example}

\subsection{Persistent Laplacian operator}

Let $\mathcal{G} : (\mathbb{Z}, \leq) \to \mathbf{Grp}$ be a persistence group, and let $\theta=\{\theta_t\}_{t \in \mathbb{Z}}$ be a persistent function on $\mathcal{G}$.
By Remark \ref{remark:persistent_function}, any function $\eta_t: \mathcal{G}_t \to \mathbb{C}$ can extend a persistent function $\eta=\{\eta_r\}_{r\leq t}$ on $\mathcal{G}$.

\begin{definition}
For any integers $s \leq t$, the $(s,t)$-persistent convolution operator is defined as
\[
T_\theta^{s,t} : \mathrm{Hom}(\mathcal{G}_t,\mathbb{C}) \to  \mathrm{Hom}(\mathcal{G}_t,\mathbb{C}), \quad T_\theta^{s,t}\eta_t = (\theta \ast_{s} \eta)_t,
\]
where $\eta_t: \mathcal{G}_t \to \mathbb{C}$ is a function and $\eta$ is the extended persistent function of $\eta_t$ on $\mathcal{G}$.
\end{definition}

\begin{definition}
Let $\eta_t: \mathcal{G}_t \to \mathbb{C}$ be any function. For any integers $s \leq t$, the \textit{$(s,t)$-persistent Laplacian} operator
\[
\mathcal{L}_{s,t} : \mathrm{Hom}(\mathcal{G}_t, \mathbb{C}) \to \mathrm{Hom}(\mathcal{G}_t, \mathbb{C})
\]
is defined by
\[
\mathcal{L}_{s,t} := d_s \cdot \mathrm{id} - T_\theta^{s,t},
\]
where the scalar $d_s$ is given by $d_s = \sum_{g \in \mathcal{G}_s} \theta_s(g)$.
\end{definition}
In particular, when $s = t$, the persistent Laplacian $\mathcal{L}_{s,t}$ reduces to the classical Laplacian associated with the function $\theta_s$ on the group $\mathcal{G}_s$.

By definition, for any $g \in \mathcal{G}_t$, the $(s,t)$-persistent Laplacian acts as
\begin{equation*}
  \mathcal{L}_{s,t}(\eta_t)(g) = d_s \cdot \eta_t(g) - \sum_{h \in \mathcal{G}_s} \theta_s(h) \cdot \eta_t(\phi_{s,t}(h^{-1})g),
\end{equation*}
where $d_s = \sum_{h \in \mathcal{G}_s} \theta_s(h)$. The operator $\mathcal{L}_{s,t}$ is linear and vanishes on constant functions. Under appropriate conditions on the function $\theta_s$, the operator $\mathcal{L}_{s,t}$ is also self-adjoint and positive semi-definite.

\begin{proposition}\label{proposition:adjoint}
Let $\mathcal{G} : (\mathbb{Z}, \leq) \to \mathbf{Grp}$ be a persistence group, and let $\theta = \{\theta_t\}_{t \in \mathbb{Z}}$ be a persistent function such that for each $s$, $\theta_s(g^{-1}) = \overline{\theta_s(g)}$ for all $g \in \mathcal{G}_s$. Then for any $s \leq t$, the $(s,t)$-persistent Laplacian operator $\mathcal{L}_{s,t}$ is self-adjoint with respect to the standard inner product
\[
\langle f, g \rangle = \sum_{x \in \mathcal{G}_t} f(x) \cdot \overline{g(x)},\quad f,g\in \mathrm{Hom}(\mathcal{G}_t, \mathbb{C}).
\]
\end{proposition}

\begin{proof}
For any $f, g \in \mathrm{Hom}(\mathcal{G}_t, \mathbb{C})$, we compute the inner product
\begin{align*}
  \langle T_\theta^{s,t} f, g \rangle &= \sum_{x \in \mathcal{G}_t} \left( \sum_{h \in \mathcal{G}_s} \theta_s(h) f(\phi_{s,t}(h^{-1})x) \right) \cdot \overline{g(x)}  \\
  &= \sum_{h \in \mathcal{G}_s} \theta_s(h) \sum_{x \in \mathcal{G}_t} f(\phi_{s,t}(h^{-1})x) \cdot \overline{g(x)}.
\end{align*}
On the other hand, we have
\begin{align*}
  \langle f, T_\theta^{s,t} g \rangle &= \sum_{x \in \mathcal{G}_t} f(x) \cdot \overline{ \sum_{h \in \mathcal{G}_s} \theta_s(h) g(\phi_{s,t}(h^{-1})x)} \\
  &= \sum_{h \in \mathcal{G}_s} \overline{\theta_s(h)} \sum_{x \in \mathcal{G}_t} f(x) \cdot \overline{g(\phi_{s,t}(h^{-1})x)} \\
  &= \sum_{h \in \mathcal{G}_s} \overline{\theta_s(h)} \sum_{y \in \mathcal{G}_t} f(\phi_{s,t}(h)y) \cdot \overline{g(y)}.
\end{align*}
Here, we use the change of variables $y = \phi_{s,t}(h^{-1})x$.

Now, using the Hermitian condition $\overline{\theta_s(h)} = \theta_s(h^{-1})$, and letting $h' = h^{-1}$, we obtain
\[
\langle f, T_\theta^{s,t} g \rangle = \sum_{h' \in \mathcal{G}_s} \theta_s(h') \sum_{y \in \mathcal{G}_t} f(\phi_{s,t}(h'^{-1})y) \cdot \overline{g(y)} = \langle T_\theta^{s,t} f, g \rangle.
\]
Hence, $T_\theta^{s,t}$ is self-adjoint. Since the $(s,t)$-persistent Laplacian is defined as
\[
\mathcal{L}_{s,t} = d_s \cdot \mathrm{id} - T_\theta^{s,t}, \quad \text{with } d_s = \sum_{h \in \mathcal{G}_s} \theta_s(h),
\]
it follows that $\mathcal{L}_{s,t}$ is also self-adjoint.
\end{proof}

\begin{proposition}
Let $\mathcal{G} : (\mathbb{Z}, \leq) \to \mathbf{Grp}$ be a persistence group, and let $\theta = \{\theta_t\}$ be a persistent function on $\mathcal{G}$. Suppose that for all $h \in \mathcal{G}_s$, the condition $\theta_s(h^{-1}) = \overline{\theta_s(h)}$ holds, and that $\theta_s(h) \geq 0$. Then the operator $\mathcal{L}_{s,t}$ is positive semi-definite. That is, for all $f \in \mathrm{Hom}(\mathcal{G}_t, \mathbb{C})$, we have
\[
\langle \mathcal{L}_{s,t} f, f \rangle \geq 0.
\]
\end{proposition}

\begin{proof}
Fix a function $f \in \mathrm{Hom}(\mathcal{G}_t, \mathbb{C})$. By definition of $\mathcal{L}_{s,t}$, we have
\[
\langle \mathcal{L}_{s,t} f, f \rangle = d_s \langle f, f \rangle - \langle T_\theta^{s,t} f, f \rangle.
\]
The convolution operator $T_\theta^{s,t}$ acts as
\[
T_\theta^{s,t} f(x) = \sum_{h \in \mathcal{G}_s} \theta_s(h) f(\phi_{s,t}(h^{-1}) x).
\]
It follows that
\[
\langle T_\theta^{s,t} f, f \rangle = \sum_{x \in \mathcal{G}_t} \left( \sum_{h \in \mathcal{G}_s} \theta_s(h) f(\phi_{s,t}(h^{-1}) x) \right) \overline{f(x)} = \sum_{h \in \mathcal{G}_s} \theta_s(h) \sum_{x \in \mathcal{G}_t} f(\phi_{s,t}(h^{-1}) x) \overline{f(x)}.
\]
Changing variables by setting $y = \phi_{s,t}(h^{-1}) x$, and using the Hermitian property $\theta_s(h^{-1}) = \overline{\theta_s(h)}$, we can rewrite the sum as
\[
\sum_{x \in \mathcal{G}_t} f(\phi_{s,t}(h^{-1}) x) \overline{f(x)} = \sum_{y \in \mathcal{G}_t} f(y) \overline{f(\phi_{s,t}(h) y)}.
\]
Thus, we obtain
\[
\langle \mathcal{L}_{s,t} f, f \rangle = \frac{1}{2} \sum_{h \in \mathcal{G}_s} \theta_s(h) \sum_{x \in \mathcal{G}_t} \left( |f(x)|^2 + |f(\phi_{s,t}(h^{-1}) x)|^2 - 2 \mathrm{Re}(f(\phi_{s,t}(h^{-1}) x) \overline{f(x)}) \right),
\]
which simplifies to
\[
\langle \mathcal{L}_{s,t} f, f \rangle = \frac{1}{2} \sum_{h \in \mathcal{G}_s} \theta_s(h) \sum_{x \in \mathcal{G}_t} |f(x) - f(\phi_{s,t}(h^{-1}) x)|^2.
\]
Each term in the sum is non-negative since $\theta_s(h) \geq 0$, hence the total sum is non-negative. This completes the proof.
\end{proof}

\subsection{An application on symmetry groups}

From Section~\ref{section:span}, we recall the pseudofunctor
\[
\mathrm{Sym} : \mathbb{S}_n(M) \to \mathrm{Span}(\mathbf{Grp}),
\]
which assigns to each $n$-configuration $X$ a symmetry group $\mathrm{Sym}(X)$, and to each morphism $f : X \to Y$ in $\mathbb{S}_n(M)$ a span of group homomorphisms
\[
\mathrm{Sym}(X) \xleftarrow{f^\flat} \mathrm{Sym}_f(X) \xrightarrow{f^\sharp} \mathrm{Sym}(Y).
\]

Consider the group embedding $f^\sharp : \mathrm{Sym}_f(X) \to \mathrm{Sym}(Y)$. Fix a point $x \in X$ and $y = f(x) \in Y$. Consider the following functions defined by
\begin{equation*}
  \theta_X(\pi) = \|\pi \cdot x - x\|^2, \qquad \pi \in \mathrm{Sym}_f(X),
\end{equation*}
and
\begin{equation*}
  \theta_Y(\pi) = \|f^{-1}(\pi \cdot y) - f^{-1}(y)\|^2, \qquad \pi \in \mathrm{Sym}(Y).
\end{equation*}
Note that for any $\pi \in \mathrm{Sym}_f(X)$ and any $y\in Y$, we have
\begin{equation*}
  f^\sharp(\pi)(y) = f(\pi(f^{-1}(y))).
\end{equation*}
It follows that
\begin{equation*}
  \theta_Y(f^\sharp(\pi)) =  \|\pi(f^{-1}(y)) - f^{-1}(y)\|^2 = \|\pi(x) - x\|^2 = \theta_X(\pi).
\end{equation*}
Therefore, we obtain the following result.
\begin{lemma}
The pair $\theta=\{\theta_X, \theta_Y\}$ defines a persistent function along the group homomorphism $f^\sharp : \mathrm{Sym}_f(X) \to \mathrm{Sym}(Y)$.
\end{lemma}

The function $\theta_X$ measures the squared Euclidean displacement of a point $x \in X$ under the group action of $\pi \in \mathrm{Sym}_f(X)$. Intuitively, this function encodes how much each group element $\pi$ perturbs the configuration around the point $x$. In this way, $\theta_X$ can be interpreted as a descriptor of local symmetry breaking or deviation from invariance at $x$. In addition, this construction is closely related to the (normalized) Laplacian operator in spectral graph theory, where such Laplacians are typically defined as sums of squared displacements over local neighborhoods.

Now, the persistent Fourier transform associated with the pair $(X, Y)$ is expressed as
\begin{equation*}
\widehat{\theta}_{X,Y}(\varrho_X)
= \sum_{\pi \in \mathrm{Sym}_f(X)} \|\pi \cdot x - x\|^2 \cdot \varrho_{X,Y}(\pi^{-1}),
\end{equation*}
where $\varrho_{X,Y}$ denotes the persistent representation of the irreducible representation $\varrho_X$ of $\mathrm{Sym}_f(X)$. As established in Proposition~\ref{proposition:fourier_note}, the persistent representation $\varrho_{X,Y}$ acts on a space $f^\sharp V_{\varrho_X}$ isomorphic to $V_{\varrho_X}$ via the linear extension of $f^\sharp$. Consequently, for any vector $v \in V_{\varrho_X}$, the action can be identified with the original representation:
\[
\varrho_{X,Y}(\pi)(f^\sharp(v)) = f^\sharp(\varrho_X(\pi)v).
\]
This implies that the form of the persistent Fourier transform remains the same for different choices of $\varrho_X$. However, since the representation spaces vary, the Fourier transform
\[
\widehat{\theta}_{X,Y}(\varrho_X) : f^{\sharp}V_{\varrho_X} \to f^{\sharp}V_{\varrho_X}
\]
is defined on different spaces depending on $\varrho_X$.

Next, we analyze the spectral information of the persistent Fourier transform $\widehat{\theta}_{X,Y}(\varrho_X)$.
By Theorem~\ref{theorem:persistent_spectral}, we have
\begin{equation*}
  \| \widehat{\theta}_{X,Y}(\varrho_X)\| = \|\widehat{\theta}_{X}(\varrho_X)\|,
\end{equation*}
where
\begin{equation*}
  \widehat{\theta}_{X}(\varrho_X) = \sum_{\pi \in \mathrm{Sym}_f(X)} \|\pi \cdot x - x\|^2 \cdot \varrho_{X}(\pi^{-1}).
\end{equation*}
By a straightforward calculation, we obtain
\[
  \|\widehat{\theta}_X(\varrho_X)\|^2 = \sum_{\pi, \pi' \in \mathrm{Sym}_f(X)} \|\pi \cdot x - x\|^2\|\pi' \cdot x - x\|^2 \cdot \tr\left( \varrho_X(\pi^{-1} \pi') \right).
\]
Let $H$ be a subgroup of $\mathrm{Sym}_f(X)$. Consider the restriction of the function
\[
\theta_X(\pi) = \|\pi \cdot x - x\|^2
\]
to the subgroup $H$, and define the corresponding restricted Fourier transform by
\[
\widehat{\theta}_{X,H}(\varrho_X) = \sum_{h \in H} \theta_X(h) \cdot \varrho_X(h^{-1}),
\]
where $\varrho_X$ denotes a unitary representation of $\mathrm{Sym}_f(X)$ restricted to $H$. This allows us to analyze the spectral contribution of specific symmetry subgroups.

\begin{proposition}
Suppose $H = \{ \pi \in \mathrm{Sym}_f(X) \mid \pi \cdot x = x \}$ is the stabilizer subgroup of the configuration point $x \in X$. Then for any $\pi \in H$, the restricted Fourier transform $\widehat{\theta}_{X,H}(\varrho_X)$ vanishes.
\end{proposition}
\begin{proof}
Since $\pi \cdot x = x$ for all $\pi \in H$, it follows that $\theta_X(\pi) = 0$. Hence, every term in the sum vanishes and $\widehat{\theta}_{X,H}(\varrho_X) = 0$.
\end{proof}

\begin{corollary}
If $\mathrm{Sym}_f(X)$ stabilizes the point $x$, then the Fourier transform $\widehat{\theta}_{X}(\varrho_X)$ vanishes. In particular, if $\mathrm{Sym}_f(X)$ is the trivial group, then $\widehat{\theta}_{X}(\varrho_X) = 0$.
\end{corollary}

We now recall the definition of the persistent Laplacian operator
\begin{equation*}
  \mathcal{L}_{X,Y} = d_X \cdot \mathrm{id} - T_{\theta}^{X,Y},
\end{equation*}
where the scalar coefficient is given by
\begin{equation*}
  d_X = \sum_{\pi \in \mathrm{Sym}_f(X)} \theta_X(\pi) = \sum_{\pi \in \mathrm{Sym}_f(X)} \|\pi \cdot x - x\|^2,
\end{equation*}
and the operator $T_{\theta}^{X,Y}$ acts on a function $\eta : \mathrm{Sym}(Y) \to \mathbb{C}$ by the formula
\begin{equation*}
  (T_{\theta}^{X,Y} \eta)(\sigma) = \sum_{\pi \in \mathrm{Sym}_f(X)} \theta_X(\pi) \cdot \eta\left( f^\sharp(\pi)^{-1} \sigma \right), \quad \sigma \in \mathrm{Sym}(Y).
\end{equation*}
Substituting the definition of $d_X$ into $\mathcal{L}_{X,Y}$, we obtain the explicit form
\begin{equation*}
  (\mathcal{L}_{X,Y} \eta)(\sigma) = \sum_{\pi \in \mathrm{Sym}_f(X)} \|\pi \cdot x - x\|^2 \cdot \left[ \eta(\sigma) - \eta\left( f^\sharp(\pi)^{-1} \sigma \right) \right].
\end{equation*}

This operator quantifies the deviation of the function $\eta$ under symmetry-induced shifts, and structurally resembles a non-normalized Laplacian operator defined via group convolution, with weights determined by the geometric displacement under the finite symmetry group $\mathrm{Sym}_f(X)$.

%% file: equivariant_persistence.tex
\section{Equivariant persistent cohomology}\label{section:equivariant_persistenc}

Equivariant cohomology is a standard topological invariant for capturing the intrinsic symmetric structures of spaces. The theoretical formulation of persistent equivariant cohomology was introduced in \cite{adams2024persistent}, which lays the algebraic foundation for the analysis of filtered $G$-spaces. The practical value of combining persistence and equivariance is also validated by topological deep learning architectures such as the TopNets model presented in \cite{verma2024topological}. In this section, we further study equivariant persistence modules, discuss the topological visualization and stability of persistent equivariant cohomology, and address the related computational aspects.

\subsection{Borel equivariant cohomology}
Let $G$ be a Hausdorff topological group, $X$ a left $G$-space, and fix a base field $\mathbb{K}$ of coefficients. There exists a contractible topological space $EG$ admitting a free right $G$-action, called the \textit{total space of the universal principal $G$-bundle}. Such a space can be constructed explicitly as follows:
\begin{itemize}
    \item For compact Lie groups $G$, the total space $EG = \varinjlim_{n} E_nG$ is the inductive limit of finite-dimensional free $G$-CW complexes;
    \item For discrete groups $G$, the total space $EG$ is the contractible simplicial complex with vertex set $G$, whose simplices are spanned by finite subsets of $G$.
\end{itemize}
The quotient space $BG := EG/G$ is the \textit{classifying space} of $G$. Up to weak homotopy equivalence, the pair $(EG,BG)$ is uniquely determined by $G$. The quotient projection $EG \twoheadrightarrow BG$ is a principal $G$-bundle, and every principal $G$-bundle over a paracompact space $Y$ is pulled back from $EG\to BG$ via a unique map $Y\to BG$ up to homotopy.

On the product $EG\times X$, there is a right $G$-action compatible with the left $G$-action on $X$ given by
\[
(e,x)\cdot g = (e\cdot g,\,g^{-1}\cdot x),\quad \forall g\in G,\,e\in EG,\,x\in X.
\]
The \textit{Borel construction} of $X$ is the orbit space
\[
X_G := EG\times_G X = (EG\times X)/G.
\]
\begin{definition}
The \textit{Borel equivariant cohomology ring} of $X$ is the ordinary singular cohomology of the Borel construction:
\[
H_G^p(X) := H^p(X_G;\mathbb{K}),\quad p\geq 0.
\]
\end{definition}
For the dual Borel equivariant homology, we use the notation $H_\ast^G(X) = H_\ast(X_G;\mathbb{K})$. If $X=\mathrm{pt}$ is the one-point $G$-space, then $X_G=BG$, so $H_G^\ast(\mathrm{pt})=H^{\ast}(BG)$. Let us denote $\mathcal{H}_G^\ast=H^{\ast}(BG)$. The natural fibration $\pi\colon X_G\to BG,\ [e,x]\mapsto [e]$ induced by projection to the $EG$-factor yields a pullback ring homomorphism
\[
\pi^*\colon \mathcal{H}_G^\ast\to H_G^\ast(X),
\]
endowing $H_G^\ast(X)$ with the structure of a graded commutative $\mathcal{H}_G^\ast$-algebra. Here, the cup product on $H_G^\ast(X)$ makes it a graded commutative ring.

\begin{remark}
The coefficient ring $\mathcal{H}_G^\ast \cong H^\ast(BG; \mathbb{K})$ is determined by the group topology and field characteristic. 
For finite groups, $\mathcal{H}_G^\ast \cong \mathbb{K}$ if $\operatorname{char}(\mathbb{K})=0$ or $\operatorname{char}(\mathbb{K}) \nmid |G|$, in which case the module structure simplifies significantly. Otherwise, it forms a polynomial algebra; for instance, $\mathcal{H}_{\mathbb{Z}/p}^\ast \cong \mathbb{F}_p[c] \otimes \Lambda[\tau]$, where $c$ denotes the first Chern class with $\deg(c)=2$ and $\tau$ is an exterior generator of degree $1$. 
For connected compact Lie groups, $\mathcal{H}_G^\ast$ is generally a polynomial ring generated by characteristic classes; for example, $\mathcal{H}_{S^1}^\ast \cong \mathbb{K}[c]$.
\end{remark}

\begin{proposition}[\cite{tu2020introductory}]
Let $f\colon X\to Y$ be a $G$-equivariant map of $G$-spaces. It descends to a continuous map $f_G\colon X_G\to Y_G,\ [e,x]\mapsto [e,f(x)]$, which yields a contravariant graded ring homomorphism
\[
f_G^*\colon H_G^\ast(Y)\to H_G^\ast(X).
\]
Thus, $H_G^\ast(-)$ is a contravariant functor from the category of $G$-spaces to the category of graded commutative $\mathcal{H}_G^\ast$-algebras.
\end{proposition}

\begin{remark}
The pullback homomorphism $\pi^* \colon \mathcal{H}_G^\ast \to H_G^\ast(X)$ equips $H_G^\ast(X)$ with the structure of a graded $\mathcal{H}_G^\ast$-module via scalar multiplication $\lambda\cdot \alpha = \pi^*(\lambda)\cup \alpha$ for all $\lambda\in\mathcal{H}_G^\ast$, $\alpha\in H_G^\ast(X)$.
Moreover, the cup product on $H_G^\ast(X)$ satisfies the graded commutativity law
\[
\alpha\cup\beta = (-1)^{\deg\alpha\cdot\deg\beta}\,\beta\cup\alpha
\]
for homogeneous elements $\alpha,\beta\in H_G^\ast(X)$, and the module multiplication commutes with the internal cup product.
Thus, $H_G^\ast(X)$ is a graded commutative $\mathcal{H}_G^\ast$-algebra.
\end{remark}

\subsection{Persistence $G$-complexes and equivariant persistent cohomology}

Let $G$ denote a Hausdorff topological group, and let $(P,\leq)$ be a totally ordered index poset, most commonly $P=\mathbb{R}$ or $P=\mathbb{Z}$.

\begin{definition}
A \textit{persistence $G$-space} is a functor
\[
\mathcal{X}\colon (P,\leq) \to \mathbf{Top}_G
\]
valued in the category $\mathbf{Top}_G$ of left $G$-spaces and $G$-equivariant continuous maps.
\end{definition}

Specifically, a \textit{$G$-filtered space} is a pair $(X,\{X_t\}_{t\in P})$, where $X$ is a left $G$-space and $\{X_t\}_{t\in P}$ is a nested family of subspaces of $X$ satisfying the following three axioms:
\begin{enumerate}[label=($\roman*$)]
    \item For all $t\in P$, $g\cdot X_t = X_t$ holds for every $g\in G$;
    \item If $s\leq t$ in $P$, then $X_s \subseteq X_t$ as subspaces of $X$;
    \item The canonical inclusion $i_{s,t}\colon X_s \hookrightarrow X_t$ is a $G$-equivariant continuous map whenever $s\leq t$.
\end{enumerate}
If each subspace $X_t$ carries the structure of a $G$-simplicial complex and every inclusion $i_{s,t}$ is a simplicial $G$-map, we call $(X,\{X_t\}_{t\in P})$ a \textit{$G$-filtered complex}.

\begin{example}
Let $X$ be a $G$-space, and suppose $f\colon X\to \mathbb{R}$ is a $G$-invariant real-valued function, meaning $f(g\cdot x)=f(x)$ for all $g\in G$ and $x\in X$. The sublevel sets
\[
X_t = f^{-1}\big((-\infty,t]\big),\quad t\in\mathbb{R},
\]
form a $G$-filtered space.
\end{example}

\begin{example}\label{example:G_filtered}
Let $E$ be a finite-dimensional Euclidean space, and $X\subset E$ a finite set of points with centroid $\bar{x}\in E$. Let $G\leq \operatorname{Iso}(E)$ be a subgroup of the Euclidean isometry group acting on $X$, and assume the centroid $\bar{x}$ is fixed by every element of $G$, i.e., $g\cdot \bar{x} = \bar{x}$ for all $g\in G$.
Define the radial distance function $f\colon X\to \mathbb{R}$ by $f(x) = \|x - \bar{x}\|$. For each real radius $r\ge 0$, the sublevel set
\[
X_r = f^{-1}\big((-\infty,r]\big) = \bigl\{x\in X \,\big|\, \|x - \bar{x}\| \le r\bigr\}
\]
yields a $G$-filtered space $(X,\{X_r\}_{r\in\mathbb{R}})$. To verify the $G$-invariance axiom: for any $g\in G$ and $x\in X$,
\[
f(g\cdot x) = \|g\cdot x - \bar{x}\| = \|g\cdot x - g\cdot \bar{x}\| = \|x - \bar{x}\| = f(x),
\]
where the third equality follows from the isometry property $\|g u - g v\|=\|u-v\|$ and the final equality uses the $G$-invariance of the centroid $\bar{x}$.
\end{example}

\begin{definition}
Let $\mathcal{H}_G^\ast = H^\ast(BG; \mathbb{K})$ be the equivariant cohomology ring of a point. A \textit{copersistence graded $\mathcal{H}_G^\ast$-module} is a functor
\[
    M: (P, \leq)^{\mathrm{op}} \longrightarrow \mathbf{Mod}_{\mathcal{H}_G^\ast}^{\mathrm{gr}},
\]
where $\mathbf{Mod}_{\mathcal{H}_G^\ast}^{\mathrm{gr}}$ denotes the category of graded $\mathcal{H}_G^\ast$-modules.
\end{definition}

Specifically, for each parameter $t \in P$, the assignment $M(t)$ is a graded $\mathcal{H}_G^\ast$-module; for $s \leq t$, the functor induces an $\mathcal{H}_G^\ast$-linear map $\varphi_{s,t}: M(t) \to M(s)$ (reversing the order) satisfying the compatibility condition $\varphi_{s,t} \circ \varphi_{t,u} = \varphi_{s,u}$ for $s \leq t \leq u$.

For a fixed integer $k$, the $k$-th component is given by
\[
    M^k(t) = M(t)^k, \quad M^k(s \leq t) = \varphi_{s,t}|_{M(t)^k}.
\]
The action of the ring $\mathcal{H}_G^\ast$ on the component $M(t)^k$ is given by
\[
    \mathcal{H}_G^p \cdot M(t)^k \subseteq M(t)^{p+k},\quad p \geq 0.
\]
Hence, $M(t)^k$ is a left $\mathcal{H}_G^0$-module. Since $BG$ is connected, we have $\mathcal{H}_G^0 \cong \mathbb{K}$. Consequently, $M^k(t)$ carries a natural $\mathbb{K}$-vector space structure, and the maps are $\mathbb{K}$-linear. Thus, the component $M^k$ forms a copersistence module over the base field $\mathbb{K}$.

\begin{definition}
Let $\mathcal{X} \colon (P, \leq) \longrightarrow \mathbf{Top}_G$ be a persistence $G$-space. The \textit{equivariant persistent cohomology} of $\mathcal{X}$ is defined as the functor
\[
    H_G^\ast(\mathcal{X}) \colon (P, \leq)^{\mathrm{op}} \longrightarrow \mathbf{Mod}_{\mathcal{H}_G^\ast}^{\mathrm{gr}},
\]
which acts on objects and morphisms as follows:
\begin{itemize}
    \item For each parameter $t \in P$, it assigns the graded $\mathcal{H}_G^\ast$-module
    \[
        H_G^\ast(\mathcal{X})(t) = H_G^\ast(X_t).
    \]
    \item For each relation $s \leq t$ in $P$, the structure map
    \[
        H_G^\ast(\mathcal{X})(s \leq t) = (i_{s,t})_G^\ast \colon H_G^\ast(X_t) \longrightarrow H_G^\ast(X_s)
    \]
    is the $\mathcal{H}_G^\ast$-linear homomorphism induced by the $G$-equivariant inclusion $i_{s,t} \colon X_s \hookrightarrow X_t$.
\end{itemize}
\end{definition}

The persistent cohomology can be interpreted as a copersistence module over $\mathcal{H}_G^\ast$, tracking the creation and propagation of topological features in reverse order with respect to the group action.

\begin{example}
Consider the unit sphere $S^2 \subset \mathbb{R}^3$ with the group $G = S^1$ acting by rotation about the $z$-axis. Define a $G$-invariant height function $f: S^2 \to \mathbb{R}$ by $f(x,y,z) = z$. This induces a $G$-filtered complex $\{X_t\}_{t \in \mathbb{R}}$ defined by sublevel sets $X_t = f^{-1}((-\infty, t])$. 

We analyze the $\mathcal{H}_{S^1}^\ast$-module structure of $H_{S^1}^\ast(X_t)$ as follows:
\begin{itemize}
    \item When $t < -1$, we have $X_t = \emptyset$. The module is trivial.
    
    \item When $-1 \le t < 1$, the sublevel set $X_t$ is a $G$-invariant polar cap that contracts equivariantly to the South Pole $S$. Since $S$ is a fixed point, the inclusion $S \hookrightarrow X_t$ induces an isomorphism
    \[
        H_{S^1}^\ast(X_t) \cong H_{S^1}^\ast(\mathrm{pt}) \cong \mathbb{K}[c],
    \]
    where $c \in H^2(BS^1)$ is the first Chern class. This is a \textit{free $\mathcal{H}_{S^1}^\ast$-module of rank 1}, generated by the fundamental class of the connected component.
    
    \item When $t \ge 1$, we obtain $X_t = S^2$. The $S^1$-action has two fixed points (North and South Poles). By the GKM description (or Mayer-Vietoris sequence), the equivariant cohomology is
    \[
        H_{S^1}^\ast(S^2) \cong \mathbb{K}[c] \oplus \mathbb{K}[c].
    \]
    This is a \textit{free $\mathcal{H}_{S^1}^\ast$-module of rank 2}.
\end{itemize}

The transition map $H_{S^1}^\ast(X_t) \to H_{S^1}^\ast(X_s)$ for $s \in (-1, 1)$ and $t \ge 1$ is surjective. Specifically, the two free generators originate at infinity and vanish at the North Pole and the South Pole, respectively.
\end{example}

\subsection{Topological visualizations of equivariant persistent cohomology}

First, we recall the duality between the persistence diagrams of persistent homology and persistent cohomology.
Let $\mathcal{X}: (\mathbb{R}, \le) \to \mathbf{Top}$ be a persistence topological space. Applying the homology and cohomology functors yields two distinct types of persistence modules. The homology $H_k(\mathcal{X})$ constitutes a persistence module, where the structure maps $\iota_{s,t}: H_k(\mathcal{X}_s) \to H_k(\mathcal{X}_t)$ are induced by the forward inclusions $\mathcal{X}_s \hookrightarrow \mathcal{X}_t$. 
Conversely, the cohomology $H^k(\mathcal{X})$ forms a copersistence module via the restriction maps $r_{t,s}: H^k(\mathcal{X}_t) \to H^k(\mathcal{X}_s)$.

Assuming field coefficients $\mathbb{K}$, the Universal Coefficient Theorem provides a natural isomorphism $H^k(\mathcal{X}_t) \cong \mathrm{Hom}(H_k(\mathcal{X}_t),\mathbb{K})$ for each parameter $t$. 
This establishes a categorical duality: the persistent cohomology module is isomorphic to the vector space dual of the persistent homology module. This duality induces a strict geometric symmetry in their respective persistence diagrams.

\begin{definition}
Let $\mathcal{D}_k^\text{hom}$ be the persistence diagram for the homology module $H_k(K_\bullet)$. Each topological feature is represented by a point $(b,d) \in \mathbb{R}^2$ with $b < d$, corresponding to an interval decomposition $[b, d)$. Geometrically, all points in $\mathcal{D}_k^\text{hom}$ lie strictly above the diagonal $\Delta = \{ (x,y) \in \mathbb{R}^2 \mid x=y \}$.

The persistence diagram for cohomology, $\mathcal{D}_k^\text{coh}$, is obtained by reflecting $\mathcal{D}_k^\text{hom}$ across the diagonal $\Delta$. Formally, for every feature $(b, d) \in \mathcal{D}_k^\text{hom}$, there exists a dual feature in the cohomology diagram:
\[
    (b, d) \in \mathcal{D}_k^\text{hom} \longleftrightarrow (d, b) \in \mathcal{D}_k^\text{coh}.
\]
Consequently, points in $\mathcal{D}_k^\text{coh}$ satisfy $d < b$ and lie strictly below the diagonal. The persistence length $\ell = |d - b|$ remains invariant under this transformation.
\end{definition}

\begin{figure}[htbp]
    \centering
    \begin{minipage}{0.48\textwidth}
        \centering
        \begin{tikzpicture}[scale=0.6]
            \draw[->] (0,0) node[below left] {$0$} -- (6,0) node[right] {$b$ (birth)};
            \draw[->] (0,0) -- (0,6) node[above] {$d$ (death)};
            
            \draw[dashed, gray] (0,0) -- (5.5,5.5) node[above right, black] {$d=b$};
            
            \fill[mainblue] (2,5) circle (2.5pt);
            \node[mainblue, above] at (2,5) {$(b,d)=(2,5)$};
            \node[font=\small] at (1.2, 3.5) {$H_k$};
            
            \fill[symred] (5,2) circle (2.5pt);
            \node[symred, below right] at (5,2) {$(d,b)=(5,2)$};
            \node[font=\small] at (3.5, 1.2) {$H^k$};
            
            \draw[dotted, gray] (2,5) -- (5,2);
        \end{tikzpicture}
    \end{minipage}
    \hfill
    \begin{minipage}{0.48\textwidth}
        \centering
        \begin{tikzpicture}[scale=0.9]
            \draw[->] (0,0) -- (6,0) node[right] {$t$};
            \draw[->] (0,-0.5) -- (0,3) node[above] {Bars};
            
            \draw (2, 0.1) -- (2, -0.1) node[below] {$2$};
            \draw (5, 0.1) -- (5, -0.1) node[below] {$5$};
            
            \draw[mainblue, line width=4pt] (2, 2) -- (5, 2);
            \node[black, left] at (0, 2) {$H_k$};
            

            \draw[symred, line width=4pt] (2, 1) -- (5, 1);
            \node[black, left] at (0, 1) {$H^k$};
            
            \draw[dotted, gray] (2, 0) -- (2, 2.5);
            \draw[dotted, gray] (5, 0) -- (5, 2.5);
            
        \end{tikzpicture}
    \end{minipage}
    \vspace{-0.5em}
    \caption{While the points in the persistence diagram are reflected across the diagonal, the barcodes for $H_k$ and $H^k$ share the identical interval $[b, d) = [2, 5)$.}
\end{figure}

Let $\mathcal{X} \colon (\mathbb{R},\leq) \to \mathbf{Top}_G$ be a persistence $G$-space. Its equivariant persistent cohomology
\[
H_G^\ast(\mathcal{X}) \colon (\mathbb{R},\leq)^{\mathrm{op}} \to \mathbf{Mod}_{\mathcal{H}_G^\ast}^{\mathrm{gr}}
\]
sends each real parameter $t\in\mathbb{R}$ to the graded $\mathcal{H}_G^\ast$-module $H_G^\ast\big(\mathcal{X}(t)\big)$. 

To understand its decomposition, we first recall the classical interval module over a field $\mathbb{K}$, defined as $\mathbb{I}[b, d) \colon (\mathbb{R},\leq)^{\mathrm{op}} \to \mathbf{Vect}_{\mathbb{K}}$, which is isomorphic to $\mathbb{K}$ for $t \in [b, d)$ and zero otherwise, with identities as restriction maps. A classical persistence module over $\mathbb{K}$ decomposes canonically as a finite direct sum of such interval modules.

However, the $\mathcal{H}_G^\ast$-module $H_G^\ast\big(\mathcal{X}(t)\big)$ can be structurally complex in general. We focus on the tractable cases where $G=S^1$ or $G=\mathbb{Z}/p$, for which $\mathcal{H}_G^\ast$ is a Principal Ideal Domain (PID).  
Note that while the slice $H_G^\ast\big(\mathcal{X}(t)\big)$ is a module over a PID at each fixed $t$, the persistent module $H_G^\ast(\mathcal{X})$ is a functor whose structure is governed by the monotonicity of annihilator ideals.
Specifically, for any $s \le t$, the restriction maps satisfy $\operatorname{Ann}(x_t) \subseteq \operatorname{Ann}(x_s)$. 
This implies that the algebraic type of a generator can only \textit{degenerate}: a free element may become torsion, or a higher-order torsion element may degenerate into a lower-order one; the reverse process is impossible.

Consequently, the structure theorem yields a decomposition into indecomposable persistent modules that capture this degradation history:
\[
H_G^\ast(\mathcal{X}) \cong 
\underbrace{\left( \bigoplus_{i} \mathbb{I}[b_i, d_i) \otimes \mathcal{H}_G^\ast \right)}_{\text{Free Intervals}} 
\oplus 
\underbrace{\left( \bigoplus_{j} \mathbb{I}[b_j, d_j) \otimes \frac{\mathcal{H}_G^\ast}{(p_j^{n_j})} \right)}_{\text{Constant Torsion Intervals}}
\oplus 
\underbrace{\left( \bigoplus_{l} \mathcal{T}_l \right)}_{\text{Transition Intervals}}.
\]
The summands are classified as follows:
\begin{enumerate}[label=($\roman*$)]
    \item \textit{Free Intervals}: The generator remains free throughout its entire lifetime $[b_k, d_k)$, faithfully encoding the group symmetry.
    \item \textit{Constant Torsion Intervals}: The generator is born as torsion and maintains a constant torsion order $n_j$ throughout its lifetime.
    \item \textit{Transition Intervals} ($\mathcal{T}_l$): These capture the phenomenon of symmetry degradation. A generator may start for free (or high-order torsion) at birth $b$ and degenerate into torsion (or lower-order torsion) at a transition time $t_{\text{trans}}$.
\end{enumerate}
Thus, the barcode of equivariant persistence provides a stratification of features based not only on their lifetime but also on the evolution of their symmetry properties.

To visualize the decomposition, we propose a \textit{stratified barcode} representation. Unlike the classical barcode which stacks all bars on a single axis, we separate bars into horizontal lanes based on their algebraic type:
\begin{itemize}
    \item The free lane contains bars for summands of type $\mathbb{I}[b, d) \otimes \mathcal{H}_G^\ast$.
    \item The torsion lane contains bars for constant torsion summands.
\end{itemize}
Transition modules, corresponding to symmetry degradation, are represented by \textit{stepped bars} that originate in the free Lane and drop to the torsion lane at the transition time $t_{\text{trans}}$. This vertical movement visually encodes the monotonic growth of the annihilator ideal and the loss of symmetry information over time.

\begin{figure}[h]
    \centering
    \begin{tikzpicture}[
        scale=0.9,
        >=Stealth,
        bar/.style={line width=4pt}, 
        free/.style={bar, blue!60!black},
        torsion/.style={bar, red!70!black},
        transition/.style={bar, purple!60!black},
        label/.style={font=\small, inner sep=2pt}
    ]
    
    \fill[blue!5] (0, 2) rectangle (11, 4.5);
    \fill[red!5] (0, 0) rectangle (11, 2);

    \draw[thick] (0, 2) -- (11.5, 2);

    \node[rotate=90] at (-0.8, 3.25) {Free};
    \node[rotate=90] at (-0.8, 1) {Torsion};

    \draw[->, thick] (0, -0.1) -- (0, 5);
    \draw[->, thick] (0, 0) -- (11.5, 0) node[right] {$t$};
    
    \foreach \x in {0, 1, 2, 3, 4, 5} {
        \draw (\x*2, 0.1) -- (\x*2, -0.1) node[below] {$\x$};
    }

    
    \draw[free] (0, 4) -- (10, 4);
    \node[label, above] at (3, 4.1) {Feature A (Global Free)};

    \draw[transition] (2, 1.2) -- (6, 1.2); 
    \draw[black!60, line width=1.5pt, ->] (6, 3) -- (6, 1.2); 
    \draw[transition] (6, 3) -- (8, 3); 
    \node[label, above] at (4, 1.3) {Feature B};

    \draw[torsion] (4, 0.6) -- (8, 0.6);
    \node[label, below] at (6, 0.5) {Feature C};
    \end{tikzpicture}
    \caption{Illustration of stratified barcode for equivariant persistence. The vertical drop indicates the transition from free to torsion.}
\end{figure}

\subsection{Equivariant Vietoris-Rips complexes}

We now extend the discussion to simplicial complexes generated by metric data, focusing on the Vietoris-Rips construction under group actions.

\begin{proposition}\label{proposition:equivariant_VR}
    Let $X$ be a finite subset of Euclidean space $E$, and let $G\leq \operatorname{Iso}(E)$ be an isometry subgroup acting on $X$. For any fixed radius $r\geq 0$, the Vietoris-Rips complex $\operatorname{VR}_r(X)$ admits a canonical simplicial $G$-action.
\end{proposition}

\begin{proof}
    Define the action vertex-wise by $g\cdot x = g(x)$ for all $g\in G$ and $x\in X$. 
    Let $\sigma = \{x_0, \dots, x_k\}$ be a simplex in $\operatorname{VR}_r(X)$; by definition, the pairwise distances satisfy $\|x_i - x_j\| \leq r$. 
    Since $G$ acts by isometries, for the image simplex $g\cdot \sigma = \{g(x_0), \dots, g(x_k)\}$, we have
    \[
        \|g(x_i) - g(x_j)\| = \|x_i - x_j\| \leq r.
    \]
    Thus, $g\cdot \sigma$ remains a simplex in $\operatorname{VR}_r(X)$. 
    This assignment respects face inclusions and group composition, defining a simplicial $G$-action.
\end{proof}

\begin{proposition}
Let $X$ be a finite subset of Euclidean space $E$, and let $G\leq \operatorname{Iso}(E)$ be an isometry subgroup acting on $X$. The family of Vietoris–Rips complexes $\{\operatorname{VR}_t(X)\}_{t\in\mathbb{R}}$ constitutes a $G$-filtered simplicial complex.
\end{proposition}

\begin{proof}
By Proposition \ref{proposition:equivariant_VR}, for every parameter $t \in \mathbb{R}$, the complex $\operatorname{VR}_t(X)$ admits a simplicial $G$-action defined vertex-wise by $g \cdot x = g(x)$.

To verify the $G$-filtered structure, it remains to show that the inclusion maps are equivariant. For any $s \le t$, consider the inclusion $\iota_{s,t}: \operatorname{VR}_s(X) \hookrightarrow \operatorname{VR}_t(X)$. For any simplex $\sigma \in \operatorname{VR}_s(X)$ and any $g \in G$, we have
\[
    \iota_{s,t}(g \cdot \sigma) = g \cdot \sigma = g \cdot \iota_{s,t}(\sigma).
\]
Since $\iota_{s,t}$ commutes with the group action, the family $\{\operatorname{VR}_t(X)\}_t$ is a $G$-filtered space.
\end{proof}

Applying the Borel construction to the filtration yields the topological summary of the data's symmetry.

\begin{definition}
    Let $X$ be a finite $G$-metric space. The \textit{equivariant Vietoris-Rips persistence module} is the functor
    \[
        H_G^\ast(\operatorname{VR}_\bullet(X)) : (\mathbb{R}, \leq)^{\mathrm{op}} \longrightarrow \mathbf{Mod}_{\mathcal{H}_G^\ast}^{\mathrm{gr}}.
    \]
    For each parameter $t$, the module is $H_G^\ast(\operatorname{VR}_t(X))$, and the structure maps are induced by the equivariant inclusions $\operatorname{VR}_s(X) \hookrightarrow \operatorname{VR}_t(X)$ for $s \leq t$.
\end{definition}

The algebraic structure of this persistence module captures the ``birth'' or ``disappearance'' of symmetric simplices.

A key distinction from non-equivariant persistence is the role of torsion. While $H^\ast(\operatorname{VR}_t(X))$ is typically a vector space over $\mathbb{K}$ (thus characterized by Betti numbers), $H_G^\ast(\operatorname{VR}_t(X))$ is a module over $\mathcal{H}_G^\ast$. If $G$ is finite and $|\mathbb{K}|$ divides $|G|$, or if $G$ is a continuous group like $S^1$, the module structure may exhibit torsion or complex relations derived from the characteristic classes of the group, reflecting the non-trivial equivariant topology of the data cloud.

\begin{example}\label{example:Z2_action}
Consider the finite group $G = \mathbb{Z}/2\mathbb{Z} = \{e, \tau\}$ acting on the Euclidean line $E = \mathbb{R}$ by reflection about the origin: $\tau \cdot x = -x$.
Let the point cloud $X = \{-1, 1\} \subset \mathbb{R}$.
The group acts freely on $X$: $\tau \cdot (-1) = 1$ and $\tau \cdot (1) = -1$.
Assume the coefficient field is $\mathbb{K} = \mathbb{F}_2$. The Borel cohomology ring is $\mathcal{H}_G^\ast \cong H^\ast(B\mathbb{Z}/2; \mathbb{F}_2) \cong \mathbb{F}_2[c]$, where $c$ is a degree 1 generator.

We analyze the Vietoris-Rips filtration $\operatorname{VR}_r(X)$ as follows:
\begin{itemize}
    \item When $0 \le r < 2$, the complex consists of two discrete vertices $X$. Since the action is free, the Borel construction $X_G$ is homotopy equivalent to the quotient $X/G$, which is a single point. Thus, $H_G^\ast(X) \cong H^\ast(\mathrm{pt}) \cong \mathbb{F}_2$. As an $\mathcal{H}_G^\ast$-module, this is the quotient $\mathbb{F}_2[c]/(c) \cong \mathbb{F}_2$. \textit{This is a torsion module}: the generator $c$ acts as zero.
        
    \item When $r \ge 2$, an edge connects $\{-1\}$ and $\{1\}$, forming a line segment $\Delta^1$. This segment contains the origin $0$ (in its interior), which is a $G$-fixed point. Since a $G$-space containing a fixed point is $G$-contractible, we have
        \[
            H_G^\ast(\Delta^1) \cong H_G^\ast(\mathrm{pt}) \cong \mathbb{F}_2[c].
        \]
        \textit{This is a free module} of rank 1.
\end{itemize}

\begin{figure}[h]
    \centering
    \begin{tikzpicture}[
        scale=0.9,
        >=Stealth,
        bar/.style={line width=4pt}, 
        free/.style={bar, blue!60!black},
        torsion/.style={bar, red!70!black},
        transition/.style={bar, purple!60!black},
        label/.style={font=\small, inner sep=2pt}
    ]
    
    \fill[blue!5] (0, 2) rectangle (11, 4.5);
    \fill[red!5] (0, 0) rectangle (11, 2);

    \draw[thick] (0, 2) -- (11.5, 2);

    \node[rotate=90] at (-0.8, 3.25) {Free};
    \node[rotate=90] at (-0.8, 1) {Torsion};

    \draw[->, thick] (0, -0.1) -- (0, 5);
    \draw[->, thick] (0, 0) -- (11.5, 0) node[right] {$r$};
    
    \node[below] at (0, -0.1) {$0$};
    \node[below] at (2, -0.1) {$2$};
    \draw (2, 0.1) -- (2, -0.1);

    
    \draw[transition] (0, 1) -- (2, 1);
    \node[label, below] at (1, 0.9) {$\mathbb{F}_2[c]/(c)$};
    
    \draw[black!60, line width=1.5pt, ->] (2, 3.5) -- (2, 1);
    
    \draw[transition] (2, 3.5) -- (10, 3.5);
    \node[label, above] at (6, 3.6) {$\mathbb{F}_2[c]$};
    \end{tikzpicture}
    \caption{Stratified barcode for equivariant cohomology in Example \ref{example:Z2_action}. The bar is born in the free lane from infinity and jumps to the torsion lane at $r=2$ when the fixed point disappears.}
\end{figure}

The persistence structure map $H_G^\ast(\operatorname{VR}_{\ge 2}(X)) \to H_G^\ast(\operatorname{VR}_{< 2}(X))$ corresponds to the projection
\[
    \mathbb{F}_2[c] \longrightarrow \mathbb{F}_2[c]/(c).
\]
In the early stage, the topological feature is ``blind'' to the group's characteristic class $c$, resulting in a torsion module where $c$ vanishes. As the complex grows to include a fixed point, the module structure ``upgrades'' from torsion to free. The characteristic class $c$, which was trivial for $r<2$, becomes a non-trivial persistent element for $r\geq 2$, reflecting the transition from a free orbit to a trivial isotropy type.
\end{example}

\subsection{Stability for the equivariant persistent cohomology}

\begin{definition}
Let $\mathcal{X}, \mathcal{Y}: (P, \leq) \longrightarrow \mathbf{Top}_G$ be two persistence $G$-spaces. A \textit{morphism of persistence spaces} $f: \mathcal{X} \to \mathcal{Y}$ consists of a family of $G$-equivariant continuous maps $\{f_t: \mathcal{X}(t) \to \mathcal{Y}(t)\}_{t \in P}$ that are compatible with the persistence structure. Specifically, for every pair $s \leq t$ in $P$, the following diagram commutes:
\[
\begin{tikzcd}
\mathcal{X}(s) \arrow[r, "\mathcal{X}(s \leq t)"] \arrow[d, "f_s"'] & \mathcal{X}(t) \arrow[d, "f_t"] \\
\mathcal{Y}(s) \arrow[r, "\mathcal{Y}(s \leq t)"'] & \mathcal{Y}(t)
\end{tikzcd}
\]
In other words, $f$ is a natural transformation between the functors $\mathcal{X}$ and $\mathcal{Y}$.
\end{definition}

Let $\mathbf{Pers}\text{-}\mathbf{Top}_G$ denote the category of persistence $G$-spaces and $\mathbf{Pers}_G$ the category of graded persistence $\mathcal{H}_G^\ast$-modules. 

\begin{proposition}
The equivariant Borel cohomology construction defines a contravariant functor
\[
H_G^\ast: \mathbf{Pers}\text{-}\mathbf{Top}_G \longrightarrow \mathbf{Pers}_G.
\]
\end{proposition}
\begin{proof}
For a persistence $G$-space $\mathcal{X}$, the object $H_G^\ast(\mathcal{X})$ is defined pointwise as the diagram of modules $\{H_G^\ast(\mathcal{X}(t))\}_{t \in P}$. The structure maps are the pullbacks induced by $\mathcal{X}(s \leq t)$, and the persistence module axioms follow directly from the functoriality of $\mathcal{X}$ and the contravariance of $H_G^\ast$.

Given a morphism $f: \mathcal{X} \to \mathcal{Y}$, we define $H_G^\ast(f): H_G^\ast(\mathcal{Y}) \to H_G^\ast(\mathcal{X})$ component-wise via $f_t^\ast$. To see that this is a morphism of persistence modules, let $s \leq t$. The commutativity of the original diagram of spaces, $f_t \circ \mathcal{X}(s \leq t) = \mathcal{Y}(s \leq t) \circ f_s$, implies the following diagram commutes in the category of modules:
\[
\begin{tikzcd}
H_G^\ast(\mathcal{Y}(t)) \arrow[r, "f_t^\ast"] \arrow[d, "\mathcal{Y}(s \leq t)^\ast"'] & H_G^\ast(\mathcal{X}(t)) \arrow[d, "\mathcal{X}(s \leq t)^\ast"] \\
H_G^\ast(\mathcal{Y}(s)) \arrow[r, "f_s^\ast"'] & H_G^\ast(\mathcal{X}(s)).
\end{tikzcd}
\]
Finally, since Borel cohomology is contravariant on $\mathbf{Top}_G$, it preserves identities and reverses compositions:
\[
H_G^\ast(\mathrm{id}_{\mathcal{X}}) = \mathrm{id}, \quad H_G^\ast(g \circ f) = H_G^\ast(f) \circ H_G^\ast(g).
\]
Thus, $H_G^\ast$ is a contravariant functor.
\end{proof}

The interleaving distance is defined for both persistence $G$-spaces and graded $\mathcal{H}_G^\ast$-modules as a metric between persistence objects.

\begin{definition}
Let $\mathfrak{C} = \mathbf{Top}_G$ be the category of $G$-topological spaces.
For two persistence $G$-spaces $\mathcal{X}, \mathcal{Y}: (\mathbb{R}, \leq) \to \mathbf{Top}_G$, their \textit{interleaving distance} is defined as
\[
d_I(\mathcal{X}, \mathcal{Y}) = \inf \left\{ \varepsilon \geq 0 \ \middle|\ \text{$\mathcal{X}$ and $\mathcal{Y}$ are $\varepsilon$-interleaved in $\mathfrak{C}^\mathbb{R}$} \right\}.
\]
\end{definition}

\begin{definition}
Let $\mathfrak{C} = \mathbf{Mod}_{\mathcal{H}_G^\ast}^{\mathrm{gr}}$ be the category of graded $\mathcal{H}_G^\ast$-modules.
For two persistence modules $\mathcal{M}, \mathcal{N}: (\mathbb{R}, \leq) \to \mathbf{Mod}_{\mathcal{H}_G^\ast}^{\mathrm{gr}}$, their \textit{interleaving distance} is defined as
\[
d_I(\mathcal{M}, \mathcal{N}) = \inf \left\{ \varepsilon \geq 0 \ \middle|\ \text{$\mathcal{M}$ and $\mathcal{N}$ are $\varepsilon$-interleaved in $\mathfrak{C}^\mathbb{R}$} \right\}.
\]
\end{definition}

\begin{proposition}
The equivariant cohomology functor $H_G^\ast: \mathbf{Pers}\text{-}\mathbf{Top}_G \to \mathbf{Pers}_G$ is $1$-Lipschitz with respect to the interleaving distances. Specifically, for any two persistence $G$-spaces $\mathcal{X}$ and $\mathcal{Y}$, we have
\[
d_I\bigl(H_G^\ast(\mathcal{X}), H_G^\ast(\mathcal{Y})\bigr) \leq d_I(\mathcal{X}, \mathcal{Y}).
\]
\end{proposition}

\begin{proof}
If $\mathcal{X}$ and $\mathcal{Y}$ are $\varepsilon$-interleaved, there exist natural transformations $\phi: \mathcal{X} \Rightarrow \Sigma^\varepsilon \mathcal{Y}$ and $\psi: \mathcal{Y} \Rightarrow \Sigma^\varepsilon \mathcal{X}$ satisfying the triangular identities.
Since $H_G^\ast$ is a functor, it preserves the shift operation (i.e., $H_G^\ast(\Sigma^\varepsilon \mathcal{X}) \cong \Sigma^\varepsilon H_G^\ast(\mathcal{X})$) and sends natural transformations to natural transformations.
Thus, the induced maps $H_G^\ast(\phi)$ and $H_G^\ast(\psi)$ constitute an $\varepsilon$-interleaving between $H_G^\ast(\mathcal{X})$ and $H_G^\ast(\mathcal{Y})$.
Taking the infimum over all such $\varepsilon$ yields the inequality.
\end{proof}

Let $X$ be a $G$-space, and suppose $f\colon X\to \mathbb{R}$ is a $G$-invariant real-valued function, meaning $f(g\cdot x)=f(x)$ for all $g\in G$ and $x\in X$. The sublevel sets
\[
X^f_t = f^{-1}\big((-\infty,t]\big),\quad t\in\mathbb{R},
\]
form a persistence $G$-space $\mathcal{X}_f: (\mathbb{R}, \leq) \to \mathbf{Top}_G$.

\begin{theorem}
Let $f, g: X \to \mathbb{R}$ be two $G$-invariant functions on a $G$-space $X$. Let $\mathcal{X}_f$ and $\mathcal{X}_g$ denote the corresponding persistence $G$-spaces defined by sublevel sets.
Then we have
\[
d_I\bigl(H_G^\ast(\mathcal{X}_f), H_G^\ast(\mathcal{X}_g)\bigr) \leq \| f - g \|_\infty.
\]
Specifically, if $\| f - g \|_\infty \leq \varepsilon$, then the persistence modules $H_G^\ast(\mathcal{X}_f)$ and $H_G^\ast(\mathcal{X}_g)$ are $\varepsilon$-interleaved.
\end{theorem}

\begin{proof}
Let $\varepsilon = \| f - g \|_\infty$. For any $t \in \mathbb{R}$, we have the inclusion of sublevel sets
\[
X^f_t = f^{-1}\big((-\infty, t]\big) \subseteq g^{-1}\big((-\infty, t+\varepsilon]\big) = X^g_{t+\varepsilon}.
\]
This inclusion holds because if $f(x) \leq t$, then $g(x) \leq f(x) + \varepsilon \leq t + \varepsilon$. Since $f$ and $g$ are $G$-invariant, these inclusion maps are $G$-equivariant. Collecting these maps for all $t$ defines a natural transformation $\phi: \mathcal{X}_f \Rightarrow \Sigma^\varepsilon \mathcal{X}_g$. Similarly, we have a natural transformation $\psi: \mathcal{X}_g \Rightarrow \Sigma^\varepsilon \mathcal{X}_f$. These natural transformations satisfy the commutativity requirements of an $\varepsilon$-interleaving in the category of persistence $G$-spaces. Applying the contravariant functor $H_G^\ast$ yields an $\varepsilon$-interleaving between $H_G^\ast(\mathcal{X}_f)$ and $H_G^\ast(\mathcal{X}_g)$, which completes the proof.
\end{proof}

\subsection{Finite computation of equivariant cohomology}

While conceptually elegant, the Borel construction is computationally difficult for TDA. The space $EG$ is typically infinite-dimensional; for instance, $ES^1 \simeq S^\infty$. Directly computing the homology of $X_G$ involves managing infinite-dimensional chain complexes, which is infeasible for large datasets like point clouds. To overcome the limitations, we turn to Bredon cohomology.

\begin{definition}
The \textit{orbit category} $\mathrm{Or}_G$ has as objects the left cosets $G/H$ for all subgroups $H \le G$, and as morphisms all $G$-equivariant maps between these coset spaces.
\end{definition}

\begin{definition}
A \textit{coefficient system} $\underline{M}$ is a contravariant functor from $\mathrm{Or}_G$ to the category of abelian groups (or graded $R$-modules).
\end{definition}

For the computation of Borel cohomology, we employ the \textit{Borel coefficient system} $\underline{M}_{\mathrm{Borel}}$, defined by
\[
\underline{M}_{\mathrm{Borel}}(G/H) = H^*(BH; R),
\]
where $BH$ denotes the classifying space of the subgroup $H\le G$, and $R$ is the base ring. The contravariant functoriality is given by the restriction maps in group cohomology induced by subgroup inclusions.

For a regular $G$-simplicial complex $K$, the Bredon cochain complex is indexed by the simplices of the orbit space $K/G$, with coefficients determined by the stabilizer subgroups.

\begin{definition}
Let $\Sigma_n(K/G)$ denote the set of $n$-simplices in the orbit space. The \textit{Bredon cochain groups} are defined as
\[
C^n_{\mathrm{br}}(K; \underline{M}) = \prod_{\bar{\sigma} \in \Sigma_n(K/G)} \underline{M}(G_{\sigma}),
\]
where $G_{\sigma}$ is the stabilizer of a representative simplex $\sigma$ of the orbit $\bar{\sigma}$.

For any $(n+1)$-simplex $\bar{\sigma}$ and its $n$-dimensional face $\bar{\tau}$, the regularity of the $G$-complex implies the inclusion of stabilizers $G_\sigma \le G_\tau$. The coboundary $\delta^n: C^n_\mathrm{br} \to C^{n+1}_\mathrm{br}$ is defined via the restriction maps $\mathrm{Res}^{G_\tau}_{G_\sigma}: \underline{M}(G_\tau) \to \underline{M}(G_\sigma)$. More precisely, for $\alpha \in C^n_{\mathrm{br}}(K; \underline{M})$, the coboundary evaluated on $\bar{\sigma}$ is given by
\[
\delta^n(\alpha)(\bar{\sigma}) = \sum_{\substack{\bar{\tau} \subset \bar{\sigma} \\ \dim\bar{\tau}=n}} [\bar{\sigma}:\bar{\tau}] \, \mathrm{Res}^{G_\tau}_{G_\sigma}\bigl(\alpha(\bar{\tau})\bigr),
\]
where $[\bar{\sigma}:\bar{\tau}]$ denotes the incidence number (sign) following the standard simplicial cochain convention.
\end{definition}

\begin{definition}
The cohomology of the cochain complex $(C^\ast_{\mathrm{br}}, \delta^\ast)$ is called the \textit{Bredon cohomology} of $K$ with coefficients in $\underline{M}$, denoted $H^*_{\mathrm{br}}(K; \underline{M})$.

When $\underline{M}$ is a graded coefficient system (such as $\underline{M}_{\mathrm{Borel}}$), the complex is bigraded by simplicial degree $p$ and internal degree $q$, and the total Bredon cohomology is graded by total degree $k = p+q$.
\end{definition}

\begin{theorem}[{\cite[Chapter IV, Section 2]{may1996equivariant}}]
Let $G$ be a finite group and $K$ be a finite-dimensional regular $G$-simplicial complex. There is a natural isomorphism of graded $R$-modules
\[
H^*_{\mathrm{br}}(K; \underline{M}_{\mathrm{Borel}}) \cong H^*_G(K; R),
\]
where the left-hand side is taken in total degree.
\end{theorem}
This theorem implies that computing the equivariant cohomology of $K$ reduces to computing the cohomology of the finite quotient complex $K/G$ enriched with coefficients from the group cohomology of the stabilizers.

\begin{remark}
The sum in the definition of $C^n_\mathrm{br}(K)$ is taken only over the orbit types that actually appear in $K$. If a subgroup $H \le G$ is not a stabilizer of any simplex in $K$, it does not contribute to the complex. In practical TDA scenarios, the set of occurring stabilizers is usually very small.
\end{remark}

\begin{example}
Consider the cyclic group $G = \mathbb{Z}_2 = \{e, \tau\}$ acting simplicially on the $2$-simplex $\sigma = [x_0, x_1, x_2]$. Let the action be defined such that the generator $\tau$ fixes the vertex $x_0$ while transposing the vertices $x_1$ and $x_2$. This action is not regular: the simplex $\sigma$ is invariant under $\tau$ (implying $G_\sigma = \mathbb{Z}_2$), yet it contains vertices with trivial stabilizers ($G_{x_1} = G_{x_2} = \{e\}$). Consequently, the face stabilizer relation $G_\sigma \le G_{\text{face}}$ fails, preventing the direct application of the standard Bredon complex definition.

To obtain a regular $G$-complex, we subdivide the irregular edge $[x_1, x_2]$ by inserting its midpoint $m$. Since the resulting two $2$-simplices have trivial stabilizers, the complex is now regular. The orbit space $K'/G$ has the following simplicial structure:
\begin{itemize}
    \item \textit{Vertices:} The orbit space has three vertices: $v_0$, the fixed vertex $x_0$ with stabilizer $G_{v_0} = \mathbb{Z}_2$; $m$, the midpoint of the edge $[x_1, x_2]$ which is fixed pointwise with stabilizer $G_m = \mathbb{Z}_2$; and $v_1$, representing the free orbit of $\{x_1, x_2\}$ with trivial stabilizer $G_{v_1} = \{e\}$.
    \item \textit{Edges:} There is one fixed edge $[v_0, m]$ with stabilizer $\mathbb{Z}_2$, and two free edges $[v_0, v_1]$ and $[v_1, m]$ with trivial stabilizers.
    \item \textit{2-Simplices:} There is a single $2$-simplex orbit $[v_0, v_1, m]$ with trivial stabilizer $\{e\}$.
\end{itemize}

Let $\underline{M}_{\mathrm{Borel}}$ be the Borel coefficient system with $R = \mathbb{F}_2$. We have $\underline{M}(\mathbb{Z}_2) \cong H^*(B\mathbb{Z}_2; \mathbb{F}_2) \cong \mathbb{F}_2[w]$ ($\deg w = 1$) and $\underline{M}(\{e\}) \cong H^*(pt; \mathbb{F}_2) \cong \mathbb{F}_2$ concentrated in degree $0$. The Bredon cochain groups $C_{\mathrm{br}}^p(K'; \underline{M})$ are given by
\begin{align*}
C^0_{\mathrm{br}} &\cong \mathbb{F}_2[w] \oplus \mathbb{F}_2[w] \oplus \mathbb{F}_2, \\
C^1_{\mathrm{br}} &\cong \mathbb{F}_2[w] \oplus \mathbb{F}_2 \oplus \mathbb{F}_2, \\
C^2_{\mathrm{br}} &\cong \mathbb{F}_2.
\end{align*}
The differentials are determined by the restriction maps. The map $\delta^0: C^0 \to C^1$ acts on the fixed components via the diagonal map $(p(w), q(w)) \mapsto p(w)+q(w)$, while the restriction to the free components extracts constant terms. Consequently, $\ker(\delta^0)$ consists of classes $(p, p, p(0))$ with $p \in \mathbb{F}_2[w]$, yielding $H^0_{\mathrm{br}} \cong \mathbb{F}_2[w]$. The map $\delta^1: C^1 \to C^2$ restricts the polynomial part via the augmentation map $p(w) \mapsto p(0)$ and sums the coefficients from the free edges. This map is surjective, resulting in $H^{\ge 1}_{\mathrm{br}} = 0$. We recover the Borel cohomology
\[
H_G^\ast(\sigma; \mathbb{F}_2) \cong H^\ast_{\mathrm{br}}(K'; \underline{M}) \cong \mathbb{F}_2[w].
\]
This calculation demonstrates the use of a finite cochain complex for computing equivariant cohomology.
\end{example}

Let $G$ be a finite group and $K$ be a finite $G$-simplicial complex. The barycentric subdivision $\mathrm{Sd}(K)$ yields a regular $G$-complex. Nevertheless, barycentric subdivision tends to produce a rather complicated combinatorial structure. To obtain a regular $G$-complex while minimizing the subdivision of simplices, we employ an inductive strategy on the dimension of the simplices.

\begin{definition}
A simplex $\sigma \in K$ is called \textit{irregular} if its stabilizer $G_\sigma$ acts non-trivially on its vertex set.
\end{definition}

The minimal regular subdivision $K_{\mathrm{reg}}$ is constructed iteratively:
\begin{enumerate}[label=($\roman*$)]
    \item The $0$-skeleton remains unchanged. Vertices are always regular.
    
    \item Assume the $(n-1)$-skeleton has been processed. For each $n$-simplex $\sigma$ in the current complex:
    \begin{itemize}
        \item If $\sigma$ is regular, leave $\sigma$ intact.
        \item If $\sigma$ is irregular, subdivide $\sigma$ by introducing a single new vertex $v_\sigma$ at its barycenter. Cone the boundary of $\sigma$ (which is already regular by the inductive hypothesis) to $v_\sigma$.
    \end{itemize}
    \item The new vertex $v_\sigma$ inherits the stabilizer $G_\sigma$. The stabilizers of the new $n$-simplices are determined by the intersection of stabilizers: $G_{[v_\sigma, \tau]} = G_\sigma \cap G_\tau$.
\end{enumerate}

\begin{proposition}
The iterative subdivision $K_{\mathrm{reg}}$ is a regular $G$-complex. Furthermore, $K_{\mathrm{reg}}$ is minimal in the sense that subdivision is performed only on simplices whose symmetry violates regularity.
\end{proposition}

\begin{proof}
The subdivision introduces a vertex $v_\sigma$ at the barycenter of $\sigma$, which is fixed pointwise by $G_\sigma$. Any new simplex $\tau'$ is of the form $[v_\sigma, \eta]$ where $\eta \subset \partial \sigma$ is a face of the boundary. By the inductive hypothesis, the action on $\partial \sigma$ is regular. 

If an element $g \in G_\sigma$ permutes the vertices of $\eta$, then $\tau'$ lies in a free orbit with trivial stabilizer, satisfying regularity trivially. If $g$ fixes $\eta$ pointwise, it fixes $\tau'$ pointwise. Consequently, for any $g \in G_{\tau'}$, $g$ fixes all vertices of $\tau'$, establishing the regularity of the subdivision. Minimality follows from the restriction of the operation to simplices with non-trivial stabilizer action on their vertices.
\end{proof}

In the case of the $\mathbb{Z}_2$ action on a triangle, the edge $[x_1, x_2]$ is irregular. Subdividing this edge (inserting a midpoint) forces the adjacent $2$-simplex to split into two free orbits. This restores regularity without requiring the addition of a barycenter for the $2$-simplex itself, exemplifying the minimality of this strategy.

Mathematical pseudocode for the minimal regular subdivision is given in Algorithm \ref{algorithm:regular_subdivision}.

\begin{algorithm}[H]
\caption{Minimal Regular Subdivision}\label{algorithm:regular_subdivision}
\begin{algorithmic}[1]
\Require{Finite $G$-simplicial complex $K$; Finite group $G$.}
\Ensure{Regular $G$-simplicial complex $K_{\mathrm{reg}}$.}

\State $K_{\mathrm{reg}} \leftarrow$ Vertex set of $K$ \Comment{Initialize with 0-skeleton}
\State $D \leftarrow \dim(K)$

\For{$n = 1$ \textbf{to} $D$}
    \State Let $\Sigma_n^{\mathrm{curr}}$ be the set of candidate $n$-simplices formed by the $(n-1)$-simplices in $K_{\mathrm{reg}}$.
    \For{\textbf{each} simplex $\sigma \in \Sigma_n^{\mathrm{curr}}$}
        \State Compute stabilizer $G_\sigma \leftarrow \{g \in G \mid g \cdot \sigma = \sigma\}$.
        \State \textbf{Check Regularity:} Does $G_\sigma$ fix all vertices of $\sigma$ pointwise?
        \If{$\sigma$ is regular} \Comment{Preserve simplex if regular}
            \State Add $\sigma$ to $K_{\mathrm{reg}}$.
        \Else \Comment{Subdivide if irregular}
            \State Introduce new vertex $v_\sigma$ at barycenter of $\sigma$.
            \State Add $v_\sigma$ to $K_{\mathrm{reg}}$.
            \For{\textbf{each} face $\eta \subset \partial \sigma$}
                \State Construct new simplex $[v_\sigma, \eta]$.
                \State Add $[v_\sigma, \eta]$ to $K_{\mathrm{reg}}$.
                \Comment{New simplices are regular}
            \EndFor
        \EndIf
    \EndFor
\EndFor
\end{algorithmic}
\end{algorithm}

In summary, the computation of Borel equivariant cohomology $H^*_G(K; R)$ is accomplished through a finite procedure that avoids the direct manipulation of infinite-dimensional classifying spaces. The method proceeds in three distinct stages. First, the input $G$-simplicial complex $K$ is transformed into a regular $G$-complex $K_{\mathrm{reg}}$ via a minimal subdivision algorithm. This regularization is essential to ensure that for any simplex and its face, the inclusion of stabilizer subgroups $G_\sigma \le G_\tau$ holds, thereby satisfying the structural requirement for the Bredon complex. Second, the Bredon cochain complex $C^*_{\mathrm{br}}(K_{\mathrm{reg}}; \underline{M}_{\mathrm{Borel}})$ is constructed. This complex is indexed by the simplices of the finite orbit space $K_{\mathrm{reg}}/G$, with coefficients given by the cohomology of the classifying spaces of the stabilizers, $H^*(BG_\sigma; R)$. The coboundary operator is explicitly defined via simplicial incidence numbers combined with restriction maps in group cohomology. Finally, the computation reduces to solving the linear algebra problem of finding the kernel and image of the coboundary operator. By the foundational theorem relating Bredon and Borel cohomology, the resulting cohomology groups are isomorphic to $H^*_G(K; R)$, yielding computationally tractable access to equivariant topological data analysis. We outline its computational procedure in Algorithm \ref{algorithm:equivariant_cohomology}.

\begin{algorithm}[H]
\caption{Computing Equivariant Cohomology}\label{algorithm:equivariant_cohomology}
\begin{algorithmic}[1]
\Require{Finite $G$-simplicial complex $K$; Finite group $G$; Coefficient ring $R$.}
\Ensure{$H^*_G(K; R)$ (graded $R$-module).}

\State \textbf{Phase 1: Minimal Regular Subdivision}
\State Build $K_{\mathrm{reg}}$ inductively on dimension.
\For{each simplex $\sigma$}
    \If{stabilizer $G_\sigma$ acts non-trivially on vertices}
        \State Subdivide $\sigma$ by coning $\partial \sigma$ to barycenter $v_\sigma$.
    \Else
        \State Retain $\sigma$.
    \EndIf
\EndFor

\State \textbf{Phase 2: Bredon Cochain Complex Construction}
\State Let $\underline{M}$ be the Borel coefficient system: $\underline{M}(H) = H^*(BH; R)$.
\State Initialize chain groups $C^n_{\mathrm{br}} \leftarrow 0$ for $n=0, \dots, D$.

\For{$n = 0$ \textbf{to} $D$}
    \State Let $\mathcal{O}_n$ be the set of orbits of $n$-simplices in $K_{\mathrm{reg}}/G$.
    \State $C^n_{\mathrm{br}} \leftarrow \prod_{\bar{\sigma} \in \mathcal{O}_n} \underline{M}(G_\sigma)$.
\EndFor
\For{\textbf{each} orbit representative $\bar{\sigma} \in \mathcal{O}_{n+1}$}
    \For{\textbf{each} face $\bar{\tau} \subset \bar{\sigma}$ with $\dim(\bar{\tau})=n$}
        \State Let $\mathrm{Res}^{G_\tau}_{G_\sigma}: \underline{M}(G_\tau) \to \underline{M}(G_\sigma)$ be the restriction map.
        \State Component of $\delta^n$: $[\bar{\sigma}:\bar{\tau}] \cdot \mathrm{Res}^{G_\tau}_{G_\sigma}(\alpha(\bar{\tau}))$.
    \EndFor
\EndFor

\State \textbf{Phase 3: Cohomology Computation}
\State Compute kernel $Z^n \leftarrow \ker(\delta^n)$ and image $B^n \leftarrow \mathrm{im}(\delta^{n-1})$.
\State \Return $H^*_{\mathrm{br}}(K; \underline{M}) \cong Z^n / B^n$.
\end{algorithmic}
\end{algorithm} 

%% file: computation_aspects.tex
\section{Computation on symmetries}\label{section:computation}

In the preceding sections, we developed theoretical aspects for understanding symmetries of finite configurations. Building on this foundation, we now turn to the computational aspects of symmetry, motivated by practical applications. Our focus will include the computation and representation of symmetries in low-dimensional spaces, the algorithmic determination of symmetry groups, the construction and analysis of persistent symmetries and their associated barcodes, the evaluation of symmetry defects, as well as the formulation of local symmetries.

\subsection{Symmetry patterns in low-dimensional spaces}

In the Euclidean plane, a finite configuration may exhibit both mirror symmetry and rotational symmetry. Let $\bar{x}$ denote the centroid of the configuration. A reflection symmetry can be represented by the pair $(\bar{x}, \alpha)$, where $\alpha$ is the angle between the reflection axis and the $x$-axis. The angle $\alpha$ can take any value in the interval $[0^\circ, 180^\circ)$.
For rotational symmetry, let $\theta$ denote the angle of rotation. Then a rotational symmetry is represented by the pair $(\bar{x}, \theta)$, where the rotation is centered at $\bar{x}$ and rotates the configuration counterclockwise by angle $\theta$. In this case, $\theta$ must be of the form $\theta = \frac{360^\circ}{m}$ for some integer $m \geq 1$. In particular, when $m=1$, the identity isometry is the only rotational symmetry in $\Sym(X)$.

By the classification of finite subgroups of the orthogonal group $O(2)$, we have the following result \cite{armstrong1997groups}.

\begin{proposition}
Let $X$ be a finite subset of the Euclidean plane $\mathbb{R}^2$. Then the symmetry group $\Sym(X)$ is isomorphic to either a finite cyclic group or a dihedral group. That is,
\[
\Sym(X) \cong C_m \quad \text{or} \quad \Sym(X) \cong D_m
\]
for some integer $m \geq 1$, where $C_m$ denotes the cyclic group of order $m$, and $D_m$ denotes the dihedral group of order $2m$.
\end{proposition}

We note that the dihedral group $D_m$ is generated by a rotation of order $m$ and a single reflection. Thus, the rotational symmetries in $\Sym(X)$ are generated by a single rotation $\rho_\theta$ about the centroid $\bar{x}$, where $\theta = \frac{360^\circ}{m}$. The mirror symmetries in $\Sym(X)$ are determined by reflections across lines passing through $\bar{x}$ and forming specific angles with the $x$-axis.

Therefore, the symmetry group of a finite configuration $X$ in the plane can be characterized by a triple $(\bar{x}, m, \{\alpha_i\})$, where $\bar{x}$ is the centroid of the point configuration, and $m$ is the order of rotational symmetry. In particular, $\Sym(X)$ contains a rotation $\rho_{\theta}$ about $\bar{x}$, with $\theta = \frac{360^\circ}{m}$, and the group of rotations forms a cyclic subgroup $\langle \rho_{\theta} \rangle \cong C_m$. Each angle $\alpha_i \in [0^\circ, 180^\circ)$ specifies the orientation of a reflection axis passing through $ \bar{x}$ and forming an angle $\alpha_i$ with the $x$-axis. The corresponding reflection across this axis is denoted $R_{\alpha_i}$, and is also contained in $\Sym(X)$.

\begin{example}
Consider the finite configuration $X = \{ x_1, x_2, x_3 \}$ in the plane, where the points are given by
\[
x_1 = \left(0, 1\right), \quad
x_2 = \left(-\frac{\sqrt{3}}{2}, -\frac{1}{2}\right), \quad
x_3 = \left(\frac{\sqrt{3}}{2}, -\frac{1}{2}\right).
\]
These three points form an equilateral triangle centered at the centroid
\[
\bar{x} = \frac{1}{3}(x_1 + x_2 + x_3) = (0,0).
\]

\begin{figure}[h]
\centering
\begin{tikzpicture}[>=Stealth]

\begin{scope}[scale=2]
\coordinate (x1) at (0,1);
\coordinate (x2) at (-0.866, -0.5);
\coordinate (x3) at (0.866, -0.5);

\coordinate (G) at (0,0);

\draw[->, color=axisgray, text=black] (-1.2,0) -- (1.2,0) node[right] {$x$};
\draw[->, color=axisgray, text=black] (0,-0.8) -- (0,1.2) node[above left] {$y$};

\draw[line width=1.2pt, color=mainblue] (x1) -- (x2) -- (x3) -- cycle;

\filldraw[color=mainblue] (x1) circle (0.02) node[above right, text=black] {$x_1$};
\filldraw[color=mainblue] (x2) circle (0.02) node[below left, text=black] {$x_2$};
\filldraw[color=mainblue] (x3) circle (0.02) node[below right, text=black] {$x_3$};

\filldraw[color=mainblue] (G) circle (0.02) node[above left, text=black] {$\bar{x}$};

\draw[dashed, color=symred, line width=0.8pt] (G) -- (0,1.3);
\draw[dashed, color=symred, line width=0.8pt] (G) -- (-1.0392,-0.6);
\draw[dashed, color=symred, line width=0.8pt] (G) -- (1.0392,-0.6);

\draw[->, color=arcangle, line width=0.8pt] (0.172,-0.1) arc[start angle=-30, end angle=90, radius=0.2];
\node[text=black] at (0.32,0.15) {\scriptsize $120^\circ$};
\end{scope}

\begin{scope}[xshift=3.5cm]
\node at (0,0) {$\Longrightarrow$};
\end{scope}

\begin{scope}[scale=0.8, xshift=7cm]
\coordinate (x1) at (0,1);
\coordinate (x2) at (-0.866, -0.5);
\coordinate (x3) at (0.866, -0.5);

\draw[line width=1.2pt, color=mainblue] (0,0) circle (1);

\filldraw[color=mainblue] (0,0) circle (0.03) node[above right, text=black] {$\bar{x}$};

\draw[line width=0.8pt, color=mainblue] (0,0) -- (x1);
\draw[line width=0.8pt, color=mainblue] (0,0) -- (x2);
\draw[line width=0.8pt, color=mainblue] (0,0) -- (x3);

\node[text=black] at (0,-1.5) {$m=3$};
\end{scope}

\end{tikzpicture}
\caption{The finite configuration $X = \{x_1, x_2, x_3\}$ forming an equilateral triangle centered at the centroid $\bar{x}=(0,0)$. The symmetry group $\Sym(X)$ is the dihedral group $D_3$, encoded by the triple $(\bar{x}, 3, \{30^\circ, 90^\circ, 150^\circ\})$, which specifies the centroid, the order of rotational symmetry ($120^\circ$), and the directions of the reflection axes.}
\label{fig:triangle_symmetry_encoded}
\end{figure}

The symmetry group $\Sym(X)$ is the dihedral group $D_3$, consisting of 6 elements: three rotations and three reflections. The rotational symmetries are generated by the rotation $\rho_{120^\circ}$ about $\bar{x}$, forming the cyclic subgroup $\langle \rho_{120^\circ} \rangle \cong C_3$. The reflection symmetries correspond to reflections across three axes passing through $\bar{x}$ and aligned with the lines connecting $\bar{x}$ to each vertex.

Thus, the symmetry group $\Sym(X)$ can be encoded by the triple
\[
(\bar{x}, 3, \{30^\circ, 90^\circ, 150^\circ\}),
\]
where $\bar{x} = (0,0)$ is the centroid; $3$ indicates the order of rotational symmetry generated by $\rho_{120^\circ}$; and the angles $30^\circ$, $90^\circ$, and $150^\circ$ denote the directions of the reflection axes relative to the $x$-axis.
\end{example}

Now, let us consider the three-dimensional case. In this setting, the symmetry group of a finite configuration is a subgroup of the orthogonal group $O(3)$. Recall that $O(3)$ consists of all real $3 \times 3$ orthogonal matrices, i.e.,
\[
O(3) = \left\{ A \in \mathrm{GL}_3(\mathbb{R}) \,\middle|\, A^T A = I \right\}.
\]
Each element $A \in O(3)$ satisfies either $\det A = 1$ or $\det A = -1$. Elements with $\det A = 1$ correspond to rotations, while elements with $\det A = -1$ represent compositions of rotations with reflections.

Every matrix $A \in O(3)$ can be represented by a tuple
\[
(\theta, \mathbf{u}, \mathbf{n}) \in [0, \pi] \times S^2 \times (S^2 \cup \{\mathbf{0}\}),
\]
where
\begin{itemize}
  \item $\theta \in [0, \pi]$ denotes the rotation angle;
  \item $\mathbf{u} \in S^2$ is the unit vector specifying the rotation axis;
  \item $\mathbf{n} \in S^2 \cup \{\mathbf{0}\}$ is the unit normal vector of the reflection plane, or the zero vector $\mathbf{0}$ if no reflection is involved.
\end{itemize}

When $\mathbf{n} = \mathbf{0}$, the matrix $A$ represents a proper rotation given by the Rodrigues formula
\[
A = R_{\mathbf{u}}(\theta) = I + \sin\theta\,[\mathbf{u}]_\times + (1 - \cos\theta)\,[\mathbf{u}]_\times^2,
\]
where
\[
[\mathbf{u}]_\times =
\begin{bmatrix}
0 & -u_3 & u_2 \\
u_3 & 0 & -u_1 \\
-u_2 & u_1 & 0
\end{bmatrix}.
\]
denotes the skew-symmetric matrix of $\mathbf{u} = (u_1, u_2, u_3)^T \in S^2$.
When $\mathbf{n} \in S^2$, the matrix $A$ corresponds to an improper rotation, expressed as a composition of a proper rotation and a reflection
\[
A = R_{\mathbf{u}}(\theta) \cdot R_f, \quad \text{where} \quad R_f = I - 2\mathbf{n} \mathbf{n}^T.
\]

Based on the above, the symmetry group of a finite point configuration $X$ can be represented by the following data
\begin{equation*}
  (\bar{x}, \{(\theta_i,\mathbf{u}_i)\},\{\mathbf{n}_j\}),
\end{equation*}
where
\begin{itemize}
  \item $\bar{x}$ is the centroid of the configuration;
  \item Each pair $(\theta_i, \mathbf{u}_i)$ corresponds to a rotational symmetry, with rotation angle $\theta_i$ and rotation axis $\mathbf{u}_i$ passing through $\bar{x}$;
  \item Each vector $\mathbf{n}_j$ specifies a reflectional symmetry, where the reflection plane is perpendicular to $\mathbf{n}_j$ and passes through $\bar{x}$.
\end{itemize}
The above data encodes the symmetry patterns of a configuration in three-dimensional space and serves as a basis for effective visualization.

\begin{example}
Consider a regular tetrahedron in $\mathbb{R}^3$ with vertices given by
\begin{equation*}
  x_1 = (1, 1, 1), \quad x_2 = (-1, -1, 1), \quad x_3 = (-1, 1, -1), \quad x_4 = (1, -1, -1).
\end{equation*}
The centroid of the configuration is located at the origin:
\[
\bar{x} = \frac{1}{4}(x_1 + x_2 + x_3 + x_4) = (0, 0, 0).
\]

The proper rotational symmetry group (orientation-preserving symmetries) of the tetrahedron is the alternating group $A_4$, which has 12 elements. The full symmetry group, denoted $T_d$, includes both rotations and reflections, isomorphic to the symmetric group $S_4$, and contains 24 elements.

The symmetry data of the configuration can be represented as
\[
(\bar{x}, \{(\theta_i, \mathbf{u}_i)\}, \{\mathbf{n}_j\}),
\]
where
\begin{itemize}
  \item Rotational symmetries $(\theta_i, \mathbf{u}_i)$ (12 elements):
  \begin{itemize}
    \item \emph{Identity transformation} (1 element).

    \item \emph{Eight rotations of order 3} ($120^\circ$ and $240^\circ$) (Four axes of 3-fold rotation symmetry, passing through each vertex and the centroid of the opposite face):
    \[
    \begin{aligned}
    &\mathbf{u} = \tfrac{1}{\sqrt{3}}(1,1,1), \quad \theta = 120^\circ, 240^\circ; \\
    &\mathbf{u} = \tfrac{1}{\sqrt{3}}(-1,-1,1), \quad \theta = 120^\circ, 240^\circ; \\
    &\mathbf{u} = \tfrac{1}{\sqrt{3}}(-1,1,-1), \quad \theta = 120^\circ, 240^\circ; \\
    &\mathbf{u} = \tfrac{1}{\sqrt{3}}(1,-1,-1), \quad \theta = 120^\circ, 240^\circ.
    \end{aligned}
    \]

    \item \emph{Three rotations of order 2} ($180^\circ$) (Three axes of 2-fold rotation symmetry, connecting midpoints of opposite edges):
    \[
    \begin{aligned}
    &\mathbf{u} = (1, 0, 0), \quad \theta = 180^\circ; \\
    &\mathbf{u} = (0, 1, 0), \quad \theta = 180^\circ; \\
    &\mathbf{u} = (0, 0, 1), \quad \theta = 180^\circ.
    \end{aligned}
    \]
  \end{itemize}

  \item Reflection symmetries $\mathbf{n}_j$: There are six reflection planes. Each plane passes through $\bar{x}$, contains one edge of the tetrahedron, and bisects the opposite edge. The unit normal vectors for these planes are given by the permutations of $(1, -1, 0)$:
  \[
  \begin{aligned}
  \mathbf{n}_1 &= \tfrac{1}{\sqrt{2}}(1, -1, 0), &
  \mathbf{n}_2 &= \tfrac{1}{\sqrt{2}}(1, 1, 0), \\
  \mathbf{n}_3 &= \tfrac{1}{\sqrt{2}}(1, 0, -1), &
  \mathbf{n}_4 &= \tfrac{1}{\sqrt{2}}(1, 0, 1), \\
  \mathbf{n}_5 &= \tfrac{1}{\sqrt{2}}(0, 1, -1), &
  \mathbf{n}_6 &= \tfrac{1}{\sqrt{2}}(0, 1, 1).
  \end{aligned}
  \]
\end{itemize}
Together, the above data encodes the full symmetry structure of the tetrahedral configuration in three-dimensional Euclidean space. The full symmetry group of the regular tetrahedron, denoted by $T_d \cong S_4$, is a finite subgroup of $O(3)$ of order 24. It consists of the following elements:
\begin{itemize}
  \item 1 identity transformation;
  \item 8 rotations by $120^\circ$ and $240^\circ$ about the four vertex-to-opposite-face axes;
  \item 3 rotations by $180^\circ$ about the three axes connecting midpoints of opposite edges;
  \item 6 reflections across planes passing through the centroid;
  \item 6 improper rotations (roto-reflections of order 4), corresponding to the compositions of $90^\circ$ rotations about axes through midpoints of opposite edges with reflections perpendicular to those axes.
\end{itemize}
These 24 orthogonal transformations form the full symmetry group of the tetrahedron and capture both its rotational and reflectional symmetries.
\end{example}

Finally, it is worth noting that even if two finite configurations $X$ and $Y$ have isomorphic symmetry groups, the specific symmetry transformations they admit are generally different. Each symmetry is inherently tied to both the position and orientation of the configuration. In our previous description of rotational and reflectional symmetries, we explicitly incorporated positional information (via the centroid) and directional information (via rotation angles and normal vectors). If we disregard positional information, we may define an equivalence relation on symmetries: two symmetries $\pi$ and $\pi'$ are considered equivalent if they differ only by a translation. Under this equivalence, the centroid component of the symmetry data can be omitted when describing the structure of the symmetry group.

\subsection{Computation of symmetry groups}

In practical computations, the metric space under consideration is typically the Euclidean space. Since the isometry group of $\mathbb{R}^k$ is the Euclidean group $E(k) = O(k) \ltimes \mathbb{R}^k$,
the symmetry group of a finite configuration can be regarded as a subgroup of the orthogonal group $O(k)$ after centering the configuration. In the case $k = 2$, the orthogonal group $O(2)$ consists of all rotations and reflections that fix the centroid in the plane. Therefore, we only need to consider rotational and mirror symmetries in this setting.

\begin{lemma}\label{lemma:computation_symmetry}
Let $X$ be an $n$-configuration in $\mathbb{R}^k$ whose centroid is located at the origin. Then the symmetry group of $X$ is given by
\[
\Sym(X) = \mathrm{Aut}(X) \cap O(k),
\]
where $\mathrm{Aut}(X)$ is the automorphism group of the configuration $X$ (viewed as a permutation of points), and $O(k)$ is the orthogonal group.
\end{lemma}

\begin{proof}
By definition, we have $\Sym(X) \leq \mathrm{Aut}(X)$, since any symmetry must permute the points of $X$. Moreover, since the centroid of $X$ lies at the origin and any symmetry must preserve the centroid, it follows that $\Sym(X) \leq O(k)$. Hence, we have
\[
\Sym(X) \leq \mathrm{Aut}(X) \cap O(k).
\]
Conversely, let $\pi \in \mathrm{Aut}(X) \cap O(k)$. Then $\pi$ is both a permutation of the configuration $X$ and an orthogonal transformation of $\mathbb{R}^k$. In particular, $\pi$ preserves all distances and maps $X$ onto itself. Therefore, we have $\pi \in \Sym(X)$. It follows that
\[
\Sym(X) = \mathrm{Aut}(X) \cap O(k),
\]
as claimed.
\end{proof}

Based on Lemma~\ref{lemma:computation_symmetry}, a naive approach to compute $\Sym(X)$ is to check which permutations of $X$ are isometries by verifying if they preserve pairwise distances. However, this method quickly becomes computationally infeasible as $n$ grows, due to the factorial growth of the number of permutations ($n!$).

To improve efficiency, we leverage the radial symmetry about the centroid. First, we compute the centroid of the point set $X$, denoted by $\bar{x}$. Next, we compute the distance from each point in $X$ to the centroid $\bar{x}$, which defines a function
\begin{equation*}
  \phi : X \to \mathbb{R}, \quad x \mapsto \phi(x) = d(x, \bar{x}).
\end{equation*}
Let $D = \operatorname{im}(\phi)$ denote the set of all distances from points in $X$ to the centroid $\bar{x}$.

\begin{proposition}\label{proposition:computation_symmetry}
For each $r \in D$, denote $X_r = \{ x \in X \mid \phi(x) = r \}$. Then we have
\begin{equation*}
  \Sym(X) = \bigcap\limits_{r \in D} \Sym(X_r).
\end{equation*}
\end{proposition}

\begin{proof}
Suppose $\pi \in \Sym(X)$. By definition, we have $\pi(X) = X$, which implies that
\begin{equation*}
  \pi(\bar{x}) = \bar{x}.
\end{equation*}
Moreover, for any $x \in X$, since $\pi$ is an isometry, we have
\begin{equation*}
  d(x, \bar{x}) = d(\pi(x), \pi(\bar{x})) = d(\pi(x), \bar{x}).
\end{equation*}
It follows that $\phi(x) = \phi(\pi(x))$, and hence $\pi(X_r) = X_r$ for every $r \in D$. This shows that $\pi \in \Sym(X_r)$ for all $r \in D$, and thus
\[
\pi \in \bigcap\limits_{r \in D} \Sym(X_r).
\]
Therefore, $\Sym(X) \subseteq \bigcap\limits_{r \in D} \Sym(X_r)$.

Conversely, suppose $\pi \in \bigcap\limits_{r \in D} \Sym(X_r)$. Note that
\begin{equation*}
  X = \bigsqcup\limits_{r \in D} X_r,
\end{equation*}
so we have
\begin{equation*}
  \pi(X) = \bigsqcup\limits_{r \in D} \pi(X_r) = \bigsqcup\limits_{r \in D} X_r = X.
\end{equation*}
Since $\pi$ preserves the partition of $X$ into the level sets $X_r$ and acts as an isometry on each $X_r$, it follows that $\pi$ is an isometry of $X$. Hence, $\pi \in \Sym(X)$. Thus, we have $\bigcap\limits_{r \in D} \Sym(X_r)\subseteq \Sym(X)$, which completes the proof.
\end{proof}

Proposition~\ref{proposition:computation_symmetry} suggests a divide-and-conquer strategy. Instead of searching for symmetries in the entire configuration $X$, we can partition $X$ into radial subsets $X_r$ based on their distance to the centroid. By computing the symmetry group for each radial subset individually and then taking their intersection, we recover the full symmetry group of $X$. This approach is particularly effective when the points are distributed across multiple concentric shells, as it reduces the complexity of the symmetry search problem for each subset.

The centroid of the configuration is first translated to the origin. Then, by computing distances from each point to the origin, the configuration is partitioned into subsets $X_r$ of points lying at the same radial distance. According to Proposition~\ref{proposition:computation_symmetry}, the full symmetry group $\Sym(X)$ can be recovered as the intersection of the symmetry groups of these radial subsets. The Algorithm \ref{alg:symmetry} formalizes this procedure.

\begin{algorithm}[h]
\caption{Compute Symmetry Group $\Sym(X)$ for a Finite Configuration $X \subset \mathbb{R}^k$}
\label{alg:symmetry}
\begin{algorithmic}[1]
\Require A finite configuration $X = \{x_1, x_2, \dots, x_n\} \subset \mathbb{R}^k$
\Ensure The symmetry group $\Sym(X)$

\State Compute the centroid $\bar{x} = \frac{1}{n} \sum_{i=1}^{n} x_i$
\State Translate all points so that $\bar{x}$ is at the origin: $X \gets \{x_i - \bar{x} \mid x_i \in X\}$
\State Compute the distance function $\phi(x) = \|x\|$ for each $x \in X$
\State Let $D = \operatorname{im}(\phi)$ be the set of distinct distances from the origin
\For{each $r \in D$}
    \State Denote $X_r = \{x \in X \mid \phi(x) = r\}$
    \State Compute the local symmetry group $\Sym(X_r)$ using exhaustive isometry testing:
    \For{each permutation $\pi$ of $X_r$}
        \If{$\pi$ preserves all pairwise distances in $X_r$}
            \State Determine if $\pi$ is a rotation or a reflection
            \State Record $\pi$ as a valid symmetry of $X_r$
        \EndIf
    \EndFor
\EndFor
\State Compute the intersection $\Sym(X) = \bigcap_{r \in D} \Sym(X_r)$
\State \Return $\Sym(X)$
\end{algorithmic}
\end{algorithm}

\subsection{Multiscale computation of symmetries}

In Section~\ref{section:categorification}, we introduced two representations of persistent symmetries: the symmetry barcode and the polybarcode.
Both constructions aim to capture the persistence of symmetries over time within a given persistence configuration, tracking when a symmetry acts invariantly on the evolving dataset. Whereas the polybarcode represents this persistence as a set of discrete time points, the symmetry barcode encodes it as a collection of intervals, enabling a more direct visualization of symmetry duration. Moreover, the polybarcode is more naturally aligned with geometric perspectives, whereas the symmetry barcode lends itself well to topological data analysis through its intuitive bar representation.

\begin{table}[htbp]
  \centering
  \small
  \renewcommand{\arraystretch}{1.3}
  \setlength{\tabcolsep}{12pt}
  
  \rowcolors{2}{gray!10}{white} 
  
  \begin{tabular}{@{}m{3.8cm} m{5.6cm} m{5.6cm}@{}}
    \toprule
    \textbf{Aspect} & \textbf{Polybarcode} & \textbf{Symmetry barcode} \\
    \midrule
    \textbf{Definition} &
    For a fixed symmetry $\pi$, records the set of times $I(\pi) = \{t \mid \pi(\mathcal{F}_t) = \mathcal{F}_t\}$ as a subset of $\mathbb{R}$ &
    Multiset of intervals $\{[a_i,b_i)\}$, where each interval represents a lifespan of a symmetry under persistence \\
    \midrule
    \textbf{Persistent type} &
    Records the time set during which a specific symmetry exists &
    Tracks the continuous duration over which a symmetry persists \\
    \midrule    
    \textbf{Data representation} &
    Set-valued function associating symmetries to time subsets; possibly disconnected sets &
    Collection of intervals (bars) \\
    \midrule
    \textbf{Visualization} &
    Abstract set-theoretic object, less intuitive for direct plotting &
    Directly visualized as barcode diagrams \\
    \midrule 
    \textbf{Distance and comparison} &
    Expansion distance, symmetric difference distance, and interleaving distance &
    Bottleneck distance, Wasserstein distance, and interleaving distance \\
    \midrule
    \textbf{Suitability for analysis} &
    Suitable for geometric studies emphasizing exact times of persistence &
    More suited for topological data analysis leveraging interval persistence \\
    \bottomrule
  \end{tabular}
  \caption{Comparison between polybarcode and symmetry barcode}
  \label{tab:poly_vs_sym_barcode}
\end{table}

Let us first review the span of symmetry groups. Given a morphism $f: X \to Y$ in the category $\mathbb{S}_n(M)$, one can associate a span of symmetry groups
\begin{equation*}
  \mathrm{Sym}(X)\xleftarrow{f^{\flat}} \mathrm{Sym}_{f}(X) \xrightarrow{f^{\sharp}} \mathrm{Sym}(Y),
\end{equation*}
where the middle term is given by
\begin{equation*}
  \mathrm{Sym}_f(X) = \left\{ \sigma \in \mathrm{Sym}(X) \mid f \circ \sigma \circ f^{-1} \in \mathrm{Sym}(Y) \right\}.
\end{equation*}

Intuitively, an element $\sigma \in \mathrm{Sym}(X)$ belongs to $\mathrm{Sym}_f(X)$ if and only if the conjugation $f \circ \sigma \circ f^{-1}$ preserves the geometric distances in $Y$. That is, $\sigma$ is compatible with the geometry of $Y$ via the morphism $f$. The basic idea of the algorithm for computing $\mathrm{Sym}_f(X)$ is as follows. We begin by computing the full symmetry group $\mathrm{Sym}(X)$, using a method such as normalizing $X$ so that its centroid lies at the origin, and then enumerating all permutations that preserve pairwise distances. For each permutation $\sigma \in \mathrm{Sym}(X)$, we construct a candidate symmetry of $Y$ by conjugating $\sigma$ through $f$. Specifically, we define the map $\tilde{\sigma} := f \circ \sigma \circ f^{-1}: Y \to Y$. We then test whether $\tilde{\sigma}$ is indeed an element of $\mathrm{Sym}(Y)$. This involves checking whether $\tilde{\sigma}$ preserves all pairwise distances, i.e., whether it is an isometry of $Y$. If this condition holds, we retain the original $\sigma$ as an element of $\mathrm{Sym}_f(X)$. Otherwise, $\sigma$ is discarded. Repeating this verification process for each $\sigma \in \mathrm{Sym}(X)$, we obtain the desired subgroup $\mathrm{Sym}_f(X)$.

Now we revisit the notion of the symmetry barcode introduced in Section~\ref{section:categorification}. Let $\mathcal{F}: (\mathbb{R}, \leq) \to \mathcal{S}_{n}(M)$ be a persistence configuration.
A symmetry bar associated to a nontrivial symmetry $\pi \in \mathrm{Sym}(\mathcal{F}_a)$ is defined as the maximal interval $[a, b)$ such that, for all $t \in [a, b)$, the element $\pi$ lies in $\mathrm{Sym}_{f_{a,t}}(\mathcal{F}_a)$. The corresponding symmetry barcode is the collection of all such symmetry bars.
Algorithm~\ref{algorithm:symmetry_barcode} presents a straightforward procedure for computing the symmetry barcode of a persistent symmetry group. We expect that more efficient and advanced algorithms will be developed in the future.

\begin{algorithm}[H]
\caption{Symmetry Barcode Extraction}\label{algorithm:symmetry_barcode}
\begin{algorithmic}[1]
\Require A discrete parameter list $A = [a_0, a_1, \dots, a_n]$
\Require A list of symmetry groups $G[i] = \mathrm{Sym}(\mathcal{F}_{a_i})$
\Require Transition maps $\varphi[i][j]$ mapping $\sigma \in G[i]$ to a candidate isometry $\tilde{\sigma}$ of $\mathcal{F}_{a_j}$
\Ensure A symmetry barcode as a multiset of intervals $[a_i, a_j)$

\State Initialize an empty list $\mathcal{B} \gets [\,]$ \Comment{barcode list}

\For{$i = 0$ to $n$}
  \ForAll{$\sigma \in G[i]$}
    \State $birth \gets a_i$
    \State $death \gets \infty$
    \For{$j = i+1$ to $n$}
      \State $\tilde{\sigma} \gets \varphi[i][j](\sigma)$
      \If{$\tilde{\sigma} \notin G[j]$}
        \State $death \gets a_j$
        \State \textbf{break}
      \EndIf
    \EndFor
    \State Append interval $[birth, death)$ to $\mathcal{B}$
  \EndFor
\EndFor

\State \Return $\mathcal{B}$
\end{algorithmic}
\end{algorithm}

A direct simplification of the algorithm for computing the symmetry barcode arises from the observation that a symmetry $\pi$ and its inverse $\pi^{-1}$ always share the same symmetry bar, see Proposition~\ref{proposition:inverse_bar}. Therefore, when computing persistence, it suffices to consider only half of the elements in each symmetry group.

We now turn to a discussion of the polybarcode of persistence configurations.
To compute the polybar $I(\pi)$ of a given symmetry $\pi \in \Iso(M)$, we proceed as follows. Let $\mathcal{F}: \mathbb{R} \to \mathcal{S}_n(M)$ be a persistence configuration, and let $T = \{t_1, t_2, \dots, t_m\} \subset \mathbb{R}$ be a finite time grid over which we examine the configuration. For each time $t_i \in T$, we evaluate the configuration $\mathcal{F}_{t_i}$, and check whether the symmetry $\pi$ preserves the configuration at that time, i.e., whether $\pi(\mathcal{F}_{t_i}) = \mathcal{F}_{t_i}$. If this condition holds, then $t_i$ is included in the polybar $I(\pi)$. Formally, we define
\[
I(\pi) = \left\{t \in T \,\middle|\, \pi(\mathcal{F}_t) = \mathcal{F}_t \right\}.
\]
This set $I(\pi)$ thus records the collection of time points at which $\pi$ acts as a true symmetry of the configuration. In practical computations, equality between configurations may be tested up to a small numerical tolerance to account for floating-point errors or noise.

In contrast to the symmetry barcode, where the symmetry itself may evolve as the filtration time progresses, reflecting the changing group structure, the polybar tracks a fixed symmetry $\pi$ throughout the filtration. In practice, it is often useful to consider certain symmetries, such as those differing by a translation, as equivalent. Such equivalence helps produce denser polybars and eliminates the interference caused by translational variability.

\subsection{Computation of symmetry defect}

In Section~\ref{section:symmetry_defect}, we introduced the notion of symmetry defect as a quantitative measure of the asymmetry of a finite configuration that may not be perfectly symmetric. When computing symmetry groups and identifying existing symmetries, the candidate symmetries are typically drawn from the automorphism group of the configuration. Since the automorphism group is finite, such computations are effectively reduced to a finite search problem.

However, in the computation of the symmetry defect, the set of candidate symmetries is generally taken to be the full group of isometries of the Euclidean space, including arbitrary translations, rotations, and reflections. This renders a direct search over candidate symmetries infeasible due to the continuum nature of the search space. As a result, we adopt certain approximation schemes to identify and evaluate candidate symmetries in practice.

\subsubsection*{Approximate symmetries in Euclidean plane}

First, since the symmetries of a finite configuration consist only of rotations and reflections, we disregard translations in the selection of candidate isometries. Second, for a given finite configuration $X$, we restrict our attention to approximate symmetries whose rotation centers and reflection axes are required to pass through the centroid of $X$. This constraint significantly reduces the number of candidate symmetries to be considered and facilitates the computation.

For a finite configuration $X$ in the Euclidean plane, we approximate symmetries by considering a discrete set of reflections and rotations centered at the centroid $O$ of $X$. We first consider a set of $m$ reflection axes passing through $O$, each making an angle of $i \cdot \frac{180^\circ}{m}$ with the $x$-axis for $i = 0, 1, \dots, m-1$. Each such line determines a reflection symmetry $\pi_i$. For each $\pi_i$, we compute the symmetry defect using distance
\[
  \mu(X,\pi_i) = \left( \inf_{\gamma \neq {\pi_i}|_X}   \sum_{x \in X} d(x, \gamma(x))^2  \right)^{1/2}.
\]
Next, we consider $n-1$ rotations about the point $O$ by angles $j \cdot \frac{360^\circ}{n}$ for $j = 1,2, \dots, n-1$, resulting in a set of approximate rotational symmetries $\sigma_j$. Similarly, we compute their symmetry defects
\[
\mu(X, \sigma_j) = \left( \inf_{\gamma \neq {\sigma_j}|_X}   \sum_{x \in X} d(x, \gamma(x))^2  \right)^{1/2}.
\]
We record all such approximate symmetries along with their corresponding symmetry defects. This yields a quantitative profile of the approximate symmetries of $X$ and their degrees of deviation from perfect symmetry. When $m$ and $n$ are sufficiently large, the set of sampled candidate symmetries becomes dense enough to approximate the true symmetry defect of $X$ within any desired level of precision. Those symmetries with smaller defect values provide meaningful approximations of the actual symmetries of $X$.

However, particular care is needed when interpreting small values of the symmetry defect $\mu(X, \sigma_j)$, especially in the case of rotations by small angles. Consider, for instance, the rotation $\sigma_1$ by a small angle. A low symmetry defect in this case does not necessarily imply that $\sigma_1$ approximates a genuine symmetry of $X$. Instead, it may arise simply because the rotation maps each point $x \in X$ to a nearby location, allowing the natural bijection $\gamma = \sigma_1|_X$ to nearly minimize the defect.
Since the definition of $\mu(X, \sigma_j)$ explicitly excludes the bijection $\gamma = \sigma_j|_X$ from the infimum, we must determine whether this excluded bijection would have been the actual minimizer. If so, the small defect value reflects only a near-identity point matching caused by the small rotation angle, rather than indicating meaningful approximate symmetry. Such cases should be carefully distinguished and excluded from further analysis.

A similar procedure can also be applied to compute symmetry measures, with the main difference lying in the specific numerical expressions used in the calculations.

\subsubsection*{Approximate symmetries in 3-dimensional Euclidean space}

To analyze the approximate symmetries of a finite configuration $X$ in 3-dimensional Euclidean space, we consider three types of orthogonal transformations centered at the centroid $O$ of $X$: reflections, rotations, and compositions thereof. Each class of transformations gives rise to a collection of candidate symmetries against which symmetry defects and measures can be computed.

\paragraph{1. Mirror Symmetries.}
Let $m \in \mathbb{N}$ be a resolution parameter. We sample unit vectors $\mathbf{n}_i \in S^2$ on the unit sphere, indexed by $i = 1, \dots, m$, using a uniform or geodesic grid on the sphere. Each vector $\mathbf{n}_i$ determines a mirror plane $H_i$ passing through the origin with normal $\mathbf{n}_i$. The corresponding reflection $\pi_i$ across $H_i$ is a linear isometry satisfying
\[
\pi_i(\mathbf{x}) = \mathbf{x} - 2 \langle \mathbf{x}, \mathbf{n}_i \rangle \mathbf{n}_i.
\]
We collect all such $\pi_i$ as candidate mirror symmetries.

\paragraph{2. Rotations.}
Let $l, n \in \mathbb{N}$. We again sample unit vectors $\mathbf{v}_j \in S^2$ for $j = 1, \dots, l$, which serve as rotation axes through the origin. For each axis $\mathbf{v}_j$, we define a family of rotation angles $\theta_k = \frac{2\pi k}{n}$ for $k = 1, \dots, n$. The rotation $\sigma_{j,k}$ about $\mathbf{v}_j$ by angle $\theta_k$ is a proper orthogonal transformation
\[
\sigma_{j,k} = \exp(\theta_k [\mathbf{v}_j]_\times),
\]
where $[\mathbf{v}_j]_\times$ denotes the skew-symmetric matrix associated to $\mathbf{v}_j$ for the cross product. We collect all such $\sigma_{j,k}$ as candidate rotational symmetries.

\paragraph{3. Compositions of Reflections and Rotations.}
Any element $T \in O(3)$ with determinant $-1$ can be written as the composition of a fixed mirror reflection and a rotation. In particular, we may fix the reflection $\overline{\pi}$ to be the mirror reflection with respect to the $xOy$-plane, given by
\[
\overline{\pi}(x, y, z) = (x, y, -z).
\]
Then for every $T \in O(3)$ with $\det(T) = -1$, there exists a rotation $\sigma \in SO(3)$ such that
\[
T = \overline{\pi} \circ \sigma.
\]
Therefore, to construct a complete set of candidate symmetries in $O(3)$, we begin with a set of candidate rotational symmetries $\{\sigma_{j,k}\}$, and then for each $\sigma_{j,k}$, we include the compositions $\overline{\pi} \circ \sigma_{j,k}$. This yields a new set of candidate symmetries.

Altogether, the candidate symmetry set is given by
\[
\mathsf{Cand}(X) = \{ \pi_i \}_{i=1}^m \cup \{ \sigma_{j,k} \}_{j=1,k=1}^{l,n} \cup \{ \overline{\pi} \circ\sigma_{j,k} \}_{j=1,k=1}^{l,n} .
\]
This collection forms a finite, symmetry-rich approximation of the full orthogonal group $O(3)$, suitable for the analysis of approximate symmetries in finite point clouds.

Finally, for each element in the candidate symmetry set $\mathsf{Cand}(X)$, we compute the corresponding symmetry defect and symmetry measure, recording these values as features for subsequent data analysis. Moreover, similar to the two-dimensional case, it is necessary to exclude those candidate symmetries whose small symmetry defects arise solely from minor perturbations corresponding to very small rotation angles. Such near-identity transformations do not represent meaningful approximate symmetries and should be filtered out to avoid spurious results.

\begin{example}
Consider the planar point set
\[
X = \{(-1, 1),\  (1, 1),\  (1, -1),\  (-1, 0)\} \subset \mathbb{R}^2.
\]
The centroid of $X$ is
\[
\bar{x} = \frac{1}{4} \sum_{i=1}^4 x_i = \left(0, \frac{1}{4} \right).
\]
We consider mirror reflections through lines passing through the centroid, with directions given by $v_\theta = (\cos\theta, \sin\theta)$, where the angle $\theta = \frac{k\pi}{180}$, $k = 0,1,\dots,179$. The mirror reflection $\pi_\theta$ is given by
\[
\pi_\theta(x) = x - 2 \cdot \langle x - \bar{x}, v_\theta^\perp \rangle v_\theta^\perp,
\quad v_\theta^\perp = (-\sin\theta, \cos\theta).
\]
We also consider rotational symmetries about the centroid with rotation angle $\phi = \frac{2k\pi}{180}$, $k = 0, \dots, 179$. Each rotation $\sigma_\phi$ is computated by
\[
\sigma_\phi(x) = R_\phi(x - \bar{x}) + \bar{x}, \quad
R_\phi =
\begin{bmatrix}
\cos\phi & -\sin\phi \\
\sin\phi & \cos\phi
\end{bmatrix}.
\]

As a concrete example, consider the reflection $\pi_y$ about the vertical axis passing through the centroid. In the coordinate system centered at $\bar{x}$ (denoted by $X' = X - \bar{x}$), this operation corresponds to $(x, y)' \to (-x, y)'$.
Applying $\pi_y$ to each point in $X'$ gives
\[
Y = \{(1, 0.75),\ (-1, 0.75),\ (-1, -1.25),\ (1, -0.25)\}.
\]
We compute the squared Euclidean distance matrix $D$ between points in $X'$ and $Y$ as
\[
D = \begin{bmatrix}
4 & 0 & 4 & 5 \\
0 & 4 & 8 & 1 \\
4 & 8 & 4 & 1 \\
5 & 1 & 1 & 4
\end{bmatrix}.
\]
Using the optimal assignment (Hungarian) algorithm to minimize the sum of matched squared distances, one optimal matching is
\[
\begin{cases}
(-1, 0.75) \leftrightarrow (-1, 0.75), \quad d^2 = 0 \\
(1, 0.75) \leftrightarrow (1, 0.75), \quad d^2 = 0 \\
(1, -1.25) \leftrightarrow (1, -0.25), \quad d^2 = 1 \\
(-1, -0.25) \leftrightarrow (-1, -1.25), \quad d^2 = 1
\end{cases}
\]
The total squared distance is
\[
0 + 0 + 1 + 1 = 2,
\]
thus the symmetry defect is
\[
\mu(X, \pi_y) = \sqrt{2} \approx 1.414.
\]

Now, consider the case of rotation by $180^{\circ}$ about the centroid. In the centered coordinate system, the image set is
\[
Y = \{ (1, -0.75),\ (-1, -0.75),\ (-1, 1.25),\ (1, 0.25) \}.
\]
To compute the symmetry defect under this rotation, we consider all possible matchings $\gamma$ between points in $X'$ and $Y$, and calculate the sum of squared distances
\[
S(\gamma) = \sum_{i=1}^4 \| x'_i - y_{\gamma(i)} \|^2.
\]
By enumerating all permutations, the minimal sum is achieved by the matching
\[
\gamma = \begin{pmatrix}0 & 1 & 2 & 3 \\ 2 & 3 & 0 & 1 \end{pmatrix},
\]
with
\[
S(\gamma) = 1.
\]
Therefore, the symmetry defect is
\[
\mu(X, \sigma_{\pi}) = \sqrt{1} = 1.
\]

In fact, through algorithmic computation, the minimum symmetry defect of $0.6130$ is attained for a mirror symmetry axis forming an angle of approximately $171^\circ$ with the $x$-axis. The smallest rotational symmetry defect is $1$, which corresponds precisely to the rotation by $180^\circ$.

\end{example}

\subsection{Parametrization of symmetries for point cloud data}

We have outlined computational approaches for studying symmetries of finite configurations. For dynamic configurations, the evolution of the configuration naturally induces a filtration parameter given by time or space. In contrast, for static data such as point clouds, the construction of an appropriate parametrization allows for the extraction of multi-scale and multi-level symmetry features, providing more targeted schemes for practical computation and applications.

\subsubsection*{Parametrization of approximate symmetry groups}

Let $X$ be a finite point set in Euclidean space $\mathbb{R}^k$. For a given isometry $\pi \in O(k)$ (assuming the centroid of $X$ is at the origin), we have the symmetry defect function
\[
\widetilde{\mu}(X, -) : O(k) \to \mathbb{R}.
\]
Recall that, for a parameter $\varepsilon > 0$, we can define the approximate symmetry group by
\begin{equation*}
  \Sym_\varepsilon(X) = \left\{ \pi \in O(k) \mid \widetilde{\mu}(X, \pi) \leq \varepsilon \right\}.
\end{equation*}
As the filtration parameter $\varepsilon$ varies, we obtain a filtration of approximate symmetry groups. Specifically, let
\begin{equation*}
  \varepsilon_{0} < \varepsilon_{1} < \cdots < \varepsilon_{m}.
\end{equation*}
Then we have a nested filtration
\begin{equation*}
  \Sym_{\varepsilon_0}(X) \hookrightarrow \Sym_{\varepsilon_1}(X) \hookrightarrow \cdots \hookrightarrow \Sym_{\varepsilon_m}(X).
\end{equation*}
This indicates that, as $\varepsilon$ increases, new approximate symmetries emerge. Therefore, for each approximate symmetry, one can record its \emph{birth time}---the smallest value of $\varepsilon$ at which it appears. These data can be recorded as symmetry features of the point set $X$.

On the other hand, we may introduce a parametrization on the point set $X$ itself. Let $\bar{x}$ denote the centroid of $X$, and let
\begin{equation*}
  X_r = \{ x \in X \mid d(x, \bar{x}) \leq r \}.
\end{equation*}
In this way, we obtain a nested filtration of points indexed by the distance parameter $r$. A direct motivation for this construction comes from Proposition~\ref{proposition:computation_symmetry}:
if $\pi$ is a symmetry of $X$, namely
\begin{equation*}
  \mu(X, \pi) = 0,
\end{equation*}
then it follows that
\begin{equation*}
  \mu(X_r, \pi) = 0, \quad \text{for any } r > 0.
\end{equation*}
However, for real-world data, $\mu(X, \pi)$ and $\mu(X_r, \pi)$ are not necessarily directly related. This observation provides additional features for characterizing the symmetries present in the data. Moreover, for given parameters $r > 0$ and $\varepsilon > 0$, we obtain the parametrized approximate symmetry group $\Sym_{\varepsilon}(X_r)$. Fixing $r$ and increasing $\varepsilon$, a symmetry $\pi$ appears in $\Sym_{\varepsilon_0}(X_r)$ at some threshold $\varepsilon_0$. Conversely, fixing $\varepsilon$ and increasing $r$, a symmetry $\pi$ may emerge in $\Sym_{\varepsilon}(X_r)$ or disappear at some stage. Such behavior yields a richer collection of symmetry features, which can be recorded for the study of the symmetry structure of the data.

The basic idea for computing symmetry features using parametrized approximate symmetry groups is as follows. First, a suitable set of candidate symmetries $\Gamma$ is selected from $O(k)$ as uniformly as possible. Next, a sequence of distance parameters
\[
r_1 < r_2 < \cdots < r_p
\]
is chosen in ascending order. For each fixed $r_i$, we compute the symmetry defect for every symmetry in $\Gamma$. These symmetry defects correspond exactly to the birth times of the symmetries at distance $r_i$. Consequently, all desired symmetry features can be extracted from this process. The algorithmic procedure is presented in Algorithm \ref{alg:param_symmetry}.
Notably, for each symmetry $\pi$, its symmetry features can be represented by a binary matrix or visualized through birth-time curves, among other methods.

\begin{algorithm}[h]
\caption{Compute Symmetry Features from Parametrized Approximate Symmetry Groups}
\label{alg:param_symmetry}
\begin{algorithmic}[1]
\Require Finite point set $X = \{x_1,\dots,x_n\} \subset \mathbb{R}^k$;
         distance parameters $R = \{r_1 < \dots < r_p\}$;
         thresholds $E = \{\varepsilon_1 < \dots < \varepsilon_q\}$
\Ensure Feature list $\mathcal{F}$ of tuples $(\pi, r_i, \varepsilon_{\text{birth}})$

\State Compute centroid $\bar{x} = \frac{1}{n} \sum_{i=1}^n x_i$ and translate $X \gets \{x_i - \bar{x}\}$
\State Select a suitable candidate set $\Gamma \subset O(k)$ uniformly
\State Initialize $\mathcal{F} \gets \varnothing$

\For{$i=1$ to $p$}
  \State $X_{r_i} \gets \{ x \in X \mid \|x\| \le r_i \}$
  \ForAll{$\pi \in \Gamma$}
    \State Initialize $\varepsilon_{\text{birth}} \gets +\infty$
    \State Compute defect $d \gets \mu(X_{r_i}, \pi)$
    \For{$j=1$ to $q$}
      \If{$d \le \varepsilon_j$}
        \State Update $\varepsilon_{\text{birth}} \gets \min(\varepsilon_{\text{birth}}, \varepsilon_j)$
      \EndIf
    \EndFor
    \If{$\varepsilon_{\text{birth}} \neq +\infty$}
      \State Append $(\pi, r_i, \varepsilon_{\text{birth}})$ to $\mathcal{F}$
    \EndIf
  \EndFor
\EndFor

\State \Return $\mathcal{F}$
\end{algorithmic}
\end{algorithm}

Moreover, another similar algorithm for computing parametrized symmetries replaces symmetry defect with a symmetry measure. The key difference lies in the range and invariance properties of these functions: symmetry defect values range from zero to positive infinity and heavily depend on the distances between data points; thus, the defect varies under similarity transformations. In contrast, the symmetry measure takes values between zero and one, where one corresponds to perfect symmetry, and is invariant under similarity transformations. Since the algorithmic procedure based on the symmetry measure is almost identical to that of symmetry defect, we omit its detailed description here.

\subsubsection*{Parametrization by point sets}

Recall from Section~\ref{section:symmetry_defect} that for a given finite point set $X$ and symmetries $\pi, \sigma \in O(k)$, we define the distance
\[
  d_{X,p}(\pi, \sigma) = \left( \sum_{x \in X} d(\pi(x), \sigma(x))^p \right)^{1/p}.
\]
If $X$ is not symmetric, i.e., its symmetry defect is strictly positive, then $d_{X,p}(\pi, \sigma) = 0$ if and only if $\pi = \sigma$. Thus, $(O(k), d_{X,p})$ forms a metric space parameterized by the point set $X$.

Based on this idea, we can compute symmetry features from a different perspective. First, fix a finite subset $\Gamma \subseteq O(k)$, which serves as a discrete \emph{probe set} for the symmetries of $X$, since $O(k)$ is uncountable and cannot be handled exhaustively. The elements of $\Gamma$, i.e., candidate symmetries, can be viewed as points in the metric space $(O(k), d_{X,p})$. Naturally, we can construct simplicial complexes over $\Gamma$, such as the Vietoris-Rips complex, the alpha complex, or the \v{C}ech complex. Here, we take the Vietoris-Rips complex as an example and consider the metric
\[
  d_X = d_{X,2} = \left( \sum_{x \in X} \|\pi(x) - \sigma(x)\|^2 \right)^{1/2}.
\]
For any $\varepsilon > 0$, we obtain the Vietoris-Rips complex $\mathcal{R}_{\varepsilon}(\Gamma)$. Persistent homology techniques can then be applied to filtration $\{\mathcal{R}_{\varepsilon}(\Gamma)\}_{\varepsilon > 0}$.

In particular, the zero-dimensional Betti numbers and their corresponding barcodes capture the clustering behavior of symmetries on the set $X$: a large number of long bars indicates that the symmetries are distributed unevenly, while very short bars correspond to a more uniform distribution of symmetry degrees. This provides an intuitive but traditionally difficult-to-quantify property of the symmetry structure.

In the following, we present several approaches for sampling the probe set $\Gamma$.

\paragraph{Uniform sampling}
A straightforward approach to construct the probe set $\Gamma \subseteq O(k)$ is to sample elements uniformly over the orthogonal group. For the special orthogonal group $SO(k)$, uniform sampling can be achieved by generating random orthogonal matrices using techniques such as the QR decomposition of matrices with entries drawn from a standard normal distribution, or by sampling uniformly over parameterizations such as Euler angles or quaternions. To include reflections and obtain samples in the full orthogonal group $O(k)$, one can augment the uniformly sampled rotations with appropriate reflection operations. This method ensures a broad and representative coverage of the symmetry space, but its computational cost increases rapidly with dimension $k$.

\paragraph{Grid discretization}
Another common method is to discretize the parameter space of $O(k)$ into a finite grid. For instance, one may discretize the space of Euler angles or other parameterizations into evenly spaced points, thereby forming a finite and structured probe set. This approach provides explicit control over the sampling resolution and can be straightforward to implement. However, the dimensionality of the parameter space causes the number of grid points to grow exponentially, leading to potential computational infeasibility in higher dimensions.

\paragraph{Randomized sampling with importance weighting}
This method involves randomly sampling candidates from $O(k)$ according to a probability distribution that favors elements that are likely to correspond to approximate symmetries. Importance weighting can be based on preliminary estimates of symmetry defects or domain-specific heuristics. By focusing computational effort on promising regions of the symmetry space, this approach can increase efficiency while maintaining effectiveness in detecting relevant symmetries.

\paragraph{Heuristic and data-driven sampling}
Heuristic or data-driven methods utilize prior knowledge or learned models to guide the selection of candidate symmetries. For example, one might use historical data or machine learning techniques to predict likely symmetry transformations. These approaches can incorporate domain expertise and adaptively refine the probe set to improve the detection of symmetries with fewer candidate elements, although they may risk missing unexpected or rare symmetries.

\subsection{Local symmetry}

Locally symmetric spaces, such as locally symmetric complex manifolds, are classical objects in complex geometry, fundamental to the study of geometric rigidity and symmetry structures \cite{jost1986strong,mostow1973strong}.
In this final subsection, we present several ideas and strategies for analyzing local symmetry in data. 
These approaches are not only useful in symmetry analysis but also extend naturally to other forms of local geometric or topological computations.

In practical applications, point cloud data is more common and typically consists of a large number of data points. This often makes the direct computation of global symmetry both computationally intensive and insufficiently informative. In contrast, studying local symmetry can provide both computational advantages and richer structural insights. Since local neighborhoods contain fewer points, symmetry analysis becomes more efficient and naturally parallelizable.

Research on local symmetry can generally be conducted from two complementary perspectives: \textit{geometric local symmetry}, which focuses on the spatial distribution and configuration of points; and \textit{semantic local symmetry}, which incorporates application-specific meaning, where local symmetric patterns correspond to functionally or structurally meaningful units, such as molecular functional groups or repetitive elements in man-made structures.

The computation and application of geometric local symmetry involve a variety of technical approaches. Here, we present two conceptual strategies that serve as the foundation for many practical methods: cover-based analysis and local multi-scale analysis.

\paragraph{Cover-based analysis.}

Let $X$ be a finite point set, and let $\{U_i\}_{i=1}^N$ be a cover of $X$. For each $i = 1, \dots, N$, we compute local symmetry features such as the local symmetry group $\mathrm{Sym}(U_i)$ or the symmetry defect of $U_i$. These local descriptors can then be aggregated and used for subsequent data analysis.
A critical step in this approach is the selection of the cover $\{U_i\}$. A natural strategy is to partition $X$ into clusters that reflect local geometric or topological features. One common approach is $k$-nearest neighbor clustering in machine learning. Another method is based on topological connectivity: one constructs a neighborhood graph $G_X = (X, E)$ by connecting pairs of points that lie within a fixed radius $\epsilon$, and then defines each $U_i$ as a connected component or a cluster derived from graph-based algorithms such as spectral clustering or modularity maximization. These strategies are designed to balance locality and redundancy. Small subsets $U_i$ ensure computational efficiency and capture local structure accurately, while allowing overlaps between them preserves global coherence.

\paragraph{Local multi-scale analysis.}

Given a finite point set $X$ in Euclidean space, we fix a reference point, for example, the centroid of $X$, denoted by $\bar{x}$. Centered at $\bar{x}$, we consider a family of neighborhoods defined by $U_r = X \cap B(\bar{x}, r)$, where $B(\bar{x}, r)$ is the closed ball of radius $r$ centered at $\bar{x}$. As the radius $r$ increases from zero, the subset $U_r$ grows monotonically, forming a multi-scale filtration of $X$. At each scale, we may compute symmetry-related features such as the local symmetry group of $U_r$ or its symmetry defect. This process captures how local symmetry evolves as the region of interest expands, providing a multi-scale description of the geometric structure around $\bar{x}$.
Alternatively, the neighborhood of $\bar{x}$ can be defined using other selection criteria, such as cubical neighborhoods, strip-shaped regions, or $k$-nearest neighbor sets. In particular, for the $k$-nearest neighbor (k-NN) method, one defines $U_k(\bar{x})$ to be the set consisting of the $k$ points in $X$ closest to $\bar{x}$ with respect to the Euclidean metric. As $k$ increases, $U_k(\bar{x})$ forms an expanding sequence of subsets that also induces a natural multi-scale structure.

Moreover, the two approaches described above can be naturally combined by introducing a hierarchy of covers that reflect different levels of locality. Specifically, starting with a coarse cover $\mathcal{U}^{(0)} = \{U_i^{(0)}\}_{i \in I_0}$ of the point set $X$, we consider successive refinements
\[
\mathcal{U}^{(0)} \succ \mathcal{U}^{(1)} \succ \cdots \succ \mathcal{U}^{(p)},
\]
where each $\mathcal{U}^{(k)} = \{U_i^{(k)}\}_{i \in I_k}$ is a cover of $X$ and satisfies the refinement condition: for every $U_i^{(k)} \in \mathcal{U}^{(k)}$, there exists $U_j^{(k-1)} \in \mathcal{U}^{(k-1)}$ such that $U_i^{(k)} \subseteq U_j^{(k-1)}$. At each level $k$, one may compute local symmetry features. This multiscale structure enables localized symmetry analysis at varying resolutions, facilitating both fine-grained geometric insight and computational efficiency.

In addition to constructing a multiscale family of covers on $X$, one may also perform multiscale analysis within each cover element. Given a fixed cover $\{U_i\}_{i=1}^N$ of $X$, we consider, for each $U_i$, its barycenter $\bar{x}_i$, and define a one-parameter family of nested neighborhoods
\[
U_{i,r} = U_i \cap B(\bar{x}_i, r),
\]
where $B(\bar{x}_i, r)$ denotes the closed ball of radius $r$ centered at $\bar{x}_i$. As the parameter $r$ increases, the subset $U_{i,r}$ expands, giving rise to a local multiscale filtration. Such intra-cluster multiscale structures complement the global hierarchy of covers and provide enhanced flexibility and resolution in symmetry analysis.

Another approach to studying local symmetry is through semantic local symmetry, which emphasizes domain-specific meaning and structural roles. In practical applications, this often requires more refined analysis. For example, in molecular biology, local symmetric patterns arise in ring structures or functional groups of large molecules; in protein analysis, specific structural units may exhibit symmetry; in neural networks, recurring motifs often reflect symmetrical organization. More generally, in complex systems, any localized region with distinct structure or function may be a candidate for symmetry analysis.